# Some series and integrals involving the Riemann zeta function, binomial coefficients and the harmonic numbers

## Volume II(a)

Donal F. Connon

18 February 2008

## Abstract

In this series of seven papers, predominantly by means of elementary analysis, we establish a number of identities related to the Riemann zeta function, including the following:

$$\psi(u) = \sum_{k=0}^{\infty} \frac{1}{k+1} \sum_{j=0}^{k} \binom{k}{j} (-1)^j \log(u+j)$$

$$\log \Gamma(u) = \sum_{n=0}^{\infty} \frac{1}{n+1} \sum_{k=0}^{n} \binom{n}{k} (-1)^k (u+k) \log(u+k) + \frac{1}{2} - u + \frac{1}{2} \log(2\pi)$$

$$\sum_{k=1}^{\infty} \left[ \frac{(2k)!}{2^{2k}(k!)^2} \right]^3 \sum_{j=0}^{k-1} \frac{1}{2j+1} = \frac{\pi}{24} \frac{6\pi - \sqrt{2}\,\Gamma(1/4)}{\Gamma^4(3/4)}$$

$$2\log G(1+u) = -\sum_{n=0}^{\infty} \frac{1}{n+1} \sum_{k=0}^{n} \binom{n}{k} (-1)^k (u+k)^2 \log(u+k) + 2u \log \Gamma(u) + 2u^2 - u$$

$$+ 2\varsigma'(-1) - \frac{11}{12} - \frac{1}{2}(u-1)(3u+2)$$

$$\frac{1}{2} \int_0^x u \log \sin \pi u \, du = [\varsigma'(-2,x) + \varsigma'(-2,1-x)] - 2x[\varsigma'(-1,x) - \varsigma'(-1,1-x)] + \frac{\varsigma(3)}{2\pi^2}$$

$$\int_0^1 \log \Gamma_3(x+1) dx = -\frac{1}{24} \log(2\pi) + \frac{3\varsigma(3)}{8\pi^2}$$

where $\varsigma'(s,x) = \dfrac{\partial}{\partial s} \varsigma(s,x)$ and $G(x)$ and $\Gamma_3(x)$ are the Barnes double and triple gamma functions.

Whilst this paper is mainly expository, some of the formulae reported in it are believed to be new, and the paper may also be of interest specifically due to the fact that most of the various identities have been derived by elementary methods.

**CONTENTS OF VOLUMES I TO VI:**               **Volume/page**

**SECTION:**

**APPENDICES (Volume VI):**

**A**. Some properties of the Bernoulli numbers and the Bernoulli polynomials

**B**. A well-known integral

**C**. Euler's reflection formula for the gamma function and related matters

**D**. A very elementary proof of $\dfrac{\pi^2}{8} = \sum\limits_{n=0}^{\infty} \dfrac{1}{(2n+1)^2}$

**E**. Some aspects of Euler's constant $\gamma$ and the gamma function

**F**. Elementary aspects of Riemann's functional equation for the zeta function

**ACKNOWLEDGEMENTS**

**REFERENCES**

## 4. SOME ASPECTS OF THE GENERALISED HARMONIC NUMBERS

### ELEMENTARY PROOFS OF THE FLAJOLET AND SEDGEWICK IDENTITIES

In this part we derive alternative proofs of the Flajolet and Sedgewick identities (3.16a), (3.16b) and (3.16c) contained in Volume I.

First of all, we consider the following integral (using L'Hôpital's rule we note that the integrand remains finite as $x \to 0$)

(4.1.1)
$$J = \int_0^t \frac{1 - (1-x)^n}{x} \, dx$$

Using the binomial expansion we have

(4.1.2)
$$(1-x)^n = \sum_{k=0}^n \binom{n}{k} (-1)^k x^k$$

and $J$ therefore becomes

(4.1.3)
$$J = \int_0^t \sum_{k=1}^n \binom{n}{k} (-1)^{k+1} x^{k-1} \, dx = \sum_{k=1}^n \binom{n}{k} (-1)^{k+1} \frac{t^k}{k}$$

Alternatively, using the substitution $u = 1 - x$ in (4.1.1), we get

(4.1.4)
$$J = \int_{1-t}^1 \frac{1 - u^n}{1 - u} \, du = \int_{1-t}^1 \sum_{k=1}^n u^{k-1} \, du = \sum_{k=1}^n \int_{1-t}^1 u^{k-1} \, du$$

(having used the geometric series in the second part).

(4.1.5)
$$= \sum_{k=1}^n \frac{1}{k} - \sum_{k=1}^n \frac{(1-t)^k}{k} = \sum_{k=1}^n \frac{1 - (1-t)^k}{k}$$

Equating (4.1.3) and (4.1.5), we obtain

(4.1.6)
$$J = \int_0^t \frac{1 - (1-x)^n}{x} \, dx = \sum_{k=1}^n \binom{n}{k} (-1)^{k+1} \frac{t^k}{k} = \sum_{k=1}^n \frac{1 - (1-t)^k}{k}$$

and if we let $t = 1$, we have

$$(4.1.7) \qquad \sum_{k=1}^{n} \binom{n}{k} \frac{(-1)^{k+1}}{k} = \sum_{k=1}^{n} \frac{1}{k} = H_n = H_n^{(1)} \text{ (by definition)}$$

This therefore proves the Flajolet and Sedgewick identity (3.16a): see also the collection of formulae at (4.4.155a) et seq.

We previously used the identity (4.1.7) in (3.5), where it was mentioned that the formula was due to Euler. An alternative proof of (4.1.7) was given by Buschman [37] in 1958 using the digamma function $\psi(x) = \dfrac{d}{dx} \log \Gamma(x) = \dfrac{\Gamma'(x)}{\Gamma(x)}$, where $\psi(x)$ satisfies the relationships:

$$(4.1.7a) \qquad \psi(n+1) - \psi(1) = H_n$$

$$(4.1.7b) \qquad \psi(x+a) = \psi(a) + \sum_{k=1}^{\infty} \frac{(-1)^{k+1}}{k} \frac{x(x-1)...(x-k+1)}{a(a+1)...(a+k-1)}$$

and then substituting $x = n$ and $a = 1$. Formula (4.1.7b) is contained in Higher Transcendental Functions by Erdélyi et al. [58b] and is convergent for $\text{Re}(x+a) > 0$. The digamma function is employed later in (4.3.32). Reference should also be made to Appendix E of Volume VI for more details of the digamma function.

We now divide (4.1.6) by $t$ and integrate over the interval $[0, x]$ to obtain

$$(4.1.8) \qquad \sum_{k=1}^{n} \binom{n}{k} (-1)^{k+1} \frac{x^k}{k^2} = \sum_{k=1}^{n} \frac{1}{k} \int_{0}^{x} \frac{1-(1-t)^k}{t} dt$$

$$(4.1.9) \qquad = \sum_{k=1}^{n} \frac{1}{k} \left\{ \sum_{l=1}^{k} \binom{k}{l} (-1)^{l+1} \frac{x^l}{l} \right\} \qquad \text{, using (4.1.3).}$$

$$(4.1.10) \qquad = \sum_{k=1}^{n} \frac{1}{k} \left\{ \sum_{l=1}^{k} \frac{1-(1-x)^l}{l} \right\} \qquad \text{, using (4.1.6).}$$

Letting $x = 1$ we obtain

$$(4.1.11) \qquad \sum_{k=1}^{n} \binom{n}{k} \frac{(-1)^{k+1}}{k^2} = \sum_{k=1}^{n} \frac{1}{k} \left\{ \sum_{l=1}^{k} \binom{k}{l} \frac{(-1)^{l+1}}{l} \right\} \qquad \text{, using (4.1.9).}$$

$$(4.1.12) \qquad = \sum_{k=1}^{n} \frac{1}{k} \left\{ \sum_{l=1}^{k} \frac{1}{l} \right\} \qquad \text{, using (4.1.10).}$$

Using (4.1.12) we can write

$$(4.1.13) \qquad \sum_{k=1}^{n} \binom{n}{k} \frac{(-1)^{k+1}}{k^2} = \sum_{k=1}^{n} \frac{H_k}{k}$$

From the formula (3.17) given by Adamchik we have

$$(4.1.14) \qquad \sum_{k=1}^{n} \frac{H_k}{k} = \frac{1}{2}\left(H_n^{(1)}\right)^2 + \frac{1}{2} H_n^{(2)}$$

and hence we have an alternative derivation of Flajolet and Sedgewick's formula (3.16b).

$$(4.1.15) \qquad S_n(2) = \sum_{k=1}^{n} \binom{n}{k} \frac{(-1)^{k+1}}{k^2} = \frac{1}{2}\left(H_n^{(1)}\right)^2 + \frac{1}{2} H_n^{(2)}$$

As a matter of interest, I also found formula (4.1.14) reported by M.E. Levenson in a 1938 volume of The American Mathematical Monthly [98] in a problem concerning the evaluation of

$$(4.1.15a) \qquad \Gamma''(1) = \int_0^\infty e^{-x} \log^2 x \, dx = \gamma^2 + \varsigma(2)$$

where $\gamma$ is Euler's constant defined by

$$(4.1.15b) \qquad \gamma = \lim_{n \to \infty}\left(H_n - \log n\right)$$

Some 20 years later Levenson [99] published another article in The American Mathematical Monthly entitled "A recursion formula for $\int_0^\infty e^{-t} \log^{n+1} t \, dt$", and this provides us with identities for $\Gamma^{(n)}(1)$. A different proof of (4.1.15a) is given in Appendix E of Volume VI.

We can repeat the exercise again by dividing (4.1.10) by $x$ and integrating over the interval $[0, t]$. This results in

$$(4.1.16) \qquad \sum_{k=1}^{n} \binom{n}{k}(-1)^{k+1} \frac{t^k}{k^3} = \sum_{k=1}^{n} \frac{1}{k} \sum_{l=1}^{k} \frac{1}{l} \int_0^t \frac{1-(1-x)^l}{x} \, dx$$

$$(4.1.17) \qquad = \sum_{k=1}^{n} \frac{1}{k} \sum_{l=1}^{k} \frac{1}{l} \left\{ \sum_{m=1}^{l} \frac{1-(1-t)^m}{m} \right\}$$

Letting $t = 1$, (4.1.17) becomes

(4.1.18)  $$\sum_{k=1}^{n}\binom{n}{k}\frac{(-1)^{k+1}}{k^3} = \sum_{k=1}^{n}\frac{1}{k}\sum_{l=1}^{k}\frac{1}{l}\sum_{m=1}^{l}\frac{1}{m} = \sum_{k=1}^{n}\frac{1}{k}\sum_{l=1}^{k}\frac{H_l}{l}$$

Quite obviously the method may be extended to produce identities for $\sum_{k=1}^{n}\binom{n}{k}\frac{(-1)^{k+1}}{k^s}$

where $s$ is an integer.

The formula (4.1.18) is not new: the following generalised identity is referred to as Dilcher's formula in the Mathworld website [134] and also see the papers by Dilcher [54], Hernández [79] and Prodinger [109].

(4.1.18a)  $$\sum_{k=1}^{n}\binom{n}{k}\frac{(-1)^{k+1}}{k^s} = \sum_{1 \le i_1 \le i_2 \le \ldots \le i_s \le n}\frac{1}{i_1 i_2 \ldots i_s}$$

Reference should also be made to (E.60) in Volume VI for the following identity given by Olds [94aa] in 1936

(4.1.18b)  $$n\sum_{k=1}^{n-1}\binom{n-1}{k}\frac{(-1)^k}{(k+1)^m} = \sum_{i_1=1}^{n}\frac{1}{i_1}\sum_{i_2=1}^{i_1}\frac{1}{i_2}\sum_{i_3=1}^{i_2}\frac{1}{i_3}\ldots\sum_{i_{m-1}=1}^{i_{m-2}}\frac{1}{i_{m-1}}$$

Dividing (4.1.9) by $x$, and integrating as before, we obtain

(4.1.19)  $$\sum_{k=1}^{n}\binom{n}{k}\frac{(-1)^{k+1}}{k^3} = \sum_{k=1}^{n}\frac{1}{k}\left\{\sum_{l=1}^{k}\binom{k}{l}\frac{(-1)^{l+1}}{l^2}\right\}$$

and, using (4.1.13) to substitute for the sum in parentheses, we get

(4.1.20)  $$\sum_{k=1}^{n}\binom{n}{k}\frac{(-1)^{k+1}}{k^3} = \sum_{k=1}^{n}\frac{1}{k}\sum_{l=1}^{k}\frac{H_l}{l}$$

(4.1.21)  $$= \frac{1}{2}\left\{\sum_{k=1}^{n}\frac{\left(H_k^{(1)}\right)^2}{k} + \sum_{k=1}^{n}\frac{H_k^{(2)}}{k}\right\}$$

The formula (4.1.21) was obtained by substituting the equivalent of (4.1.14) in (4.1.20). Formula (4.1.21) is the same as Adamchik's result (3.20c) and hence we have an alternative proof of the Flajolet and Sedgewick equation (3.16c) for $S_n(3)$ (which is rather good news because the author has still not familiarised himself with Mellin transforms, Rice/Nörlund integrals and Bell polynomials!).

We therefore have the following useful identities:

(4.1.22a)
$$\sum_{k=1}^{n}\binom{n}{k}\frac{(-1)^{k+1}}{k} = \sum_{k=1}^{n}\frac{1}{k} = H_n = H_n^{(1)}$$

(4.1.22b)
$$\sum_{k=1}^{n}\binom{n}{k}\frac{(-1)^{k+1}}{k^2} = \sum_{k=1}^{n}\frac{1}{k}\left\{\sum_{l=1}^{k}\binom{k}{l}\frac{(-1)^{l+1}}{l}\right\}$$

$$= \sum_{k=1}^{n}\frac{H_k^{(1)}}{k}$$

$$= \frac{1}{2}\left(H_n^{(1)}\right)^2 + \frac{1}{2}H_n^{(2)}$$

$$= \sum_{k=1}^{n}\frac{1}{k}\left\{\sum_{l=1}^{k}\frac{1}{l}\right\}$$

(4.1.22c)
$$\sum_{k=1}^{n}\binom{n}{k}\frac{(-1)^{k+1}}{k^3} = \sum_{k=1}^{n}\frac{1}{k}\left\{\sum_{l=1}^{k}\binom{k}{l}\frac{(-1)^{l+1}}{l^2}\right\}$$

$$= \sum_{k=1}^{n}\frac{1}{k}\sum_{l=1}^{k}\frac{H_l^{(1)}}{l}$$

$$= \frac{1}{2}\left\{\sum_{k=1}^{n}\frac{\left(H_k^{(1)}\right)^2}{k} + \sum_{k=1}^{n}\frac{H_k^{(2)}}{k}\right\}$$

$$= \sum_{k=1}^{n}\frac{1}{k}\sum_{l=1}^{k}\frac{1}{l}\sum_{m=1}^{l}\frac{1}{m}$$

We also note that Spieß [123bi] has derived the following identities

(4.1.22d)
$$\sum_{k=1}^{n}\frac{1}{k}\frac{1}{n-k+1} = \frac{2}{n+1}H_n^{(1)}$$

$$\sum_{k=2}^{n}\frac{2}{k}\frac{1}{n-k+1}H_{k-1}^{(1)} = \frac{3}{n+1}\left[\left(H_n^{(1)}\right)^2 - H_n^{(2)}\right]$$

$$\sum_{k=2}^{n} \frac{2}{k} \frac{2}{n-k+1} H_{k-1}^{(1)} H_{n-k}^{(1)} = \frac{4}{n+1} \left[ \left( H_n^{(1)} \right)^3 - 3 H_n^{(1)} H_n^{(2)} + 2 H_n^{(3)} \right]$$

Note that the coefficients in the final identity are the same as (4.3.26).

Some additional identities involving the harmonic numbers are considered below.

Let us revisit equation (4.1.6)

$$\sum_{k=1}^{n} \binom{n}{k} (-1)^{k+1} \frac{t^k}{k} = \sum_{k=1}^{n} \frac{1-(1-t)^k}{k}$$

and if we let $t \to 1-t$, we have

$$\sum_{k=1}^{n} \binom{n}{k} (-1)^{k+1} \frac{(1-t)^k}{k} = \sum_{k=1}^{n} \frac{1}{k} - \sum_{k=1}^{n} \frac{t^k}{k}$$

With $t=1$ we have as before

$$\sum_{k=1}^{n} \binom{n}{k} \frac{(-1)^{k+1}}{k} = \sum_{k=1}^{n} \frac{1}{k} = H_n = H_n^{(1)}$$

Combining these we obtain

$$\sum_{k=1}^{n} \binom{n}{k} (-1)^{k+1} \frac{(1-t)^k}{k} - \sum_{k=1}^{n} \binom{n}{k} \frac{(-1)^{k+1}}{k} = -\sum_{k=1}^{n} \frac{t^k}{k}$$

and this can be written as

$$\sum_{k=1}^{n} \binom{n}{k} (-1)^{k+1} \frac{(1-t)^k - 1}{k} = -\sum_{k=1}^{n} \frac{t^k}{k}$$

Dividing by $t$ we get

$$\sum_{k=1}^{n} \binom{n}{k} (-1)^{k+1} \frac{(1-t)^k - 1}{kt} = -\sum_{k=1}^{n} \frac{t^{k-1}}{k}$$

The binomial theorem gives us

(4.1.23)      $$\frac{(1-t)^k - 1}{t} = \sum_{j=1}^{k} (1-t)^{j-1}$$

and hence we have

$$\sum_{k=1}^{n}\binom{n}{k}(-1)^{k+1}\frac{1}{k}\sum_{j=1}^{k}(1-t)^{j-1}=-\sum_{k=1}^{n}\frac{t^{k-1}}{k}$$

A simple integration gives us

$$\int_0^x\sum_{k=1}^{n}\binom{n}{k}(-1)^{k+1}\frac{1}{k}\sum_{j=1}^{k}(1-t)^{j-1}dt=-\sum_{k=1}^{n}\frac{x^k}{k^2}$$

and therefore we have

(4.1.23a) $$\sum_{k=1}^{n}\binom{n}{k}(-1)^{k+1}\frac{1}{k}\sum_{j=1}^{k}\frac{(1-x)^j-1}{j}=-\sum_{k=1}^{n}\frac{x^k}{k^2}$$

We now divide by $x$ and obtain

$$\sum_{k=1}^{n}\binom{n}{k}(-1)^{k+1}\frac{1}{k}\sum_{j=1}^{k}\frac{(1-x)^j-1}{jx}=-\sum_{k=1}^{n}\frac{x^{k-1}}{k^2}$$

and after substituting (4.1.23) we have

$$\sum_{k=1}^{n}\binom{n}{k}(-1)^{k+1}\frac{1}{k}\sum_{j=1}^{k}\frac{1}{j}\sum_{i=1}^{j}(1-x)^{i-1}=-\sum_{k=1}^{n}\frac{x^{k-1}}{k^2}$$

With $x=1$ in (4.1.23a) we obtain

(4.1.23b) $$\sum_{k=1}^{n}\binom{n}{k}(-1)^{k+1}\frac{H_k}{k}=H_n^{(2)}$$

Dividing by $t$ and integrating in the range $[1,t]$ we obtain

(4.1.24) $$\sum_{k=1}^{n}\binom{n}{k}(-1)^{k+1}\frac{1}{k}\sum_{j=1}^{k}\frac{1}{j}\sum_{i=1}^{j}\frac{(1-t)^i}{i}=-\sum_{k=1}^{n}\frac{t^k}{k^3}+H_n^{(3)}$$

Letting $t=0$ we have a known result

(4.1.25) $$\sum_{k=1}^{n}\binom{n}{k}(-1)^{k+1}\frac{1}{k}\sum_{j=1}^{k}\frac{H_j^{(1)}}{j}=H_n^{(3)}$$

Four other very different proofs of (4.1.25) are contained in [17a] (and (4.1.24) corrects the minor misprint contained in the solution provided by O.P. Lossers).

In one of the solutions to (4.1.25) in [17a], A.N.'t Woord employs the relation

(4.1.26)
$$b_n = \sum_{k=1}^{n} \binom{n}{k} a_k \;\Rightarrow \sum_{k=1}^{n} \frac{b_k}{k} = \sum_{k=1}^{n} \binom{n}{k} \frac{a_k}{k}$$

which is easily proved by mathematical induction. It is clear that (4.1.26) is true for $n = 1$ and for $n > 1$ we have

$$\sum_{k=1}^{n} \binom{n}{k} \frac{a_k}{k} = \sum_{k=1}^{n} \binom{n-1}{k} \frac{a_k}{k} + \sum_{k=1}^{n} \binom{n-1}{k-1} \frac{a_k}{k}$$

$$= \sum_{k=1}^{n-1} \frac{b_k}{k} + \sum_{k=1}^{n} \binom{n}{k} \frac{a_k}{n} = \sum_{k=1}^{n} \frac{b_k}{k}$$

Applying this identity to (3.16a) we have

$$\sum_{k=1}^{n} \binom{n}{k} \frac{(-1)^{k+1}}{k} = H_n^{(1)}$$

and therefore

$$\sum_{k=1}^{n} \frac{H_k^{(1)}}{k} = \sum_{k=1}^{n} \binom{n}{k} \frac{(-1)^{k+1}}{k^2}$$

(as was obtained previously in (4.1.12)). Repeating the process we may easily derive the identity in (4.1.18)

$$\sum_{k=1}^{n} \frac{1}{k} \; \sum_{l=1}^{k} \frac{1}{l} \; \sum_{m=1}^{l} \frac{1}{m} = \sum_{k=1}^{n} \binom{n}{k} \frac{(-1)^{k+1}}{k^3}$$

Using (4.1.23b) and (4.1.26) we see that

(4.1.27)
$$\sum_{k=1}^{n} \frac{H_k^{(2)}}{k} = \sum_{k=1}^{n} \binom{n}{k} (-1)^{k+1} \frac{H_k^{(1)}}{k^2}$$

We now refer back to (4.1.25)

$$\sum_{k=1}^{n}\binom{n}{k}(-1)^{k+1}\frac{1}{k}\sum_{j=1}^{k}\frac{H_j^{(1)}}{j}=H_n^{(3)}$$

and note that

$$\sum_{j=1}^{k}\frac{H_j^{(1)}}{j}=-\sum_{j=1}^{k}\int_0^1(1-x)^{j-1}\log x\,dx$$

$$=-\int_0^1\frac{1-(1-x)^k}{x}\log x\,dx$$

We then have

$$-\sum_{k=1}^{n}\binom{n}{k}(-1)^{k+1}\frac{1}{k}\int_0^1\frac{1-(1-x)^k}{x}\log x\,dx=H_n^{(3)}$$

$$\int_0^1\sum_{k=1}^{n}\binom{n}{k}(-1)^{k+1}\frac{1}{k}\frac{(1-x)^k-1}{x}\log x\,dx=H_n^{(3)}$$

and, using (4.1.6), this may be written as

$$\int_0^1\left[\sum_{k=1}^{n}\frac{1-(x)^k}{k}-H_n^{(1)}\right]\frac{\log x}{x}\,dx=H_n^{(3)}$$

We then have

$$\int_0^1\sum_{k=1}^{n}\frac{x^k}{k}\frac{\log x}{x}\,dx=H_n^{(3)}$$

$$\int_0^1\sum_{k=1}^{n}\int_0^x y^{k-1}dy\,\frac{\log x}{x}\,dx=H_n^{(3)}$$

$$\int_0^1 dx\int_0^x\frac{1-y^n}{1-y}\frac{\log x}{x}\,dy=H_n^{(3)}$$

and summation gives us

$$\int_0^1 dx\int_0^x\frac{\varsigma(s)-Li_s(y)}{1-y}\frac{\log x}{x}\,dy=\sum_{n=1}^{\infty}\frac{H_n^{(3)}}{n^s}$$

We note that the Wolfram Integrator is able to evaluate $\int\limits_0^x \dfrac{\varsigma(s) - Li_s(y)}{1-y} dy$ for $s = 2, 3$

and for $s = 3$ we have

$$\int\limits_0^x \frac{\varsigma(3) - Li_3(y)}{1-y} dy = \frac{1}{2}\big[Li_2(x)\big]^2 + \log(1-x)[Li_3(x) - \varsigma(3)]$$

This then gives us

(4.1.28)    $\int\limits_0^1 \left(\dfrac{1}{2}\big[Li_2(x)\big]^2 + \log(1-x)[Li_3(x) - \varsigma(3)]\right)\dfrac{\log x}{x} dx = \sum\limits_{n=1}^{\infty} \dfrac{H_n^{(3)}}{n^3}$

SOME IDENTITIES DERIVED FROM THE HASSE/SONDOW EQUATIONS

In this Part we derive identities for $\sum\limits_{k=0}^{n} \binom{n}{k} \dfrac{(-1)^k}{(k+x)^s}$ , noting that similar terms appear in

the Hasse/Sondow formula given in (3.11)

$$\varsigma_a(s) = \sum_{n=0}^{\infty} \frac{1}{2^{n+1}} \sum_{k=0}^{n} \binom{n}{k} \frac{(-1)^k}{(k+1)^s}$$

We know from [75, p.188] that

(4.2.1)    $\sum\limits_{k=0}^{n} \binom{n}{k} \dfrac{(-1)^k}{k+x} = \dfrac{n!}{x(1+x)...(n+x)} = g(x)$

(4.2.2)    $= x^{-1}\binom{n+x}{n}^{-1}$    $, x \notin \{0, -1, ..., -n\}$

(An alternative proof of (4.2.1) is contained in (4.4.1). I subsequently discovered that yet another proof, using partial fractions, is contained in Further Mathematics [108, p.43]: I must have read, and since forgotten, that proof when I was 17!). The connection between $g(x)$ and the gamma function $\Gamma(x)$ is referred to in more detail in Appendix E of Volume VI.

In passing, we note that an application of (4.1.26) to (4.2.1) results in

$$\sum_{k=1}^{n}\left(\frac{k!}{kx(1+x)...(k+x)}-\frac{1}{kx}\right)=\sum_{k=1}^{n}\binom{n}{k}\frac{(-1)^k}{k(k+x)}$$

and, upon multiplying across by $x$, we see that

$$\sum_{k=1}^{n}\left(\frac{k!}{k(1+x)...(k+x)}\right)-H_n^{(1)}=\sum_{k=1}^{n}\binom{n}{k}\frac{(-1)^k}{k}-\sum_{k=1}^{n}\binom{n}{k}\frac{(-1)^k}{k+x}=-H_n^{(1)}-\frac{n!}{x(1+x)...(n+x)}+\frac{1}{x}$$

We therefore obtain

(4.2.2a)     $$\sum_{k=1}^{n}\left(\frac{k!}{k(1+x)...(k+x)}\right)=\frac{1}{x}-\frac{n!}{x(1+x)...(n+x)}$$

which may be contrasted with Sondow's result (4.4.92) in Volume III

$$\gamma=\int_{1}^{\infty}\sum_{n=1}^{\infty}\frac{n!}{(n+1)x(x+1)...(x+n)}\,dx$$

Equation (4.2.2a) is a particular case of the identity employed by Apéry [131a] in the proof of the irrationality of $\varsigma(3)$

$$\sum_{k=1}^{n}\frac{a_1 a_2 ... a_{k-1}}{(x+a_1)(x+a_2)...(x+a_k)}=\frac{1}{x}$$

Now, returning to the main story, with $x=1$ in (4.2.1) we have

(4.2.3)     $$g(1)=\sum_{k=0}^{n}\binom{n}{k}\frac{(-1)^k}{k+1}=\frac{1}{(n+1)}$$

and it is also seen for example that

$$g(2)=\sum_{k=0}^{n}\binom{n}{k}\frac{(-1)^k}{k+2}=\frac{1}{(n+1)(n+2)}$$

$$g(3)=\sum_{k=0}^{n}\binom{n}{k}\frac{(-1)^k}{k+3}=\frac{2!}{(n+1)(n+2)(n+3)}$$

$$g(4)=\sum_{k=0}^{n}\binom{n}{k}\frac{(-1)^k}{k+4}=\frac{3!}{(n+1)(n+2)(n+3)(n+4)}$$

Using (4.2.3), and letting $s = 1$ in the Hasse/Sondow identity (3.11), results in the simple expression

(4.2.4) $$\varsigma_a(1) = \sum_{n=0}^{\infty} \frac{1}{2^{n+1}} \frac{1}{n+1} = \sum_{n=1}^{\infty} \frac{1}{n 2^n} = \log 2$$

and, using the definition of $\varsigma_a(s)$ in (1.1), this therefore implies that

(4.2.5) $$\lim_{s \to 1} \left[ \left( 1 - 2^{1-s} \right) \varsigma(s) \right] = \log 2.$$

An alternative proof of (4.2.5) is contained in [122]. This limit is usefully employed in Volume II(b).

We have from (4.2.1)

(4.2.6) $$\log g(x) = \log n! - \sum_{k=0}^{n} \log(k+x)$$

Hence, by taking the logarithmic derivative of (4.2.6) we have

(4.2.7) $$\frac{g'(x)}{g(x)} = -\left\{ \frac{1}{x} + \frac{1}{1+x} \cdots + \frac{1}{n+x} \right\} = -h(x)$$

From (4.2.7) we have

(4.2.8) $$h^{(p)}(x) = (-1)^p \, p! \sum_{k=0}^{n} \frac{1}{(k+x)^{p+1}}$$

and specifically

(4.2.9) $$h^{(p)}(1) = (-1)^p \, p! \, H_{n+1}^{(p+1)}$$

We therefore have from (4.2.7)

(4.2.10) $$\frac{g'(1)}{g(1)} = -h(1) = -H_{n+1}$$

and accordingly

(4.2.11) $$g'(1) = -\frac{H_{n+1}}{n+1}$$

The derivative of the left-hand side of (4.2.1) is elementary and, by differentiating it $p$ times, we obtain

$$(4.2.12a) \qquad g^{(p)}(x) = (-1)^p\, p! \sum_{k=0}^{n} \binom{n}{k} \frac{(-1)^k}{(k+x)^{p+1}}$$

$$(4.2.13b) \qquad g^{(p)}(1) = (-1)^p\, p! \sum_{k=0}^{n} \binom{n}{k} \frac{(-1)^k}{(k+1)^{p+1}}$$

Therefore we have

$$(4.2.14) \qquad \sum_{k=0}^{n} \binom{n}{k} \frac{(-1)^k}{(k+1)^{p+1}} = \frac{(-1)^p\, g^{(p)}(1)}{p!}$$

In particular

$$(4.2.15a) \qquad g'(x) = -\sum_{k=0}^{n} \binom{n}{k} \frac{(-1)^k}{(x+k)^2}$$

and

$$(4.2.15b) \qquad g'(1) = -\sum_{k=0}^{n} \binom{n}{k} \frac{(-1)^k}{(1+k)^2}$$

Equating (4.2.12a) and (4.2.15b) gives us the identity

$$(4.2.16) \qquad \sum_{k=0}^{n} \binom{n}{k} \frac{(-1)^k}{(1+k)^2} = \frac{H_{n+1}}{n+1} = -g'(1)$$

Substituting the above in the Hasse/Sondow identity (3.11)

$$(4.2.17) \qquad \varsigma_a(s) = \sum_{n=0}^{\infty} \frac{1}{2^{n+1}} \sum_{k=0}^{n} \binom{n}{k} \frac{(-1)^k}{(k+1)^s}$$

$$= \sum_{n=0}^{\infty} \frac{1}{2^{n+1}} \frac{(-1)^{s-1} g^{(s-1)}(1)}{(s-1)!}$$

and this gives us for $s = 2$

$$(4.2.18a) \qquad \varsigma_a(2) = \sum_{n=0}^{\infty} \frac{1}{2^{n+1}} \frac{H_{n+1}}{n+1}$$

(4.2.18b)
$$= \sum_{n=1}^{\infty} \frac{H_n}{n2^n}$$

Therefore we have

$$\frac{\pi^2}{12} = \sum_{n=1}^{\infty} \frac{H_n}{n2^n}$$

This is in agreement with (3.31) where we showed that

(4.2.19)
$$\frac{1}{2}\log^2(1-x) + Li_2(x) = \sum_{n=1}^{\infty} \frac{H_n}{n}x^n$$

Letting $x = 1/2$ we have

(4.2.20)
$$\frac{1}{2}\log^2 2 + Li_2(1/2) = \sum_{n=1}^{\infty} \frac{H_n}{n2^n}$$

and therefore, using Euler's identity (3.43a) for $Li_2(1/2)$, we again deduce

(4.2.21)
$$\frac{\pi^2}{12} = \sum_{n=1}^{\infty} \frac{H_n}{n2^n}$$

which may be compared with the formula for $\log 2$ in (3.7)

$$2\log 2 = \sum_{n=1}^{\infty} \frac{H_n}{2^n}$$

The result (4.2.21) is well-known: see for example [46].

Using (3.16a) and (4.2.18b) we have

(4.2.21a)
$$\varsigma_a(2) = \sum_{n=0}^{\infty} \frac{1}{n2^n} \sum_{k=1}^{n} \binom{n}{k} \frac{(-1)^{k+1}}{k}$$

Employing the second Hasse identity (3.12),

(4.2.22)
$$\varsigma(s) = \frac{1}{s-1} \sum_{n=0}^{\infty} \frac{1}{n+1} \sum_{k=0}^{n} \binom{n}{k} \frac{(-1)^k}{(k+1)^{s-1}}$$

$$= \frac{1}{s-1} \sum_{n=0}^{\infty} \frac{1}{n+1} \frac{(-1)^{s-2} g^{(s-2)}(1)}{(s-2)!}$$

we obtain for $s = 2$

(4.2.23)
$$\varsigma(2) = \sum_{n=0}^{\infty} \frac{1}{n+1} \frac{1}{n+1} = \sum_{n=1}^{\infty} \frac{1}{n^2}$$

and this at least proves that my arithmetic was correct!!

Let us now consider (4.2.8) again: we have

(4.2.24)
$$\frac{g'(x)}{g(x)} = -h(x)$$

and differentiating (4.2.24) gives

(4.2.25)
$$\frac{gg'' - (g')^2}{g^2} = -h'$$

The left-hand side of (4.2.25) may be written as

(4.2.26)
$$= \frac{g''}{g} - \left\{ \frac{g'}{g} \right\}^2$$

$$= \frac{g''}{g} - h^2$$

Therefore we have

(4.2.27)
$$\frac{g''}{g} = h^2 - h'$$

and hence

(4.2.28)
$$g''(1) = \frac{1}{n+1} \left\{ (H_{n+1}^{(1)})^2 + H_{n+1}^{(2)} \right\}$$

Accordingly we have

(4.2.29)
$$g''(1) = \frac{1}{n+1} \left\{ (H_{n+1}^{(1)})^2 + H_{n+1}^{(2)} \right\} = 2 \sum_{k=0}^{n} \binom{n}{k} \frac{(-1)^k}{(1+k)^3}$$

Using the Hasse/Sondow identity (3.11) we obtain

$$(4.2.30a) \qquad \varsigma_a(3) = \frac{1}{2}\sum_{n=0}^{\infty}\frac{1}{2^{n+1}}\,\frac{1}{n+1}\left\{(H_{n+1}^{(1)})^2 + H_{n+1}^{(2)}\right\}$$

which is equivalent to

$$(4.2.30b) \qquad \varsigma_a(3) = \frac{1}{2}\sum_{n=1}^{\infty}\frac{1}{n2^n}\left\{(H_n^{(1)})^2 + H_n^{(2)}\right\}$$

An independent proof of this is set out below. From (4.2.30b) we have

$$(4.2.31) \qquad S = \frac{1}{2}\sum_{n=1}^{\infty}\frac{1}{2^n}\,\frac{(H_n^{(1)})^2}{n} + \frac{1}{2}\sum_{n=1}^{\infty}\frac{1}{2^n}\,\frac{H_n^{(2)}}{n}$$

$$= \frac{1}{4}D + \frac{1}{4}E$$

$$(4.2.32) \qquad = \frac{1}{4}\left\{\frac{7}{4}\varsigma(3)\right\} + \frac{1}{4}\left\{\frac{5}{4}\varsigma(3)\right\} = \frac{3}{4}\varsigma(3)$$

where we have employed the formulae for $D$ and $E$ given in (3.22) and (3.49) respectively. We may also derive (4.2.30b) by combining (3.106b) and (3.106b) from Volume I.

The second Hasse identity (3.12) gives us

$$(4.2.33) \qquad \varsigma(3) = \frac{1}{2}\sum_{n=0}^{\infty}\frac{H_{n+1}}{(n+1)^2} = \frac{1}{2}\sum_{n=1}^{\infty}\frac{H_n}{n^2}$$

The above formula is not new: Shen [120] derived the following result in 1995 (see also formula (3.105e) and [126, p.253] and formula (3.105b) in this paper)

$$(4.2.34) \qquad \varsigma(3) = \sum_{n=1}^{\infty}\frac{H_{n-1}}{n^2}$$

and using the substitution $H_{n-1} = H_n - \frac{1}{n}$ it is clear that (4.2.33) and (4.2.34) are equivalent (defining $H_0 = 0$). This result was also obtained by Briggs et al [35] in 1955 and previously by several others, including Euler and Ramanujan [21, Part I, p.252].

Using (3.16a) to substitute for $H_n$, equation (4.2.33) may be written as

(4.2.34a)
$$\varsigma(3) = \frac{1}{2} \sum_{n=1}^{\infty} \frac{1}{n^2} \sum_{k=1}^{n} \binom{n}{k} \frac{(-1)^{k+1}}{k}$$

Similarly, using (3.16b) we may write (4.2.30b) as

(4.2.34b)
$$\varsigma_a(3) = \sum_{n=0}^{\infty} \frac{1}{n2^n} \sum_{k=1}^{n} \binom{n}{k} \frac{(-1)^{k+1}}{k^2}$$

Continuing in the same manner as before we obtain from (4.2.27)

(4.2.35)
$$\frac{gg''' - g''g'}{g^2} = 2hh' - h''$$

The left hand side of (4.2.35) may be expressed as

$$= \frac{g'''}{g} - \frac{g''}{g} \frac{g'}{g}$$

$$= \frac{g'''}{g} - (h^2 + h')(-h)$$

Therefore

$$\frac{g'''}{g} = -h^3 + 3h'h - h''$$

and

(4.2.36a)
$$\frac{g'''(1)}{g(1)} = -\left(H_{n+1}^{(1)}\right)^3 - 3H_{n+1}^{(1)} H_{n+1}^{(2)} - 2H_{n+1}^{(3)}$$

(4.2.36b)
$$g'''(1) = -3! \sum_{k=0}^{n} \binom{n}{k} \frac{(-1)^k}{(1+k)^4} = -\frac{1}{n+1} \left\{ \left(H_{n+1}^{(1)}\right)^3 + 3H_{n+1}^{(1)} H_{n+1}^{(2)} + 2H_{n+1}^{(3)} \right\}$$

The Hasse/Sondow identity (3.11) then gives us

(4.2.37)
$$\varsigma_a(4) = \frac{1}{3!} \sum_{n=0}^{\infty} \frac{1}{2^{n+1}} \frac{1}{n+1} \left\{ \left(H_n^{(1)}\right)^3 + 3H_n^{(1)} H_n^{(2)} + 2H_n^{(3)} \right\}$$

which may be written as

$$\varsigma_a(4) = \sum_{n=1}^{\infty} \frac{1}{n 2^n} \left\{ \frac{1}{6} \left( H_n^{(1)} \right)^3 + \frac{1}{2} H_n^{(1)} H_n^{(2)} + \frac{1}{3} H_n^{(3)} \right\}$$

We immediately note that the terms contained in parentheses are identical to those contained in the Flajolet and Sedgewick formula (3.16c) for $S_n(3)$ which was obtained using Bell polynomials (and also note the connection between (4.2.30) and $S_n(2)$ ). As explained in more detail in (4.4.177) in Volume IV, this similarity is not a mere coincidence.

Using (3.16c) the equation (4.2.37) becomes

(4.2.38) $$\varsigma_a(4) = \sum_{n=1}^{\infty} \frac{1}{n 2^n} \sum_{k=1}^{n} \binom{n}{k} \frac{(-1)^{k+1}}{k^3}$$

The second Hasse identity (3.12) gives us for $s = 4$

(4.2.39) $$\varsigma(4) = \frac{1}{6} \sum_{n=0}^{\infty} \frac{1}{(n+1)^2} \left\{ (H_{n+1}^{(1)})^2 + H_{n+1}^{(2)} \right\}$$

(4.2.40) $$= \frac{1}{6} \sum_{n=1}^{\infty} \frac{1}{n^2} \left\{ (H_n^{(1)})^2 + H_n^{(2)} \right\}$$

By using (3.16b) equation (4.2.40) can be written as

(4.2.41) $$\varsigma(4) = \frac{1}{3} \sum_{n=1}^{\infty} \frac{1}{n^2} \sum_{k=1}^{n} \binom{n}{k} \frac{(-1)^{k+1}}{k^2}$$

In their 1995 paper, "On an Intriguing Integral and Some Series Related to $\varsigma(4)$" [27], D. Borwein and J. M. Borwein relate that Enrico Au-Yeung (an undergraduate student in the Faculty of Mathematics in Waterloo), on the basis of numerical approximation, conjectured that:

(4.2.42) $$\sum_{n=1}^{\infty} \frac{\left( H_n^{(1)} \right)^2}{n^2} = \frac{17\pi^4}{360} = \frac{17}{4} \varsigma(4)$$

The Borweins then provide us with a most enjoyable tour through "diverse topics in classical analysis" and they end up with an analytical proof of (4.2.42). They also pointed out that a similar result had been obtained some four years earlier by P.J. De Doelder [55] but they had not previously realised that this work was of relevance because of the rather arcane title of his paper, "On some series containing $\psi(x) - \psi(y)$ and $(\psi(x) - \psi(y))^2$ for certain values of $x$ and $y$". (With respect, the title does indeed lack some transparency! Its relevance does however become more apparent after reading (4.3.32) et seq).

I then recalled that I had seen a reference to the infinite series (4.2.42) in a 1988 article by Castellanos, entitled "The Ubiquitous $\pi$". This paper [42, p.86] stated that "A series where the terms are the squares of the arithmetic means of the partial sums of the harmonic series was discovered by H.F. Sandham" and unfortunately it then incorrectly referenced the source to J.L.F. Bertrand's book, "Traité de Calcul Différentiel et de Calcul Intégral" (1864) . I thumbed through quite a lot of this useful treatise (on the internet [22]) before finally realising that it was not particularly germane to my quest. Eventually I discovered the real Sandham and found out that he received his PhD from Trinity College Dublin in 1958 and that the Castellanos reference should have been to a problem proposed and solved by Sandham in 1949 [118].

Sandham was latterly at the Institute for Advanced Studies in Dublin as was Erwin Schrödinger of quantum mechanics fame. In 1946, Sandham [117] presented a remarkably concise proof of the integral identity

$$(4.2.42a) \qquad\qquad \int\limits_0^\infty e^{-x^2}dx = \frac{\sqrt{\pi}}{2}$$

Whilst this is not immediately relevant to the present work, the proof is so elegant that I decided to include it in Appendix B to this paper (in the hope that it would become standard fare in undergraduate classes in the not too distant future). As it turns out, some time after I drafted the previous sentence, I subsequently found the need to employ (4.2.42a) several times in this series of papers.

It may be noted that another elegant proof of (4.2.42) was recently given by Freitas [69a] using a simple logarithmic double integral. Another different proof is contained at (4.4.156e).

Let us now return to the main story.

Using (4.2.40) and (4.2.42) enables us to show that

$$(4.2.43) \qquad\qquad \sum_{n=1}^\infty \frac{H_n^{(2)}}{n^2} = \frac{7\pi^4}{360} = \frac{7}{4}\varsigma(4)$$

Equation (4.2.43) is not new: see for example the 1997 paper by Flajolet and Salvy, "Euler Sums and Contour Integral Representations" [69, p.23] (this paper is considered further below).

Combining (4.2.42) and (4.2.43) we obtain

$$(4.2.44) \qquad\qquad \sum_{n=1}^\infty \frac{\left(H_n^{(1)}\right)^2}{n^2} + \sum_{n=1}^\infty \frac{H_n^{(2)}}{n^2} = \frac{\pi^4}{15} = 6\varsigma(4)$$

and using Adamchik's formula (3.17) we have

$$(4.2.45) \qquad \sum_{n=1}^{\infty} \frac{1}{n^2} \sum_{k=1}^{n} \frac{H_k^{(1)}}{k} = 3\varsigma(4)$$

Changing the order of summation of (4.2.45) using (3.23) we have

$$(4.2.46) \qquad \sum_{n=1}^{\infty} \frac{H_n^{(1)}}{n} \sum_{k=n}^{\infty} \frac{1}{k^2} = \sum_{n=1}^{\infty} \frac{H_n^{(1)}}{n} \left( \varsigma(2) - H_{n-1}^{(2)} \right) = 3\varsigma(4)$$

It will be noted that is not possible to dispense with the brackets in the above series because $\sum_{n=1}^{\infty} \frac{H_n^{(1)}}{n}$ is clearly divergent due to the fact that its $n$ th term exceeds $1/n$.

Continuing as before we see that

$$(4.2.47) \qquad \varsigma(5) = \frac{1}{4!} \sum_{n=0}^{\infty} \frac{1}{(n+1)^2} \left\{ \left( H_{n+1}^{(1)} \right)^3 + 3 H_{n+1}^{(1)} H_{n+1}^{(2)} + 2 H_{n+1}^{(3)} \right\}$$

or alternatively

$$(4.2.47a) \qquad \varsigma(5) = \frac{1}{4!} \sum_{n=1}^{\infty} \frac{1}{n^2} \left\{ \left( H_n^{(1)} \right)^3 + 3 H_n^{(1)} H_n^{(2)} + 2 H_n^{(3)} \right\}$$

The derivation of an expression for $\varsigma_a(5)$ is shown in (4.3.29).

By using (3.16c) equation (4.2.47a) can be written as

$$(4.2.48) \qquad \varsigma(5) = \frac{1}{4} \sum_{n=1}^{\infty} \frac{1}{n^2} \sum_{k=1}^{n} \binom{n}{k} \frac{(-1)^{k+1}}{k^3}$$

As indicated by their name, Euler was the first person to consider Euler Sums represented by

$$(4.2.49) \qquad S_{p,q} = \sum_{n=1}^{\infty} \frac{H_n^{(p)}}{n^q}$$

Euler commenced this investigation in 1742 during his correspondence with Goldbach (1690-1764) [21, p.253].

In 1995 Borwein et al. [28] proved that for an odd weight $m = p + q$, the linear sums are reducible to zeta values. The following formula was also obtained by Flajolet and Salvy [69] using impressive contour integration techniques (see also the papers by Boyadzhiev [31] and [32]).

$$(4.2.50) \qquad \sum_{n=1}^{\infty} \frac{H_n^{(p)}}{n^q} = \varsigma(m)\left\{ \frac{1}{2} - \frac{(-1)^p}{2}\binom{m-1}{p} - \frac{(-1)^p}{2}\binom{m-1}{q} \right\} + \frac{1-(-1)^p}{2}\varsigma(p)\varsigma(q)$$

$$+(-1)^p \sum_{k=1}^{\lfloor p/2 \rfloor} \binom{m-2k-1}{q-1} \varsigma(2k)\varsigma(m-2k)$$

$$+(-1)^p \sum_{k=1}^{\lfloor q/2 \rfloor} \binom{m-2k-1}{p-1} \varsigma(2k)\varsigma(m-2k)$$

where $\varsigma(1)$ should be interpreted as equal to nil wherever it occurs in the summation.

Therefore, using (4.2.47a) we should be able to derive an identity for $\sum_{n=1}^{\infty} \frac{H_n^{(1)} H_n^{(2)}}{n^2}$ and indeed similar formulae with higher weights (noting however that Flajolet and Salvy consider that no finite formula in terms of zeta values is likely to exist for the cubic sums $S_{1^3,q} = \sum_{n=1}^{\infty} \frac{\left(H_n^{(1)}\right)^3}{n^q}$ of odd weight exceeding 10 [69, p.17] ). See also (4.4.168a).

I undoubtedly require a wet towel in order to fully understand the machinations of the Flajolet and Salvy paper: however, every time I don one, the fabric dries out before I complete the task! Fully cognisant of my status as an amateur mathematician, I leave it to others more experienced than I to explore these aspects further. A good working definition of the expression "amateur mathematician" is contained in [6b].

Some of the identities which were proved above are collected below for ease of reference.

$$\varsigma_a(s) = \sum_{n=0}^{\infty} \frac{1}{2^{n+1}} \sum_{k=0}^{n} \binom{n}{k} \frac{(-1)^k}{(k+1)^s}$$

$$\varsigma_a(1) = \sum_{n=1}^{\infty} \frac{1}{n2^n} = \log 2$$

$$\varsigma_a(2) = \sum_{n=1}^{\infty} \frac{H_n}{n2^n}$$

$$\varsigma_a(3) = \frac{1}{2}\sum_{n=1}^{\infty}\frac{1}{n2^n}\left\{(H_n^{(1)})^2 + H_n^{(2)}\right\} = \sum_{n=0}^{\infty}\frac{1}{n2^n}\sum_{k=1}^{n}\binom{n}{k}\frac{(-1)^{k+1}}{k^2}$$

$$\varsigma_a(4) = \sum_{n=1}^{\infty}\frac{1}{n2^n}\left\{\frac{1}{6}\left(H_n^{(1)}\right)^3 + \frac{1}{2}H_n^{(1)}H_n^{(2)} + \frac{1}{3}H_n^{(3)}\right\} = \sum_{n=1}^{\infty}\frac{1}{n2^n}\sum_{k=1}^{n}\binom{n}{k}\frac{(-1)^{k+1}}{k^3}$$

$$\varsigma_a(5) = \sum_{n=1}^{\infty}\frac{1}{n2^n}\left\{\left(H_n^{(1)}\right)^4 + 6\left(H_n^{(1)}\right)^2 H_n^{(2)} + 8H_n^{(1)}H_n^{(3)} + 3\left(H_n^{(2)}\right)^2 + 6H_n^{(4)}\right\}$$

$$= \sum_{n=1}^{\infty}\frac{1}{n2^n}\sum_{k=1}^{n}\binom{n}{k}\frac{(-1)^{k+1}}{k^4}$$

$$\varsigma(s) = \frac{1}{s-1}\sum_{n=0}^{\infty}\frac{1}{n+1}\sum_{k=0}^{n}\binom{n}{k}\frac{(-1)^k}{(k+1)^{s-1}}$$

$$\varsigma(2) = \sum_{n=1}^{\infty}\frac{1}{n^2}$$

$$\varsigma(3) = \frac{1}{2!}\sum_{n=1}^{\infty}\frac{H_n}{n^2} = \frac{1}{2}\sum_{n=1}^{\infty}\frac{1}{n^2}\sum_{k=1}^{n}\binom{n}{k}\frac{(-1)^{k+1}}{k}$$

$$\varsigma(4) = \frac{1}{3!}\sum_{n=1}^{\infty}\frac{1}{n^2}\left\{(H_n^{(1)})^2 + H_n^{(2)}\right\} = \frac{1}{3}\sum_{n=1}^{\infty}\frac{1}{n^2}\sum_{k=1}^{n}\binom{n}{k}\frac{(-1)^{k+1}}{k^2}$$

$$\varsigma(5) = \frac{1}{4!}\sum_{n=1}^{\infty}\frac{1}{n^2}\left\{\left(H_n^{(1)}\right)^3 + 3H_n^{(1)}H_n^{(2)} + 2H_n^{(3)}\right\} = \frac{1}{4}\sum_{n=1}^{\infty}\frac{1}{n^2}\sum_{k=1}^{n}\binom{n}{k}\frac{(-1)^{k+1}}{k^3}$$

On the basis of this limited data one may conjecture that

(4.2.51) $$\varsigma(s+2) = \frac{1}{s+1}\sum_{n=1}^{\infty}\frac{1}{n^2}\sum_{k=1}^{n}\binom{n}{k}\frac{(-1)^{k+1}}{k^s} \qquad\qquad , s \geq 1$$

(4.2.52) $$\varsigma_a(s) = \sum_{n=1}^{\infty}\frac{1}{n2^n}\sum_{k=1}^{n}\binom{n}{k}\frac{(-1)^{k+1}}{k^{s-1}} \qquad\qquad , s \geq 1$$

These conjectures are proved (and in fact generalised) in equations (4.4.58) and (4.4.45) respectively: the generalised identities are

$$\sum_{n=1}^{\infty}\frac{1}{n^2}\sum_{k=1}^{n}\binom{n}{k}\frac{(-1)^k x^k}{k^s}=-(s+1)Li_{s+2}(x)+\log x\,Li_{s+1}(x)$$

$$\sum_{n=1}^{\infty}\frac{1}{n2^n}\sum_{k=1}^{n}\binom{n}{k}\frac{x^k}{k^s}=Li_{s+1}(x)$$

Equation (4.2.52) may be compared with (3.65) which is shown below

$$\varsigma_a(s)=\sum_{n=1}^{\infty}\frac{1}{2^{n+1}}\sum_{k=1}^{n}\binom{n}{k}\frac{(-1)^{k+1}}{k^s}$$

From (4.2.1) we have

$$\frac{1}{n!}\sum_{k=0}^{n}\binom{n}{k}\frac{(-1)^k}{k+x}=\frac{1}{x(1+x)...(n+x)}$$

and hence

$$\frac{y^{n+x}}{n!}\sum_{k=0}^{n}\binom{n}{k}\frac{(-1)^k}{k+x}=\frac{y^{n+x}}{x(1+x)...(n+x)}$$

The summation gives us

$$\sum_{n=0}^{\infty}\frac{y^{n+x}}{n!}\sum_{k=0}^{n}\binom{n}{k}\frac{(-1)^k}{k+x}=\sum_{n=0}^{\infty}\frac{y^{n+x}}{x(1+x)...(n+x)}$$

I noted in Nielsen's book [104a, p.9] that

$$\sum_{n=0}^{\infty}\frac{y^{n+x}}{x(1+x)...(n+x)}=e^y\int_0^y e^{-t}t^{x-1}dt$$

and we therefore obtain

$$e^y\int_0^y e^{-t}t^{x-1}dt=\sum_{n=0}^{\infty}\frac{y^{n+x}}{n!}\sum_{k=0}^{n}\binom{n}{k}\frac{(-1)^k}{k+x}=y^x\sum_{n=0}^{\infty}\frac{y^n}{n!}\sum_{k=0}^{n}\binom{n}{k}\frac{(-1)^k}{k+x}$$

I'm not sure if this will lead us anywhere but it does look interesting. The function $\gamma(x,y)=\int_0^y e^{-t}t^{x-1}dt$ is known as the incomplete gamma function [126, p.11] and this rather loose connection leads us nicely into the next part.

## A CONNECTION WITH THE GAMMA, BETA AND PSI FUNCTIONS

As is well-known, the preceding analysis can also be carried out in a concise manner by using the gamma function (whose birth was evidenced in a letter from Euler to Goldbach dated 13 October 1729 [126, p.1]). The gamma function is defined by

$$(4.3.1) \qquad \Gamma(x) = \int_0^\infty t^{x-1} e^{-t} dt \qquad , \operatorname{Re}(x) > 0$$

and, using integration by parts, it is easily shown that

$$(4.3.2) \qquad \Gamma(x+1) = x\Gamma(x)$$

Since $\Gamma(1) = 1$, it is easily shown that $\Gamma(n+1) = n!$.

A summary of some of the properties of the gamma function is contained in Appendix E of Volume VI for ease of reference.

We have from (4.3.2)

$$(4.3.3) \qquad \Gamma(x) = \frac{\Gamma(x+1)}{x} = \frac{\Gamma(x+2)}{x(x+1)} = ... = \frac{\Gamma(x+n+1)}{x(x+1)...(x+n)} \qquad , \; x > 0$$

Using (4.3.3) we may therefore write (4.2.1) as follows

$$(4.3.4) \qquad g(x) = \frac{n!}{x(1+x)...(n+x)} = \frac{n!\Gamma(x)}{\Gamma(x+n+1)}$$

and upon differentiating we obtain

$$(4.3.5) \qquad \frac{g'(x)}{n!} = \frac{\Gamma(x+n+1)\Gamma'(x) - \Gamma(x)\Gamma'(x+n+1)}{\Gamma^2(x+n+1)}$$

$$= \frac{\Gamma(x)}{\Gamma(x+n+1)} \left\{ \frac{\Gamma'(x)}{\Gamma(x)} - \frac{\Gamma'(x+n+1)}{\Gamma(x+n+1)} \right\}$$

$$(4.3.6) \qquad = \frac{g(x)}{n!} \left\{ \psi(x) - \psi(x+n+1) \right\}$$

where the digamma function $\psi(x)$ is the logarithmic derivative of $\Gamma(x)$.

$$(4.3.7) \qquad \psi(x) = \frac{d}{dx}\log\Gamma(x) = \frac{\Gamma'(x)}{\Gamma(x)}$$

Therefore we have

$$(4.3.8) \qquad g'(x) = g(x)\left\{\psi(x) - \psi(x+n+1)\right\}$$

Differentiating (4.3.8) we have

$$(4.3.9) \qquad g''(x) = g(x)\left\{\psi(x) - \psi(x+n+1)\right\}^2 + g(x)\left\{\psi'(x) - \psi'(x+n+1)\right\}$$

By taking the logarithm of (4.3.2) we get

$$(4.3.10) \qquad \log\Gamma(x+1) = \log x + \log\Gamma(x)$$

and upon differentiating (4.3.10) we obtain

$$(4.3.11) \qquad \frac{\Gamma'(x+1)}{\Gamma(x+1)} = \frac{1}{x} + \frac{\Gamma'(x)}{\Gamma(x)}$$

Equation (4.3.11) may therefore be written as

$$(4.3.12) \qquad \psi(x+1) = \frac{1}{x} + \psi(x)$$

$$(4.3.13) \qquad \psi(x+2) = \frac{1}{x+1} + \psi(x+1)$$

$$= \frac{1}{x} + \frac{1}{x+1} + \psi(x)$$

and (4.3.13) is easily generalised to

$$(4.3.14) \qquad \psi(x+n+1) = \frac{1}{x} + \frac{1}{x+1} + ... \frac{1}{x+n} + \psi(x)$$

For completeness we mention that $\psi(1) = -\gamma$ where $\gamma$ is Euler's constant defined in (4.1.15b) (see also Appendix E of Volume VI and [78, p.58]). We may also note from that

(4.3.14a) $$\lim_{x \to 0} \left( \frac{1}{x} + \psi(x) \right) = \lim_{x \to 0} \left[ \psi(x+n+1) - \left( \frac{1}{x+1} + \ldots \frac{1}{x+n} \right) \right]$$

$$= \psi(n+1) - H_n = -\gamma$$

Letting $x = 1$ we obtain

(4.3.15) $\psi(n+2) = H_{n+1} + \psi(1)$

(4.3.16) $\psi(n+1) = H_n + \psi(1) = H_n - \gamma$

And substituting (4.3.15) in (4.3.8) and (4.3.9) respectively we obtain

(4.3.17) $g'(1) = -\dfrac{H_{n+1}}{n+1}$

(4.3.18) $g''(1) = \dfrac{1}{n+1} \left\{ (H_{n+1}^{(1)})^2 + H_{n+1}^{(2)} \right\}$

and these accord with the formulae previously found in (4.2.16) and (4.2.29).

We note from (4.3.14) that

(4.3.19) $\psi^{(p)}(x+n+1) - \psi^{(p)}(x) = (-1)^p \, p! \left\{ \dfrac{1}{x^{p+1}} + \dfrac{1}{(x+1)^{p+1}} + \ldots \dfrac{1}{(x+n)^{p+1}} \right\}$

(4.3.20) $\psi^{(p)}(n+2) - \psi^{(p)}(1) = (-1)^p \, p! \, H_{n+1}^{(p+1)}$

For convenience, let us define $\Delta(x) = -\left\{ \psi(x+n+1) - \psi(x) \right\}$ so that (4.3.8), (4.3.9) and (4.3.19) read respectively

(4.3.21) $g'(x) = g(x)\Delta(x)$

(4.3.22) $g''(x) = g(x)\Delta^2(x) + g(x)\Delta'(x)$

(4.3.23) $\Delta^{(p)}(1) = (-1)^p \, p! \, H_{n+1}^{(p+1)}$

By differentiating (4.3.22) we obtain

(4.3.24) $g^{(3)}(x) = 2g(x)\Delta(x)\Delta'(x) + g'(x)\Delta^2(x) - g(x)\Delta''(x) - g'(x)\Delta'(x)$

Therefore

(4.3.25)
$$\frac{g^{(3)}(x)}{g(x)} = 3\Delta(x)\Delta'(x) - \Delta^3(x) - \Delta''(x)$$

and we obtain for $x = 1$

(4.3.26)
$$(n+1)g^{(3)}(1) = -\left\{3H_{n+1}^{(1)}H_{n+1}^{(2)} + \left(H_{n+1}^{(1)}\right)^3 + 2H_{n+1}^{(3)}\right\}$$

which is of course identical to the formula (4.2.36a) obtained by the previous method.

Similarly, by differentiating (4.3.25) we obtain the following expression for $g^{(4)}(1)$

$$\frac{g(x)g^{(4)}(x) - g^{(3)}(x)g^{(1)}}{g^2(x)} = 3\Delta(x)\Delta''(x) + 3\left(\Delta'(x)\right)^2 - 3\Delta^2(x)\Delta'(x) - \Delta^{(3)}(x)$$

and we then have

(4.3.27)
$$\frac{g^{(4)}(x)}{g(x)} = 2\Delta(x)\Delta''(x) + 3\left(\Delta'(x)\right)^2 - 6\Delta^2(x)\Delta'(x) - \Delta^{(3)}(x) + \Delta^4(x)$$

Substituting $x = 1$ we have

(4.3.28)
$$\frac{g^{(4)}(1)}{g(1)} = \left(H_{n+1}^{(1)}\right)^4 + 6\left(H_{n+1}^{(1)}\right)^2 H_{n+1}^{(2)} + 8H_{n+1}^{(1)}H_{n+1}^{(3)} + 3\left(H_{n+1}^{(2)}\right)^2 + 6H_{n+1}^{(4)}$$

and hence we have

(4.3.29)
$$\varsigma_a(5) = \sum_{n=1}^{\infty}\frac{1}{n2^n}\left\{\left(H_n^{(1)}\right)^4 + 6\left(H_n^{(1)}\right)^2 H_n^{(2)} + 8H_n^{(1)}H_n^{(3)} + 3\left(H_n^{(2)}\right)^2 + 6H_n^{(4)}\right\}$$

The result (4.3.28) was also obtained by Vermaseren [133, p.20] in 1999 who commented that the factors of (4.3.28) are related to the cycle structure of the permutation group (and, from my limited knowledge of combinatorics, I assume that this implies a reference to the Stirling numbers). See also the recent paper by Blümhein [24] where he used shuffle products to obtain these relations for use in Quantum Field Theory. The formula can of course also be derived using the general Flajolet and Sedgewick expression in (3.14).

It is clear that this method can be readily extended to determine $g^{(p)}(1)$ (albeit the algebra becomes rather tedious). As a matter of interest, we have the following relationship:

(4.3.30)
$$\frac{g^{(p)}(x)}{g(x)} = -\frac{g^{(p-1)}(x)}{g(x)} + \frac{d}{dx}\left\{\frac{g^{(p-1)}(x)}{g(x)}\right\} = e^x\frac{d}{dx}\left\{e^{-x}\frac{g^{(p-1)}(x)}{g(x)}\right\}$$

According to Ginsburg's note [70aa] this relationship hints at a connection with the Stirling numbers of the second kind because we have

(4.3.31)  $\qquad \varphi_0 = e^{-x}$

$$\varphi_1 = \frac{d}{dx}(x\varphi_0) = e^{-x}(1-x)$$

$$\varphi_2 = \frac{d}{dx}(x\varphi_1) = e^{-x}(1-3x+x^2)$$

$$\varphi_3 = \frac{d}{dx}(x\varphi_2) = e^{-x}(1-7x+6x^2+x^3)$$

and the coefficients are indeed the Stirling numbers of the second kind $S(n,k)$.

## AN APPLICATION OF THE DIGAMMA FUNCTION TO THE DERIVATION OF EULER SUMS

We now recall (4.1.7b)

(4.3.32)  $\qquad \psi(x+a) - \psi(a) = \sum_{k=1}^{\infty} \frac{(-1)^{k+1}}{k} \frac{x(x-1)...(x-k+1)}{a(a+1)...(a+k-1)}$

which converges for $\mathrm{Re}\,(x+a) > 0$. According to Raina and Ladda [109a], this summation formula is due to Nörlund (see [105], [105(i)] and also Ruben's note [111a]).

We know from [126, p.22] that (a proof is also given in Appendix E of Volume VI)

(4.3.33)  $\qquad \psi^{(n)}(z) = (-1)^{n+1} n! \sum_{k=0}^{\infty} \frac{1}{(k+z)^{n+1}} = (-1)^{n+1} n! \varsigma(n+1, z)$

and hence

(4.3.33a)  $\qquad \psi^{(n)}(1) = (-1)^{n+1} n! \varsigma(n+1, 1) = (-1)^{n+1} n! \varsigma(n+1)$

Differentiating (4.3.32) with respect to $x$ we get

$$(4.3.34) \qquad \psi'(x+a) = \sum_{k=1}^{\infty} \frac{(-1)^{k+1}}{k} \frac{s_k'(x)}{a(a+1)...(a+k-1)}$$

where

$$(4.3.34a) \qquad s_k(x) = x(x-1)...(x-k+1) = \sum_{j=0}^{k} s(k,j)x^j$$

and $s(k,j)$ are the Stirling numbers of the first kind which are briefly referred to in (3.105) in Volume 1.

We know from (3.105a) that $s_k'(0) = s(k,1) = (-1)^{k+1}(k-1)!$ and we therefore get for $x = 0$

$$(4.3.35) \qquad \psi'(a) = \sum_{k=1}^{\infty} \frac{1}{k} \frac{(k-1)!}{a(a+1)...(a+k-1)}$$

$$= \sum_{k=1}^{\infty} \frac{\Gamma(k)}{k} \frac{\Gamma(a)}{\Gamma(a+k)}$$

This result was also reported by Ruben [111a] in 1976.

Combining (4.3.35) and (4.2.1) we get

$$\psi'(a) = \sum_{k=1}^{\infty} \frac{1}{k} \sum_{j=0}^{k-1} \binom{k-1}{j} \frac{(-1)^j}{j+a}$$

which we shall also see in (4.3.70).

We also see that

$$\psi'(a) = \sum_{k=1}^{\infty} \frac{1}{k} \frac{\Gamma(k)\Gamma(a)}{\Gamma(a+k)} = \sum_{k=1}^{\infty} \frac{1}{k} B(k,a)$$

$$= \sum_{k=1}^{\infty} \frac{1}{k} \int_{0}^{1} x^{k-1}(1-x)^{a-1} dx$$

$$= -\int_{0}^{1} \frac{\log(1-x)}{x}(1-x)^{a-1} dx$$

which is equivalent to (4.3.68). More generally we have

$$\psi^{(p)}(a) = -\int\limits_0^1 \frac{\log^{p+1}(1-x)}{x}(1-x)^{a-1}dx$$

With $a = 1$, in accordance with (4.3.33a) this becomes

(4.3.36) $$\psi'(1) = \sum_{k=1}^\infty \frac{1}{k^2}$$

and more generally we see that

(4.3.37) $$\psi'(n+1) = \sum_{k=1}^\infty \frac{1}{k^2}\frac{k!(n)!}{(n+k)!} = \sum_{k=1}^\infty \frac{1}{k^2}\binom{n+k}{k}^{-1}$$

A further differentiation of (4.3.34) with respect to $x$ results in

(4.3.38) $$\psi''(x+a) = \sum_{k=1}^\infty \frac{(-1)^{k+1}}{k}\frac{s_k''(x)}{a(a+1)...(a+k-1)}$$

We know from (4.3.34a) that $s_k''(0) = 2s(k,2)$ and from (3.105iii) that $s(k,2) = (-1)^k(k-1)!H_{k-1}^{(1)}$ , and we therefore get from (4.3.38)

(4.3.39) $$\psi''(a) = -2\sum_{k=1}^\infty \frac{1}{k}\frac{(k-1)!H_{k-1}^{(1)}}{a(a+1)...(a+k-1)}$$

and with $a = 1$ this becomes

(4.3.40) $$\psi''(1) = -2\sum_{k=1}^\infty \frac{H_{k-1}^{(1)}}{k^2} = -2\sum_{k=1}^\infty \frac{H_k^{(1)}}{k^2} + 2\varsigma(3) = -2\varsigma(3)$$

Hence we have shown by a new method that

(4.3.40a) $$\sum_{k=1}^\infty \frac{H_k^{(1)}}{k^2} = 2\varsigma(3)$$

Differentiating (4.3.35) we obtain

(4.3.40b) $$\psi''(a) = \sum_{k=1}^\infty \frac{\Gamma(k)}{k}\frac{\Gamma(a)}{\Gamma(a+k)}[\psi(a) - \psi(a+k)]$$

From (4.3.34a) it is apparent that $s_k^{(3)}(0) = 3!s(k,3)$ where (3.104) gives us

$(4.3.41)$ $$s(k,3) = (-1)^{k+1} \frac{(k-1)!}{2}\left\{\left(H_{k-1}^{(1)}\right)^2 - H_{k-1}^{(2)}\right\}$$

and hence we get by differentiating $(4.3.38)$ and letting $x=0$

$(4.3.42)$ $$\psi^{(3)}(a) = 3\sum_{k=1}^{\infty} \frac{1}{k} \frac{(k-1)!\left\{\left(H_{k-1}^{(1)}\right)^2 - H_{k-1}^{(2)}\right\}}{a(a+1)...(a+k-1)}$$

Letting $a=1$ this becomes

$$\psi^{(3)}(1) = 3\sum_{k=1}^{\infty} \frac{\left\{\left(H_{k-1}^{(1)}\right)^2 - H_{k-1}^{(2)}\right\}}{k^2}$$

Therefore we obtain using $(4.3.33a)$

$(4.3.43)$ $$\psi^{(3)}(1) = 3\sum_{k=1}^{\infty} \frac{\left\{\left(H_{k-1}^{(1)}\right)^2 - H_{k-1}^{(2)}\right\}}{k^2} = 3!\varsigma(4)$$

We may write

$(4.3.43a)$ $$\left(H_{k-1}^{(1)}\right)^2 - H_{k-1}^{(2)} = \left(H_k^{(1)} - \frac{1}{k}\right)^2 - H_k^{(2)} + \frac{1}{k^2} = \left(H_k^{(1)}\right)^2 - 2\frac{H_k^{(1)}}{k} - H_k^{(2)} + 2\frac{1}{k^2}$$

and hence

$(4.3.44)$ $$\sum_{k=1}^{\infty} \frac{\left\{\left(H_{k-1}^{(1)}\right)^2 - H_{k-1}^{(2)}\right\}}{k^2} = \sum_{k=1}^{\infty} \frac{\left[H_k^{(1)}\right]^2}{k^2} - 2\sum_{k=1}^{\infty} \frac{H_k^{(1)}}{k^3} - \sum_{k=1}^{\infty} \frac{H_k^{(2)}}{k^2} + 2\sum_{k=1}^{\infty} \frac{1}{k^4} = 2\varsigma(4)$$

Therefore we obtain

$(4.3.44a)$ $$\sum_{k=1}^{\infty} \frac{\left[H_k^{(1)}\right]^2}{k^2} - 2\sum_{k=1}^{\infty} \frac{H_k^{(1)}}{k^3} - \sum_{k=1}^{\infty} \frac{H_k^{(2)}}{k^2} = 0$$

From $(4.3.35)$ we have

$$\psi'(a) = \sum_{k=1}^{\infty} \frac{1}{k} \frac{(k-1)!}{a(a+1)...(a+k-1)}$$

and, using $(4.3.55)$, differentiation gives us

$$(4.3.45) \qquad \psi''(a) = -\sum_{k=1}^{\infty} \frac{1}{k}(k-1)! A(a)\big[\psi(a+k)-\psi(a)\big]$$

where $A(a) = \dfrac{1}{a(a+1)...(a+k-1)} = \dfrac{\Gamma(a)}{\Gamma(a+k)}$

Letting $a=1$ we get the familiar result

$$(4.3.45a) \qquad \psi''(1) = -\sum_{k=1}^{\infty} \frac{1}{k^2}\big[\psi(1+k)-\psi(1)\big] = -\sum_{k=1}^{\infty} \frac{H_k^{(1)}}{k^2}$$

A further differentiation results in

$$(4.3.45b) \qquad \psi^{(3)}(a) = -\sum_{k=1}^{\infty} \frac{1}{k}(k-1)! A(a)\Big\{\big[\psi'(a+k)-\psi'(a)\big]-\big[\psi(a+k)-\psi(a)\big]^2\Big\}$$

and with $a=1$ we get

$$(4.3.45c) \qquad \psi^{(3)}(1) = -\sum_{k=1}^{\infty} \frac{1}{k^2}\Big\{\big[\psi'(1+k)-\psi'(1)\big]-\big[H_k^{(1)}\big]^2\Big\}$$

Since $\psi(z+k)-\psi(z) = \sum_{j=1}^{k} \dfrac{1}{z+j-1}$ we see that [126, p.22]

$$(4.3.45d) \qquad \psi^{(n)}(z+k)-\psi^{(n)}(z) = (-1)^n n! \sum_{j=1}^{k} \frac{1}{(z+j-1)^{n+1}}$$

and therefore

$$(4.3.45e) \qquad \psi^{(n)}(1+k)-\psi^{(n)}(1) = (-1)^n n! \sum_{j=1}^{k} \frac{1}{j^{n+1}} = (-1)^n n! H_k^{(n+1)}$$

Hence we obtain from (4.3.45c) with $n=1$

$$\psi^{(3)}(1) = \sum_{k=1}^{\infty} \frac{H_k^{(2)}}{k^2} + \sum_{k=1}^{\infty} \frac{\big[H_k^{(1)}\big]^2}{k^2}$$

and therefore we have another proof of (4.2.40)

$$(4.3.45f) \qquad 6\zeta(4) = \sum_{k=1}^{\infty} \frac{H_k^{(2)}}{k^2} + \sum_{k=1}^{\infty} \frac{\big[H_k^{(1)}\big]^2}{k^2}$$

Let us now employ a slightly different approach by differentiating (4.3.39) with respect to $a$

$$\psi''(a) = -2\sum_{k=1}^{\infty}\frac{1}{k}\frac{(k-1)!H_{k-1}^{(1)}}{a(a+1)...(a+k-1)}$$

to obtain

(4.3.46)     $$\psi^{(3)}(a) = 2\sum_{k=1}^{\infty}\frac{1}{k}\frac{(k-1)!H_{k-1}^{(1)}\left[\psi(a+k)-\psi(a)\right]}{a(a+1)...(a+k-1)}$$

We get with $a=1$

$$\psi^{(3)}(1) = 2\sum_{k=1}^{\infty}\frac{H_{k-1}^{(1)}H_k^{(1)}}{k^2} = 2\sum_{k=1}^{\infty}\frac{\left[H_k^{(1)}\right]^2}{k^2} - 2\sum_{k=1}^{\infty}\frac{H_k^{(1)}}{k^3}$$

Since $\psi^{(3)}(1) = 3!\varsigma(4)$ we get

(4.3.46a)     $$3\varsigma(4) = \sum_{k=1}^{\infty}\frac{\left[H_k^{(1)}\right]^2}{k^2} - \sum_{k=1}^{\infty}\frac{H_k^{(1)}}{k^3}$$

Hence we have three simultaneous equations (4.3.44a), (4.3.45f) and (4.3.46a)

$$0 = \sum_{k=1}^{\infty}\frac{\left[H_k^{(1)}\right]^2}{k^2} - 2\sum_{k=1}^{\infty}\frac{H_k^{(1)}}{k^3} - \sum_{k=1}^{\infty}\frac{H_k^{(2)}}{k^2}$$

$$6\varsigma(4) = \sum_{k=1}^{\infty}\frac{H_k^{(2)}}{k^2} + \sum_{k=1}^{\infty}\frac{\left[H_k^{(1)}\right]^2}{k^2}$$

$$3\varsigma(4) = \sum_{k=1}^{\infty}\frac{\left[H_k^{(1)}\right]^2}{k^2} - \sum_{k=1}^{\infty}\frac{H_k^{(1)}}{k^3}$$

from which we can determine the various Euler sums (which we have previously obtained by alternative means). We then obtain in a more unified manner

(4.3.46b)     $$\sum_{k=1}^{\infty}\frac{\left(H_k^{(1)}\right)^2}{k^2} = \frac{17}{4}\varsigma(4) \qquad \sum_{k=1}^{\infty}\frac{H_k^{(1)}}{k^3} = \frac{5}{4}\varsigma(4) \qquad \sum_{k=1}^{\infty}\frac{H_k^{(2)}}{k^2} = \frac{7}{4}\varsigma(4)$$

It is clear that this method may be easily extended by using further expressions for the Stirling numbers of the first kind. In fact we have

$$(4.3.47) \qquad \psi^{(p)}(x+a) = \sum_{k=1}^{\infty} \frac{(-1)^{k+1}}{k} \frac{s_k^{(p)}(x)}{a(a+1)...(a+k-1)}$$

$$\psi^{(p)}(a) = \sum_{k=1}^{\infty} \frac{(-1)^{k+1}}{k} \frac{s_k^{(p)}(0)}{a(a+1)...(a+k-1)}$$

Since

$$(4.3.47a) \qquad s_k^{(p)}(0) = p!s(k,p)$$

we have

$$(4.3.48) \qquad \psi^{(p)}(a) = \sum_{k=1}^{\infty} \frac{(-1)^{k+1}}{k} \frac{p!s(k,p)}{a(a+1)...(a+k-1)}$$

and with $a = 1$ we obtain

$$(4.3.49) \qquad \psi^{(p)}(1) = p! \sum_{k=1}^{\infty} \frac{(-1)^{k+1}}{k.k!} s(k,p)$$

Reference to (4.3.33a) then shows that

$$(4.3.50) \qquad \varsigma(p+1) = (-1)^p \sum_{k=1}^{\infty} \frac{(-1)^k}{k.k!} s(k,p)$$

This result was previously obtained in 1995 by Shen [120] by employing a different method. Another proof was given in [30a]. In [120] Shen also referred to Morley's identity [135, p.301] which is valid for $\operatorname{Re}(a) < 2/3$

$$(4.3.51) \qquad 1 + \sum_{k=1}^{\infty} \left[ \frac{a(a+1)...(a+k-1)}{k!} \right]^3 = \cos(\pi a/2) \frac{\Gamma\left(1 - \frac{3}{2}a\right)}{\Gamma^3\left(1 - \frac{1}{2}a\right)}$$

and used it to show that

$$\sum_{k=2}^{\infty} \frac{H_k^{(1)}}{(k+1)^3} = \frac{1}{4}\varsigma(4)$$

We have

$$S_k(a) = a(a+1)...(a+k-1) = \frac{\Gamma(a+k)}{\Gamma(a)}$$

$$S_k'(a) = s_k(a)\big[\psi(a+k) - \psi(a)\big]$$

and hence we get by differentiation of (4.3.51)

(4.3.52)
$$3\sum_{k=1}^{\infty}\left[\frac{S_k(a)}{k!}\right]^3\big[\psi(a+k) - \psi(a)\big] =$$

$$\frac{3}{2}\cos(\pi a/2)\frac{-\Gamma^3\left(1-\frac{1}{2}a\right)\Gamma'\left(1-\frac{3}{2}a\right) + \Gamma\left(1-\frac{3}{2}a\right)\Gamma^2\left(1-\frac{1}{2}a\right)\Gamma'\left(1-\frac{1}{2}a\right)}{\Gamma^6\left(1-\frac{1}{2}a\right)}$$

$$-\frac{\pi}{2}\sin(\pi a/2)\frac{\Gamma\left(1-\frac{3}{2}a\right)}{\Gamma^3\left(1-\frac{1}{2}a\right)}$$

We have Legendre's duplication formula for the gamma function [8a, p.22]

$$\sqrt{\pi}\,\Gamma(2z) = 2^{2z-1}\Gamma(z)\Gamma\left(z+\frac{1}{2}\right) \quad , z \neq 0, -\frac{1}{2}, -1, -\frac{3}{2},...$$

and therefore we have [126, p.4]

$$S_k(1/2) = \frac{\Gamma(k+1/2)}{\Gamma(1/2)} = \frac{\sqrt{\pi}\,\Gamma(2k)}{\Gamma(1/2)2^{2k-1}\Gamma(k)} = \frac{(2k)!}{2^{2k}k!}$$

Using [126, p.20]

(4.3.53)
$$\big[\psi(k+1/2) - \psi(1/2)\big] = 2\sum_{j=0}^{k-1}\frac{1}{2j+1}$$

we get

$$6\sum_{k=1}^{\infty}\left[\frac{(2k)!}{2^{2k}(k!)^2}\right]^3\sum_{j=0}^{k-1}\frac{1}{2j+1} =$$

$$-\frac{3\sqrt{2}}{4}\frac{\Gamma'\left(\frac{1}{4}\right)}{\Gamma^3\left(\frac{3}{4}\right)}+\frac{3\sqrt{2}}{4}\frac{\Gamma\left(\frac{1}{4}\right)\Gamma'\left(\frac{3}{4}\right)}{\Gamma^4\left(\frac{3}{4}\right)}-\frac{\pi\sqrt{2}}{4}\frac{\Gamma\left(\frac{1}{4}\right)}{\Gamma^3\left(\frac{3}{4}\right)}$$

$$=\frac{3\sqrt{2}}{4}\frac{1}{\Gamma^4\left(\frac{3}{4}\right)}\left[\Gamma\left(\frac{1}{4}\right)\Gamma'\left(\frac{3}{4}\right)-\Gamma\left(\frac{3}{4}\right)\Gamma'\left(\frac{1}{4}\right)\right]-\frac{\pi\sqrt{2}}{4}\frac{\Gamma\left(\frac{1}{4}\right)}{\Gamma^3\left(\frac{3}{4}\right)}$$

Differentiating Euler's reflection identity we obtain

$$\Gamma'(x)\Gamma(1-x)-\Gamma(x)\Gamma'(1-x)=-\frac{\pi^2\cos\pi x}{\sin^2\pi x}$$

and with $x=1/4$ we obtain

$$\Gamma\left(\frac{1}{4}\right)\Gamma'\left(\frac{3}{4}\right)-\Gamma\left(\frac{3}{4}\right)\Gamma'\left(\frac{1}{4}\right)=\sqrt{2}\pi^2$$

Therefore we obtain

$$(4.3.54)\qquad \sum_{k=1}^{\infty}\left[\frac{(2k)!}{2^{2k}(k!)^2}\right]^3\sum_{j=0}^{k-1}\frac{1}{2j+1}=\frac{\pi}{24}\frac{6\pi-\sqrt{2}\Gamma(1/4)}{\Gamma^4(3/4)}$$

Similar identities are given by Srivastava and Choi [126, p.251ff].

Let us now differentiate (4.3.32) with respect to $a$ this time. This gives us

$$\frac{\partial}{\partial a}\left[\psi(x+a)-\psi(a)\right]=\sum_{k=1}^{\infty}\frac{(-1)^{k+1}}{k}A'(a)x(x-1)...(x-k+1)$$

where $A(a)=\dfrac{1}{a(a+1)...(a+k-1)}=\dfrac{\Gamma(a)}{\Gamma(a+k)}$ and $A(1)=\dfrac{1}{k!}$.

We have

$$(4.3.55)\qquad A'(a)=-A(a)\left[\psi(a+k)-\psi(a)\right]$$

and therefore

(4.3.56)
$$\frac{\partial}{\partial a}\big[\psi(x+a)-\psi(a)\big]=\sum_{k=1}^{\infty}\frac{(-1)^k}{k}A(a)\big[\psi(a+k)-\psi(a)\big]x(x-1)...(x-k+1)$$

We now differentiate (4.3.56) with respect to $x$ and obtain

(4.3.56a)
$$\frac{\partial}{\partial x}\frac{\partial}{\partial a}\big[\psi(x+a)\big]=\sum_{k=1}^{\infty}\frac{(-1)^k}{k}A(a)\big[\psi(a+k)-\psi(a)\big]s'(x)$$

Therefore we have

(4.3.56b)
$$\frac{\partial}{\partial x}\frac{\partial}{\partial a}\big[\psi(x+a)\big]\bigg|_{x=0}=-\sum_{k=1}^{\infty}\frac{(k-1)!}{k}\frac{\Gamma(a)}{\Gamma(a+k)}\big[\psi(a+k)-\psi(a)\big]$$

and this is equivalent to (4.3.45).

Another differentiation of (4.4.45b) gives us

(4.3.57)
$$\psi^{(4)}(a)=-\sum_{k=1}^{\infty}\frac{1}{k}(k-1)!A(a)\left\{\begin{array}{l}\big[\psi''(a+k)-\psi''(a)\big]-3\big[\psi(a+k)-\psi(a)\big]\big[\psi'(a+k)-\psi'(a)\big]\\[2mm]+\big[\psi(a+k)-\psi(a)\big]^3\end{array}\right\}$$

and $a=1$ gives us

(4.3.57a)
$$\psi^{(4)}(1)=-\sum_{k=1}^{\infty}\frac{1}{k^2}\left\{2H_k^{(3)}+3H_k^{(1)}H_k^{(2)}+\Big[H_k^{(1)}\Big]^3\right\}$$

We therefore have

(4.3.57b)
$$24\varsigma(5)=2\sum_{k=1}^{\infty}\frac{H_k^{(3)}}{k^2}+3\sum_{k=1}^{\infty}\frac{H_k^{(1)}H_k^{(2)}}{k^2}+\sum_{k=1}^{\infty}\frac{\Big[H_k^{(1)}\Big]^3}{k^2}$$

which is the same as (4.2.27a), which was derived by reference to the Hasse identity (3.12).

We now differentiate (4.3.46) to get

(4.3.58)
$$\psi^{(4)}(a)=2\sum_{k=1}^{\infty}\frac{1}{k}\frac{(k-1)!H_{k-1}^{(1)}\big[\psi'(a+k)-\psi'(a)\big]-\big[\psi(a+k)-\psi(a)\big]^2}{a(a+1)...(a+k-1)}$$

and therefore we have

(4.3.58a) $$\psi^{(4)}(1) = -4!\varsigma(5) = -2\sum_{k=1}^{\infty}\frac{H_{k-1}^{(1)}\left(H_k^{(2)}+\left[H_k^{(1)}\right]^2\right)}{k^2}$$

We have

$$\sum_{k=1}^{\infty}\frac{H_{k-1}^{(1)}\left(H_k^{(2)}+\left[H_k^{(1)}\right]^2\right)}{k^2} = \sum_{k=1}^{\infty}\frac{H_k^{(1)}\left(H_k^{(2)}+\left[H_k^{(1)}\right]^2\right)}{k^2} - \sum_{k=1}^{\infty}\frac{\left(H_k^{(2)}+\left[H_k^{(1)}\right]^2\right)}{k^3}$$

$$= \sum_{k=1}^{\infty}\frac{H_k^{(1)}H_k^{(2)}}{k^2} + \sum_{k=1}^{\infty}\frac{\left[H_k^{(1)}\right]^3}{k^2} - \sum_{k=1}^{\infty}\frac{H_k^{(2)}}{k^3} - \sum_{k=1}^{\infty}\frac{\left[H_k^{(1)}\right]^2}{k^3}$$

Therefore we get

(4.3.58b) $$12\varsigma(5) = \sum_{k=1}^{\infty}\frac{H_k^{(1)}H_k^{(2)}}{k^2} + \sum_{k=1}^{\infty}\frac{\left[H_k^{(1)}\right]^3}{k^2} - \sum_{k=1}^{\infty}\frac{H_k^{(2)}}{k^3} - \sum_{k=1}^{\infty}\frac{\left[H_k^{(1)}\right]^2}{k^3}$$

The next procedure is to differentiate (4.3.42)

$$\psi^{(3)}(a) = 3\sum_{k=1}^{\infty}\frac{1}{k}\frac{(k-1)!\left\{\left(H_{k-1}^{(1)}\right)^2 - H_{k-1}^{(2)}\right\}}{a(a+1)...(a+k-1)}$$

and we get

(4.3.59) $$\psi^{(4)}(a) = -3\sum_{k=1}^{\infty}\frac{1}{k}\frac{(k-1)!\left\{\left(H_{k-1}^{(1)}\right)^2 - H_{k-1}^{(2)}\right\}\left[\psi(a+k)-\psi(a)\right]}{a(a+1)...(a+k-1)}$$

Hence we have

$$\psi^{(4)}(1) = -3\sum_{k=1}^{\infty}\frac{\left\{\left(H_{k-1}^{(1)}\right)^2 - H_{k-1}^{(2)}\right\}H_k^{(1)}}{k^2}$$

and, using (4.3.43a) we have

$$= -3\sum_{k=1}^{\infty}\frac{\left\{\left(H_k^{(1)}\right)^2 - 2\frac{H_k^{(1)}}{k} - H_k^{(2)} + 2\frac{1}{k^2}\right\}H_k^{(1)}}{k^2}$$

$$= -3\sum_{k=1}^{\infty}\frac{\left[H_k^{(1)}\right]^3}{k^2} + 6\sum_{k=1}^{\infty}\frac{\left[H_k^{(1)}\right]^2}{k^3} + 3\sum_{k=1}^{\infty}\frac{H_k^{(1)}H_k^{(2)}}{k^2} - 6\sum_{k=1}^{\infty}\frac{H_k^{(1)}}{k^4}$$

Hence we obtain

(4.3.59a) $$8\zeta(5) = \sum_{k=1}^{\infty}\frac{\left[H_k^{(1)}\right]^3}{k^2} - 2\sum_{k=1}^{\infty}\frac{\left[H_k^{(1)}\right]^2}{k^3} - \sum_{k=1}^{\infty}\frac{H_k^{(1)}H_k^{(2)}}{k^2} + 2\sum_{k=1}^{\infty}\frac{H_k^{(1)}}{k^4}$$

Similarly, we have from (3.105i)

$$s(k,4) = (-1)^k\frac{(k-1)!}{6}\left\{\left(H_{k-1}^{(1)}\right)^3 - 3H_{k-1}^{(1)}H_{k-1}^{(2)} + 2H_{k-1}^{(3)}\right\}$$

and hence we get from (4.3.47a)

(4.3.60) $$\psi^{(4)}(a) = -4\sum_{k=1}^{\infty}\frac{1}{k}\frac{(k-1)!\left\{\left(H_{k-1}\right)^3 - 3H_{k-1}H_{k-1}^{(2)} + 2H_{k-1}^{(3)}\right\}}{a(a+1)...(a+k-1)}$$

and

(4.3.60a) $$\psi^{(4)}(1) = -4\sum_{k=1}^{\infty}\frac{\left\{\left(H_{k-1}^{(1)}\right)^3 - 3H_{k-1}^{(1)}H_{k-1}^{(2)} + 2H_{k-1}^{(3)}\right\}}{k^2} = -4!\zeta(5)$$

We see that

$$\left(H_{k-1}^{(1)}\right)^3 = \left(H_k^{(1)} - \frac{1}{k}\right)^3 = \left(H_k^{(1)}\right)^3 - 3\frac{\left[H_k^{(1)}\right]^2}{k} + 3\frac{H_k^{(1)}}{k^2} - \frac{1}{k^3}$$

$$H_{k-1}^{(1)}H_{k-1}^{(2)} = \left(H_k^{(1)} - \frac{1}{k}\right)\left(H_k^{(2)} - \frac{1}{k^2}\right) = H_k^{(1)}H_k^{(2)} - \frac{H_k^{(2)}}{k} - \frac{H_k^{(1)}}{k^2} + \frac{1}{k^3}$$

and therefore we have

$$\sum_{k=1}^{\infty}\frac{\left\{\left(H_{k-1}^{(1)}\right)^3 - 3H_{k-1}^{(1)}H_{k-1}^{(2)} + 2H_{k-1}^{(3)}\right\}}{k^2} =$$

$$\sum_{k=1}^{\infty}\frac{\left(H_k^{(1)}\right)^3}{k^2} - 3\sum_{k=1}^{\infty}\frac{\left(H_k^{(1)}\right)^2}{k^3} + 6\sum_{k=1}^{\infty}\frac{H_k^{(1)}}{k^4} - 3\sum_{k=1}^{\infty}\frac{H_k^{(1)}H_k^{(2)}}{k^2} + 2\sum_{k=1}^{\infty}\frac{H_k^{(3)}}{k^2} - 6\zeta(5)$$

Accordingly we see that

(4.3.60b)

$$\sum_{k=1}^{\infty} \frac{\left(H_k^{(1)}\right)^3}{k^2} - 3\sum_{k=1}^{\infty} \frac{\left(H_k^{(1)}\right)^2}{k^3} + 6\sum_{k=1}^{\infty} \frac{H_k^{(1)}}{k^4} - 3\sum_{k=1}^{\infty} \frac{H_k^{(1)} H_k^{(2)}}{k^2} + 3\sum_{k=1}^{\infty} \frac{H_k^{(2)}}{k^3} + 2\sum_{k=1}^{\infty} \frac{H_k^{(3)}}{k^2} = 12\varsigma(5)$$

Therefore, we have four simultaneous equations (4.3.57b), (4.3.58b), (4.3.59a) and (4.3.60b) involving six unknowns

$$\sum_{k=1}^{\infty} \frac{H_k^{(1)}}{k^4}, \quad \sum_{k=1}^{\infty} \frac{H_k^{(1)} H_k^{(2)}}{k^2}, \quad \sum_{k=1}^{\infty} \frac{\left[H_k^{(1)}\right]^2}{k^3}, \quad \sum_{k=1}^{\infty} \frac{\left[H_k^{(1)}\right]^3}{k^2}, \quad \sum_{k=1}^{\infty} \frac{H_k^{(2)}}{k^3} \text{ and } \sum_{k=1}^{\infty} \frac{H_k^{(3)}}{k^2}$$

$$24\varsigma(5) = 3\sum_{k=1}^{\infty} \frac{H_k^{(1)} H_k^{(2)}}{k^2} + \sum_{k=1}^{\infty} \frac{\left[H_k^{(1)}\right]^3}{k^2} + 2\sum_{k=1}^{\infty} \frac{H_k^{(3)}}{k^2}$$

$$12\varsigma(5) = \sum_{k=1}^{\infty} \frac{H_k^{(1)} H_k^{(2)}}{k^2} - \sum_{k=1}^{\infty} \frac{\left[H_k^{(1)}\right]^2}{k^3} + \sum_{k=1}^{\infty} \frac{\left[H_k^{(1)}\right]^3}{k^2} - \sum_{k=1}^{\infty} \frac{H_k^{(2)}}{k^3}$$

$$8\varsigma(5) = 2\sum_{k=1}^{\infty} \frac{H_k^{(1)}}{k^4} - \sum_{k=1}^{\infty} \frac{H_k^{(1)} H_k^{(2)}}{k^2} - 2\sum_{k=1}^{\infty} \frac{\left[H_k^{(1)}\right]^2}{k^3} + \sum_{k=1}^{\infty} \frac{\left[H_k^{(1)}\right]^3}{k^2}$$

$$12\varsigma(5) = 6\sum_{k=1}^{\infty} \frac{H_k^{(1)}}{k^4} - 3\sum_{k=1}^{\infty} \frac{H_k^{(1)} H_k^{(2)}}{k^2} - 3\sum_{k=1}^{\infty} \frac{\left(H_k^{(1)}\right)^2}{k^3} + \sum_{k=1}^{\infty} \frac{\left(H_k^{(1)}\right)^3}{k^2} + 3\sum_{k=1}^{\infty} \frac{H_k^{(2)}}{k^3} + 2\sum_{k=1}^{\infty} \frac{H_k^{(3)}}{k^2}$$

However, we already know from Flajolet and Salvy [69] and Borwein et al. [28] that (see also (3.110m))

$$\sum_{k=1}^{\infty} \frac{H_k^{(1)}}{k^4} = 3\varsigma(5) - \varsigma(2)\varsigma(3)$$

$$\sum_{k=1}^{\infty} \frac{H_k^{(2)}}{k^3} = -\frac{9}{2}\varsigma(5) + 3\varsigma(2)\varsigma(3)$$

$$\sum_{k=1}^{\infty} \frac{H_k^{(3)}}{k^2} = \frac{11}{2}\varsigma(5) - 2\varsigma(2)\varsigma(3)$$

and therefore we may solve for the other three unknowns. The solutions are:

(4.3.60c)
$$\sum_{k=1}^{\infty} \frac{H_k^{(1)} H_k^{(2)}}{k^2} = \varsigma(5) + \varsigma(2)\varsigma(3)$$

(4.3.60d)
$$\sum_{k=1}^{\infty} \frac{\left[H_k^{(1)}\right]^3}{k^2} = 10\varsigma(5) + \varsigma(2)\varsigma(3)$$

(4.3.60e)
$$\sum_{k=1}^{\infty} \frac{\left[H_k^{(1)}\right]^2}{k^3} = \frac{7}{2}\varsigma(5) - \varsigma(2)\varsigma(3)$$

some of which we also saw in Volume I.

Equations (4.3.60c), (4.3.60d) and (4.3.60e) were also derived by Panholzer and Prodinger in [105a]. Equation (4.3.60c) was also derived by Choi and Srivastava [45aci] in 2005. Flajolet and Salvy [69] reported that

(4.3.60f)
$$\sum_{k=1}^{\infty} \frac{\left[H_k^{(1)}\right]^3}{(k+1)^2} = \frac{15}{2}\varsigma(5) + \varsigma(2)\varsigma(3)$$

We will now consider various identities for $\psi^{(5)}(a)$ in a more systematic manner. First of all we have

(4.3.61i)
$$\psi^{(5)}(a) = \frac{\partial}{\partial a} \frac{\partial^4}{\partial x^4} F(a,x) \Big|_{x=0}$$

where $F(a,x) = \psi(x+a) - \psi(a)$ and this is followed by four other related identities

(4.3.61ii)
$$\psi^{(5)}(a) = \frac{\partial^2}{\partial a^2} \frac{\partial^3}{\partial x^3} F(a,x) \Big|_{x=0}$$

(4.3.61iii)
$$\psi^{(5)}(a) = \frac{\partial^3}{\partial a^3} \frac{\partial^2}{\partial x^2} F(a,x) \Big|_{x=0}$$

(4.3.61iv)
$$\psi^{(5)}(a) = \frac{\partial^4}{\partial a^4} \frac{\partial}{\partial x} F(a,x) \Big|_{x=0}$$

(4.3.61v)
$$\psi^{(5)}(a) = \frac{\partial^5}{\partial x^5} F(a,x) \Big|_{x=0}$$

More generally, it is easily seen that there will be $n$ corresponding identities for $\psi^{(n)}(a)$.

We have

$$A(a) = \frac{1}{a(a+1)...(a+k-1)} = \frac{\Gamma(a)}{\Gamma(a+k)} \quad \text{and} \quad A(1) = \frac{1}{k!}.$$

and therefore successive derivatives are

$$A^{(1)}(a) = -A(a)\big[\psi(a+k) - \psi(a)\big] = -A(a)B(a)$$

$$A^{(2)}(a) = A(a)\big[-B^{(1)}(a) + B^2(a)\big]$$

$$A^{(3)}(a) = A(a)\big[-B^{(2)}(a) + 3B(a)B^{(1)}(a) - B^3(a)\big]$$

$$A^{(4)}(a) = A(a)\big[-B^{(3)}(a) + 4B(a)B^{(2)}(a) + 3\big[B^{(1)}(a)\big]^2 - 6B^2(a)B^{(1)}(a) + B^4(a)\big]$$

$$A^{(5)}(a) = A(a)\begin{bmatrix} -B^{(4)}(a) + 5B(a)B^{(3)}(a) + 10\big[B^{(1)}(a)\big]B^{(2)}(a) \\ -10B^2(a)B^{(2)}(a) - 15B(a)\big[B^{(1)}(a)\big]^2 + 10B^3(a)B^{(1)}(a) - B^5(a) \end{bmatrix}$$

We have from (4.3.45), (4.3.42), (4.3.39) and (4.3.35) respectively

(4.3.62i)
$$\left.\frac{\partial^4}{\partial x^4}\psi(x+a)\right|_{x=0} = \psi^{(4)}(a) = -6\sum_{k=1}^{\infty}\frac{1}{k}\frac{(k-1)!\left\{\big(H_{k-1}\big)^3 - 3H_{k-1}H_{k-1}^{(2)} + 2H_{k-1}^{(3)}\right\}}{a(a+1)...(a+k-1)}$$

(4.3.62ii)
$$\left.\frac{\partial^3}{\partial x^3}\psi(x+a)\right|_{x=0} = \psi^{(3)}(a) = 3\sum_{k=1}^{\infty}\frac{1}{k}\frac{(k-1)!\left\{\big(H_{k-1}^{(1)}\big)^2 - H_{k-1}^{(2)}\right\}}{a(a+1)...(a+k-1)}$$

(4.3.62iii)
$$\left.\frac{\partial^2}{\partial x^2}\psi(x+a)\right|_{x=0} = \psi^{(2)}(a) = -2\sum_{k=1}^{\infty}\frac{1}{k}\frac{(k-1)!H_{k-1}^{(1)}}{a(a+1)...(a+k-1)}$$

(4.3.62iv)
$$\left.\frac{\partial}{\partial x}\psi(x+a)\right|_{x=0} = \psi^{(1)}(a) = \sum_{k=1}^{\infty}\frac{1}{k}\frac{(k-1)!}{a(a+1)...(a+k-1)}$$

For example we get

(4.3.63i)

$$\psi^{(5)}(a) = \frac{\partial}{\partial a}\frac{\partial^4}{\partial x^4}\psi(x+a)\bigg|_{x=0} = -6\sum_{k=1}^{\infty}\frac{1}{k}\frac{(k-1)!\left\{\left(H_{k-1}\right)^3 - 3H_{k-1}H_{k-1}^{(2)} + 2H_{k-1}^{(3)}\right\}\left[\psi(a+k) - \psi(a)\right]}{a(a+1)...(a+k-1)}$$

(4.3.63ii)

$$\psi^{(5)}(a) = \frac{\partial^2}{\partial a^2}\frac{\partial^3}{\partial x^3}\psi(x+a)\bigg|_{x=0}$$

$$= 3\sum_{k=1}^{\infty}\frac{1}{k}\frac{(k-1)!\left\{\left(H_{k-1}^{(1)}\right)^2 - H_{k-1}^{(2)}\right\}}{a(a+1)...(a+k-1)}\left\{-\left[\psi'(a+k) - \psi'(a)\right] + \left[\psi(a+k) - \psi(a)\right]^2\right\}$$

The corresponding expressions for the other identities are too lengthy to write out explicitly but they may of course be easily derived by simple algebra.

We have with $a = 1$

(4.3.64) $A^{(1)}(1) = -\frac{1}{k!}H_k^{(1)}$

$$A^{(2)}(1) = \frac{1}{k!}\left(H_k^{(2)} + \left[H_k^{(1)}\right]^2\right)$$

$$A^{(3)}(1) = \frac{1}{k!}\left(-2H_k^{(3)} + 3H_k^{(1)}H_k^{(2)} - \left[H_k^{(1)}\right]^3\right)$$

$$A^{(4)}(1) = \frac{1}{k!}\left(6H_k^{(4)} + 8H_k^{(1)}H_k^{(3)} + 3\left[H_k^{(2)}\right]^2 + 6\left[H_k^{(1)}\right]^2 H_k^{(2)} + \left[H_k^{(1)}\right]^4\right)$$

$$A^{(5)}(1) = \frac{1}{k!}\left(\begin{array}{l}-24H_k^{(5)} - 30H_k^{(1)}H_k^{(4)} - 20H_k^{(2)}H_k^{(3)} - 20\left[H_k^{(1)}\right]^2 H_k^{(3)}\\[2mm] -15H_k^{(1)}\left[H_k^{(2)}\right]^2 - 10\left[H_k^{(1)}\right]^3 H_k^{(2)} - \left[H_k^{(1)}\right]^5\end{array}\right)$$

We then obtain

(4.3.65i) $\quad \psi^{(5)}(1) = 6\sum_{k=1}^{\infty}\frac{1}{k^2}\left\{\left(H_{k-1}\right)^3 - 3H_{k-1}H_{k-1}^{(2)} + 2H_{k-1}^{(3)}\right\}H_k$

$$= 6\sum_{k=1}^{\infty} \frac{\left(H_k^{(1)}\right)^4}{k^2} - 18\sum_{k=1}^{\infty} \frac{\left(H_k^{(1)}\right)^3}{k^3} + 36\sum_{k=1}^{\infty} \frac{\left(H_k^{(1)}\right)^2}{k^4} - 18\sum_{k=1}^{\infty} \frac{\left(H_k^{(1)}\right)^2 H_k^{(2)}}{k^2} + 12\sum_{k=1}^{\infty} \frac{H_k^{(1)} H_k^{(3)}}{k^2} - 36\sum_{k=1}^{\infty} \frac{H_k^{(1)}}{k^5}$$

$$\psi^{(5)}(1) = 3\sum_{k=1}^{\infty} \frac{1}{k^2}\left\{\left(H_{k-1}^{(1)}\right)^2 - H_{k-1}^{(2)}\right\}\left\{\left(H_k^{(1)}\right)^2 + H_k^{(2)}\right\}$$

$$= 3\sum_{k=1}^{\infty} \frac{\left\{\left(H_k^{(1)}\right)^2 - 2\frac{H_k^{(1)}}{k} - H_k^{(2)} + 2\frac{1}{k^2}\right\}\left\{\left(H_k^{(1)}\right)^2 + H_k^{(2)}\right\}}{k^2}$$

$$= 3\sum_{k=1}^{\infty} \frac{\left(H_k^{(1)}\right)^4}{k^2} - 6\sum_{k=1}^{\infty} \frac{\left(H_k^{(1)}\right)^3}{k^3} - 3\sum_{k=1}^{\infty} \frac{\left(H_k^{(1)}\right)^2 H_k^{(2)}}{k^2} + 6\sum_{k=1}^{\infty} \frac{\left(H_k^{(1)}\right)^2}{k^4}$$

(4.3.65ii)

$$+ 3\sum_{k=1}^{\infty} \frac{\left(H_k^{(1)}\right)^2 H_k^{(2)}}{k^2} - 6\sum_{k=1}^{\infty} \frac{H_k^{(1)} H_k^{(2)}}{k^3} - 3\sum_{k=1}^{\infty} \frac{\left(H_k^{(2)}\right)^2}{k^2} + 6\sum_{k=1}^{\infty} \frac{H_k^{(2)}}{k^4}$$

$$\psi^{(5)}(1) = -2\sum_{k=1}^{\infty} \frac{1}{k^2}\left(-2H_k^{(3)} + 3H_k^{(1)} H_k^{(2)} - \left[H_k^{(1)}\right]^3\right)H_{k-1}^{(1)}$$

(4.3.65iii) $\quad = 4\sum_{k=1}^{\infty} \frac{H_k^{(1)} H_k^{(3)}}{k^2} - 6\sum_{k=1}^{\infty} \frac{\left(H_k^{(1)}\right)^2 H_k^{(2)}}{k^2} + 2\sum_{k=1}^{\infty} \frac{\left(H_k^{(1)}\right)^4}{k^2}$

$$+ 4\sum_{k=1}^{\infty} \frac{H_k^{(3)}}{k^3} - 6\sum_{k=1}^{\infty} \frac{H_k^{(1)} H_k^{(2)}}{k^3} + 2\sum_{k=1}^{\infty} \frac{\left(H_k^{(1)}\right)^3}{k^3}$$

(4.3.65iv) $\quad \psi^{(5)}(1) = \sum_{k=1}^{\infty} \frac{1}{k^2}\left(6H_k^{(4)} + 8H_k^{(1)} H_k^{(3)} + 3\left[H_k^{(2)}\right]^2 + 6\left[H_k^{(1)}\right]^2 H_k^{(2)} + \left[H_k^{(1)}\right]^4\right)$

Using Adamchik's formula (3.105vi) we may obtain $s(n,5)$ and hence another linear equation from (4.3.61v). It should also be noted that Coffey [45b], using the contour integral approach of Flajolet and Salvy [69], has previously computed that

$$-20\varsigma^2(3) = 5\sum_{k=1}^{\infty} \frac{\left[H_k^{(1)}\right]^4}{k^2} - 20\sum_{k=1}^{\infty} \frac{\left[H_k^{(1)}\right]^3}{k^3} + 30\sum_{k=1}^{\infty} \frac{\left[H_k^{(1)}\right]^2}{k^4} - 20\sum_{k=1}^{\infty} \frac{H_k^{(1)}}{k^5} - 30\sum_{k=1}^{\infty} \frac{H_k^{(2)} H_k^{(2)}}{k^2}$$

$$+40\sum_{k=1}^{\infty}\frac{H_k^{(1)}H_k^{(2)}}{k^3}+20\sum_{k=1}^{\infty}\frac{H_k^{(1)}H_k^{(3)}}{k^2}-15\sum_{k=1}^{\infty}\frac{H_k^{(2)}}{k^4}+10\sum_{k=1}^{\infty}\frac{\left[H_k^{(2)}\right]^2}{k^2}-10\sum_{k=1}^{\infty}\frac{H_k^{(3)}}{k^3}-5\sum_{k=1}^{\infty}\frac{H_k^{(4)}}{k^2}$$

$$-30\varsigma(2)\sum_{k=1}^{\infty}\frac{\left[H_k^{(1)}\right]^2}{k^2}+40\varsigma(2)\sum_{k=1}^{\infty}\frac{H_k^{(1)}}{k^3}+20\varsigma(2)\sum_{k=1}^{\infty}\frac{H_k^{(2)}}{k^2}-20\varsigma(3)\sum_{k=1}^{\infty}\frac{H_k^{(1)}}{k^2}$$

Some of the other unknowns may also be obtained from [69, p.25] and also from Zheng's recent paper [142b].

## GAUSS HYPERGEOMETRIC SUMMATION

Instead of using Nörlund's identity (4.3.32) we may employ the well-known Gauss hypergeometric summation formula [8a, p.64] (a reference to Gauss is appropriate because Gauss died on 23 February 1855, exactly 98 years before the birth of the author of this paper!)

$$(4.3.66) \qquad \sum_{n=0}^{\infty}\frac{(a)_n(b)_n}{n!(c)_n}=\frac{\Gamma(c)\Gamma(c-a-b)}{\Gamma(c-a)\Gamma(c-b)} \qquad \mathrm{Re}\,(c-a-b)>0$$

The familiar hypergeometric series is

$$_2F_1(a,b;c;z)=\sum_{n=0}^{\infty}\frac{(a)_n(b)_n}{(c)_n}\frac{z^n}{n!}$$

By definition $(a)_n=a(a+1)...(a+n-1)=(-1)^n(-a)_n$ for $n\geq 1$ and $(a)_0=1$; therefore letting $a\to -x$ and $b\to -y$ we get

$$(4.3.66a) \qquad \sum_{n=0}^{\infty}\frac{s_n(x)s_n(y)}{n!(c)_n}=\frac{\Gamma(c)\Gamma(c+x+y)}{\Gamma(c+x)\Gamma(c+y)} \qquad \mathrm{Re}\,(c+x+y)>0$$

where $s_n(x)$ is defined by (4.3.34a). For convenience we designate

$$S(x,y;c)=\frac{\Gamma(c)\Gamma(c+x+y)}{\Gamma(c+x)\Gamma(c+y)}$$

With differentiation we obtain

$$\frac{\partial}{\partial x}S(x,y;c) = S(x,y;c)\big[\psi(c+x+y)-\psi(c+x)\big]$$

$$\frac{\partial}{\partial y}\frac{\partial}{\partial x}S(x,y;c) =$$

$$S(x,y;c)\Big(\psi'(c+x+y)+\big[\psi(c+x+y)-\psi(c+y)\big]\big[\psi(c+x+y)-\psi(c+x)\big]\Big)$$

and therefore we have

$$\frac{\partial}{\partial y}\frac{\partial}{\partial x}S(x,y;c)\bigg|_{(0,0;c)} = \psi'(c)$$

We have seen before at (4.3.35) that $s_n'(0) = s(n,1) = (-1)^{n+1}(n-1)!$ and hence we obtain for $x=y=0$ and $c=1$

$$\sum_{n=1}^{\infty}\frac{1}{n^2} = \psi'(1) = \varsigma(2)$$

We have

$$\frac{\partial^2}{\partial x^2}S(x,y;c) = S(x,y;c)\big[\psi'(c+x+y)-\psi'(c+x)\big]+S(x,y;c)\big[\psi(c+x+y)-\psi(c+x)\big]^2$$

$$\frac{\partial}{\partial y}\frac{\partial^2}{\partial x^2}S(x,y;c) =$$

$$S(x,y;c)\begin{cases}\psi''(c+x+y)+\big[\psi'(c+x+y)-\psi'(c+x)\big]\big[\psi(c+x+y)-\psi(c+y)\big]\\[2mm] +2\big[\psi(c+x+y)-\psi(c+x)\big]\psi'(c+x+y)\\[2mm] +\big[\psi(c+x+y)-\psi(c+x)\big]^2\big[\psi(c+x+y)-\psi(c+y)\big]\end{cases}$$

and hence $\left. \dfrac{\partial}{\partial y}\dfrac{\partial^2}{\partial x^2}S(x,y;c)\right|_{(0,0;c)}=\psi''(c)$. This then gives us (4.3.40) when $c=1$. A further differentiation results in

$$\dfrac{\partial^2}{\partial y^2}\dfrac{\partial^2}{\partial x^2}S(x,y;c)=$$

$$S(x,y;c)\left\{\begin{array}{l}\psi'''(c+x+y)+\big[\psi'(c+x+y)-\psi'(c+x)\big]\big[\psi'(c+x+y)-\psi'(c+y)\big]\\[1.5em]+\psi''(c+x+y)\big[\psi(c+x+y)-\psi(c+y)\big]\\[1.5em]+2\psi''(c+x+y)\big[\psi(c+x+y)-\psi(c+x)\big]+2\big[\psi'(c+x+y)\big]^2\\[1.5em]+\big[\psi(c+x+y)-\psi(c+x)\big]^2\big[\psi'(c+x+y)-\psi'(c+y)\big]\\[1.5em]+2\big[\psi(c+x+y)-\psi(c+x)\big]\big[\psi(c+x+y)-\psi(c+y)\big]\psi'(c+x+y)\end{array}\right\}$$

$$+S(x,y;c)\left\{\begin{array}{l}\psi''(c+x+y)\big[\psi(c+x+y)-\psi(c+y)\big]\\[1.5em]+\big[\psi'(c+x+y)-\psi'(c+x)\big]\big[\psi(c+x+y)-\psi(c+y)\big]^2\\[1.5em]+2\big[\psi(c+x+y)-\psi(c+x)\big]\big[\psi(c+x+y)-\psi(c+y)\big]\psi'(c+x+y)\\[1.5em]+\big[\psi(c+x+y)-\psi(c+x)\big]^2\big[\psi(c+x+y)-\psi(c+y)\big]^2\end{array}\right\}$$

With $x=y=0$ we get

$$\left.\dfrac{\partial^2}{\partial y^2}\dfrac{\partial^2}{\partial x^2}S(x,y;c)\right|_{(0,0,c)}=\psi'''(c)+2\big[\psi'(c)\big]^2$$

and with $c=1$ we obtain

$$\left.\dfrac{\partial^2}{\partial y^2}\dfrac{\partial^2}{\partial x^2}S(x,y;c)\right|_{(0,0,1)}=6\varsigma(4)+2\varsigma^2(2)$$

Since $s_n''(0)=2s(n,2)=2(-1)^n(n-1)!H_{n-1}^{(1)}$ we deduce that

$$4\sum_{n=1}^{\infty}\frac{\left[H_{n-1}^{(1)}\right]^2}{n^2}=6\varsigma(4)+2\varsigma^2(2)$$

We see that $\displaystyle\sum_{n=1}^{\infty}\frac{\left[H_{n-1}^{(1)}\right]^2}{n^2}=\sum_{n=1}^{\infty}\frac{\left[H_{n}^{(1)}\right]^2}{n^2}-2\sum_{n=1}^{\infty}\frac{H_{n}^{(1)}}{n^3}+\varsigma(4)$ and hence we deduce

$$\sum_{n=1}^{\infty}\frac{\left[H_{n}^{(1)}\right]^2}{n^2}=\frac{17}{4}\varsigma(4)$$

which we previously obtained in (4.3.46b).

It is easily seen from (4.3.41) that $\dfrac{\partial^2}{\partial y^2}\dfrac{\partial^3}{\partial x^3}S(x,y;c)$ will result in the following series

$\displaystyle\sum_{n=1}^{\infty}\frac{H_{n-1}^{(1)}\left\{\left(H_{n-1}^{(1)}\right)^2-H_{n-1}^{(2)}\right\}}{n^2}$ containing the product of the terms $2(-1)^n(n-1)!H_{n-1}^{(1)}$ and

$(-1)^{n+1}\dfrac{(n-1)!}{2}\left\{\left(H_{n-1}^{(1)}\right)^2-H_{n-1}^{(2)}\right\}$. The procedure may then be extended by taking higher order derivatives. It remains to be checked whether the equations thus obtained are linearly independent from the results previously obtained using Nörlund's identity.

Note that we may also differentiate with respect to $c$ to obtain for example

$$\frac{\partial}{\partial c}\frac{\partial}{\partial y}\frac{\partial}{\partial x}S(x,y;c)=$$

$$S(x,y;c)\begin{pmatrix}\psi''(c+x+y)+\left[\psi(c+x+y)-\psi(c+y)\right]\left[\psi'(c+x+y)-\psi'(c+x)\right]\\[2mm]+\left[\psi'(c+x+y)-\psi'(c+y)\right]\left[\psi(c+x+y)-\psi(c+x)\right]\end{pmatrix}$$

$$+S(x,y;c)\left(\psi'(c+x+y)+\left[\psi(c+x+y)-\psi(c+y)\right]\left[\psi(c+x+y)-\psi(c+x)\right]\right)\times$$

$$\left(\psi(c+x+y)+\psi(c)-\psi(c+x)-\psi(c+y)\right)$$

STIRLING NUMBERS REVISITED

Let $f(x) = \log\Gamma(n+x) - \log\Gamma(1+x)$. Then we have

$$f^{(k)}(0) = \psi^{(k-1)}(n) - \psi^{(k-1)}(1) = (-1)^{k+1}(k-1)!H_{n-1}^{(k)}$$

and we therefore have the Maclaurin expansion

(4.3.66b)   $f(x) = \log\Gamma(n+x) - \log\Gamma(1+x) = \log\Gamma(n) + \sum_{k=1}^{\infty}\frac{(-1)^{k+1}}{k}H_{n-1}^{(k)}x^k$

(it may be worthwhile integrating this).This may be written as

(4.3.66c)   $\dfrac{\Gamma(n+x)}{\Gamma(1+x)\Gamma(n)} = \exp\left[\sum_{k=1}^{\infty}\frac{(-1)^{k+1}}{k}H_{n-1}^{(k)}x^k\right]$

and expansion of the exponential function leads to

(4.3.66d)

$$\frac{\Gamma(n+x)}{\Gamma(1+x)\Gamma(n)} = 1 + H_{n-1}^{(1)}x + \frac{1}{2}\left(\left[H_{n-1}^{(1)}\right]^2 - H_{n-1}^{(2)}\right)x^2 + \frac{1}{6}\left(\left[H_{n-1}^{(1)}\right]^3 - 3H_{n-1}^{(1)}H_{n-1}^{(2)} + 2H_{n-1}^{(3)}\right)x^3 + O(x^4)$$

and

(4.3.66e)

$$\frac{\Gamma(1+x)\Gamma(n)}{\Gamma(n+x)} = 1 - H_{n-1}^{(1)}x + \frac{1}{2}\left(\left[H_{n-1}^{(1)}\right]^2 + H_{n-1}^{(2)}\right)x^2 - \frac{1}{6}\left(\left[H_{n-1}^{(1)}\right]^3 + 3H_{n-1}^{(1)}H_{n-1}^{(2)} + 2H_{n-1}^{(3)}\right)x^3 + O(x^4)$$

The terms involving the generalised harmonic numbers are ubiquitous! I came across (4.3.66c) in a 1999 paper entitled "Analytic two-loop results for self energy- and vertex-type diagrams with one non-zero mass" by Fleisher et al [69aa].

Letting $n = 2$ in (4.3.66b) results in the familiar Maclaurin expansion for $\log(1+x)$.

It is easily seen that

$$\Gamma(n+x) = x(x+1)...(x+n-1)\Gamma(x)$$

and therefore we have

$$\log \Gamma(n+x) - \log \Gamma(1+x) = \log[(x+1)...(x+n-1)]$$

We see that

$$\sum_{k=1}^{n-1} \log(x+k) = \sum_{k=1}^{n-1} \log\left(1+\frac{x}{k}\right) + \sum_{k=1}^{n-1} \log k = \sum_{k=1}^{n-1} \log\left(1+\frac{x}{k}\right) + \log \Gamma(n)$$

and hence we have

$$\log \Gamma(n+x) - \log \Gamma(1+x) - \log \Gamma(n) = \sum_{k=1}^{n-1} \log\left(1+\frac{x}{k}\right)$$

Using (4.3.80a) this becomes

$$= \sum_{k=1}^{n-1} \sum_{j=1}^{m-1} \frac{(-1)^{j+1}}{j}\left(\frac{x}{k}\right)^j + \sum_{k=1}^{n-1} \frac{(-1)^{m-1}}{k^{m-1}} \int_0^x \frac{t^{m-1}}{k+t} dt$$

$$= \sum_{k=1}^{n-1} \frac{1}{k^j} \sum_{j=1}^{m-1} \frac{(-1)^{j+1} x^j}{j} + \sum_{k=1}^{n-1} \frac{(-1)^{m-1}}{k^{m-1}} \int_0^x \frac{t^{m-1}}{k+t} dt$$

$$= \sum_{j=1}^{m-1} \frac{(-1)^{j+1}}{j} H_{n-1}^{(j)} x^j + \sum_{k=1}^{n-1} \frac{(-1)^{m-1}}{k^{m-1}} \int_0^x \frac{t^{m-1}}{k+t} dt$$

As $m \to \infty$ we again see that

$$\log \Gamma(n+x) - \log \Gamma(1+x) - \log \Gamma(n) = \sum_{j=1}^{\infty} \frac{(-1)^{j+1}}{j} H_{n-1}^{(j)} x^j$$

We also have a form equivalent to (4.3.66c)

$$\frac{\Gamma(n+x)}{\Gamma(1+x)\Gamma(n)} = \exp\left[\sum_{j=1}^{m-1} \frac{(-1)^{j+1}}{j} H_{n-1}^{(j)} x^j\right] \exp\left[\sum_{k=1}^{n-1} \frac{(-1)^{m-1}}{k^{m-1}} \int_0^x \frac{t^{m-1}}{k+t} dt\right]$$

We have

$$\frac{\Gamma(1+x)\Gamma(n)}{\Gamma(n+x)} = \frac{x\Gamma(x)\Gamma(n)}{\Gamma(n+x)} = xB(x,n)$$

and from (E.35) in Volume VI we have

$$\log B(x,n) = \log\left(\frac{x+n}{xn}\right) + \int_0^1 \frac{(1-v^x)(1-v^n)}{(1-v)\log v}\,dv \qquad x,n > 0$$

Therefore we have

(4.3.66f)
$$-\sum_{k=1}^{\infty} \frac{(-1)^{k+1}}{k} H_{n-1}^{(k)} x^k = \log\left(1+\frac{x}{n}\right) + \int_0^1 \frac{(1-v^x)(1-v^n)}{(1-v)\log v}\,dv$$

An alternative proof is shown below.

We see that $1-v^x = 1-e^{x\log v} = -\sum_{k=1}^{\infty} \dfrac{x^k \log^k v}{k!}$ and therefore we have

$$\int_0^1 \frac{(1-v^x)(1-v^n)}{(1-v)\log v}\,dv = -\sum_{k=1}^{\infty} \frac{x^k}{k!} \int_0^1 \frac{(1-v^n)\log^{k-1} v}{1-v}\,dv$$

$$= -\sum_{k=1}^{\infty} \frac{x^k}{k!} \sum_{j=1}^{n-1} \int_0^1 v^j \log^{k-1} v\,dv$$

We have $\int_0^1 v^j v^a\,dv = \dfrac{1}{j+a+1}$ and therefore differentiation with respect to $a$ results in

$$\int_0^1 v^j v^a \log^{k-1} v\,dv = \frac{(-1)^{k-1}(k-1)!}{(j+a+1)^k}.$$

Hence we have

$$\int_0^1 v^j \log^{k-1} v\,dv = \frac{(-1)^{k-1}(k-1)!}{(j+1)^k}$$

$$\sum_{j=0}^{n-1} \int_0^1 v^j \log^{k-1} v\,dv = (-1)^{k-1}(k-1)!H_n^{(k)}$$

and thus

(4.3.66fi)
$$\int_0^1 \frac{(1-v^x)(1-v^n)}{(1-v)\log v}\,dv = -\sum_{k=1}^{\infty} \frac{(-1)^{k+1}}{k} H_n^{(k)} x^k$$

$$= -\sum_{k=1}^{\infty} \frac{(-1)^{k+1}}{k} \left( H_{n-1}^{(k)} + \frac{1}{n^k} \right) x^k$$

$$= -\sum_{k=1}^{\infty} \frac{(-1)^{k+1}}{k} H_{n-1}^{(k)} x^k - \sum_{k=1}^{\infty} \frac{(-1)^{k+1}}{k} \frac{x^k}{n^k}$$

$$= -\sum_{k=1}^{\infty} \frac{(-1)^{k+1}}{k} H_{n-1}^{(k)} x^k - \log\left( 1 + \frac{x}{n} \right)$$

Differentiating (4.3.66fi) results in

$$\int_0^1 \frac{v^x(1-v^n)}{1-v} dv = \sum_{k=1}^{\infty} (-1)^{k+1} H_n^{(k)} x^{k-1}$$

and as $x \to 0$ we rediscover (4.1.4)

$$\int_0^1 \frac{(1-v^n)}{1-v} dv = H_n^{(1)}$$

Successive differentiations give us

$$\int_0^1 \frac{v^x(1-v^n)\log^p v}{1-v} dv = \sum_{k=p+1}^{\infty} (-1)^{k+1} H_n^{(k)} (k-1)(k-2)...(k-p) x^{k-p-1}$$

and hence we see that

(4.3.66fii)     $$\int_0^1 \frac{(1-v^n)\log^p v}{1-v} dv = (-1)^p p! H_n^{(p+1)}$$

As $n \to \infty$ we see that

$$\int_0^1 \frac{\log^p v}{1-v} dv = (-1)^p p! \varsigma(p+1)$$

See also (E.35c) in Volume VI.

Completing the summation of (4.3.66fii) gives us for $|x| < 1$

$$(-1)^p p! \sum_{n=1}^{\infty} \frac{H_n^{(p+1)}}{n} x^n = \sum_{n=1}^{\infty} \frac{1}{n} x^n \int_0^1 \frac{(1-v^n)\log^p v}{1-v} dv$$

$$= \sum_{n=1}^{\infty} \int_0^1 u^{n-1} du \int_0^1 \frac{(1-v^n)x^n \log^p v}{1-v} dv$$

We therefore obtain

(4.3.66fiii) $\qquad (-1)^p p! \sum_{n=1}^{\infty} \frac{H_n^{(p+1)}}{n} x^n = \int_0^1 \frac{[\log(1-xv) - \log(1-x)] \log^p v}{(1-v)} dv$

With $p = 1$ we have

$$\sum_{n=1}^{\infty} \frac{H_n^{(2)}}{n} x^n = -\int_0^1 \frac{[\log(1-xv) - \log(1-x)] \log v}{(1-v)} dv$$

We note from (3.106f) that

$$\sum_{n=1}^{\infty} \frac{H_n^{(2)}}{n} x^n = -2 Li_3\left(-\frac{x}{1-x}\right)$$

Integrating (4.3.66fiii) will give us an expression for $\displaystyle \sum_{n=1}^{\infty} \frac{H_n^{(p+1)}}{n^2} x^n$.

$$(-1)^p p! \sum_{n=1}^{\infty} \frac{H_n^{(p+1)}}{n^2} u^n = \int_0^u \frac{dx}{x} \int_0^1 \frac{[\log(1-xv) - \log(1-x)] \log^p v}{1-v} dv$$

We then have

$$\int_0^u \frac{dx}{x} \int_0^1 \frac{[\log(1-xv) - \log(1-x)] \log^p v}{1-v} dv = \int_0^1 \frac{[Li_2(u) - Li_2(uv)] \log^p v}{1-v} dv$$

and therefore we obtain

(4.3.66fiv) $\qquad (-1)^p p! \sum_{n=1}^{\infty} \frac{H_n^{(p+1)}}{n^2} u^n = \int_0^1 \frac{[Li_2(u) - Li_2(uv)] \log^p v}{1-v} dv$

With $u = 1$ we have

$$(-1)^p p! \sum_{n=1}^{\infty} \frac{H_n^{(p+1)}}{n^2} = \int_0^1 \frac{[\varsigma(2) - Li_2(v)] \log^p v}{1-v} dv$$

and with $u = -1$ we have since $Li_s(-1) = (2^{1-s} - 1)\varsigma(s)$

$$(-1)^p \, p! \sum_{n=1}^{\infty} (-1)^n \frac{H_n^{(p+1)}}{n^2} = -\int_0^1 \frac{\left[\frac{1}{2}\varsigma(2) + Li_2(-v)\right]\log^p v}{1-v} dv$$

This is easily extended by further integration to

(4.3.66fv) $\qquad (-1)^p \, p! \sum_{n=1}^{\infty} \frac{H_n^{(p+1)}}{n^r} u^n = \int_0^1 \frac{\left[Li_r(u) - Li_r(uv)\right]\log^p v}{1-v} dv$

and with $p = 0$ we have

(4.3.66fvi) $\qquad \sum_{n=1}^{\infty} \frac{H_n^{(1)}}{n^r} u^n = \int_0^1 \frac{\left[Li_r(u) - Li_r(uv)\right]}{1-v} dv$

In particular, we see that

(4.3.66fvii) $\qquad \sum_{n=1}^{\infty} \frac{H_n^{(1)}}{n^2} = \int_0^1 \frac{\left[\varsigma(2) - Li_2(v)\right]}{1-v} dv$

and therefore we have using (4.2.33)

(4.3.66fviii) $\qquad 2\varsigma(3) = \int_0^1 \frac{\left[\varsigma(2) - Li_2(v)\right]}{1-v} dv$

$\square$

From (4.3.66d) and using the formulae (3.105i) for the Stirling numbers of the first kind we see that

$$(-1)^{n+1} \frac{\Gamma(n+x)}{\Gamma(1+x)} = s(n,1) - s(n,2)x + s(n,3)x^2 - s(n,4)x^3 + O(x^4)$$

and this may be written as

(4.3.66g) $\qquad (x)_n = \frac{\Gamma(n+x)}{\Gamma(x)} = \frac{x\,\Gamma(n+x)}{\Gamma(1+x)} = \sum_{k=0}^{\infty} (-1)^{n+k} s(n,k)\, x^k$

which is simply the generating function for the Stirling numbers of the first kind [26, p.56]

$$(x)_n = x(x+1)...(x+n-1) = \sum_{k=0}^{\infty} (-1)^{n+k} s(n,k) x^k = \sum_{k=1}^{\infty} (-1)^{n+k} s(n,k) x^k$$

since $s(n,k) = 0$ for $k \geq n+1$. Using (4.3.66g) is probably the simplest way of evaluating the Stirling numbers (by successive differentiation); for example

(4.3.66h)
$$\frac{\Gamma(n+x)}{\Gamma(1+x)} \big[ x\psi(n+x) + 1 - x\psi(1+x) \big] = \sum_{k=1}^{\infty} (-1)^{n+k} k\, s(n,k) x^{k-1}$$

and with $x = 0$ we obtain $s(n,1) = (-1)^{n+1}(n-1)!$.

A further differentiation results in

$$\frac{\Gamma(n+x)}{\Gamma(1+x)} \big[ \psi(n+x) - \psi(1+x) \big] \big[ x\psi(n+x) + 1 - x\psi(1+x) \big]$$

$$+ \frac{\Gamma(n+x)}{\Gamma(1+x)} \Big( x\big[ \psi'(n+x) - \psi'(1+x) \big] + \big[ \psi(n+x) - \psi(1+x) \big] \Big) = \sum_{k=2}^{\infty} (-1)^{n+k} k(k-1) s(n,k) x^{k-2}$$

and with $x = 0$ we easily obtain $s(n,2) = (-1)^n H_{n-1}^{(1)}$.

The above equation may be written as

$$x L(x) \Delta^2(x) + 2 L(x) \Delta(x) + x L(x) \Delta^{(1)}(x) = \sum_{k=2}^{\infty} (-1)^{n+k} k(k-1) s(n,k) x^{k-2}$$

where for convenience

$$L(x) = \frac{\Gamma(n+x)}{\Gamma(1+x)}$$

$$\Delta(x) = \psi(n+x) - \psi(1+x)$$

and $\qquad L'(x) = L(x)\Delta(x)$.

The next derivative is then easily determined using this notation and we note the connection with (4.3.31).

It is easily seen that

$$x L''(x) + 2L'(x) = x L(x)\Delta^2(x) + 2L(x)\Delta(x) + x L(x)\Delta^{(1)}(x)$$

and this format is eminently suitable for Leibniz's rule for differentiation. It may also be noted that

$$e^{-x}\frac{d}{dx}\left(xe^{x}L(x)\right)=x\,L''(x)+2L'(x)$$

Hence we have

$$\frac{d}{dx}\left(xe^{x}L(x)\right)=\sum_{k=2}^{\infty}(-1)^{n+k}k(k-1)s(n,k)\,x^{k-2}e^{x}$$

Upon integration we obtain

$$te^{t}L(t)=\sum_{k=2}^{\infty}(-1)^{n+k}k(k-1)s(n,k)\int_{0}^{t}x^{k-2}e^{x}dx$$

We have from G&R [74, p.104]

$$\int x^{p}e^{x}dx=e^{x}\left(x^{p}+\sum_{j=1}^{p}(-1)^{j}\,p(p-1)...(p-j+1)x^{p-j}\right)$$

and we may then obtain a formula for $L(x)=\dfrac{\Gamma(n+x)}{\Gamma(1+x)}$ in the form of an infinite series involving $s(n,k)$ (cf (4.3.66c)).  □

From (4.4.1) we will see that

(4.3.67)  $$\frac{1}{a(a+1)...(a+k-1)}=\frac{1}{(k-1)!}\int_{0}^{1}t^{a-1}(1-t)^{k-1}dt$$

and substituting this in (4.3.35) we get

(4.3.67a)  $$\psi'(a)=\sum_{k=1}^{\infty}\frac{1}{k}\int_{0}^{1}t^{a-1}(1-t)^{k-1}dt$$

$$=\int_{0}^{1}t^{a-1}\sum_{k=1}^{\infty}\frac{1}{k}(1-t)^{k-1}dt=\int_{0}^{1}\frac{t^{a-1}}{1-t}\sum_{k=1}^{\infty}\frac{(1-t)^{k}}{k}dt$$

Hence we have the well-known result [126, p.16]

(4.3.68)  $$\psi'(a)=-\int_{0}^{1}\frac{t^{a-1}}{1-t}\log t\,dt$$

alternative proofs of which are shown in (4.4.99k) and Appendix E of Volume VI.

Integration of (4.3.67a) gives us

$$\int\limits_{1}^{u}\psi'(a)\,da = \int\limits_{1}^{u}\sum_{k=1}^{\infty}\frac{1}{k}\int\limits_{0}^{1}t^{a-1}(1-t)^{k-1}\,dt\,da$$

and since

(4.3.68a) $$\int\limits_{1}^{u}t^{a-1}\,da = \frac{t^{u-1}-1}{\log t}$$

we obtain (assuming that the interchange of summation and integration is valid)

(4.3.68b) $$\psi(u)-\psi(1) = \sum_{k=1}^{\infty}\frac{1}{k}\int\limits_{0}^{1}\frac{t^{u-1}-1}{\log t}(1-t)^{k-1}\,dt = \int\limits_{0}^{1}\frac{1-t^{u-1}}{1-t}\,dt$$

where in the last part we have used $\log t = -\sum_{k=1}^{\infty}\frac{(1-t)^k}{k}$ (which is valid in the particular range of integration).

We shall see in (4.4.94a) that

(4.3.68c) $$\int\limits_{0}^{1}\frac{(1-y)^n}{\log y}\,dy = \sum_{j=0}^{n}\binom{n}{j}(-1)^j\log(1+j)$$

(4.3.68d) $$\int\limits_{0}^{1}\frac{y^{u-1}(1-y)^n}{\log^r y}\,dy = \frac{1}{(r-1)!}\sum_{j=0}^{n}\binom{n}{j}(-1)^j(u+j)^{r-1}\log(u+j)$$

and we therefore get

(4.3.69) $$\psi(u)-\psi(1) = \sum_{k=1}^{\infty}\frac{1}{k}\sum_{j=0}^{k-1}\binom{k-1}{j}(-1)^j\log\frac{(u+j)}{(1+j)}$$

With reference to (4.4.1) we also see that

$$\frac{1}{a(a+1)...(a+k-1)} = \frac{1}{(k-1)!}\sum_{j=0}^{k-1}\binom{k-1}{j}\frac{(-1)^j}{j+a}$$

and hence by reference to (4.3.35) we have the expression for $a > 0$ (which I suspect is in fact valid $\forall a \neq \mathbf{Z}_0^-$)

$$(4.3.70) \qquad \psi'(a) = \sum_{k=1}^{\infty} \frac{1}{k} \sum_{j=0}^{k-1} \binom{k-1}{j} \frac{(-1)^j}{j+a} = \sum_{k=0}^{\infty} \frac{1}{k+1} \sum_{j=0}^{k} \binom{k}{j} \frac{(-1)^j}{j+a}$$

Since $\psi'(1) = \varsigma(2)$ we obtain

$$\varsigma(2) = \sum_{k=0}^{\infty} \frac{1}{k+1} \sum_{j=0}^{k} \binom{k}{j} \frac{(-1)^j}{j+1}$$

Differentiation of (4.3.70) results in the identity

$$(4.3.70a) \qquad \psi^{(p+1)}(a) = (-1)^p (p)! \sum_{k=0}^{\infty} \frac{1}{k+1} \sum_{j=0}^{k} \binom{k}{j} \frac{(-1)^j}{(j+a)^{p+1}}$$

and with $a = 1$ we obtain

$$\psi^{(p+1)}(1) = (-1)^p (p)! \sum_{k=0}^{\infty} \frac{1}{k+1} \sum_{j=0}^{k} \binom{k}{j} \frac{(-1)^j}{(j+1)^{p+1}}$$

Therefore, using (4.3.33a) we have

$$(4.3.71) \qquad \varsigma(p+2) = \frac{1}{p} \sum_{k=0}^{\infty} \frac{1}{k+1} \sum_{j=0}^{k} \binom{k}{j} \frac{(-1)^j}{(j+1)^{p+1}}$$

which is a subset of the Hasse identity (3.12) for $p \in \mathbf{Z}^+$.

We also note from (4.3.67) that for $a > 0$

$$\sum_{k=0}^{\infty} \frac{1}{a(a+1)...(a+k)} = \sum_{k=0}^{\infty} \int_0^1 t^{a-1} \frac{(1-t)^k}{k!} dt$$

and accordingly we have

$$(4.3.71a) \qquad \sum_{k=0}^{\infty} \frac{1}{a(a+1)...(a+k)} = \int_0^1 t^{a-1} e^{1-t} dt = e\Gamma(a) - e\Gamma(a,1)$$

where $\Gamma(a,x)$ is the incomplete gamma function defined by

$$\Gamma(a,x) = \int_x^\infty t^{a-1} e^{-t} dt$$

If $a$ is an integer then we have

$$\Gamma(n,x) = (n-1)! e^{-x} \sum_{k=0}^{n-1} \frac{x^k}{k!}$$

and hence $\Gamma(a,1) = (a-1)! e^{-1}$.

Further series may be obtained using (4.4.1). The result when $a = 1$ is very familiar and we have

$$\sum_{k=0}^\infty \frac{1}{(k+1)!} = e - 1$$

with $a = 1/2$ we get

$$(4.3.71b) \qquad \sum_{k=0}^\infty \frac{2^{k+1}}{1.3.5..(2k+1)} = e\sqrt{\pi} - e\Gamma\left(\frac{1}{2},1\right)$$

Differentiating (4.3.71a) gives us

$$\sum_{k=0}^\infty \frac{1}{k!} \sum_{j=0}^k \binom{k}{j} \frac{(-1)^j}{(j+a)^2} = -e\Gamma'(a) + e\Gamma'(a,1)$$

where we have used (4.4.1). With $a = 1$ we get

$$\sum_{k=0}^\infty \frac{1}{k!} \sum_{j=0}^k \binom{k}{j} \frac{(-1)^j}{(j+1)^2} = e\gamma + e\Gamma'(1,1)$$

and from (4.2.16) we recall that $\sum_{j=0}^k \binom{k}{j} \frac{(-1)^j}{(j+1)^2} = \frac{H_{k+1}}{k+1}$ which gives us

$$(4.3.71c) \qquad \sum_{k=0}^\infty \frac{H_{k+1}}{(k+1)!} = \sum_{k=1}^\infty \frac{H_k}{k!} e\gamma + e\Gamma'(1,1)$$

The Mathworld website for harmonic numbers reports Gosper's formula

$$\sum_{k=0}^\infty \frac{H_k}{k!} x^k = -e^x \sum_{k=1}^\infty \frac{(-1)^k}{k.k!} x^k = e^x [\log x + \Gamma(0,x) + \gamma]$$

and a derivation is given below.

Using the known integral expression for the harmonic numbers $H_n$ (see for example Volume I)

$$H_k = -k \int_0^1 (1-u)^{k-1} \log u \, du$$

we obtain

$$\sum_{k=0}^{\infty} \frac{H_k}{k!} x^k = -\sum_{k=0}^{\infty} \frac{k \int_0^1 (1-u)^{k-1} \log u \, du}{k!} x^k$$

$$= -\int_0^1 \sum_{k=0}^{\infty} \frac{k}{k!(1-u)} x^k (1-u)^k \log u \, du$$

$$= -\int_0^1 \frac{x(1-u)e^{(1-u)x} \log u \, du}{(1-u)}$$

$$= -x e^x \int_0^1 e^{-ux} \log u \, du$$

We note that

$$\Gamma'(1) = -\gamma = \int_0^{\infty} e^{-u} \log u \, du = x \int_0^{\infty} e^{-ux} \log(ux) du = x \int_0^{\infty} e^{-ux} \log u du + x \log x \int_0^{\infty} e^{-ux} du$$

$$= x \int_0^{\infty} e^{-ux} \log u du + \log x$$

and we therefore have

$$x \int_0^{\infty} e^{-ux} \log u du = -\log x - \gamma$$

It is seen that

$$\int\limits_0^\infty e^{-ux}\log u\,du = \int\limits_0^1 e^{-ux}\log u\,du + \int\limits_1^\infty e^{-ux}\log u\,du$$

and integration by parts gives us

$$\int\limits_1^\infty e^{-ux}\log u\,du = \frac{1}{x}\left[Ei(-ux)-e^{-ux}\log u\right]\Big|_1^\infty = -\frac{1}{x}Ei(-x)$$

Noting that $\Gamma(0,x) = -Ei(-x)$ for $x>0$ completes the proof.

The formula $\sum\limits_{k=0}^\infty \dfrac{H_k}{k!}x^k = -e^x\sum\limits_{k=1}^\infty \dfrac{(-1)^k}{k.k!}x^k$ is also contained in Ramanujan's Notebooks

(see Berndt [21], Part I, p.46).

Havil shows a proof of $-\sum\limits_{k=1}^\infty \dfrac{(-1)^k}{k.k!}x^k = \log x + \Gamma(0,x) + \gamma$ in [78, p.108].

For $x=1$ we get

$$\sum\limits_{k=0}^\infty \frac{H_k}{k!} = -e\sum\limits_{k=1}^\infty \frac{(-1)^k}{k.k!} = e[\Gamma(0,1)+\gamma]$$

□

We have from (4.3.39)

$$\psi''(a) = -2\sum\limits_{k=1}^\infty \frac{1}{k}\frac{(k-1)!H_{k-1}^{(1)}}{a(a+1)...(a+k-1)}$$

and hence using (4.3.67) we see that

$$\psi''(a) = -2\sum\limits_{k=1}^\infty \frac{H_{k-1}^{(1)}}{k}\int\limits_0^1 t^{a-1}(1-t)^{k-1}dt$$

Using (4.1.6) we then have

$$\psi''(a) = -2\sum\limits_{k=1}^\infty \frac{1}{k}\int\limits_0^1 \frac{1-(1-x)^{k-1}}{x}dx\int\limits_0^1 t^{a-1}(1-t)^{k-1}dt$$

$$= -2\int\limits_0^1\int\limits_0^1 \sum\limits_{k=1}^\infty \frac{t^{a-1}(1-t)^{k-1}-t^{a-1}\left[(1-x)(1-t)\right]^{k-1}}{kx}dxdt$$

Completing the summation we obtain

$$(4.3.71d) \qquad \psi''(a) = 2\int_0^1\int_0^1 \left[\frac{t^{a-1}\log t}{x(1-t)} - \frac{t^{a-1}\log\left[1-(1-x)(1-t)\right]}{x(1-x)(1-t)}\right] dx dt$$

Since $\psi''(1) = -2\varsigma(3)$ we have

$$(4.3.71e) \qquad \varsigma(3) = 2\int_0^1\int_0^1 \left[\frac{\log\left[1-(1-x)(1-t)\right]}{x(1-x)(1-t)} - \frac{\log t}{x(1-t)}\right] dx dt$$

With the benefit of the Wolfram Integrator we have where $A = 1-t$

$$\int_0^1 \left[\frac{\log\left[1-A(1-x)\right]}{x(1-x)A} - \frac{\log t}{xA}\right] dx =$$

$$\frac{1}{A}\left[\begin{array}{l}\log\left[A(1-x)\right]\log\left[1-A(1-x)\right] + \log x\log(1-A) - Li_2\left(-\frac{Ax}{1-A}\right) \\ -Li_2\left[1-A(1-x)\right] - \log t\log x\end{array}\right]\Bigg|_0^1$$

It will be noted that the terms involving $\log x$ cancel out and thus

$$\int_0^1 \left[\frac{\log\left[1-A(1-x)\right]}{x(1-x)A} - \frac{\log t}{xA}\right] dx = -\frac{1}{A}\left(Li_2\left(-\frac{A}{1-A}\right) + \varsigma(2) + \log A\log\left[1-A\right] - Li_2\left[1-A\right]\right)$$

$$= -\frac{1}{1-t}\left(Li_2\left(-\frac{1-t}{t}\right) + \varsigma(2) + \log(1-t)\log t - Li_2(t)\right)$$

We then get

$$(4.3.71f) \qquad \varsigma(3) = -2\int_0^1 \frac{1}{1-t}\left(Li_2\left(-\frac{1-t}{t}\right) + \varsigma(2) + \log(1-t)\log t - Li_2(t)\right) dt$$

The Wolfram Integrator reports a fairly complex output for the above integral which may be worthy of further analysis.

We have

$$\int_0^1 \frac{\log(1-t)\log t}{1-t} dt = \log(1-x)Li_2(1-x) - Li_3(1-x)\Big|_0^1 = -\varsigma(3)$$

and we recall (4.3.66fviii)

$$2\varsigma(3) = \int_0^1 \frac{\left[\varsigma(2) - Li_2(t)\right]}{1-t} dt$$

This then gives us

$$3\varsigma(3) = -2\int_0^1 Li_2\left(-\frac{1-t}{t}\right)\frac{dt}{1-t}$$

More generally, upon differentiation of (4.4.71d), we obtain

$$(4.3.71g) \qquad \psi^{(p+2)}(a) = 2\int_0^1\int_0^1 \left[\frac{t^{a-1}\log^{p+1} t}{x(1-t)} - \frac{t^{a-1}\log^{p+1}\left[1-(1-x)(1-t)\right]}{x(1-x)(1-t)}\right] dx\, dt$$

## LOGARITHMIC SERIES FOR THE DIGAMMA AND GAMMA FUNCTIONS

We have from (4.3.69)

$$(4.3.72a) \qquad \psi(u) - \psi(1) = \sum_{k=0}^{\infty} \frac{1}{k+1}\sum_{j=0}^{k}\binom{k}{j}(-1)^j \log\frac{(u+j)}{(1+j)}$$

Also refer to (4.4.24g). Letting $u = n+1$ we see that

$$(4.3.73) \qquad H_n^{(1)} = \psi(n+1) - \psi(1) = \sum_{k=0}^{\infty} \frac{1}{k+1}\sum_{j=0}^{k}\binom{k}{j}(-1)^j \log\frac{(n+1+j)}{(1+j)}$$

We will see in (4.4.92a) that

$$(4.3.73a) \qquad \psi(1) = -\gamma = \sum_{k=0}^{\infty} \frac{1}{k+1}\sum_{j=0}^{k}\binom{k}{j}(-1)^j \log(1+j)$$

and we therefore obtain an infinite series for the digamma function for $u \notin \mathbf{Z}_0^-$ (see (4.4.99aix) of Volume III for a different derivation).

$$(4.3.74) \qquad \psi(u) = \sum_{k=0}^{\infty} \frac{1}{k+1}\sum_{j=0}^{k}\binom{k}{j}(-1)^j \log(u+j)$$

This result was also recently obtained in a different way by Guillera and Sondow [75aa]. It is interesting to note that (4.3.73) always produces a rational number whilst the closely related formula (4.3.73a) is believed to be irrational, and is probably transcendental.

With $u = 1/2$ we obtain

$$\psi(1/2) = \sum_{k=0}^{\infty} \frac{1}{k+1} \sum_{j=0}^{k} \binom{k}{j} (-1)^j \log(1+2j) - \log 2 \sum_{k=0}^{\infty} \frac{1}{k+1} \sum_{j=0}^{k} \binom{k}{j} (-1)^j$$

$$= \sum_{k=0}^{\infty} \frac{1}{k+1} \sum_{j=0}^{k} \binom{k}{j} (-1)^j \log(1+2j) - \log 2$$

Since [126, p.20] $\psi(1/2) = -\gamma - 2\log 2$, we obtain

(4.3.74a)
$$\sum_{k=0}^{\infty} \frac{1}{k+1} \sum_{j=0}^{k} \binom{k}{j} (-1)^j \log(1+2j) = -\gamma - \log 2$$

By an alternative method we will show in (4.4.24i) that

(4.3.74b)
$$\sum_{k=0}^{\infty} \frac{1}{k+1} \sum_{j=0}^{k} \binom{k}{j} (-1)^j \log(n+1+j) = H_n^{(1)} - \gamma$$

and therefore we see that

$$\sum_{k=0}^{\infty} \frac{1}{k+1} \sum_{j=0}^{k} \binom{k}{j} (-1)^j \log(2+j) = 1 - \gamma$$

This may also be obtained by noting that

$$\int_0^1 \psi(u+1) du = \log \Gamma(u+1) \big|_0^1 = 0$$

Another derivation of (4.3.74) is shown below. From the Hasse identity we see that

$$\varsigma(2,u) = \sum_{k=0}^{\infty} \frac{1}{k+1} \sum_{j=0}^{k} \binom{k}{j} \frac{(-1)^j}{u+j}$$

and we know that $\varsigma(2,u) = \psi'(u)$. Hence, upon integration, we have shown in a fairly elementary manner that

$$\psi(u) = \sum_{k=0}^{\infty} \frac{1}{k+1} \sum_{j=0}^{k} \binom{k}{j} (-1)^j \log(u+j) + c$$

It remains for us to determine the integration constant. We have

$$\sum_{j=0}^{k} \binom{k}{j} (-1)^j \log(u+j) = \log\left[ \prod_{j=0}^{k} (u+j)^{(-1)^j \binom{k}{j}} \right]$$

and as $u \to \infty$ we see that

$$\sum_{j=0}^{k} \binom{k}{j} (-1)^j \log(u+j) \approx \log\left[ u^{\sum_{j=0}^{k} \binom{k}{j}(-1)^j} \right] = \delta_{k,0} \log u$$

Therefore we have

$$\sum_{k=0}^{\infty} \frac{1}{k+1} \sum_{j=0}^{k} \binom{k}{j} (-1)^j \log(u+j) \to \log u \text{ as } u \to \infty$$

We know that (see, for example, (E.66a) in Volume VI) and (4.3.121aa) of this paper

$$\lim_{u \to \infty} [\psi(u) - \log u] = 0$$

and hence we may deduce that $c = 0$.

Since $\psi'(u) = \sum_{k=0}^{\infty} \frac{1}{(u+k)^2} > 0$, we note that $\sum_{k=0}^{\infty} \frac{1}{k+1} \sum_{j=0}^{k} \binom{k}{j} (-1)^j \log(u+j)$ is a strictly monotonic increasing function of $u$.

We may note from (4.3.74b) that

(4.3.74c) $\qquad \sum_{n=1}^{\infty} \frac{1}{n^s} \sum_{k=0}^{\infty} \frac{1}{k+1} \sum_{j=0}^{k} \binom{k}{j} (-1)^j \log(n+1+j) = \sum_{n=1}^{\infty} \frac{H_n^{(1)}}{n^s} - \gamma \varsigma(s)$

and with $s = 2$ we get

(4.3.74d) $\qquad \sum_{n=1}^{\infty} \frac{1}{n^2} \sum_{k=0}^{\infty} \frac{1}{k+1} \sum_{j=0}^{k} \binom{k}{j} (-1)^j \log(n+1+j) = 2\varsigma(3) - \gamma \varsigma(2)$

This is equivalent to the well-known result

$$\sum_{n=1}^{\infty}\frac{\psi(n+1)}{n^2}=\sum_{n=1}^{\infty}\frac{H_n^{(1)}-\gamma}{n^2}=2\varsigma(3)-\gamma\varsigma(2)$$

With Legendre's duplication formula [126, p.7]

$$\sqrt{\pi}\,\Gamma(2a)=2^{2a-1}\Gamma(a)\Gamma\left(a+\frac{1}{2}\right)$$

we see that

$$\frac{1}{2}\log\pi+\log\Gamma(2a)=(2a-1)\log 2+\log\Gamma(a)+\log\Gamma\left(a+\frac{1}{2}\right)$$

and differentiation results in

$$2\psi(2a)=2\log 2+\psi(a)+\psi\left(a+\frac{1}{2}\right)$$

Then using (4.3.74)

$$\psi(a)=\sum_{k=0}^{\infty}\frac{1}{k+1}\sum_{j=0}^{k}\binom{k}{j}(-1)^j\log(a+j)$$

we have

$$2\sum_{k=0}^{\infty}\frac{1}{k+1}\sum_{j=0}^{k}\binom{k}{j}(-1)^j\log(2a+j)=2\log 2+\sum_{k=0}^{\infty}\frac{1}{k+1}\sum_{j=0}^{k}\binom{k}{j}(-1)^j\log(a+j)$$

$$+\sum_{k=0}^{\infty}\frac{1}{k+1}\sum_{j=0}^{k}\binom{k}{j}(-1)^j\log\frac{(2a+1+2j)}{2}$$

or equivalently

$$2\sum_{k=0}^{\infty}\frac{1}{k+1}\sum_{j=0}^{k}\binom{k}{j}(-1)^j\log(2a+j)=\log 2+\sum_{k=0}^{\infty}\frac{1}{k+1}\sum_{j=0}^{k}\binom{k}{j}(-1)^j\log(a+j)$$

$$+\sum_{k=0}^{\infty}\frac{1}{k+1}\sum_{j=0}^{k}\binom{k}{j}(-1)^j\log(2a+1+2j)$$

This may be written as

$$\sum_{k=0}^{\infty} \frac{1}{k+1} \sum_{j=0}^{k} \binom{k}{j} (-1)^j \log \frac{(2a+j)^2}{(a+j)(2a+1+2j)} = \log 2$$

or as

$$\sum_{k=0}^{\infty} \frac{1}{k+1} \sum_{j=0}^{k} \binom{k}{j} (-1)^j \log \frac{(2a+j)^2}{2(a+j)(2a+1+2j)} = 0$$

As an aside, I noted from [16a] that

$$\sum_{j=k}^{\infty} \log \frac{(1+j)^2}{j(j+2)} = \log \left(1 + \frac{1}{k}\right)$$

and wondered if there was a connection via the Euler series acceleration technique.

Differentiation results in

$$4 \sum_{k=0}^{\infty} \frac{1}{k+1} \sum_{j=0}^{k} \binom{k}{j} \frac{(-1)^j}{(j+2a)} = \sum_{k=0}^{\infty} \frac{1}{k+1} \sum_{j=0}^{k} \binom{k}{j} \frac{(-1)^j}{(j+a)} + \sum_{k=0}^{\infty} \frac{1}{k+1} \sum_{j=0}^{k} \binom{k}{j} \frac{(-1)^j}{(j+(2a+1)/2)}$$

and more generally we have

$$2^s \sum_{k=0}^{\infty} \frac{1}{k+1} \sum_{j=0}^{k} \binom{k}{j} \frac{(-1)^j}{(j+2a)^{s-1}} = \sum_{k=0}^{\infty} \frac{1}{k+1} \sum_{j=0}^{k} \binom{k}{j} \frac{(-1)^j}{(j+a)^{s-1}} + \sum_{k=0}^{\infty} \frac{1}{k+1} \sum_{j=0}^{k} \binom{k}{j} \frac{(-1)^j}{(j+(2a+1)/2)^{s-1}}$$

and, using the Hasse formula, this is equivalent to

(4.3.74di)     $$2^s \varsigma(s, 2a) = \varsigma(s, a) + \varsigma\left(s, a + \frac{1}{2}\right)$$

With $a \rightarrow a/2$ this becomes

$$2^s \varsigma(s, a) = \varsigma\left(s, \frac{a}{2}\right) + \varsigma\left(s, \frac{a+1}{2}\right)$$

It easily follows from the definition of $\varsigma(s, a)$ that [126, p.89]

$$\varsigma\left(s, \frac{a}{2}\right) - \varsigma\left(s, \frac{a+1}{2}\right) = 2^s \varsigma_a(s, a)$$

and upon addition of the above two formula we obtain

$$2\varsigma\left(s,\frac{a}{2}\right)=2^s\left[\varsigma(s,a)+\varsigma_a(s,a)\right]$$

With $a=1$ we get

$$2\varsigma\left(s,\frac{1}{2}\right)=2^s\left[\varsigma(s)+\varsigma_a(s)\right]=2^s\left[\varsigma(s)+(1-2^{1-s})\varsigma(s)\right]$$

and this gives us the well-known formula

(4.3.74dii) $\qquad \varsigma\left(s,\frac{1}{2}\right)=\left[2^s-1\right]\varsigma(s)$

We also have from Nörlund's book [105, p.100]

$$(-1)^p\, p!\,\varsigma_a(p+1,a)=\frac{1}{2^p}\left[\psi^{(p)}\left(\frac{a+1}{2}\right)-\psi^{(p)}\left(\frac{a}{2}\right)\right]$$

$\square$

Since $\psi(1+u)-\psi(1-u)=\dfrac{1}{u}-\pi\cot\pi u$ we have

(4.3.74e) $\qquad \dfrac{1}{u}-\pi\cot\pi u=\sum_{k=0}^{\infty}\dfrac{1}{k+1}\sum_{j=0}^{k}\binom{k}{j}(-1)^j\log\dfrac{(1+u+j)}{(1-u+j)}$

and since $\psi(1-u)-\psi(u)=\pi\cot\pi u$ we also have

(4.3.74ei) $\qquad \pi\cot\pi u=\sum_{k=0}^{\infty}\dfrac{1}{k+1}\sum_{j=0}^{k}\binom{k}{j}(-1)^j\log\dfrac{(1+j-u)}{(j+u)}$

(4.3.74eiii) $\qquad \pi\cot\pi u+\psi(u)=\sum_{k=0}^{\infty}\dfrac{1}{k+1}\sum_{j=0}^{k}\binom{k}{j}(-1)^j\log(1+j-u)$

From (4.3.74ei) we may obtain a series for $\varsigma(3)$ via the integral (1.13) $\displaystyle\int_{0}^{\frac{1}{2}}u^2\cot\pi u\,du$ (see also (4.3.168a)).

$\square$

We now integrate the identity (4.3.74) to obtain

$$\log \Gamma(x) = \sum_{k=0}^{\infty} \frac{1}{k+1} \sum_{j=0}^{k} \binom{k}{j} (-1)^j \int_1^x \log(u+j)\,du$$

(4.3.74f)
$$= \sum_{k=0}^{\infty} \frac{1}{k+1} \sum_{j=0}^{k} \binom{k}{j} (-1)^j \left[ (x+j)\log(x+j) - (1+j)\log(1+j) - (x-1) \right]$$

$$= \sum_{k=0}^{\infty} \frac{1}{k+1} \sum_{j=0}^{k} \binom{k}{j} (-1)^j \, j \log \frac{(x+j)}{(1+j)} + x\psi(x) + \gamma - (x-1) \sum_{k=0}^{\infty} \frac{1}{k+1} \sum_{j=0}^{k} \binom{k}{j} (-1)^j$$

and since $\displaystyle \sum_{j=0}^{k} \binom{k}{j} (-1)^j = \delta_{0,k}$ we end up with

(4.3.75)
$$\log \Gamma(x) = \sum_{k=0}^{\infty} \frac{1}{k+1} \sum_{j=0}^{k} \binom{k}{j} (-1)^j \, j \log \frac{(x+j)}{(1+j)} + x\psi(x) + \gamma - (x-1)$$

It is curious that the right-hand side of (4.3.75) appears to remain finite as $x \to 0$ whereas the left-hand side $\to \infty$ (note however that we did start with (4.3.74) which is not valid at $x = 0$).

Differentiating (4.3.75) gives us

$$\psi(x) = \sum_{k=0}^{\infty} \frac{1}{k+1} \sum_{j=0}^{k} \binom{k}{j} (-1)^j \frac{j}{x+j} + x\psi'(x) + \psi(x) - 1$$

and we note that

$$\frac{j}{x+j} = 1 - \frac{x}{x+j}$$

which gives us

$$\sum_{k=0}^{\infty} \frac{1}{k+1} \sum_{j=0}^{k} \binom{k}{j} (-1)^j \frac{j}{x+j} = \sum_{k=0}^{\infty} \frac{1}{k+1} \sum_{j=0}^{k} \binom{k}{j} (-1)^j - \sum_{k=0}^{\infty} \frac{1}{k+1} \sum_{j=0}^{k} \binom{k}{j} (-1)^j \frac{x}{x+j}$$

We then have

$$\psi(x) = 1 - x\zeta(2,x) + x\psi'(x) + \psi(x) - 1$$

$$= 1 - x\psi'(x) + x\psi'(x) + \psi(x) - 1 = \psi(x)$$

which is where we started from.

Since $\log\Gamma(2) = 0$ we see from (4.3.75) that

$$\sum_{k=0}^{\infty}\frac{1}{k+1}\sum_{j=0}^{k}\binom{k}{j}(-1)^j j\log\frac{(2+j)}{(1+j)}+2\psi(2)+\gamma-1=0$$

and hence we get

(4.3.76)     $$\gamma-1=\sum_{k=0}^{\infty}\frac{1}{k+1}\sum_{j=0}^{k}\binom{k}{j}(-1)^j j\log\frac{(2+j)}{(1+j)}$$

Similarly we have

(4.3.76a)   $$\log(n!)=\sum_{k=0}^{\infty}\frac{1}{k+1}\sum_{j=0}^{k}\binom{k}{j}(-1)^j j\log\frac{(n-1+j)}{(1+j)}+(n-1)\psi(n-1)+\gamma-n+2$$

Letting $x=1/2$ in (4.3.75) we obtain

(4.3.76b)   $$\frac{1}{2}\big[\log\pi-1-\gamma+\log 2\big]=\sum_{k=0}^{\infty}\frac{1}{k+1}\sum_{j=0}^{k}\binom{k}{j}(-1)^j j\log\frac{(1+2j)}{2(1+j)}$$

and this may be simplified to

(4.3.76c)   $$\log\sqrt{\frac{2\pi}{e^{\gamma+1}}}=\sum_{k=0}^{\infty}\frac{1}{k+1}\sum_{j=0}^{k}\binom{k}{j}(-1)^j j\log\frac{(1+2j)}{(1+j)}$$

Using (4.3.73a) we may write (4.3.75) as

(4.3.76d)   $$\log\Gamma(x)=\sum_{k=0}^{\infty}\frac{1}{k+1}\sum_{j=0}^{k}\binom{k}{j}(-1)^j j\log(x+j)+x\psi(x)+\gamma-(x-1)$$

$$-\sum_{k=0}^{\infty}\frac{1}{k+1}\sum_{j=0}^{k}\binom{k}{j}(-1)^j j\log(1+j)$$

A more concise expression for $\log\Gamma(x)$ is given in (4.3.119).

Integration of (4.3.75) results in

$$\int_{1}^{u}\log\Gamma(x)dx=\sum_{k=0}^{\infty}\frac{1}{k+1}\sum_{j=0}^{k}\binom{k}{j}(-1)^j j\int_{1}^{u}\log\frac{(x+j)}{(1+j)}dx+\int_{1}^{u}x\psi(x)dx+(u-1)(\gamma+1)-\frac{1}{2}(u^2-1)$$

With integration by parts we have

$$\int\limits_1^u x\psi(x)\,dx = u\log\Gamma(u) - \int\limits_1^u \log\Gamma(x)\,dx$$

We also have

$$\int\limits_1^u \log\frac{(x+j)}{(1+j)}\,dx = (u+j)\log(u+j) - (1+j)\log(1+j) - (u-1) - (u-1)\log(1+j)$$

$$= j\log\frac{(u+j)}{(1+j)} + u\log\frac{(u+j)}{(1+j)} - (u-1) = (u+j)\log\frac{(u+j)}{(1+j)} - (u-1)$$

Differentiating the binomial identity we see that

$$\sum_{j=0}^k \binom{k}{j}(-1)^j j = -k\delta_{1,k}$$

and therefore we get

(4.3.76e) $$\sum_{k=0}^\infty \frac{1}{k+1}\sum_{j=0}^k \binom{k}{j}(-1)^j j = -\frac{1}{2}$$

and we accordingly obtain

$$2\int\limits_1^u \log\Gamma(x)\,dx = \sum_{k=0}^\infty \frac{1}{k+1}\sum_{j=0}^k \binom{k}{j}(-1)^j j^2 \log\frac{(u+j)}{(1+j)} + u\sum_{k=0}^\infty \frac{1}{k+1}\sum_{j=0}^k \binom{k}{j}(-1)^j j\log\frac{(u+j)}{(1+j)}$$

$$+ u\log\Gamma(u) + (u-1)\gamma + \frac{3}{2}(u-1) - \frac{1}{2}(u^2-1)$$

and substituting (4.3.75) we get (at least for $u > 0$)

(4.3.77)

$$2\int\limits_1^u \log\Gamma(x)\,dx = \sum_{k=0}^\infty \frac{1}{k+1}\sum_{j=0}^k \binom{k}{j}(-1)^j j^2 \log\frac{(u+j)}{(1+j)} + 2u\log\Gamma(u) - u^2\psi(u) - \gamma + \frac{1}{2}(u-1)(u+2)$$

Differentiating this gives us

$$2\log\Gamma(u) = \sum_{k=0}^{\infty}\frac{1}{k+1}\sum_{j=0}^{k}\binom{k}{j}(-1)^{j}\frac{j^{2}}{u+j} + 2u\psi(u) + 2\log\Gamma(u) - 2u\psi(u) + u^{2}\psi'(u) + \frac{1}{2}(2u+1)$$

and, after several cancellations, this simplifies to

$$\sum_{k=0}^{\infty}\frac{1}{k+1}\sum_{j=0}^{k}\binom{k}{j}(-1)^{j}\frac{j^{2}-u^{2}}{u+j} + \frac{1}{2}(2u+1) = 0$$

or equivalently

$$\sum_{k=0}^{\infty}\frac{1}{k+1}\sum_{j=0}^{k}\binom{k}{j}(-1)^{j}(u-j) + \frac{1}{2}(2u+1) = 0$$

whose validity is obvious.

Initially I thought that (4.3.77) was not valid for $u = 0$ because of the singularity arising from the term $\log\frac{(u+j)}{(1+j)}$ when $j = 0$; however, it is clear that the multiplicative factor of $j^{2}$ cures the problem. We therefore have (4.3.77) for $u \geq 0$ and with $u = 0$ we obtain

$$-2\int_{0}^{1}\log\Gamma(x)dx = \sum_{k=0}^{\infty}\frac{1}{k+1}\sum_{j=0}^{k}\binom{k}{j}(-1)^{j}j^{2}\log\frac{j}{1+j} - \gamma - 1$$

Since $\int_{0}^{1}\log\Gamma(x)dx = \frac{1}{2}\log(2\pi)$ we see that

(4.3.77a) $$\sum_{k=0}^{\infty}\frac{1}{k+1}\sum_{j=0}^{k}\binom{k}{j}(-1)^{j}j^{2}\log\frac{j}{1+j} = \gamma + 1 - \log(2\pi)$$

A proof of the following well-known identity is given in (C.43a) in Appendix C of Volume VI

$$\int_{t}^{t+1}\log\Gamma(x)dx = t\log t - t + \frac{1}{2}\log(2\pi)$$

and, since $\int_{1}^{2}\log\Gamma(x)dx = -1 + \frac{1}{2}\log(2\pi)$, by letting $u = 2$ in (4.3.77) we obtain

(4.3.78) $$\sum_{k=0}^{\infty}\frac{1}{k+1}\sum_{j=0}^{k}\binom{k}{j}(-1)^{j}j^{2}\log\frac{(2+j)}{(1+j)} = \log(2\pi) - 3\gamma$$

Integration by parts gives us

$$\int_1^u \log\Gamma(x)\,dx = u\log\Gamma(u) - \int_1^u x\psi(x)\,dx$$

and from (4.3.77) we therefore obtain

$$(4.3.79) \quad 2\int_1^u x\psi(x)\,dx = -\sum_{k=0}^\infty \frac{1}{k+1}\sum_{j=0}^k \binom{k}{j}(-1)^j j^2 \log\frac{(u+j)}{(1+j)} + u^2\psi(u) + \gamma - \frac{1}{2}(u-1)(u+2)$$

## AN ALTERNATIVE PROOF OF ALEXEIEWSKY'S THEOREM

We now consider moments of the digamma function (this part was inspired by Snowden's paper [120aa] which was written while he was an undergraduate at Princeton).

Let $\phi_n(x)$ be defined by

$$(4.3.79a) \quad \phi_n(x) = \int_0^x t^{n-1}\psi(1+t)\,dt$$

It may be seen that

$$\phi_1(x) = \log\Gamma(1+x)$$

and since $\psi(1+t) = \psi(t) + 1/t$ we have for $n \geq 2$

$$\phi_n(x) = \frac{x^{n-1}}{n-1} + \int_0^x t^{n-1}\psi(t)\,dt$$

Using the series for the digamma function (E.14)

$$(4.3.79b) \quad \psi(1+t) = -\gamma + \sum_{k=1}^\infty \left(\frac{1}{k} - \frac{1}{k+t}\right)$$

we get

$$\phi_n(x) = -\gamma\frac{x^n}{n} + \int_0^x t^{n-1}\sum_{k=1}^\infty \left(\frac{1}{k} - \frac{1}{k+t}\right)dt$$

$$= -\gamma \frac{x^n}{n} + \sum_{k=1}^{\infty} \left( \frac{x^n}{nk} - \int_0^x \frac{t^{n-1}}{k+t} dt \right)$$

We easily see that

$$\int_0^x \frac{t^{n-1}}{k+t} dt = \int_0^x \left( t^{n-2} + \frac{t^{n-1}}{k+t} - t^{n-2} \right) dt = \int_0^x \left( t^{n-2} - \frac{kt^{n-2}}{k+t} \right) dt$$

$$= \frac{x^{n-1}}{n-1} - k \int_0^x \frac{t^{n-2}}{k+t} dt$$

and, applying this reduction formula $n-1$ times, we obtain

(4.3.80)    $$\int_0^x \frac{t^{n-1}}{k+t} dt = P_n(x,k) + (-1)^{n-1} k^{n-1} \log\left( 1 + \frac{x}{k} \right)$$

where the homogenous polynomial $P_n(x,k)$ is defined by

$$P_n(x,k) = \sum_{j=1}^{n-1} (-1)^{j+1} \frac{k^{j-1} x^{n-j}}{n-j}$$

We may also represent $P_n(x,k)$ as

$$P_n(x,k) = (-1)^{n+1} \sum_{m=1}^{n-1} (-1)^m \frac{k^{n-m-1} x^m}{m}$$

$$= (-1)^{n+1} k^{n-1} \sum_{m=1}^{n-1} \frac{(-1)^m}{m} \left( \frac{x}{k} \right)^m$$

$$= \frac{k^{n-1}}{n} \left( \frac{x}{k} \right)^n + (-1)^{n+1} k^{n-1} \sum_{m=1}^{n} \frac{(-1)^m}{m} \left( \frac{x}{k} \right)^m$$

$$= \frac{x^n}{nk} + (-1)^{n+1} k^{n-1} \sum_{m=1}^{n} \frac{(-1)^m}{m} \left( \frac{x}{k} \right)^m$$

As noted in Adamchik's paper [6a], equation (4.3.80) is in fact just the Taylor series expansion of $\log\left( 1 + \frac{x}{k} \right)$ with the Lagrange form [17a, p.243] of the remainder

(4.3.80a) $\qquad \log\left(1+\dfrac{x}{k}\right) = \sum_{j=1}^{n-1} \dfrac{(-1)^{j+1}}{j}\left(\dfrac{x}{k}\right)^j + \dfrac{(-1)^{n-1}}{k^{n-1}} \int_0^x \dfrac{t^{n-1}}{k+t}\,dt$

This identity may also be obtained more directly by integrating the geometric sum over the interval $[0, x/k]$

$$\frac{1}{1+y} = \sum_{j=1}^{n-1} (-1)^{j-1} y^{j-1} + (-1)^{n-1} \frac{y^{n-1}}{1+y}$$

We therefore obtain Snowden's result

(4.3.81) $\qquad \phi_n(x) = -\gamma \dfrac{x^n}{n} + \sum_{k=1}^{\infty}\left(\dfrac{x^n}{nk} - P_n(x,k) + (-1)^n k^{n-1}\log\left(1+\dfrac{x}{k}\right)\right)$

and see that

$$\exp[\phi_n(x)] = \exp\left[-\gamma \frac{x^n}{n}\right] \prod_{k=1}^{\infty}\left\{\left(1+\frac{x}{k}\right)^{(-1)^n k^{n-1}} \exp\left(\frac{x^n}{nk} - P_n(x,k)\right)\right\}$$

With $n=2$ we have

$$\exp[\phi_2(x)] = \exp\left[-\gamma \frac{x^2}{2}\right] \prod_{k=1}^{\infty}\left\{\left(1+\frac{x}{k}\right)^k \exp\left(\frac{x^2}{2k} - x\right)\right\}$$

and this is immediately reminiscent of the Barnes double gamma function $\Gamma_2(x) = 1/G(x)$ defined, inter alia, by [126, p.25]

(4.3.81a) $\qquad G(1+x) = (2\pi)^{x/2} \exp\left[-\dfrac{1}{2}(\gamma x^2 + x^2 + x)\right] \prod_{k=1}^{\infty}\left\{\left(1+\dfrac{x}{k}\right)^k \exp\left(\dfrac{x^2}{2k} - x\right)\right\}$

from which we immediately see that $G(1) = 1$.
Hence we have

$$G(1+x) = (2\pi)^{x/2} \exp[\phi_2(x)] \exp\left[-\frac{1}{2}x(x+1)\right]$$

or alternatively

(4.3.82) $\qquad \phi_2(x) = \log G(1+x) - \dfrac{x}{2}\log(2\pi) + \dfrac{1}{2}x(x+1)$

This then gives us

$$(4.3.83) \qquad \phi_2(x) = \int_0^x t\,\psi(1+t)\,dt = \log G(1+x) - \frac{x}{2}\log(2\pi) + \frac{1}{2}x(x+1)$$

It is easily seen that

$$\int_1^u x\psi(x)\,dx = \int_0^{u-1}(1+t)\psi(1+t)\,dt = \int_0^{u-1}\psi(1+t)\,dt + \int_0^{u-1}t\psi(1+t)\,dt$$

Therefore we get

$$(4.3.83a) \qquad \int_1^u t\psi(t)\,dt = \log\Gamma(u) + \int_0^{u-1}t\psi(1+t)\,dt$$

and we then have from (4.3.83)

$$(4.3.83b) \qquad \int_0^{u-1}t\psi(1+t)\,dt = \log G(u) - \frac{(u-1)}{2}\log(2\pi) + \frac{1}{2}u(u-1)$$

This formula is contained in Adamchik's paper "On the Barnes function" [5a]. $G(x)$ is often referred to as the Barnes $G$-function. The Right Rev. Ernest William Barnes D.D. (1874-1953) was formerly Bishop of Birmingham and died at the age of 79 in the year of my birth. He wrote his last mathematical paper in 1910 at the relatively young age of 36, presumably due to his increasing involvement in more ethereal matters.

We therefore deduce from (4.3.79) that

$$(4.3.83c) \qquad 2\log\Gamma(u) + 2\log G(u) - (u-1)\log(2\pi) + u(u-1) =$$

$$-\sum_{k=0}^{\infty}\frac{1}{k+1}\sum_{j=0}^{k}\binom{k}{j}(-1)^j j^2 \log\frac{(u+j)}{(1+j)} + u^2\psi(u) + \gamma - \frac{1}{2}(u-1)(u+2)$$

and since [126, p.25] $G(1+u) = \Gamma(u)G(u)$ we may write this as

$$2\log G(1+u) = (u-1)\log(2\pi) - \sum_{k=0}^{\infty}\frac{1}{k+1}\sum_{j=0}^{k}\binom{k}{j}(-1)^j j^2 \log\frac{(u+j)}{(1+j)}$$

$$(4.3.84)$$

$$+ u^2\psi(u) + \gamma - \frac{1}{2}(u-1)(3u+2)$$

or as

$$\sum_{k=0}^{\infty} \frac{1}{k+1} \sum_{j=0}^{k} \binom{k}{j} (-1)^j j^2 \log \frac{(u+j)}{(1+j)} = (u-1)\log(2\pi) - 2\log G(1+u)$$

(4.3.84a)

$$+ u^2 \psi(u) + \gamma - \frac{1}{2}(u-1)(3u+2)$$

This formula is used later in (4.3.126) to derive the Gosper/Vardi functional equation for $\log G(1+u)$.

Letting $u = 0$ gives us another derivation of (4.3.77a). With $u = 2$ we have

$$2\log G(3) = \log(2\pi) - \sum_{k=0}^{\infty} \frac{1}{k+1} \sum_{j=0}^{k} \binom{k}{j} (-1)^j j^2 \log \frac{(2+j)}{(1+j)} + 4\psi(2) + \gamma - 4$$

Using $G(1+u) = \Gamma(u)G(u)$ we see that $G(3) = 1$ and, since $\psi(2) = 1 - \gamma$, we obtain another derivation of (4.3.78).

We note that

$$\int_1^u \log \Gamma(x)\,dx = \int_0^u \log \Gamma(x)\,dx - \int_0^1 \log \Gamma(x)\,dx = \int_0^u \log \Gamma(x)\,dx - \frac{1}{2}\log(2\pi)$$

and combining (4.3.77) and (4.3.84) we see that

(4.3.85) $\quad \displaystyle\int_0^u \log \Gamma(x)\,dx = \frac{u(1-u)}{2} + \frac{u}{2}\log(2\pi) - \log G(u) - (1-u)\log \Gamma(u)$

$$= \frac{u(1-u)}{2} + \frac{u}{2}\log(2\pi) - \log G(u+1) + u\log \Gamma(u)$$

and this is Alexeiewsky's theorem [45]. In fact, this result may be obtained more directly from (4.3.83) by using integration by parts since

$$\int_0^x t\,\psi(1+t)\,dt = x\log \Gamma(1+x) - \int_0^x \log \Gamma(1+t)\,dt$$

Another derivation is contained in (4.3.126e).

We may also note that

$$(4.3.85a) \quad \int_1^u \log \Gamma(x)\, dx = \frac{u(1-u)}{2} + \frac{u}{2}\log(2\pi) - \log G(u) - (1-u)\log\Gamma(u) - \frac{1}{2}\log(2\pi)$$

$$(4.3.85b) \quad \int_1^u \log\Gamma(1+x)\, dx = \frac{u(1-u)}{2} + \frac{u}{2}\log(2\pi) - \log G(u+1) + u\log\Gamma(u) - \frac{1}{2}\log(2\pi)$$

$$+ u\log u - u + 1$$

$$\int_1^u \log\Gamma(1+x)\, dx = -\frac{1}{2}(u-1)(u+2) + \frac{1}{2}(u-1)\log(2\pi) - \log G(u+1) + u\log\Gamma(1+u)$$

Upon differentiating (4.3.84) we obtain

$$2\frac{G'(1+u)}{G(1+u)} = \log(2\pi) - \sum_{k=0}^{\infty}\frac{1}{k+1}\sum_{j=0}^{k}\binom{k}{j}(-1)^j\frac{j^2}{u+j} + u^2\psi'(u) + 2u\psi(u) - 3u + \frac{1}{2}$$

$$= \log(2\pi) - \sum_{k=0}^{\infty}\frac{1}{k+1}\sum_{j=0}^{k}\binom{k}{j}(-1)^j j + u\sum_{k=0}^{\infty}\frac{1}{k+1}\sum_{j=0}^{k}\binom{k}{j}(-1)^j - u^2\sum_{k=0}^{\infty}\frac{1}{k+1}\sum_{j=0}^{k}\binom{k}{j}\frac{(-1)^j}{u+j}$$

$$+ u^2\psi'(u) + 2u\psi(u) - 3u + \frac{1}{2}$$

$$= \log(2\pi) + \frac{1}{2} + u - u^2\varsigma(2,u) + u^2\psi'(u) + 2u\psi(u) - 3u + \frac{1}{2}$$

$$= \log(2\pi) + 1 - u^2\psi'(u) + u^2\psi'(u) + 2u\psi(u) - 2u$$

Hence we have

$$(4.3.85c) \qquad 2\frac{G'(1+u)}{G(1+u)} = 1 + \log(2\pi) + 2u\psi(u) - 2u$$

and integrating this gives us another derivation of Alexeiewsky's theorem (4.3.85)

$$2\log G(x+1) = [1 + \log(2\pi)]x + 2x\log\Gamma(x) - 2\int_0^x \log\Gamma(u)\, du - x^2$$

$\square$

**Lemma:**

(4.3.85d)

$$\frac{1}{2}I(2x) = I(x) + I\left(x + \frac{1}{2}\right) - I\left(\frac{1}{2}\right) + x^2\log 2 - \left[\log\Gamma\left(\frac{1}{2}\right) + \log 2\right]x + \frac{1}{2}\log\Gamma\left(\frac{1}{2}\right) + \frac{1}{4}\log 2$$

where $\quad I(x) = \int\limits_{1}^{x} \log\Gamma(t)\,dt$

**Proof:**

The following proof is taken from a 1997 paper by Billingham and King [23aa] (and this corrects a minor misprint in their proof).

It is clear from the definition that $I'(x) = \log\Gamma(x)$ and $I''(x) = \psi(x)$. Now the digamma function satisfies the duplication formula ([1] and [126, p.15])

$$\psi(mx) = \log m + \frac{1}{m}\sum_{k=0}^{m-1}\psi\left(x + \frac{k}{m}\right) \quad , \ m \in \mathbf{N}$$

and with $m = 2$ we get

$$\psi(2x) = \frac{1}{2}\psi(x) + \frac{1}{2}\psi\left(x + \frac{1}{2}\right) + \log 2$$

Hence we have

$$I''(2x) = \frac{1}{2}I''(x) + \frac{1}{2}I''\left(x + \frac{1}{2}\right) + \log 2$$

Integration of the above results in

$$\frac{1}{2}I'(2x) = \frac{1}{2}I'(x) + \frac{1}{2}I'\left(x + \frac{1}{2}\right) + x\log 2 + c$$

and the integration constant may be determined by putting $x = 1/2$ in the above equation. This gives us

$$\frac{1}{2}\log\Gamma(1) = \frac{1}{2}\log\Gamma(1/2) + \frac{1}{2}\log\Gamma(1) + \frac{1}{2}\log 2 + c$$

A further integration results in

$$\frac{1}{4}I(2x) = \frac{1}{2}I(x) + \frac{1}{2}I\left(x + \frac{1}{2}\right) + \frac{1}{2}x^2\log 2 + cx + b$$

and $b$ is also determined by letting $x = 1/2$. Since by definition $I(1) = 0$ we then obtain a proof of the lemma.

From Alexeiewsky's theorem (4.3.85)

$$\int_0^x \log\Gamma(t)\,dt = \frac{x(1-x)}{2} + \frac{x}{2}\log(2\pi) - \log G(x) - (1-x)\log\Gamma(x)$$

and noting that

$$I(x) = \int_1^x \log\Gamma(t)\,dt = \int_0^x \log\Gamma(t)\,dt - \int_0^1 \log\Gamma(t)\,dt$$

$$= \int_0^x \log\Gamma(t)\,dt - \frac{1}{2}\log(2\pi)$$

we see that

$$I(x) = \int_1^x \log\Gamma(t)\,dt = \frac{x(1-x)}{2} + \frac{1}{2}(x-1)\log(2\pi) - \log G(x) - (1-x)\log\Gamma(x)$$

Using the above lemma we obtain the functional equation

$$\frac{x(1-2x)}{2} + \frac{1}{4}(2x-1)\log(2\pi) - \frac{1}{2}\log G(2x) - \frac{1}{2}(1-2x)\log\Gamma(2x)$$

$$= \frac{x(1-x)}{2} + \frac{1}{2}(x-1)\log(2\pi) - \log G(x) - (1-x)\log\Gamma(x)$$

$$+ \frac{(1+2x)(1-2x)}{8} + \frac{1}{4}(2x-1)\log(2\pi) - \log G\left(x + \frac{1}{2}\right) - \frac{1}{2}(1-2x)\log\Gamma\left(x + \frac{1}{2}\right)$$

$$+ x^2\log 2 - \left[\log\Gamma\left(\frac{1}{2}\right) + \log 2\right]x - I\left(\frac{1}{2}\right) + \frac{1}{2}\log\Gamma\left(\frac{1}{2}\right) + \frac{1}{4}\log 2$$

Letting $x = 1/2$, since $G(1) = 1$, we simply recover (4.3.85) evaluated at $x = 1/2$.

$$I(1/2) = \int\limits_1^{1/2} \log \Gamma(t) dt = \frac{1}{8} - \frac{1}{4} \log 2 - \frac{3}{4} \log \pi - \log G(1/2)$$

The above equation may be simplified to

(4.3.85e)

$$\log G\left(x + \frac{1}{2}\right) + \log G(x) - \frac{1}{2} \log G(2x) =$$

$$\frac{1}{2}(x-1)\log(2\pi) + \frac{1}{2}x^2 + \frac{(1+2x)(1-2x)}{8}$$

$$-\frac{1}{2}(1-2x)\left[\log \Gamma\left(x + \frac{1}{2}\right) - \log \Gamma(2x)\right] - (1-x)\log \Gamma(x)$$

$$+x^2 \log 2 - \left[\log \Gamma\left(\frac{1}{2}\right) + \log 2\right]x - I\left(\frac{1}{2}\right) + \frac{1}{2}\log \Gamma\left(\frac{1}{2}\right) + \frac{1}{4}\log 2$$

and using Legendre's duplication formula [126, p.7]

(4.3.85f)
$$\frac{1}{2}\log \pi + \log \Gamma(2x) = (2x-1)\log 2 + \log \Gamma(x) + \log \Gamma\left(x + \frac{1}{2}\right)$$

we may eliminate the terms involving the gamma function and obtain the Barnes multiplication formula.

In 1899 Barnes developed a multiplication formula for $G(nx)$ (see [126, p.30]) and a particular case is set out below [126, p.29]

(4.3.85g)
$$G^2(x)G^2\left(x + \frac{1}{2}\right)\Gamma(x) = J(x)G(2x)$$

where for convenience

$$\log J(x) = \frac{1}{4} - 3\log A + \left(-2x^2 + 3x - \frac{11}{12}\right)\log 2 + \left(x - \frac{1}{2}\right)\log \pi$$

□

Noting that [126, p.14]

(4.3.86) $$\psi(x) = \psi(-x) - \pi \cot \pi x - 1/x$$

we have

$$\int_0^u x\psi(x)\,dx - \int_0^u x\psi(-x)\,dx = -\pi \int_0^u x \cot \pi x\,dx - u$$

With integration by parts we see that

$$\int_0^u x\psi(x)\,dx = u \log \Gamma(u) - \int_0^u \log \Gamma(x)\,dx$$

$$\int_0^u x\psi(-x)\,dx = \int_0^{-u} x\psi(x)\,dx = -u \log \Gamma(-u) - \int_0^{-u} \log \Gamma(x)\,dx$$

and therefore

$$\int_0^u x\psi(x)\,dx - \int_0^u x\psi(-x)\,dx = u \log \Gamma(u) + u \log \Gamma(-u) + \int_0^{-u} \log \Gamma(x)\,dx + \int_0^u \log \Gamma(x)\,dx$$

Then, using Alexeiewsky's theorem, we obtain the well-known formula originally due to Kinkelin in 1860 [126, p.30] (an electronic version of Kinkelin's paper may be obtained from the website referred to in [52])

(4.3.87) $$\log \frac{G(1+u)}{G(1-u)} = u \log 2\pi - \pi \int_0^u x \cot \pi x\,dx$$

which we will prove in a completely different manner in (6.69b) in Volume V. It is seen that substituting $u = 1/2$ in (4.3.87) results in another proof of (3.2) in Volume I.

With regard to (4.3.87), in 1992 Freund and Zabrodin [69aa] reported the more general formula for $n \geq 2$

(4.3.87i) $$\Gamma_n(u) \left[ \Gamma_n(-u) \right]^{(-1)^{n+1}} = \exp\left[ -\pi \int_0^u x^{n-1} \cot \pi x\,dx \right]$$

where we have followed the Vignéras notation [126, p.39]

$$\Gamma_n(x) = \left[ G_n(x) \right]^{(-1)^{n-1}}$$

and $$\Gamma_1(x) = G_1(x) = \Gamma(x) \qquad \Gamma_2(x) = 1/G_2(x)$$

$\square$

From (4.3.81) we see that

$$\phi_1(x) = \log \Gamma(1+x) = -\gamma x + \sum_{k=1}^{\infty} \left[ \frac{x}{k} - \log\left(1 + \frac{x}{k}\right) \right]$$

$$\phi_2(x) = \log G(1+x) - \frac{x}{2}\log(2\pi) + \frac{1}{2}x(x+1)$$

$$= -\frac{\gamma}{2}x^2 + \sum_{k=1}^{\infty} \left[ \frac{x^2}{2k} - x + k\log\left(1 + \frac{x}{k}\right) \right]$$

Differentiation results in

$$\psi(1+x) = -\gamma + \sum_{k=1}^{\infty} \left[ \frac{1}{k} - \frac{1}{k+x} \right]$$

which gets us back to where we started from in (4.3.79b). We also see that

$$\phi_2'(x) = -\gamma x + \sum_{k=1}^{\infty} \left[ \frac{x}{k} - 1 + \frac{1}{1 + x/k} \right]$$

$$= -\gamma x + \sum_{k=1}^{\infty} \left[ \frac{x}{k} - \frac{x}{k+x} \right]$$

$$= x\psi(1+x)$$

Hence we see that

$$\frac{G'(1+x)}{G(1+x)} - \frac{1}{2}[\log(2\pi) - 1] + x = x\psi(1+x)$$

We then have upon multiplication by $G(1+x)$ and subsequent integration

$$G(1+x) - \frac{1}{2}[\log(2\pi) - 1]\int_0^u G(1+x)dx + \int_0^u xG(1+x)dx = \int_0^u x\psi(1+x)G(1+x)dx$$

More generally we have

$$\phi_n'(x) = -\gamma x^{n-1} + \sum_{k=1}^{\infty} \left( \frac{x^{n-1}}{k} - \sum_{j=1}^{n-1} (-1)^{j+1} k^{j-1} x^{n-j-1} + \frac{(-1)^n k^{n-1}}{k+x} \right)$$

We now refer back to Snowden's formulation in (4.3.79a)

$$\exp\left[\phi_n(x)\right] = \exp\left[-\gamma\frac{x^n}{n}\right]\prod_{k=1}^{\infty}\left\{\left(1+\frac{x}{k}\right)^{(-1)^n k^{n-1}}\exp\left(\frac{x^n}{nk}-P_n(x,k)\right)\right\}$$

where for $n=3$ we have

$$\exp\left[\phi_3(x)\right] = \exp\left[-\frac{1}{3}\gamma x^3\right]\prod_{k=1}^{\infty}\left\{\left(1+\frac{x}{k}\right)^{-k^2}\exp\left(\frac{1}{3}\frac{x^3}{k}-\frac{1}{2}x^2+kx\right)\right\}$$

and taking the square root gives us

$$\exp\frac{1}{2}\left[\phi_3(x)\right] = \exp\left[-\frac{1}{6}\gamma x^3\right]\prod_{k=1}^{\infty}\left\{\left(1+\frac{x}{k}\right)^{-\frac{1}{2}k^2}\exp\left(\frac{1}{6}\frac{x^3}{k}-\frac{1}{4}x^2+\frac{1}{2}kx\right)\right\}$$

Referring to (4.3.81a) we see that

$$\sqrt{\Gamma_2(1+x)} = (2\pi)^{-x/4}\exp\left[\frac{1}{4}(\gamma x^2+x^2+x)\right]\prod_{k=1}^{\infty}\left\{\left(1+\frac{x}{k}\right)^{-\frac{1}{2}k}\exp\left(-\frac{x^2}{4k}+\frac{1}{2}x\right)\right\}$$

and hence we have

$$\exp\frac{1}{2}\left[\phi_3(x)\right]\sqrt{\Gamma_2(1+x)} =$$

$$(2\pi)^{-x/4}\exp\left[\frac{1}{4}(\gamma x^2+x^2+x)-\frac{1}{6}\gamma x^3\right]\times$$

$$\prod_{k=1}^{\infty}\left\{\left(1+\frac{x}{k}\right)^{-\frac{1}{2}k(k+1)}\exp\left(-\frac{x^2}{4k}+\frac{1}{2}x\right)\exp\left(\frac{1}{6}\frac{x^3}{k}-\frac{1}{4}x^2+\frac{1}{2}kx\right)\right\}$$

$$= (2\pi)^{-x/4}\exp\left[\frac{1}{4}(\gamma x^2+x^2+x)-\frac{1}{6}\gamma x^3-\frac{1}{6}\varsigma(2)x^3\right].F(x)$$

where

$$F(x) = \prod_{k=1}^{\infty} \left\{ \left(1 + \frac{x}{k}\right)^{-\frac{1}{2}k(k+1)} \exp\left(\frac{1}{2}(k+1)x - \frac{1}{4}\left(1+\frac{1}{k}\right)x^2 + \frac{1}{6}\frac{x^3}{k} + \frac{1}{6}\frac{x^3}{k^2}\right) \right\}$$

$$= \prod_{k=1}^{\infty} \left\{ \left(1 + \frac{x}{k}\right)^{-\frac{1}{2}k(k+1)} \exp\left[\left(1+\frac{1}{k}\right)\left(\frac{1}{2}kx - \frac{1}{4}x^2 + \frac{1}{6k}x^3\right)\right] \right\}$$

From [45ac] and also [126, p.42] we find that

$$\Gamma_3(1+x) = \exp(c_1 x + c_2 x^2 + c_3 x^3) F(x)$$

where

$$c_1 = \frac{3}{8} - \frac{1}{4}\log(2\pi) - \log A \quad, \qquad c_2 = \frac{1}{4}\left[\gamma + \log(2\pi) + \frac{1}{2}\right]$$

$$c_3 = -\frac{1}{6}\left[\gamma + \varsigma(2) + \frac{3}{2}\right]$$

and hence we see that

$$\Gamma_3(1+x) =$$

$$(2\pi)^{x/4} \exp\frac{1}{2}[\phi_3(x)] \sqrt{\Gamma_2(1+x)} \exp\left[-\frac{1}{4}(\gamma x^2 + x^2 + x) + \frac{1}{6}\gamma x^3 + \frac{1}{6}\varsigma(2)x^3 + c_1 x + c_2 x^2 + c_3 x^3\right]$$

$$= (2\pi)^{x/4} \exp\frac{1}{2}[\phi_3(x)] \sqrt{\Gamma_2(1+x)} \exp\left[\left(c_1 - \frac{1}{4}\right)x + \left(c_2 - \frac{1}{4} - \frac{1}{4}\gamma\right)x^2 + \left(c_3 + \frac{1}{6}\gamma + \frac{1}{6}\varsigma(2)\right)x^3\right]$$

Therefore we have

$$\log \Gamma_3(1+x) = \frac{x}{4}\log(2\pi) + \frac{1}{2}[\phi_3(x)] + \frac{1}{2}\log\Gamma_2(1+x)$$

$$+ \left[\left(c_1 - \frac{1}{4}\right)x + \left(c_2 - \frac{1}{4} - \frac{1}{4}\gamma\right)x^2 + \left(c_3 + \frac{1}{6}\gamma + \frac{1}{6}\varsigma(2)\right)x^3\right]$$

which implies that

$$\phi_3(x) = 2\log\Gamma_3(1+x) - \frac{x}{2}\log(2\pi) - \log\Gamma_2(1+x)$$

$$-2\left[\left(c_1 - \frac{1}{4}\right)x + \left(c_2 - \frac{1}{4} - \frac{1}{4}\gamma\right)x^2 + \left(c_3 + \frac{1}{6}\gamma + \frac{1}{6}\varsigma(2)\right)x^3\right]$$

We then have

$$\phi_3(x) = \int_0^x t^2\psi(1+t)dt = 2\log\Gamma_3(1+x) - \frac{x}{2}\log(2\pi) - \log\Gamma_2(1+x)$$

$$-2\left[\left(c_1 - \frac{1}{4}\right)x + \left(c_2 - \frac{1}{4} - \frac{1}{4}\gamma\right)x^2 + \left(c_3 + \frac{1}{6}\gamma + \frac{1}{6}\varsigma(2)\right)x^3\right]$$

which may be written as

(4.3.87a)

$$\phi_3(x) = \int_0^x t^2\psi(1+t)dt = 2\log\Gamma_3(1+x) - \log\Gamma_2(1+x)$$

$$+\left(-\frac{1}{4} + \log A\right)x + \left(-\frac{1}{2}\log(2\pi) + \frac{1}{4}\right)x^2 + \frac{1}{2}x^3$$

It is therefore clear that Snowden's insight has resulted in a very compact and elegant analysis. Similarly, working with $\exp[\phi_4(x)]$ will give us an expression involving $\log\Gamma_4(1+x)$.

With integration by parts we have

$$\int_0^x t^2\psi(1+t)dt = x^2\log\Gamma(1+x) - 2\int_0^x t\log\Gamma(1+t)dt$$

and we therefore obtain

(4.3.87ai)  $\quad 2\int_0^x t\log\Gamma(1+t)dt =$

$$\left(\frac{1}{4} - 2\log A\right)x + \left(\frac{1}{2}\log(2\pi) - \frac{1}{4}\right)x^2 - \frac{1}{2}x^3 + x^2\log\Gamma(1+x) - \log G(1+x) - 2\log\Gamma_3(1+x)$$

which concurs with the more general integral $\int\limits_0^x t\log\Gamma(a+t)dt$ given in [126, p.209]. This general case can be determined in the same manner as above by starting with

$$\psi(a+t) = -\gamma + \sum_{k=1}^\infty \left( \frac{1}{k} - \frac{1}{k+a-1+t} \right)$$

instead of (4.3.79b).

We see from (4.3.85b) that

$$\int\limits_0^x \log\Gamma(1+t)dt = \frac{1}{2}[\log(2\pi)-1]x - \frac{1}{2}x^2 + x\log\Gamma(1+x) - \log G(1+x)$$

and integration by parts gives us

$$\int\limits_0^x t\log\Gamma(1+t)dt = x\left[ \frac{1}{2}[\log(2\pi)-1]x - \frac{1}{2}x^2 + x\log\Gamma(1+x) - \log G(1+x) \right]$$

$$-\int\limits_0^x \left[ \frac{1}{2}[\log(2\pi)-1]u - \frac{1}{2}u^2 + u\log\Gamma(1+u) - \log G(1+u) \right]du$$

Therefore we have

$$2\int\limits_0^x t\log\Gamma(1+t)dt = x\left[ \frac{1}{2}[\log(2\pi)-1]x - \frac{1}{2}x^2 + x\log\Gamma(1+x) - \log G(1+x) \right]$$

$$-\int\limits_0^x \left[ \frac{1}{2}[\log(2\pi)-1]u - \frac{1}{2}u^2 - \log G(1+u) \right]du$$

$$= x\left[ \frac{1}{2}[\log(2\pi)-1]x - \frac{1}{2}x^2 + x\log\Gamma(1+x) - \log G(1+x) \right]$$

$$-\frac{1}{4}[\log(2\pi)-1]x^2 + \frac{1}{6}x^3 + \int\limits_0^x \log G(1+t)\,dt$$

We then see that

(4.3.87b)

$$2\int_0^x t\log\Gamma(1+t)dt = \frac{1}{4}[\log(2\pi)-1]x^2 - \frac{1}{3}x^3 + x^2\log\Gamma(1+x) - x\log G(1+x) + \int_0^x \log G(1+t)\, dt$$

and hence we may determine the integral $\int_0^x \log G(1+t)\, dt$ using (4.3.87a). It is easily seen that

(4.3.87c) $\quad \int_0^x \log G(1+t)dt =$

$$\left(\frac{1}{4}-2\log A\right)x + \frac{1}{4}\log(2\pi)x^2 - \frac{1}{6}x^3 + (x-1)\log G(1+x) - 2\log\Gamma_3(1+x)$$

in accordance with [126, p.207]. Using integration by parts will then enable us to evaluate $\int_0^x t\log G(1+t)dt$.

We also note that in [5a] Adamchik shows the following relationship

$$\log G_3(1+u) = \int_0^u t(1-t)\psi(t)dt + \frac{1}{12}u(6u^2 - 9u + 5) - u(1-u)\varsigma'(0) - 2u\varsigma'(-1)$$

We have the well-known Maclaurin expansion [126, p.160] for $|u| < 1$ (see also (E.34d) in Volume VI)

(4.3.87d) $\quad\quad\quad \psi(1+u) = \sum_{k=1}^\infty (-1)^k \varsigma(k) u^{k-1}$

where for convenience we have designated $\varsigma(1) = \gamma$. Hence reference to (4.3.79a) gives us

$$\phi_2(x) = \int_0^x t\psi(1+t)dt = \sum_{k=1}^\infty (-1)^k \frac{\varsigma(k)}{k+1} x^{k+1}$$

and from (4.3.82) we see that [126, p.210] for $|x| < 1$ (a derivation is contained in (E.12aii) in Volume VI)

(4.3.87e) $\quad \log G(1+x) = \frac{1}{2}x\log(2\pi) - \frac{1}{2}\gamma x^2 - \frac{1}{2}x(x+1) + \sum_{k=2}^\infty (-1)^k \frac{\varsigma(k)}{k+1} x^{k+1}$

which may be compared with (4.4.117i) in Volume III

(4.3.87f) $\qquad \log \Gamma(1+x) = -\gamma x + \sum_{k=2}^{\infty} (-1)^k \frac{\varsigma(k)}{k} x^k$

(which may be obtained by integrating (4.3.87d)).

Using (4.3.87a) we obtain

$$\int_0^x t^2 \psi(1+t)\,dt = 2\log\Gamma_3(1+x) - \log\Gamma_2(1+x) + \left(-\frac{1}{4} + \log A\right)x + \left(-\frac{1}{2}\log(2\pi) + \frac{1}{4}\right)x^2 + \frac{1}{2}x^3$$

$$= \sum_{k=1}^{\infty} (-1)^k \frac{\varsigma(k)}{k+3} x^{k+3}$$

Combining (4.3.32) and (4.3.34a) we get

$$\psi(x+a) - \psi(a) = \sum_{k=1}^{\infty} \frac{(-1)^{k+1}}{k} \frac{1}{a(a+1)...(a+k-1)} \sum_{j=0}^{k} s(k,j)x^j$$

and integration results in

$$\log\Gamma(x+a) - \log\Gamma(a) - x\psi(a) = \sum_{k=1}^{\infty} \frac{(-1)^{k+1}}{k} \frac{1}{a(a+1)...(a+k-1)} \sum_{j=0}^{k} \frac{s(k,j)x^{j+1}}{j+1}$$

With $a = 1$ this becomes

(4.3.87g) $\qquad \log\Gamma(x+1) + \gamma x = \sum_{k=1}^{\infty} \frac{(-1)^{k+1}}{k.k!} \sum_{j=0}^{k} \frac{s(k,j)x^{j+1}}{j+1}$

and comparing coefficients of $x^p$ in (4.3.87f) we obtain another derivation of (4.3.50)

$$\varsigma(p+1) = (-1)^p \sum_{k=1}^{\infty} \frac{(-1)^k}{k.k!} s(k,p)$$

Integrating (4.3.87g) we obtain

$$\int_0^u \log\Gamma(x+1)\,dx + \frac{1}{2}\gamma u^2 = \sum_{k=1}^{\infty} \frac{(-1)^{k+1}}{k.k!} \sum_{j=0}^{k} \frac{s(k,j)u^{j+2}}{(j+1)(j+2)}$$

and thus using (4.3.85) we get

$$\frac{1}{2}[\log(2\pi) - 1]u - \frac{1}{2}u^2 + u\log\Gamma(u+1) - \log G(u+1) + \frac{1}{2}\gamma u^2 = \sum_{k=1}^{\infty}\frac{(-1)^{k+1}}{k.k!}\sum_{j=0}^{k}\frac{s(k,j)u^{j+2}}{(j+1)(j+2)}$$

Then using (4.3.87e) and (4.3.87f) we see after a little simplification that

(4.3.87h)     $$\sum_{k=2}^{\infty}(-1)^k\frac{\varsigma(k)}{k(k+1)}u^{k+1} = \sum_{k=1}^{\infty}\frac{(-1)^{k+1}}{k.k!}\sum_{j=0}^{k}\frac{s(k,j)u^{j+2}}{(j+1)(j+2)}$$

and unfortunately the interesting terms involving $\log(2\pi)$ and $\gamma$ cancel out.

$\square$

We see from (4.3.79a) that

$$\phi_n(1+x) = \int_0^{1+x}t^{n-1}\psi(1+t)dt$$

$$= \int_0^1 t^{n-1}\psi(1+t)dt + \int_1^{1+x}t^{n-1}\psi(1+t)dt$$

$$= \int_0^1 t^{n-1}\psi(1+t)dt + \int_0^x(t+1)^{n-1}\psi(2+t)dt$$

$$= \phi_n(1) + \frac{(x+1)^{n-1}}{n-1} + \int_0^x(t+1)^{n-1}\psi(1+t)dt$$

and we end up with Snowden's difference equation for $n \geq 2$

(4.3.87i)     $$\phi_n(1+x) = \phi_n(1) + \frac{(x+1)^{n-1}}{n-1} + \sum_{k=0}^{n-1}\binom{n-1}{k}\phi_{k+1}(x)$$

## A DIFFERENT VIEW OF THE CRIME SCENE

Looking ahead to (4.4.11d) we see that for $n \geq 1$

$$\int_0^1 t^{x-1}(1-wyt)^n dt = \sum_{k=0}^{n}\binom{n}{k}(-1)^k\frac{w^k y^k}{k+x}$$

Differentiating with respect to $w$ gives us

(4.3.88) $\qquad -ny\int_0^1 t^x(1-wyt)^{n-1}\,dt = \sum_{k=0}^n \binom{n}{k}(-1)^k \dfrac{k\,w^{k-1}y^k}{k+x}$

Letting $w=y=1$ we get

(4.3.89) $\qquad -n\int_0^1 t^x(1-t)^{n-1}\,dt = \sum_{k=0}^n \binom{n}{k}(-1)^k \dfrac{k}{k+x}$

$$= \sum_{k=0}^n \binom{n}{k}(-1)^k\left[1-\dfrac{x}{k+x}\right]$$

We have $\sum_{k=0}^n \binom{n}{k}(-1)^k = \delta_{0,n}$ and for $n\geq 1$ we have

$$= -\sum_{k=0}^n \binom{n}{k}(-1)^k \dfrac{x}{k+x} = -x\sum_{k=0}^n \binom{n}{k}\dfrac{(-1)^k}{k+x}$$

and, making reference to (4.4.1), this becomes

(4.3.89a) $\qquad \sum_{k=0}^n \binom{n}{k}(-1)^k \dfrac{k}{k+x} = \delta_{0,n} - x\int_0^1 t^{x-1}(1-t)^n\,dt$

$$= -xB(x,n+1)\ \text{ for } n\geq 1$$

This is equivalent to the Beta function identity [8a, p.47]

(4.3.90) $\qquad B(x+1,n) = \dfrac{x}{n}B(x,n+1)$

an alternative proof of which is shown below

$$B(p,q+1) = \dfrac{\Gamma(p)\Gamma(q+1)}{\Gamma(p+q+1)}$$

$$= \dfrac{q\Gamma(p)\Gamma(q)}{p\Gamma([p+1]q)}$$

$$= \dfrac{q}{p}B(p+1,q)$$

Therefore we have

(4.3.91)
$$\sum_{k=0}^{n} \binom{n}{k} (-1)^k \frac{k}{k+x} = \delta_{0,n} - x \sum_{k=0}^{n} \binom{n}{k} \frac{(-1)^k}{k+x}$$

Having regard to (4.4.8) we also see that

(4.3.92)
$$\sum_{k=0}^{n} \binom{n}{k} (-1)^k \frac{k}{k+x} = \frac{n\Gamma(x+1)\Gamma(n)}{\Gamma(x+1+n)} = \frac{\Gamma(x+1)\Gamma(n+1)}{\Gamma(x+1+n)}$$

Differentiating (4.3.88) with respect to $y$ gives us

$$-n\int_0^1 t^x (1-wyt)^{n-1} dt + (n-1)wy \int_0^1 t^{x+1} (1-wyt)^{n-2} dt = \sum_{k=0}^{n} \binom{n}{k} (-1)^k \frac{k^2 w^{k-1} y^{k-1}}{k+x}$$

Letting $w = y = 1$ then results in the combinatorial identity

(4.3.93) $\quad -n\int_0^1 t^x (1-t)^{n-1} dt + (n-1)\int_0^1 t^{x+1} (1-t)^{n-2} dt = \sum_{k=0}^{n} \binom{n}{k} (-1)^k \frac{k^2}{k+x}$

or equivalently for $n \geq 2$

(4.3.94)
$$n(n-1)B(x+1,n) - nB(x+2,n-1) = \sum_{k=0}^{n} \binom{n}{k} (-1)^k \frac{k^2}{k+x}$$

We see that

$$\frac{k^2}{k+x} = k - x + \frac{x^2}{k+x}$$

and we therefore have

$$\sum_{k=0}^{n} \binom{n}{k} (-1)^k \frac{k^2}{k+x} = \sum_{k=0}^{n} \binom{n}{k} (-1)^k k - x \sum_{k=0}^{n} \binom{n}{k} (-1)^k + x^2 \sum_{k=0}^{n} \binom{n}{k} \frac{(-1)^k}{k+x}$$

Using the binomial theorem we see, as in (3.11aa), that

$$\sum_{k=0}^{n} \binom{n}{k} (-1)^k k = -n\delta_{1,n}$$

and hence we get

(4.3.95) $\quad \sum_{k=0}^{n}\binom{n}{k}(-1)^k \frac{k^2}{k+x} = -n\delta_{1,n} - x\delta_{0,n} + x^2 \sum_{k=0}^{n}\binom{n}{k}\frac{(-1)^k}{k+x}$

$$= -n\delta_{1,n} - x\delta_{0,n} + x^2 \int_0^1 t^{x-1}(1-t)^n dt$$

Integrating (4.3.89a) gives us

$$\sum_{k=0}^{n}\binom{n}{k}(-1)^k k \log\frac{k+u}{k+1} = \delta_{0,n}(u-1) - \int_1^u x dx \int_0^1 t^{x-1}(1-t)^n dt$$

We have

$$\int_1^u x t^{x-1} dx = \frac{t^{x-1}(x\log t - 1)}{\log^2 t}\bigg|_1^u = \frac{(ut^{u-1}-1)}{\log t} - \frac{(t^{u-1}-1)}{\log^2 t}$$

and hence we get

(4.3.95a) $\quad \sum_{k=0}^{n}\binom{n}{k}(-1)^k k \log\frac{k+u}{k+1} = \delta_{0,n}(u-1) + \int_0^1\left[\frac{(t^{u-1}-1)}{\log^2 t} - \frac{(ut^{u-1}-1)}{\log t}\right](1-t)^n dt$

Having regard to (4.3.74) we obtain

$$\sum_{n=0}^{\infty}\frac{1}{n+1}\sum_{k=0}^{n}\binom{n}{k}(-1)^k k \log\frac{k+u}{k+1} = u-1+\sum_{n=0}^{\infty}\frac{1}{n+1}\int_0^1\left[\frac{(t^{u-1}-1)}{\log^2 t} - \frac{(ut^{u-1}-1)}{\log t}\right](1-t)^n dt$$

$$= u-1+\sum_{n=0}^{\infty}\int_0^1 y^n dy \int_0^1\left[\frac{(t^{u-1}-1)}{\log^2 t} - \frac{(ut^{u-1}-1)}{\log t}\right](1-t)^n dt$$

We see that

$$\sum_{n=0}^{\infty}\int_0^1 y^n (1-t)^n dy = \int_0^1 \frac{dy}{1-(1-t)y} = -\frac{\log t}{1-t}$$

and thus we get

$$\sum_{n=0}^{\infty}\int_0^1 y^n dy \int_0^1\left[\frac{(t^{u-1}-1)}{\log^2 t} - \frac{(ut^{u-1}-1)}{\log t}\right](1-t)^n dt = -\int_0^1\left[\frac{(t^{u-1}-1)}{\log t} - (ut^{u-1}-1)\right]\frac{dt}{1-t}$$

Accordingly we deduce that

(4.3.96) $\sum_{n=0}^{\infty} \frac{1}{n+1} \sum_{k=0}^{n} \binom{n}{k} (-1)^k k \log \frac{k+u}{k+1} = u - 1 - \int_0^1 \left[ \frac{(t^{u-1}-1)}{\log t} - (ut^{u-1}-1) \right] \frac{dt}{1-t}$

and reference to (4.3.75) gives us

$$\log \Gamma(u) - u\psi(u) - \gamma - 1 + u = \sum_{k=0}^{\infty} \frac{1}{k+1} \sum_{j=0}^{k} \binom{k}{j} (-1)^j j \log \frac{(u+j)}{(1+j)}$$

and hence we get

(4.3.97) $\log \Gamma(u) - u\psi(u) - \gamma = -\int_0^1 \left[ \frac{(t^{u-1}-1)}{\log t} - (ut^{u-1}-1) \right] \frac{dt}{1-t}$

We note from (E.57) in Volume VI that

$$\log \Gamma(u) - \log \Gamma(a) + \gamma(u-a) = \int_0^1 \left[ \frac{u-a}{1-t} - \frac{t^{u-1}-t^{a-1}}{(1-t)\log t} \right] dt$$

and with $a = 1$

(4.3.97a) $\log \Gamma(u) + \gamma(u-1) = \int_0^1 \left[ \frac{u-1}{1-t} - \frac{t^{u-1}-1}{(1-t)\log t} \right] dt$

Subtracting (4.3.97) from (4.3.97a) gives us back (4.3.68b)

$$\psi(u) + \gamma = \int_0^1 \frac{1-t^{u-1}}{1-t} dt$$

Differentiating (4.3.97) results in another proof of (4.3.68).

(4.3.98) $\psi'(u) = -\int_0^1 \frac{t^{u-1}\log t}{1-t} dt$

The Wolfram website for the gamma function gives the following integral which is attributed to Radović

$$\log \Gamma(u) = \int_0^1 \frac{t^u - u(t-1) - 1}{(t-1)\log t} dt$$

However, it appears that the correct representation is

(4.3.98a)
$$\log \Gamma(u) = \int_0^1 \frac{t^{u-1} - u(t-1) - 1}{(t-1)\log t} dt$$

and differentiation gives us

(4.3.98b)
$$\psi(u) = \int_0^1 \frac{t^{u-1}\log t - (t-1)}{(t-1)\log t} dt$$

(4.3.98c)
$$\psi(1) = -\gamma = \int_0^1 \frac{\log t - (t-1)}{(t-1)\log t} dt$$

Inserting (4.3.98a), (4.3.98b) and (4.3.98c) in (4.3.97) we get

$$\int_0^1 t \left[ \frac{(t^{u-1} - 1)}{\log t} - (ut^{u-1} - 1) \right] \frac{dt}{1-t} = -\int_0^1 \left[ \frac{(t^{u-1} - 1)}{\log t} - (ut^{u-1} - 1) \right] \frac{dt}{1-t}$$

and this is easily verified by integration by parts.

We recall (4.3.74)

$$\psi(u) = \sum_{k=0}^{\infty} \frac{1}{k+1} \sum_{j=0}^{k} \binom{k}{j}(-1)^j \log(u+j)$$

and the well-known integral

$$\log x = \int_0^1 \frac{t^{x-1} - 1}{\log t} dt$$

Letting $x = u + j$ we obtain

$$\psi(u) = \sum_{k=0}^{\infty} \frac{1}{k+1} \sum_{j=0}^{k} \binom{k}{j}(-1)^j \int_0^1 \frac{t^j t^{u-1} - 1}{\log t} dt$$

and we have

$$\sum_{j=0}^{k} \binom{k}{j}(-1)^j \int_0^1 \frac{t^j t^{u-1} - 1}{\log t} dt = \int_0^1 \sum_{j=0}^{k} \binom{k}{j}(-1)^j \frac{t^j t^{u-1} - 1}{\log t} dt$$

We see that

$$\sum_{j=0}^{k} \binom{k}{j} (-1)^j \left( t^j t^{u-1} - 1 \right) = t^{u-1} \left( 1 - t \right)^k - \delta_{k,0}$$

and hence we get

$$\psi(u) = \int_0^1 \sum_{k=0}^{\infty} \frac{1}{k+1} \left( \frac{t^{u-1} \left( 1 - t \right)^k - \delta_{k,0}}{\log t} \right) dt$$

which concurs with (4.3.98b).

Integrating (4.3.97) gives us

(4.3.98d) $\qquad 2 \int_0^x \log \Gamma(u) \, du - x \log \Gamma(x) - \gamma x = -\int_0^1 \left[ 2 \frac{(t^{x-1} - t^{-1})}{\log^2 t} - \frac{x}{\log t} - \frac{x t^{x-1}}{\log t} + x \right] \frac{dt}{1-t}$

and with $x = 1$ we get

(4.3.98e) $\qquad \int_0^1 \log \Gamma(u) \, du - \gamma = -\int_0^1 \left[ \frac{(1 - t^{-1})}{\log^2 t} - \frac{1}{\log t} + \frac{1}{2} \right] \frac{dt}{1-t}$

We may also consider the following summation with (4.3.95a)

$$\sum_{n=0}^{\infty} v^n \sum_{k=0}^{n} \binom{n}{k} (-1)^k k \log \frac{k+u}{k+1} = \sum_{n=0}^{\infty} v^n \delta_{0,n}(u-1) + \sum_{n=0}^{\infty} v^n \int_0^1 \left[ \frac{(t^{u-1} - 1)}{\log^2 t} - \frac{(ut^{u-1} - 1)}{\log t} \right] (1-t)^n \, dt$$

which gives us via the geometric series

(4.3.99) $\qquad \sum_{n=0}^{\infty} v^n \sum_{k=0}^{n} \binom{n}{k} (-1)^k k \log \frac{k+u}{k+1} = u - 1 + \int_0^1 \left[ \frac{(t^{u-1} - 1)}{\log^2 t} - \frac{(ut^{u-1} - 1)}{\log t} \right] \frac{dt}{1 - v(1-t)}$

and with $v = 1/2$ we get (see also (4.4.112z) in Volume III)

(4.3.100) $\qquad \sum_{n=0}^{\infty} \frac{1}{2^{n+1}} \sum_{k=0}^{n} \binom{n}{k} (-1)^k k \log \frac{k+u}{k+1} = u - 1 + \int_0^1 \left[ \frac{(t^{u-1} - 1)}{\log^2 t} - \frac{(ut^{u-1} - 1)}{\log t} \right] \frac{dt}{1+t}$

Differentiating (4.3.100) results in

(4.3.101) $\qquad \sum_{n=0}^{\infty} \frac{1}{2^{n+1}} \sum_{k=0}^{n} \binom{n}{k} (-1)^k \frac{k}{k+u} = 1 - u \int_0^1 \frac{t^{u-1}}{1+t} \, dt$

and thus

$$\sum_{n=0}^{\infty}\frac{1}{2^{n+1}}\sum_{k=0}^{n}\binom{n}{k}(-1)^k - u\sum_{n=0}^{\infty}\frac{1}{2^{n+1}}\sum_{k=0}^{n}\binom{n}{k}\frac{(-1)^k}{k+u} = 1 - u\int_0^1\frac{t^{u-1}}{1+t}\,dt$$

Hence we have

$$(4.3.101a) \qquad \sum_{n=0}^{\infty}\frac{1}{2^{n+1}}\sum_{k=0}^{n}\binom{n}{k}\frac{(-1)^k}{k+u} = \int_0^1\frac{t^{u-1}}{1+t}\,dt$$

which we shall also see in (4.4.117l). Further differentiations will result in

$$(4.3.102) \qquad (-1)^p\,p!\sum_{n=0}^{\infty}\frac{1}{2^{n+1}}\sum_{k=0}^{n}\binom{n}{k}\frac{(-1)^k}{(k+u)^{p+1}} = \int_0^1\frac{t^{u-1}\log^p t}{1+t}\,dt$$

which we will also see in (4.4.13c). Using the Hasse/Sondow formula (3.11), we see that

$$(4.3.102a) \qquad \varsigma_a(p+1) = \sum_{n=0}^{\infty}\frac{1}{2^{n+1}}\sum_{k=0}^{n}\binom{n}{k}\frac{(-1)^k}{(k+1)^{p+1}} = \frac{(-1)^p}{p!}\int_0^1\frac{\log^p t}{1+t}\,dt$$

Integrating (4.3.95) gives us

$$\sum_{k=0}^{n}\binom{n}{k}(-1)^k\,k^2\log\frac{k+u}{k+1} = -n(u-1)\delta_{1,n} - \frac{1}{2}(u^2-1)\delta_{0,n} + \int_1^u x^2\,dx\int_0^1 t^{x-1}(1-t)^n\,dt$$

We have using integration by parts

$$\int_1^u x^2 t^{x-1}\,dx = \frac{x^2 t^{x-1}}{\log t} - \frac{2x t^{x-1}}{\log^2 t} + \frac{2 t^{x-1}}{\log^3 t}\bigg|_1^u = \frac{(u^2 t^{u-1}-1)}{\log t} - \frac{2(u t^{u-1}-1)}{\log^2 t} + \frac{2(t^{u-1}-1)}{\log^3 t}$$

and hence we get

$$\sum_{k=0}^{n}\binom{n}{k}(-1)^k\,k^2\log\frac{k+u}{k+1} =$$

$$-n(u-1)\delta_{1,n} - \frac{1}{2}(u^2-1)\delta_{0,n} + \int_0^1\left[\frac{(u^2 t^{u-1}-1)}{\log t} - \frac{2(u t^{u-1}-1)}{\log^2 t} + \frac{2(t^{u-1}-1)}{\log^3 t}\right](1-t)^n\,dt$$

As before we obtain on summation

$$\sum_{n=0}^{\infty} \frac{1}{n+1} \sum_{k=0}^{n} \binom{n}{k} (-1)^k k^2 \log \frac{k+u}{k+1} = -\sum_{n=0}^{\infty} \frac{1}{n+1} \left[ n(u-1)\delta_{1,n} + \frac{1}{2}(u^2-1)\delta_{0,n} \right]$$

$$-\int_0^1 \left[ (u^2 t^{u-1}-1) - \frac{2(ut^{u-1}-1)}{\log t} + \frac{2(t^{u-1}-1)}{\log^2 t} \right] \frac{dt}{1-t}$$

Since

$$-\sum_{n=0}^{\infty} \frac{1}{n+1} \left[ n(u-1)\delta_{1,n} + \frac{1}{2}(u^2-1)\delta_{0,n} \right] = -\frac{1}{2}(u^2-1) - \frac{1}{2}(u-1) = -\frac{1}{2}(u-1)(u+2)$$

we then determine that

(4.3.103)
$$\sum_{n=0}^{\infty} \frac{1}{n+1} \sum_{k=0}^{n} \binom{n}{k} (-1)^k k^2 \log \frac{k+u}{k+1} =$$

$$-\frac{1}{2}(u-1)(u+2) - \int_0^1 \left[ (u^2 t^{u-1}-1) - \frac{2(ut^{u-1}-1)}{\log t} + \frac{2(t^{u-1}-1)}{\log^2 t} \right] \frac{dt}{1-t}$$

(and we may also compare this with (4.3.84)).

A differentiation of (4.3.103) gives us

(4.3.104)
$$\sum_{n=0}^{\infty} \frac{1}{n+1} \sum_{k=0}^{n} \binom{n}{k} (-1)^k \frac{k^2}{k+u} = -u - \frac{1}{2} - u^2 \int_0^1 \frac{t^{u-1} \log t \, dt}{1-t}$$

and letting $u=1$ produces

(4.3.104a)
$$\sum_{n=0}^{\infty} \frac{1}{n+1} \sum_{k=0}^{n} \binom{n}{k} (-1)^k \frac{k^2}{k+1} = -\frac{3}{2} - \int_0^1 \frac{\log t \, dt}{1-t} = \varsigma(2) - \frac{3}{2}$$

We see that

$$\sum_{n=0}^{\infty} \frac{1}{n+1} \sum_{k=0}^{n} \binom{n}{k} (-1)^k \frac{k^2}{k+1} = \sum_{n=0}^{\infty} \frac{1}{n+1} \sum_{k=0}^{n} \binom{n}{k} (-1)^k (k-1) + \sum_{n=0}^{\infty} \frac{1}{n+1} \sum_{k=0}^{n} \binom{n}{k} \frac{(-1)^k}{k+1}$$

and, employing the Hasse formula (3.12), we see that this concurs with (4.3.104a).

Substituting (4.3.103) in (4.3.77) gives us

(4.3.105)

$$2\int_1^u \log\Gamma(x)dx = -\int_0^1 \left[ (u^2t^{u-1}-1) - \frac{2(ut^{u-1}-1)}{\log t} + \frac{2(t^{u-1}-1)}{\log^2 t} \right] \frac{dt}{1-t} + 2u\log\Gamma(u) - u^2\psi(u) - \gamma$$

Differentiating the above gives us (4.3.68).

Combining (4.3.84) and (4.3.103) we see that

(4.3.105a) $(u-1)\log(2\pi) - 2\log G(1+u) + u^2\psi(u) + \gamma - u(u-1) =$

$$-\int_0^1 \left[ (u^2t^{u-1}-1) - \frac{2(ut^{u-1}-1)}{\log t} + \frac{2(t^{u-1}-1)}{\log^2 t} \right] \frac{dt}{1-t}$$

One could also integrate the above to obtain another integral expression for $\int_0^x \log G(1+u)\,du$ .

It is easy to see that (4.3.105a) works with $u=1$. When $u=2$ we have

$$\log(2\pi) - 3\gamma + 2 = -\int_0^1 \left[ (4t-1) - \frac{2(2t-1)}{\log t} + \frac{2(t-1)}{\log^2 t} \right] \frac{dt}{1-t}$$

With reference to (4.3.68b) and (4.3.97)

$$\psi(u) + \gamma = -\int_0^1 \frac{t^{u-1}-1}{1-t}dt$$

$$\log\Gamma(u) - u\psi(u) - \gamma = -\int_0^1 \left[ \frac{(t^{u-1}-1)}{\log t} - (ut^{u-1}-1) \right] \frac{dt}{1-t}$$

we may write them as

$$u^2[\psi(u) + \gamma] = -\int_0^1 \frac{u^2[t^{u-1}-1]}{1-t}dt$$

$$-2u\log\Gamma(u) + 2u^2\psi(u) + 2\gamma u = \int_0^1 \left[ \frac{2u(t^{u-1}-1)}{\log t} - 2u(ut^{u-1}-1) \right] \frac{dt}{1-t}$$

Referring to the right-hand side of (4.3.105a), simple algebra gives us

$$-\int_0^1 \left[ (u^2 t^{u-1} - 1) - \frac{2(ut^{u-1} - 1)}{\log t} + \frac{2(t^{u-1} - 1)}{\log^2 t} \right] \frac{dt}{1-t} =$$

$$-\int_0^1 \left[ u^2[t^{u-1} - 1] + u^2 - 1 - 2\frac{u(t^{u-1} - 1)}{\log t} - 2\frac{u-1}{\log t} + \frac{2(t^{u-1} - 1)}{\log^2 t} \right] \frac{dt}{1-t}$$

$$= -\int_0^1 \left[ u^2[t^{u-1} - 1] + u^2 - 1 - 2u(ut^{u-1} - 1) - 2\frac{u(t^{u-1} - 1)}{\log t} + 2u(ut^{u-1} - 1) - 2\frac{u-1}{\log t} + \frac{2(t^{u-1} - 1)}{\log^2 t} \right] \frac{dt}{1-t}$$

$$= u^2[\psi(u) + \gamma] - 2u \log \Gamma(u) + 2u^2 \psi(u) + 2\gamma u$$

$$-\int_0^1 \left[ u^2 - 1 - 2u(ut^{u-1} - 1) - 2\frac{u-1}{\log t} + \frac{2(t^{u-1} - 1)}{\log^2 t} \right] \frac{dt}{1-t}$$

Therefore we have from the left-hand side of (4.3.105a)

$$(u-1)\log(2\pi) - 2\log G(1+u) + u^2 \psi(u) + \gamma - u(u-1) =$$

$$u^2[\psi(u) + \gamma] - 2u \log \Gamma(u) + 2u^2 \psi(u) + 2\gamma u - \int_0^1 \left[ u^2 - 1 - 2u(ut^{u-1} - 1) - 2\frac{u-1}{\log t} + \frac{2(t^{u-1} - 1)}{\log^2 t} \right] \frac{dt}{1-t}$$

and hence

$$(u-1)\log(2\pi) - 2\log G(1+u) + 2u \log \Gamma(u) - 2u^2 \psi(u) - \gamma(u^2 + 2u - 1) - \frac{1}{2}u(u-1) =$$

$$-\int_0^1 \left[ u^2 - 1 - 2u(ut^{u-1} - 1) - 2\frac{u-1}{\log t} + \frac{2(t^{u-1} - 1)}{\log^2 t} \right] \frac{dt}{1-t}$$

We see that

$$-\int_0^1 \left[ u^2 - 1 - 2u(ut^{u-1} - 1) - 2\frac{u-1}{\log t} + \frac{2(t^{u-1} - 1)}{\log^2 t} \right] \frac{dt}{1-t} =$$

$$-\int_0^1 \left[ u^2 - 1 - 2u^2(t^{u-1} - 1) + 2u - 2u^2 - 2\frac{u-1}{\log t} + \frac{2(t^{u-1} - 1)}{\log^2 t} \right] \frac{dt}{1-t}$$

$$= -2u^2[\psi(u) + \gamma] - \int_0^1 \left[ -(u-1)^2 - 2\frac{u-1}{\log t} + \frac{2(t^{u-1}-1)}{\log^2 t} \right] \frac{dt}{1-t}$$

Finally, we obtain

$$(u-1)\log(2\pi) - 2\log G(1+u) + 2u\log\Gamma(u) + \gamma(u-1)^2 - u(u-1) =$$

$$= -\int_0^1 \left[ -(u-1)^2 - 2\frac{u-1}{\log t} + \frac{2(t^{u-1}-1)}{\log^2 t} \right] \frac{dt}{1-t}$$

Write the above equation as

$$F(u) - 2\log G(1+u) = -\int_0^1 \left[ -(u-1)^2 - 2\frac{u-1}{\log t} + \frac{2(t^{u-1}-1)}{\log^2 t} \right] \frac{dt}{1-t}$$

and differentiate to obtain

$$F'(u) - 2\frac{d}{du}\log G(1+u) = -2\int_0^1 \left[ -(u-1) + \frac{t^{u-1}-1}{\log t} \right] \frac{dt}{1-t}$$

$$F''(u) - 2\frac{d^2}{du^2}\log G(1+u) = 2\int_0^1 \frac{1-t^{u-1}}{1-t} dt$$

$$= 2[\psi(u) + \gamma]$$

We see that

$$F^{(2)}(u) - 2\frac{d^2}{du^2}\log G(1+u) = 2[\psi(u) + \gamma]$$

and integration over the interval $[1, x]$ results in

$$(u-1)\log(2\pi) - 2\log G(1+u) + 2u\log\Gamma(u) + \gamma(u-1)^2 - u(u-1)$$

$$= 2\int_1^u \log\Gamma(x)dx + \gamma(u^2-1) + c_1(u-1) + c_2$$

The integration constants are easily found by letting $u = 0,1$ and we get

$$= 2\int_1^u \log\Gamma(x)\,dx + \gamma(u-1)^2$$

This in fact comprises another proof of Alexeiewsky's theorem (4.3.85).

$$\frac{1}{2}(u-1)\log(2\pi) - \log G(u) - \log\Gamma(u) + u\log\Gamma(u) - u(u-1) = \int_1^u \log\Gamma(x)\,dx$$

In their paper, "Massive-Scalar Effective Actions on Anti-de Sitter Spacetime", Kamela and Burgess [82aa] prove that

$$(4.3.105b) \quad \log G_n(x+1) = (-1)^n \int_0^\infty \frac{e^{-t}}{t}\left[\frac{1-e^{-xt}}{(1-e^{-t})^n} + \sum_{m=1}^n \frac{(-1)^m}{(1-e^{-t})^{n-m}}\binom{x}{m}\right]dt$$

which, on a change of variables, becomes

$$(4.3.105c) \quad \log G_n(x+1) = (-1)^n \int_0^1 \frac{1}{\log u}\left[\frac{1-u^x}{(1-u)^n} + \sum_{m=1}^n \frac{(-1)^m}{(1-u)^{n-m}}\binom{x}{m}\right]du$$

$\square$

We now recall the Hasse formula (3.12a) from Volume I

$$(4.3.106a) \qquad \varsigma(s,u) = \frac{1}{s-1}\sum_{n=0}^\infty \frac{1}{n+1}\sum_{k=0}^n \binom{n}{k}\frac{(-1)^k}{(u+k)^{s-1}}$$

where $\varsigma(s,u)$ is the Hurwitz zeta function. The above formula is valid for all $s$ except $s=1$. Note that $\varsigma(s,u)$ is not defined at $u\in\mathbf{Z}_0^-$ if $s>1$.

Differentiation of the Hasse formula (4.3.106a) with respect to $s$ results in

(4.3.107)

$$\frac{\partial}{\partial s}\varsigma(s,u) = \varsigma'(s,u)$$

$$= -\frac{1}{(s-1)^2}\sum_{n=0}^\infty \frac{1}{n+1}\sum_{k=0}^n \binom{n}{k}\frac{(-1)^k}{(u+k)^{s-1}} - \frac{1}{s-1}\sum_{n=0}^\infty \frac{1}{n+1}\sum_{k=0}^n \binom{n}{k}(-1)^k \frac{\log(u+k)}{(u+k)^{s-1}}$$

$$= -\frac{\varsigma(s,u)}{s-1} - \frac{1}{s-1}\sum_{n=0}^\infty \frac{1}{n+1}\sum_{k=0}^n \binom{n}{k}(-1)^k \frac{\log(u+k)}{(u+k)^{s-1}}$$

This may be written as

(4.3.107a) $\qquad (s-1)\varsigma'(s,u) + \varsigma(s,u) = -\sum_{n=0}^{\infty} \frac{1}{n+1} \sum_{k=0}^{n} \binom{n}{k} (-1)^k \frac{\log(u+k)}{(u+k)^{s-1}}$

Evaluations at $s=0$ and $s=-1$ produce

(4.3.108) $\qquad \varsigma'(0,u) = \varsigma(0,u) + \sum_{n=0}^{\infty} \frac{1}{n+1} \sum_{k=0}^{n} \binom{n}{k} (-1)^k (u+k) \log(u+k)$

(4.3.109) $\qquad \varsigma'(-1,u) = \frac{1}{2} \varsigma(-1,u) + \frac{1}{2} \sum_{n=0}^{\infty} \frac{1}{n+1} \sum_{k=0}^{n} \binom{n}{k} (-1)^k (u+k)^2 \log(u+k)$

and when $u=0$ these become

(4.3.109a) $\qquad \lim_{u \to 0} [\varsigma'(0,u) - \varsigma(0,u)] = \sum_{n=0}^{\infty} \frac{1}{n+1} \sum_{k=0}^{n} \binom{n}{k} (-1)^k k \log k \quad (?)$

(4.3.109b) $\qquad \varsigma'(-1,0) = \frac{1}{2} \varsigma(-1,0) + \frac{1}{2} \sum_{n=0}^{\infty} \frac{1}{n+1} \sum_{k=0}^{n} \binom{n}{k} (-1)^k k^2 \log k$

With regard to (4.3.109b), see also (4.3.137). Note that the validity of (4.3.109a) must be highly questionable because we shall see in (4.3.116) that $\varsigma'(0,u) \to \infty$ as $u \to 0$.

We have a well known relationship between the Hurwitz zeta function and the Bernoulli polynomials $B_n(u)$ (for example, see Apostol's book [14, pp. 264-266] and Appendix E of Volume VI). This corrects the misprint in [59].

(4.3.110) $\qquad \varsigma(-m,u) = -\frac{B_{m+1}(u)}{m+1}$ for $m \in \mathbf{N}_o$

$\qquad\qquad\qquad = -\frac{1}{m+1} \sum_{n=0}^{\infty} \frac{1}{n+1} \sum_{k=0}^{n} \binom{n}{k} (-1)^k (u+k)^{m+1}$

where we have employed the Hasse representation. We also see that

$\qquad\qquad \varsigma(-m,0) = -\frac{B_{m+1}}{m+1} = -\frac{1}{m+1} \sum_{n=0}^{\infty} \frac{1}{n+1} \sum_{k=0}^{n} \binom{n}{k} (-1)^k k^{m+1}$

and thus

$$B_m = \sum_{n=0}^{\infty} \frac{1}{n+1} \sum_{k=0}^{n} \binom{n}{k} (-1)^k k^m$$

We see that $\varsigma(-2m, u) = -\dfrac{B_{2m+1}(u)}{2m+1}$ and hence $\varsigma(-2m, 0) = \varsigma(-2m, 1) = 0$.

We have the particular cases using (A.15)

(4.3.110a)   $\varsigma(0, u) = \dfrac{1}{2} - u$       $\varsigma(-1, u) = -\dfrac{1}{2}\left(u^2 - u + \dfrac{1}{6}\right)$

It is easily seen from (4.3.110) that

(4.3.110b)   $\dfrac{d}{du} \varsigma(-m, u) = -B_m(u) = m\varsigma(1-m, u)$

From (4.3.108) we have

(4.3.111)   $\varsigma'(0, u) - \varsigma(0, u) = \sum_{n=0}^{\infty} \dfrac{1}{n+1} \sum_{k=0}^{n} \binom{n}{k} (-1)^k (k+u) \log(k+u)$

and letting $u = 1$ we get

$$\varsigma'(0,1) - \varsigma(0,1) = \sum_{n=0}^{\infty} \frac{1}{n+1} \sum_{k=0}^{n} \binom{n}{k} (-1)^k \log(k+1) + \sum_{n=0}^{\infty} \frac{1}{n+1} \sum_{k=0}^{n} \binom{n}{k} (-1)^k k \log(k+1)$$

Referring to (4.3.73a) this becomes

(4.3.112)   $\varsigma'(0,1) - \varsigma(0,1) = -\gamma + \sum_{n=0}^{\infty} \dfrac{1}{n+1} \sum_{k=0}^{n} \binom{n}{k} (-1)^k k \log(k+1)$

or equivalently since $\varsigma(s, 1) = \varsigma(s)$

(4.3.112)   $\varsigma'(0) - \varsigma(0) = -\gamma + \sum_{n=0}^{\infty} \dfrac{1}{n+1} \sum_{k=0}^{n} \binom{n}{k} (-1)^k k \log(k+1) = -\dfrac{1}{2} \log(2\pi) + \dfrac{1}{2}$

where we have used (F.6) in Volume VI.

We may also note from [75aa] that the Bernoulli and Euler polynomials are given by

(4.3.112a)   $B_p(u) = \sum_{n=0}^{\infty} \dfrac{1}{n+1} \sum_{k=0}^{n} \binom{n}{k} (-1)^k (u+k)^p$

(4.3.112b)
$$E_p(u) = \sum_{n=0}^{\infty} \frac{1}{2^n} \sum_{k=0}^{n} \binom{n}{k} (-1)^k (u+k)^p$$

Since $B_p(0) = B_p(1)$ we have

(4.3.112c)
$$B_p = \sum_{n=0}^{\infty} \frac{1}{n+1} \sum_{k=0}^{n} \binom{n}{k} (-1)^k k^p = \sum_{n=0}^{\infty} \frac{1}{n+1} \sum_{k=0}^{n} \binom{n}{k} (-1)^k (1+k)^p$$

With the definition

(4.3.112d)
$$C(p,u,v) = \sum_{n=0}^{\infty} v^n \sum_{k=0}^{n} \binom{n}{k} (-1)^k (u+k)^p$$

we see that $E_p(u) = C(p, u, 1/2)$ and $B_p(u) = \int_0^1 C(p, u, v) dv$. We also have

$$\frac{\partial^{m+1}}{\partial p^{m+1}} C(p,u,v) = \sum_{n=0}^{\infty} v^n \sum_{k=0}^{n} \binom{n}{k} (-1)^k (u+k)^p \log^{m+1}(u+k)$$

$$\int_0^1 \frac{\partial^{m+1}}{\partial p^{m+1}} C(p,u,v) dv = \sum_{n=0}^{\infty} \frac{1}{n+1} \sum_{k=0}^{n} \binom{n}{k} (-1)^k (u+k)^p \log^{m+1}(u+k)$$

Then letting $p = 0$ we have the Stieltjes constants (see Volume II(b) for more details of the Stieltjes constants)

(4.3.112e)
$$\gamma_m(u) = -\frac{1}{m+1} \int_0^1 \frac{\partial^{m+1}}{\partial p^{m+1}} C(p,u,v) dv \Bigg|_{p=0} = -\frac{1}{m+1} \sum_{n=0}^{\infty} \frac{1}{n+1} \sum_{k=0}^{n} \binom{n}{k} (-1)^k \log^{m+1}(u+k)$$

It therefore appears that the function $C(p,u,v)$ merits further detailed analysis.

□

In a curious paper published in 2000 Musès [103ab] strongly argued that the discrete Bernoulli numbers should be extended to a Bernoulli function (real or complex) which he denoted by $\beta(s)$. Musès contended that the Bernoulli function could be represented by

(4.3.112f)
$$\beta(s) = 2s! \varsigma(s)(2\pi)^{-s} \cos\left(\frac{\pi s}{2}\right) \cos[\pi(1-s)]$$

and therefore if $s$ is an even integer we have

$$B_{2n} = \beta(2n) = 2(2n)!\varsigma(2n)(2\pi)^{-2n}\cos(n\pi)\cos[\pi(1-2n)]$$

which is equivalent to Euler's formula (1.7)

$$B_{2n} = \beta(2n) = 2(-1)^{n+1}(2n)!\varsigma(2n)(2\pi)^{-2n}$$

With $s = 2n+1$ Muses reports that

(4.3.112g) $\qquad B_{2n+1} = \beta(2n+1) = 2(2n+1)!\varsigma(2n+1)(2\pi)^{-2n-1}\cos\left(\frac{[2n+1]\pi}{2}\right)$

whereupon we immediately see that $B_{2n+1} = 0$ for all $n \geq 1$. Furthermore, if (4.3.112g) is valid as $n \to 0$ we have

$$B_1 = \beta(1) = \frac{1}{\pi}\lim_{n \to 0}\varsigma(2n+1)\cos\left(\frac{[2n+1]\pi}{2}\right)$$

and from Riemann's functional equation (F.4c) we may determine that

$$\lim_{s \to 1}\varsigma(s)\cos\left(\frac{s\pi}{2}\right) = -\frac{\pi}{2}$$

Hence we may determine that $B_1 = -\frac{1}{2}$ as in (A.8).

We may write (4.3.112g) in the indeterminate form

(4.3.112h) $\qquad \varsigma(2n+1) = \dfrac{(2\pi)^{2n+1}B_{2n+1}}{2(2n+1)!\cos\left(\dfrac{[2n+1]\pi}{2}\right)}$

Similarly, we may write (4.3.112f) as

(4.3.112i) $\qquad \varsigma(s) = \dfrac{(2\pi)^s \beta(s)}{2s!\cos\left(\dfrac{\pi s}{2}\right)\cos[\pi(1-s)]}$

and hence we have

$$\varsigma(2n+1) = \frac{(-1)^n (2\pi)^{2n+1}}{2(2n+1)!} \lim_{s \to 2n+1} \frac{\beta(s)}{\cos\left(\frac{\pi s}{2}\right)}$$

Using L'Hôpital's rule we have

$$\varsigma(2n+1) = \frac{(-1)^n (2\pi)^{2n+1}}{2(2n+1)!} \lim_{s \to 2n+1} \frac{\beta'(s)}{\frac{\pi}{2}\sin\left(\frac{\pi s}{2}\right)} = \frac{2(-1)^n (2\pi)^{2n}}{(2n+1)!} \lim_{s \to 2n+1} \beta'(s)$$

Having regard to (4.3.112a) let us consider

(4.3.112j) $\qquad \beta(s) = B_s(1) = \sum_{m=0}^{\infty} \frac{1}{m+1} \sum_{k=0}^{m} \binom{m}{k} (-1)^k (1+k)^s$

whereupon we have

$$\beta'(s) = \sum_{m=0}^{\infty} \frac{1}{m+1} \sum_{k=0}^{m} \binom{m}{k} (-1)^k (1+k)^s \log(1+k)$$

We therefore have

$$\varsigma(2n+1) = \frac{(-1)^n 2(2\pi)^{2n}}{(2n+1)!} \sum_{m=0}^{\infty} \frac{1}{m+1} \sum_{k=0}^{m} \binom{m}{k} (-1)^k (1+k)^s \log(1+k)$$

From (4.3.107a) we see that

$$(s-1)\varsigma'(s,u) + \varsigma(s,u) = -\sum_{m=0}^{\infty} \frac{1}{m+1} \sum_{k=0}^{m} \binom{m}{k} (-1)^k \frac{\log(u+k)}{(u+k)^{s-1}}$$

and therefore we get with $s = -2n$ and $u = 1$

$$(2n+1)\varsigma'(-2n) = \sum_{m=0}^{\infty} \frac{1}{m+1} \sum_{k=0}^{m} \binom{m}{k} (-1)^k (1+k)^{2n+1} \log(1+k)$$

since $\varsigma(-2n) = 0$. This then results in

(4.3.112k) $\qquad \varsigma(2n+1) = \frac{(-1)^n 2(2\pi)^{2n}}{(2n)!} \varsigma'(-2n)$

and this is also derived in (F.8a) in Volume VI.

Musès also contended that for negative integers

$$B_{-2n} = 2n\varsigma(2n+1)$$

and

$$B_{1-2n} = (1-2n)\varsigma(2n)$$

These results may also be formally determined by using my definition of the Bernoulli function

$$\beta(s) = B_s(1) = \sum_{m=0}^{\infty}\frac{1}{m+1}\sum_{k=0}^{m}\binom{m}{k}(-1)^k(1+k)^s$$

in conjunction with the Hasse identity (3.12)

$$\varsigma(s) = \frac{1}{s-1}\sum_{m=0}^{\infty}\frac{1}{m+1}\sum_{k=0}^{m}\binom{m}{k}\frac{(-1)^k}{(k+1)^{s-1}}$$

We may write (4.3.112f) as

$$\beta(s) = 2\Gamma(s+1)\varsigma(s)(2\pi)^{-s}\cos\left(\frac{\pi s}{2}\right)\cos[\pi(1-s)]$$

and we then have

$$\log\beta(s) = \log 2 + \log\Gamma(s+1) + \log\varsigma(s) - s\log(2\pi) + \log\cos\left(\frac{\pi s}{2}\right) + \log\cos[\pi(1-s)]$$

Logarithmic differentiation gives us

$$\frac{\beta'(s)}{\beta(s)} = \psi(s+1) + \frac{\varsigma'(s)}{\varsigma(s)} - \log(2\pi) - \frac{\pi}{2}\tan\left(\frac{\pi s}{2}\right) + \pi\tan[\pi(1-s)]$$

and we therefore have

$$\beta'(s) = \left[\psi(s+1) + \frac{\varsigma'(s)}{\varsigma(s)} - \log(2\pi) - \frac{\pi}{2}\tan\left(\frac{\pi s}{2}\right) + \pi\tan[\pi(1-s)]\right]2\Gamma(s+1)\varsigma(s)(2\pi)^{-s}\cos\left(\frac{\pi s}{2}\right)\cos[\pi(1-s)]$$

Letting $s \to 1-s$ in (4.3.112f) gives us

$$\beta(1-s) = 2\Gamma(2-s)\varsigma(1-s)(2\pi)^{s-1}\sin\left(\frac{\pi s}{2}\right)\cos(\pi s)$$

and using the Riemann functional equation (F.1) we get

$$\beta(1-s) = 2\Gamma(2-s)2(2\pi)^{-s}\Gamma(s)\cos\left(\frac{\pi s}{2}\right)\varsigma(s)(2\pi)^{s-1}\sin\left(\frac{\pi s}{2}\right)\cos(\pi s)$$

$$= \frac{1}{\pi}\Gamma(2-s)\Gamma(s)\varsigma(s)2\cos\left(\frac{\pi s}{2}\right)\sin\left(\frac{\pi s}{2}\right)\cos(\pi s)$$

$$= \frac{1}{\pi}\Gamma(2-s)\Gamma(s)\varsigma(s)\sin(\pi s)\cos(\pi s)$$

$$= \frac{1}{2\pi}\Gamma(2-s)\Gamma(s)\varsigma(s)\sin(2\pi s)$$

Since $\Gamma(2-s) = \Gamma(1+1-s) = (1-s)\Gamma(1-s)$ we have

$$\beta(1-s) = \frac{1}{2\pi}(1-s)\Gamma(s)\Gamma(1-s)\varsigma(s)\sin(2\pi s)$$

$$= \frac{1}{2}(1-s)\varsigma(s)\frac{\sin(2\pi s)}{\sin(\pi s)}$$

and therefore

(4.3.112l)            $$\beta(1-s) = (1-s)\varsigma(s)\cos(\pi s)$$

With $s = -1$ we get in agreement with (A.8)

$$\beta(2) = -2\varsigma(-1) = \frac{1}{6}$$

With $s = 1/2$ we get

$$\beta(1/2) = 0$$

which suggests that

$$\sum_{m=0}^{\infty}\frac{1}{m+1}\sum_{k=0}^{m}\binom{m}{k}(-1)^k\sqrt{1+k} = 0$$

However, it is clear that (4.3.112j) and (4.3.112l) are not consistent for all values of $s$.

It should also be noted that the paper by Musès also claims an elementary proof of the Riemann Hypothesis.

<div align="center">AN EASY PROOF OF LERCH'S IDENTITY</div>

Let us define a function $H(u)$ by

(4.3.113) $$H(u) = \varsigma'(0, u) + \frac{1}{2}\log(2\pi)$$

and, using (4.3.111), this becomes

(4.3.113a) $$H(u) = \frac{1}{2} - u + \frac{1}{2}\log(2\pi) + \sum_{n=0}^{\infty}\frac{1}{n+1}\sum_{k=0}^{n}\binom{n}{k}(-1)^k(k+u)\log(k+u)$$

Differentiation immediately gives us

$$H'(u) = -1 + \sum_{n=0}^{\infty}\frac{1}{n+1}\sum_{k=0}^{n}\binom{n}{k}(-1)^k + \sum_{n=0}^{\infty}\frac{1}{n+1}\sum_{k=0}^{n}\binom{n}{k}(-1)^k\log(k+u)$$

$$= \sum_{n=0}^{\infty}\frac{1}{n+1}\sum_{k=0}^{n}\binom{n}{k}(-1)^k\log(k+u)$$

Using (4.3.74) we then obtain

$$H'(u) = \psi(u)$$

and integration over the interval $[1, x]$ gives us

$$H(x) - H(1) = \log\Gamma(x)$$

Then, using (4.3.113) and (F.6), we see that $H(1) = 0$

My original derivation was rather more convoluted and the Gaussian scaffolding is left in place below! Differentiating the Hasse formula (3.12) we obtain

(4.3.114)
$$\varsigma'(s) = -\frac{1}{s-1}\sum_{n=0}^{\infty}\frac{1}{n+1}\sum_{k=0}^{n}\binom{n}{k}\frac{(-1)^k\log(k+1)}{(k+1)^{s-1}} - \frac{1}{(s-1)^2}\sum_{n=0}^{\infty}\frac{1}{n+1}\sum_{k=0}^{n}\binom{n}{k}\frac{(-1)^k}{(k+1)^{s-1}}$$

and with $s = 0$ we have using (F.6)

$$\varsigma'(0) = -\frac{1}{2}\log(2\pi) = \sum_{n=0}^{\infty}\frac{1}{n+1}\sum_{k=0}^{n}\binom{n}{k}(-1)^{k}(k+1)\log(k+1) - \sum_{n=0}^{\infty}\frac{1}{n+1}\sum_{k=0}^{n}\binom{n}{k}(-1)^{k}(k+1)$$

$$= \sum_{n=0}^{\infty}\frac{1}{n+1}\sum_{k=0}^{n}\binom{n}{k}(-1)^{k}(k+1)\log(k+1) - \sum_{n=0}^{\infty}\frac{1}{n+1}\sum_{k=0}^{n}\binom{n}{k}(-1)^{k}k - \sum_{n=0}^{\infty}\frac{1}{n+1}\sum_{k=0}^{n}\binom{n}{k}(-1)^{k}$$

$$= \sum_{n=0}^{\infty}\frac{1}{n+1}\sum_{k=0}^{n}\binom{n}{k}(-1)^{k}(k+1)\log(k+1) - \frac{1}{2}$$

after using (4.3.76e). Hence we see that

(4.3.115) $$\sum_{n=0}^{\infty}\frac{1}{n+1}\sum_{k=0}^{n}\binom{n}{k}(-1)^{k}(k+1)\log(k+1) = \frac{1}{2} - \frac{1}{2}\log(2\pi)$$

Then, using (4.3.74) we get

(4.3.115a) $$\sum_{n=0}^{\infty}\frac{1}{n+1}\sum_{k=0}^{n}\binom{n}{k}(-1)^{k}k\log(k+1) = \gamma + \frac{1}{2} - \frac{1}{2}\log(2\pi)$$

From (4.3.113a) we see that

(4.3.115b) $$H(1) = -\frac{1}{2} + \frac{1}{2}\log(2\pi) + \sum_{n=0}^{\infty}\frac{1}{n+1}\sum_{k=0}^{n}\binom{n}{k}(-1)^{k}(k+1)\log(k+1)$$

and we accordingly obtain $H(1) = 0$. Therefore we obtain a simple derivation of Lerch's identity

(4.3.116) $$\log\Gamma(u) = \varsigma'(0,u) - \varsigma'(0) = \varsigma'(0,u) + \frac{1}{2}\log(2\pi)$$

Lerch established the above relationship between the gamma function and the Hurwitz zeta function in 1894 (see, for example, Berndt's paper [20]). We also have

$$\lim_{u\to 0}\left[\varsigma'(0,u) + \log u\right] = \lim_{u\to 0}\left[\log\frac{\Gamma(u+1)}{\sqrt{2\pi}}\right] = -\frac{1}{2}\log(2\pi)$$

When $u = 1$ in (4.3.116) we see that

(4.3.116a) $$\varsigma'(0,1) = \varsigma'(0) = -\frac{1}{2}\log(2\pi)$$

and with $u = 2$

(4.3.116b) $\qquad \varsigma'(0,2) = \varsigma'(0) = -\frac{1}{2}\log(2\pi)$

There is in fact a much more direct proof of Lerch's identity. Using (4.3.74)

$$\psi(u) = \sum_{k=0}^{\infty} \frac{1}{k+1} \sum_{j=0}^{k} \binom{k}{j} (-1)^j \log(u+j)$$

we see from (4.3.111) that

$$\frac{d}{du} \varsigma'(0,u) = -1 + \sum_{n=0}^{\infty} \frac{1}{n+1} \sum_{k=0}^{n} \binom{n}{k} (-1)^k + \sum_{n=0}^{\infty} \frac{1}{n+1} \sum_{k=0}^{n} \binom{n}{k} (-1)^k \log(u+k)$$

$$= \sum_{n=0}^{\infty} \frac{1}{n+1} \sum_{k=0}^{n} \binom{n}{k} (-1)^k \log(u+k)$$

Therefore we see that

$$\frac{d}{du} \varsigma'(0,u) = \psi(u)$$

Integration results in

$$\varsigma'(0,u) - \varsigma'(0,1) = \log \Gamma(u)$$

Since $\varsigma(s,1) = \varsigma(s)$ we have $\varsigma'(s,1) = \varsigma'(s)$ and hence the integration constant is found to be

(4.3.117) $\qquad \varsigma'(0,1) = \varsigma'(0) = -\frac{1}{2}\log(2\pi)$

Indeed, having regard to (4.3.74a) (and this is probably the most useful format)

$$\log \Gamma(u) = \sum_{k=0}^{\infty} \frac{1}{k+1} \sum_{j=0}^{k} \binom{k}{j} (-1)^j \left[ (u+j)\log(u+j) - (1+j)\log(1+j) + 1 - u \right]$$

we see immediately from (4.3.108) and (4.3.109) that

$$\log \Gamma(u) = \varsigma'(0,u) - \varsigma'(0,1)$$

This may also be readily seen from (4.3.119).

We see from Lerch's identity that

$$\frac{d}{du}\log\Gamma(u) = \psi(u) = \frac{d}{du}\varsigma'(0,u) = \lim_{s\to 0}\frac{\partial}{\partial u}\frac{\partial}{\partial s}\varsigma(s,u)$$

$$= \lim_{s\to 0}\frac{\partial}{\partial s}\frac{\partial}{\partial u}\varsigma(s,u) = -\lim_{s\to 0}\frac{\partial}{\partial s}s\varsigma(s+1,u)$$

$$= -\lim_{s\to 0}[s\varsigma'(s+1,u)+\varsigma(s+1,u)]$$

Therefore we obtain

$$\psi(u) = -\lim_{s\to 0}[s\varsigma'(s+1,u)+\varsigma(s+1,u)]$$

or equivalently

(4.3.117a) $\qquad \psi(u) = -\lim_{s\to 1}[(s-1)\varsigma'(s,u)+\varsigma(s,u)]$

We shall also see this limit in (4.3.202). It is readily seen that

$$\psi'(u) = \lim_{u\to 0}\frac{\partial^2}{\partial u^2}\frac{\partial}{\partial s}\varsigma(s'u)$$

$$= \lim_{u\to 0}\frac{\partial}{\partial s}s(s+1)\varsigma(s+2,u)$$

$$= \lim_{u\to 0}[s(s+1)\varsigma(s+2,u)+(2s+1)\varsigma(s+2,u)]$$

and hence we get

$$\psi'(u) = \varsigma(2,u)$$

From (4.3.116) we see that

$$\lim_{u\to 0}[\log\Gamma(u)-\varsigma'(0,u)-\frac{1}{2}\log(2\pi)] = 0$$

and using L'Hôpital's rule we have

$$\lim_{u \to 0} \left[ \frac{\log \Gamma(u) - \varsigma'(0, u) - \frac{1}{2} \log(2\pi)}{u} \right] = \psi(u) - \frac{d}{du} \varsigma'(0, u) = 0$$

We may deduce from (4.3.115) and (4.3.74) that (see also (4.3.112))

$$(4.3.118) \qquad \sum_{n=0}^{\infty} \frac{1}{n+1} \sum_{k=0}^{n} \binom{n}{k} (-1)^k k \log(1+k) = \frac{1}{2} + \gamma - \frac{1}{2} \log(2\pi)$$

Hence we may write (4.3.76d) as

$$(4.3.119) \qquad \log \Gamma(u) = \sum_{n=0}^{\infty} \frac{1}{n+1} \sum_{k=0}^{n} \binom{n}{k} (-1)^k k \log(u+k) + u\psi(u) + \frac{1}{2} - u + \frac{1}{2} \log(2\pi)$$

$$= \sum_{n=0}^{\infty} \frac{1}{n+1} \sum_{k=0}^{n} \binom{n}{k} (-1)^k (u+k) \log(u+k) + \frac{1}{2} - u + \frac{1}{2} \log(2\pi)$$

and, with reference to (4.3.108) and (4.3.110a), it is easily seen that the above equation is equivalent to Lerch's identity.

It is easily seen from the above that

$$\frac{d^2}{du^2} \log \Gamma(u) = \sum_{n=0}^{\infty} \frac{1}{n+1} \sum_{k=0}^{n} \binom{n}{k} \frac{(-1)^k}{u+k} = \varsigma(2, u)$$

and therefore we have $\dfrac{d^2}{du^2} \log \Gamma(u) \geq 0$ (see [126, p.40]).

<div align="center">STIRLING'S APPROXIMATION FOR $\log \Gamma(u)$</div>

Many different derivations exist for Stirling's approximation for the gamma function and, what I believe is a new one, is presented below.

With reference to (4.3.119) we see that

$$\sum_{k=0}^{n} \binom{n}{k} (-1)^k (u+k) \log(u+k) = \log \left[ \prod_{k=0}^{n} (u+k)^{(-1)^k \binom{n}{k}(u+k)} \right]$$

and as $u \to \infty$ we see that

$$\sum_{k=0}^{n}\binom{n}{k}(-1)^k(u+k)\log(u+k) \approx \log\left[u^{\sum_{k=0}^{n}\binom{n}{k}(-1)^j(u+k)}\right]$$

Using (4.3.76e) we have

$$\sum_{k=0}^{n}\binom{n}{k}(-1)^k(u+k) = u\delta_{n,0} - n\delta_{n,1}$$

and therefore we have as $u \to \infty$

$$\sum_{n=0}^{\infty}\frac{1}{n+1}\sum_{k=0}^{n}\binom{n}{k}(-1)^k(u+k)\log(u+k) \to \sum_{n=0}^{\infty}\frac{1}{n+1}\left[u\delta_{n,0}-n\delta_{n,1}\right]\log u$$

$$= \left[u - \frac{1}{2}\right]\log u$$

Hence as $u \to \infty$ we obtain the following approximation

$$(4.3.119i) \qquad \log\Gamma(u) \approx \left[u-\frac{1}{2}\right]\log u - u + \frac{1}{2} + \frac{1}{2}\log(2\pi)$$

which, except for the curious additional factor of 1/2, is equivalent to Stirling's asymptotic approximation [126, p.8]

$$(4.3.119ii) \qquad \log\Gamma(u+1) \approx u\log u - u + \frac{1}{2}\log(2\pi u)$$

The superfluous "1/2" in (4.3.119i) niggled me for some time and many months later I noted the following asymptotic formula quoted in [6c] as $u \to \infty$

$$(4.3.119iii) \qquad \log(u+k) = \log u + \frac{k}{u} + O\left(\frac{1}{u^2}\right)$$

We then have as $u \to \infty$

$$\sum_{n=0}^{\infty}\frac{1}{n+1}\sum_{k=0}^{n}\binom{n}{k}(-1)^k(u+k)\log(u+k) \to \sum_{n=0}^{\infty}\frac{1}{n+1}\sum_{k=0}^{n}\binom{n}{k}(-1)^k\left[u\log u + k + k\log u + \frac{k^2}{u}\right]$$

and the right-hand side is equal to

$$u\log u\sum_{n=0}^{\infty}\frac{1}{n+1}\sum_{k=0}^{n}\binom{n}{k}(-1)^k + (1+\log u)\sum_{n=0}^{\infty}\frac{1}{n+1}\sum_{k=0}^{n}\binom{n}{k}(-1)^k k + \frac{1}{u}\sum_{n=0}^{\infty}\frac{1}{n+1}\sum_{k=0}^{n}\binom{n}{k}(-1)^k k^2$$

$$=u\log u - \frac{1}{2}(1+\log u)+\frac{1}{6u}=\left[u-\frac{1}{2}\right]\log u-\frac{1}{2}+\frac{1}{6u}$$

We then obtain

$$\log\Gamma(u)\approx\left[u-\frac{1}{2}\right]\log u-u+\frac{1}{2}\log(2\pi)+\frac{1}{6u}$$

and the mysterious "1/2" has fortunately disappeared. By including another term in the expansion (4.3.119iii) we may derive the more familiar form

$$\log\Gamma(u)\approx\left[u-\frac{1}{2}\right]\log u-u+\frac{1}{2}\log(2\pi)+\frac{1}{12u}$$

The asymptotic formula for $\log(u+k)$ referred to above in (4.3.119iii) is in fact simply the Maclaurin expansion which is valid for $k/u\le 1$

(4.3.119iv) $\qquad\log(u+k)-\log u=\log\left(1+\frac{k}{u}\right)=-\sum_{j=1}^{\infty}(-1)^j\frac{k^j}{ju^j}$

and using (4.3.80a) this may also be represented by

$$=-\sum_{j=1}^{N}(-1)^j\frac{k^j}{ju^j}+\frac{(-1)^N}{u^N}\int_0^k\frac{t^N}{u+t}dt$$

With the substitution $t=ky$ this becomes

$$=-\sum_{j=1}^{N}(-1)^j\frac{k^j}{ju^j}+\frac{(-1)^N k^{N+1}}{u^N}\int_0^1\frac{y^N}{u+ky}dy$$

and we have the bounds

$$\int_0^1\frac{y^N}{u+ky}dy\le\int_0^1\left|\frac{y^N}{u+ky}\right|dy\le\int_0^1\left|\frac{y^N}{ky}\right|dy=\frac{1}{kN}$$

$$\int_0^1 \frac{y^N}{u+ky} \, dy \le \int_0^1 \left| \frac{y^N}{u+ky} \right| dy \le \int_0^1 \left| \frac{y^N}{u} \right| dy = \frac{1}{u(N+1)}$$

We then obtain (it should be noted that there are major issues of convergence associated with the following presentation; rigour may be provided through an appropriate limiting procedure)

$$\sum_{n=0}^{\infty} \frac{1}{n+1} \sum_{k=0}^{n} \binom{n}{k} (-1)^k (u+k) \log(u+k) = \sum_{n=0}^{\infty} \frac{1}{n+1} \sum_{k=0}^{n} \binom{n}{k} (-1)^k (u+k) \left[ \log u - \sum_{j=1}^{\infty} (-1)^j \frac{k^j}{ju^j} \right]$$

$$= B_1(u) \log u - \sum_{n=0}^{\infty} \frac{1}{n+1} \sum_{k=0}^{n} \binom{n}{k} (-1)^k \sum_{j=1}^{\infty} (-1)^j \frac{k^{j+1}}{ju^j} - u \sum_{n=0}^{\infty} \frac{1}{n+1} \sum_{k=0}^{n} \binom{n}{k} (-1)^k \sum_{j=1}^{\infty} (-1)^j \frac{k^j}{ju^j}$$

$$= \left[ u - \frac{1}{2} \right] \log u - \sum_{n=0}^{\infty} \frac{1}{n+1} \sum_{k=0}^{n} \binom{n}{k} (-1)^k \sum_{j=1}^{\infty} (-1)^j \frac{k^{j+1}}{ju^j} - u \sum_{n=0}^{\infty} \frac{1}{n+1} \sum_{k=0}^{n} \binom{n}{k} (-1)^k \sum_{j=1}^{\infty} (-1)^j \frac{k^j}{ju^j}$$

We have

$$\sum_{n=0}^{\infty} \frac{1}{n+1} \sum_{k=0}^{n} \binom{n}{k} (-1)^k \sum_{j=1}^{\infty} (-1)^j \frac{k^{j+1}}{ju^j} = \sum_{j=1}^{\infty} \frac{(-1)^j}{ju^j} \sum_{n=0}^{\infty} \frac{1}{n+1} \sum_{k=0}^{n} \binom{n}{k} (-1)^k k^{j+1}$$

and we note from (A.23aa) in Volume VI that

$$B_j = \sum_{n=0}^{\infty} \frac{1}{n+1} \sum_{k=0}^{n} \binom{n}{k} (-1)^k k^j$$

We then have

$$\sum_{j=1}^{\infty} \frac{(-1)^j}{ju^j} \sum_{n=0}^{\infty} \frac{1}{n+1} \sum_{k=0}^{n} \binom{n}{k} (-1)^k k^{j+1} = \sum_{j=1}^{\infty} \frac{(-1)^j B_{j+1}}{ju^j}$$

and

$$u \sum_{j=1}^{\infty} \frac{(-1)^j}{ju^j} \sum_{n=0}^{\infty} \frac{1}{n+1} \sum_{k=0}^{n} \binom{n}{k} (-1)^k k^j = \sum_{j=1}^{\infty} \frac{(-1)^j B_j}{ju^{j-1}}$$

We have using even and odd Bernoulli numbers

$$\sum_{j=1}^{\infty} \frac{(-1)^j B_{j+1}}{ju^j} = \sum_{m=1}^{\infty} \frac{B_{2m+1}}{2mu^{2m}} - \sum_{m=0}^{\infty} \frac{B_{2m+2}}{(2m+1)u^{2m+1}} = -\sum_{m=0}^{\infty} \frac{B_{2m+2}}{(2m+1)u^{2m+1}}$$

$$= -\sum_{l=1}^{\infty} \frac{B_{2l}}{(2l-1)u^{2l-1}}$$

and similarly we have

$$\sum_{j=1}^{\infty} \frac{(-1)^j B_j}{ju^{j-1}} = \sum_{m=1}^{\infty} \frac{B_{2m}}{2mu^{2m-1}} - \sum_{m=0}^{\infty} \frac{B_{2m+1}}{(2m+1)u^{2m}} = \sum_{m=1}^{\infty} \frac{B_{2m}}{2mu^{2m-1}} + \frac{1}{2}$$

We then get

$$\sum_{j=1}^{\infty} \frac{(-1)^j B_{j+1}}{ju^j} + \sum_{j=1}^{\infty} \frac{(-1)^j B_j}{ju^{j-1}} = -\sum_{m=1}^{\infty} \frac{B_{2m}}{(2m-1)u^{2m-1}} + \sum_{m=1}^{\infty} \frac{B_{2m}}{2mu^{2m-1}} + \frac{1}{2}$$

$$= -\sum_{m=1}^{\infty} \frac{B_{2m}}{2m(2m-1)u^{2m-1}} + \frac{1}{2}$$

This gives us

$$\sum_{n=0}^{\infty} \frac{1}{n+1} \sum_{k=0}^{n} \binom{n}{k} (-1)^k (u+k) \log(u+k) = \left[ u - \frac{1}{2} \right] \log u + \sum_{m=1}^{\infty} \frac{B_{2m}}{2m(2m-1)u^{2m-1}} - \frac{1}{2}$$

and hence we obtain using (4.3.119)

$$(4.3.119v) \qquad \log \Gamma(u) \approx \frac{1}{2} \log(2\pi) + \left[ u - \frac{1}{2} \right] \log u - u + \sum_{m=1}^{\infty} \frac{B_{2m}}{2m(2m-1)u^{2m-1}}$$

which is also reported in [126, p.8], as corrected. Since from (1.7)

$$\varsigma(2m) = (-1)^{m+1} \frac{2^{2m-1} \pi^{2m} B_{2m}}{(2m)!}$$

we have

$$\log \Gamma(u) \approx \frac{1}{2} \log(2\pi) + \left[ u - \frac{1}{2} \right] \log u - u + \sum_{m=1}^{\infty} (-1)^{m+1} \frac{(2m)! \varsigma(2m)}{2^{2m} \pi^{2m} m(2m-1)u^{2m-1}}$$

In a similar way, it is easily shown using (4.3.74) that an asymptotic expansion for $\psi(u)$ is $\psi(u) \approx \log u$. We have [126, p.22] (as corrected)

$$(4.3.119\text{vi}) \qquad \psi(u) \approx \log u - \sum_{j=1}^{\infty} \frac{(-1)^j B_j}{ju^j}$$

$$= \log u - \frac{1}{2u} - \frac{1}{2}\sum_{j=1}^{\infty} \frac{B_{2j}}{ju^{2j}}$$

It is also possible to obtain an asymptotic expansion for $\log G(u+1)$ in the same manner as above – see (4.3.128d).

Since $\varsigma(s)$ is monotonic decreasing for $s > 1$, we see from (4.3.114) that

$$\sum_{n=0}^{\infty} \frac{1}{n+1} \sum_{k=0}^{n} \binom{n}{k} \frac{(-1)^k}{(k+1)^{s-1}} > -(s-1)\sum_{n=0}^{\infty} \frac{1}{n+1} \sum_{k=0}^{n} \binom{n}{k} \frac{(-1)^k \log(k+1)}{(k+1)^{s-1}}$$

Using (4.3.119) in the form

$$\log \Gamma(u+1) = \sum_{n=0}^{\infty} \frac{1}{n+1} \sum_{k=0}^{n} \binom{n}{k} (-1)^k (u+1+k) \log(u+1+k) - \frac{1}{2} - u + \frac{1}{2}\log(2\pi)$$

will enable us to evaluate the integrals $\int_0^x \log \Gamma(u+1) \sin 2p\pi u\, du$ and

$\int_0^x \log \Gamma(u+1) \cos 2p\pi u\, du$ as series involving the cosine and sine integrals (see (6.111) in Volume V).

It may also be noted that Dasbach [50a] has recently shown that for $(u+k) \geq 4$

$$\log(u+k) = \sum_{n=1}^{\infty} \frac{1}{n} \sum_{j=0}^{n} \binom{n}{j}\binom{2j}{j} \frac{(-1)^j}{(u+k)^j}$$

and, using suitable combinatorial identities, it may be possible to derive some interesting series from (4.3.119) etc.

With $u \to -u$ in (4.3.119) we have

$$(4.3.119\text{a}) \quad \log \Gamma(-u) = \sum_{n=0}^{\infty} \frac{1}{n+1} \sum_{k=0}^{n} \binom{n}{k} (-1)^k k \log(k-u) - u\psi(-u) + \frac{1}{2} + u + \frac{1}{2}\log(2\pi)$$

Euler's reflection formula (C.28) gives us

$$\Gamma(u)\Gamma(1-u) = \frac{\pi}{\sin \pi u}$$

$$\Gamma(1+u) = u\Gamma(u) \quad \Rightarrow \quad \Gamma(1-u) = -u\Gamma(-u)$$

and hence we have for $0 < u < 1$

(4.3.119b) $\quad \Gamma(u)\Gamma(-u) = -\dfrac{\pi}{u \sin \pi u}$ $\quad$ and $\quad$ $\log\left|\Gamma(u)\Gamma(-u)\right| = \log \dfrac{\pi}{u \sin \pi u}$

Adding the two equations (4.3.119) and (4.3.119a) together results in (provided $k^2 \neq u^2$)

(4.3.119c) $\quad \log \dfrac{\pi}{u \sin \pi u} = \displaystyle\sum_{n=0}^{\infty} \frac{1}{n+1} \sum_{k=0}^{n} \binom{n}{k} (-1)^k \, k \log(k^2 - u^2) + u\left[\psi(u) - \psi(-u)\right] + 2\log(2\pi)$

and differentiation gives us

(4.3.120)

$$-\frac{1}{u} - \pi \cot \pi u = -\sum_{n=0}^{\infty} \frac{1}{n+1} \sum_{k=0}^{n} \binom{n}{k} (-1)^k \frac{2ku}{(k^2 - u^2)} + u\left[\psi'(u) + \psi'(-u)\right] + \left[\psi(u) - \psi(-u)\right]$$

We have already seen in (4.3.86) that $-\dfrac{1}{u} - \pi \cot \pi u = \left[\psi(u) - \psi(-u)\right]$ and we have by partial fractions

$$\sum_{n=0}^{\infty} \frac{1}{n+1} \sum_{k=0}^{n} \binom{n}{k} (-1)^k \frac{2ku}{(k^2 - u^2)} = \sum_{n=0}^{\infty} \frac{1}{n+1} \sum_{k=0}^{n} \binom{n}{k} \frac{(-1)^k k}{k+u} - \sum_{n=0}^{\infty} \frac{1}{n+1} \sum_{k=0}^{n} \binom{n}{k} \frac{(-1)^k k}{k-u}$$

$$= -\sum_{n=0}^{\infty} \frac{1}{n+1} \sum_{k=0}^{n} \binom{n}{k} \frac{(-1)^k u}{k+u} + \sum_{n=0}^{\infty} \frac{1}{n+1} \sum_{k=0}^{n} \binom{n}{k} \frac{(-1)^k u}{k-u}$$

since $\dfrac{k}{k+u} = 1 - \dfrac{u}{k+u}$

Referring back to (4.3.74) we see that

(4.3.121) $$\psi'(\pm u) = \sum_{n=0}^{\infty} \frac{1}{n+1} \sum_{k=0}^{n} \binom{n}{k} \frac{(-1)^k}{k \pm u}$$

and hence it is evident that both sides of (4.3.120) are indeed equal. It is a simple exercise to reverse the order of this proof.

Multiplying (4.3.119c) by $u$, using (4.3.86) and then integrating results in a series involving $\sum_{n=0}^{\infty} \frac{1}{n+1} \sum_{k=0}^{n} \binom{n}{k} (-1)^k k(k^2 - u^2) \log(k^2 - u^2)$ which I suspect relates to the Barnes double gamma function.

With reference to (4.3.119) we see that

$$\log u = \log \Gamma(u+1) - \log \Gamma(u) = \sum_{n=0}^{\infty} \frac{1}{n+1} \sum_{k=0}^{n} \binom{n}{k} (-1)^k k \log(u+k+1)$$

$$- \sum_{n=0}^{\infty} \frac{1}{n+1} \sum_{k=0}^{n} \binom{n}{k} (-1)^k k \log(u+k) + (u+1)\psi(u+1) - 1 - u\psi(u)$$

$$= \sum_{n=0}^{\infty} \frac{1}{n+1} \sum_{k=0}^{n} \binom{n}{k} (-1)^k k \log(u+k+1)$$

$$- \sum_{n=0}^{\infty} \frac{1}{n+1} \sum_{k=0}^{n} \binom{n}{k} (-1)^k k \log(u+k) + \psi(u) + \frac{1}{u}$$

This gives us

$$\log u - \psi(u) - \frac{1}{u} = \sum_{n=0}^{\infty} \frac{1}{n+1} \sum_{k=0}^{n} \binom{n}{k} (-1)^k k \log \frac{(u+k+1)}{(u+k)}$$

and letting $u = 1$ we obtain another derivation of (4.3.76). We note that the right-hand side vanishes as $u \to \infty$ and hence we have the well-known result (E.66a)

(4.3.121aa)          $\lim_{u \to \infty} [\log u - \psi(u)] = 0$

We now consider

$$k \log \frac{(u+k+1)}{(u+k)} = k \log(u+k+1) - k \log(u+k)$$

$$= (u+k+1)\log(u+k+1) - (u+1)\log(u+k+1) - (u+k)\log(u+k) + u\log(u+k)$$

$$= (u+k+1)\log(u+k+1) - (u+k)\log(u+k) - (u+1)\log(u+k+1) + u\log(u+k)$$

and thus we see from (4.3.107a) that

$$\sum_{n=0}^{\infty}\frac{1}{n+1}\sum_{k=0}^{n}\binom{n}{k}(-1)^{k}k\log\frac{(u+k+1)}{(u+k)}=-\lim_{s\to 0}[(s-1)\varsigma'(s,u+1)+\varsigma(s,u+1)]$$

$$+\lim_{s\to 0}[(s-1)\varsigma'(s,u)+\varsigma(s,u)]-(u+1)\psi(u+1)+u\psi(u)$$

$$=\varsigma'(0,u+1)-\varsigma(0,u+1)-\varsigma'(0,u)+\varsigma(0,u)-(u+1)\psi(u+1)+u\psi(u)$$

Hence we deduce that

$$\log u-\psi(u)-\frac{1}{u}=\varsigma'(0,u+1)-\varsigma(0,u+1)-\varsigma'(0,u)+\varsigma(0,u)-(u+1)\psi(u+1)+u\psi(u)$$

which results in

$$\log u+(u+1)[\psi(u+1)-\psi(u)]-\frac{1}{u}=\varsigma'(0,u+1)-\varsigma'(0,u)-\varsigma(0,u+1)+\varsigma(0,u)$$

$$=\varsigma'(0,u+1)-\varsigma'(0,u)+1$$

Since $\psi(u+1)-\psi(u)=\frac{1}{u}$ we obtain for $u>0$

(4.3.121a)        $$\log u=\varsigma'(0,u+1)-\varsigma'(0,u)$$

With $u=1$ we see that

(4.3.121b)        $$\varsigma'(0,2)=\varsigma'(0,1)=\varsigma'(0)=-\frac{1}{2}\log(2\pi)$$

and we also have

(4.3.121c)        $$\log(u+a)=\varsigma'(0,u+a+1)-\varsigma'(0,u+a)$$

In fact, by subtracting the following two equations, we note that (4.3.121a) is merely a disguised form of Lerch's identity (4.3.116)

$$\log\Gamma(u)=\varsigma'(0,u)+\frac{1}{2}\log(2\pi)$$

$$\log \Gamma(u+1) = \log u + \log \Gamma(u) = \varsigma'(0, u+1) + \frac{1}{2}\log(2\pi)$$

From (4.3.121a) we see that

$$\lim_{u \to 0}[\log u + \varsigma'(0, u)] = \varsigma'(0, 1) = \varsigma'(0) = -\frac{1}{2}\log(2\pi)$$

and this implies that $\varsigma'(0, u) \to \infty$ as $u \to 0$ (as may also be seen from Lerch's identity).

We now let $u = 1/2$ in (4.3.119) to obtain

$$\log \Gamma\left(\frac{1}{2}\right) = \sum_{n=0}^{\infty}\frac{1}{n+1}\sum_{k=0}^{n}\binom{n}{k}(-1)^k k \log\left(k + \frac{1}{2}\right) + \frac{1}{2}\psi\left(\frac{1}{2}\right) + \frac{1}{2}\log(2\pi)$$

and this may be written as

$$\frac{1}{2}\log \pi = \sum_{n=0}^{\infty}\frac{1}{n+1}\sum_{k=0}^{n}\binom{n}{k}(-1)^k k \log(2k+1) - \log 2 \sum_{n=0}^{\infty}\frac{1}{n+1}\sum_{k=0}^{n}\binom{n}{k}(-1)^k k$$

$$-\frac{1}{2}\gamma - \log 2 + \frac{1}{2}\log(2\pi)$$

This then simplifies to

(4.3.121d) $$\sum_{n=0}^{\infty}\frac{1}{n+1}\sum_{k=0}^{n}\binom{n}{k}(-1)^k k \log(2k+1) = \frac{1}{2}\gamma$$

and we may note from (4.3.118) that

$$\sum_{n=0}^{\infty}\frac{1}{n+1}\sum_{k=0}^{n}\binom{n}{k}(-1)^k k \log(k+1) = \frac{1}{2} + \gamma - \frac{1}{2}\log(2\pi)$$

Subtraction results in

(4.3.121e) $$\sum_{n=0}^{\infty}\frac{1}{n+1}\sum_{k=0}^{n}\binom{n}{k}(-1)^k k \log\frac{2k+1}{k+1} = \frac{1}{2}[\log(2\pi) - \gamma - 1]$$

A further proof of (4.3.121d) is given below. We see from the Hasse representation that

$$\varsigma(s, 1/2) = \frac{2^{s-1}}{s-1} \sum_{n=0}^{\infty} \frac{1}{n+1} \sum_{k=0}^{n} \binom{n}{k} \frac{(-1)^k}{(2k+1)^{s-1}}$$

Differentiation results in

$$\varsigma'(s, 1/2) = -\frac{2^{s-1}}{s-1} \sum_{n=0}^{\infty} \frac{1}{n+1} \sum_{k=0}^{n} \binom{n}{k} \frac{(-1)^k \log(2k+1)}{(2k+1)^{s-1}} + \frac{2^{s-1}[(s-1)\log 2 - 1]}{(s-1)^2} \sum_{n=0}^{\infty} \frac{1}{n+1} \sum_{k=0}^{n} \binom{n}{k} \frac{(-1)^k}{(2k+1)^{s-1}}$$

and, with $s = 0$, we obtain

$$\varsigma'(0, 1/2) = \frac{1}{2} \sum_{n=0}^{\infty} \frac{1}{n+1} \sum_{k=0}^{n} \binom{n}{k} (-1)^k (2k+1) \log(2k+1) - \frac{1}{2}[1 + \log 2] \sum_{n=0}^{\infty} \frac{1}{n+1} \sum_{k=0}^{n} \binom{n}{k} (-1)^k (2k+1)$$

Using (4.3.76e) this simplifies to

$$\varsigma'(0, 1/2) = \frac{1}{2} \sum_{n=0}^{\infty} \frac{1}{n+1} \sum_{k=0}^{n} \binom{n}{k} (-1)^k (2k+1) \log(2k+1)$$

$$= \sum_{n=0}^{\infty} \frac{1}{n+1} \sum_{k=0}^{n} \binom{n}{k} (-1)^k k \log(2k+1) + \frac{1}{2} \sum_{n=0}^{\infty} \frac{1}{n+1} \sum_{k=0}^{n} \binom{n}{k} (-1)^k \log(2k+1)$$

Reference to (4.3.74) with $u = 1/2$ gives us

$$\psi\left(\frac{1}{2}\right) = \sum_{n=0}^{\infty} \frac{1}{n+1} \sum_{k=0}^{n} \binom{n}{k} (-1)^k \log(2k+1) - \log 2 \sum_{n=0}^{\infty} \frac{1}{n+1} \sum_{k=0}^{n} \binom{n}{k} (-1)^k$$

and we therefore have

(4.3.122)
$$\sum_{n=0}^{\infty} \frac{1}{n+1} \sum_{k=0}^{n} \binom{n}{k} (-1)^k \log(2k+1) = -\gamma - \log 2$$

Hence we obtain

$$\varsigma'(0, 1/2) = \sum_{n=0}^{\infty} \frac{1}{n+1} \sum_{k=0}^{n} \binom{n}{k} (-1)^k k \log(2k+1) - \frac{1}{2}\gamma - \frac{1}{2}\log 2$$

We already know from the Lerch identity (4.3.116) that

$$\log \Gamma(1/2) = \varsigma'(0, 1/2) + \frac{1}{2}\log(2\pi)$$

and hence we have

(4.3.122a) $\qquad \varsigma'(0,1/2) = -\frac{1}{2}\log 2$

which may also be seen by differentiating (4.3.74dii).

This then gives us in accordance with (4.3.121d)

$$\sum_{n=0}^{\infty} \frac{1}{n+1} \sum_{k=0}^{n} \binom{n}{k} (-1)^k k \log(2k+1) = \frac{1}{2}\gamma$$

We have by definition for $s > 1$

$$\varsigma(s,u) = \sum_{n=0}^{\infty} \frac{1}{(n+u)^s}$$

and thus

(4.3.123) $\varsigma(s,1/2) = 2^s \sum_{n=0}^{\infty} \frac{1}{(2n+1)^s} = 2^s \sum_{n=0}^{\infty} \left( \frac{1}{n^s} - \frac{1}{(2n)^s} \right) = (2^s - 1)\varsigma(s)$

and we deduce that for $s > 1$

(4.3.124) $\qquad \sum_{n=0}^{\infty} \frac{1}{n+1} \sum_{k=0}^{n} \binom{n}{k} \frac{(-1)^k}{(2k+1)^{s-1}} = (s-1)\left(2 - 2^{1-s}\right)\varsigma(s)$

Differentiation of (4.3.124) results in

$$-\sum_{n=0}^{\infty} \frac{1}{n+1} \sum_{k=0}^{n} \binom{n}{k} \frac{(-1)^k \log(2k+1)}{(2k+1)^{s-1}} = \left(2 - 2^{1-s}\right)\varsigma(s) + (s-1)\left(2 - 2^{1-s}\right)\varsigma'(s) + (s-1)2^{1-s}\varsigma(s)\log 2$$

and, with $s = 0$, we obtain with the use of (3.11a)

(4.3.124a) $\qquad \sum_{n=0}^{\infty} \frac{1}{n+1} \sum_{k=0}^{n} \binom{n}{k} (-1)^k (2k+1) \log(2k+1) = 2\varsigma(0)\log 2 = -\log 2$

Reference to (4.3.74) with $u = 1/2$ gives us

$$\psi\left(\frac{1}{2}\right) = \sum_{n=0}^{\infty} \frac{1}{n+1} \sum_{k=0}^{n} \binom{n}{k} (-1)^k \log(2k+1) - \log 2 \sum_{n=0}^{\infty} \frac{1}{n+1} \sum_{k=0}^{n} \binom{n}{k} (-1)^k$$

and we therefore have

$$(4.3.125) \qquad \sum_{n=0}^{\infty} \frac{1}{n+1} \sum_{k=0}^{n} \binom{n}{k} (-1)^k \log(2k+1) = -\gamma - \log 2$$

This concurs with the above analysis.

<div align="center">THE GOSPER/VARDI FUNCTIONAL EQUATION</div>

We now consider the function

$$j^2 \log \frac{(u+j)}{(1+j)} = j^2 \log(u+j) - j^2 \log(1+j)$$

Since $j^p = [(j+u)-u]^p$ we have

$$j^2 = (j+u)^2 - 2u(j+u) + u^2$$

$$j^2 = (j+1)^2 - 2(j+1) + 1$$

and see that

$$j^2 \log \frac{(u+j)}{(1+j)} = [(j+u)^2 - 2u(j+u) + u^2] \log(u+j) - [(j+1)^2 - 2(j+1) + 1] \log(1+j)$$

Therefore we get upon making the summation

$$\sum_{k=0}^{\infty} \frac{1}{k+1} \sum_{j=0}^{k} \binom{k}{j} (-1)^j j^2 \log \frac{(u+j)}{(1+j)}$$

$$= \sum_{k=0}^{\infty} \frac{1}{k+1} \sum_{j=0}^{k} \binom{k}{j} (-1)^j [(u+j)^2 \log(u+j) - (1+j)^2 \log(1+j)]$$

$$-2 \sum_{k=0}^{\infty} \frac{1}{k+1} \sum_{j=0}^{k} \binom{k}{j} (-1)^j [u(u+j) \log(u+j) - (1+j) \log(1+j)]$$

$$+ \sum_{k=0}^{\infty} \frac{1}{k+1} \sum_{j=0}^{k} \binom{k}{j} (-1)^j [u^2 \log(u+j) - \log(1+j)]$$

Referring to (4.3.108) and (4.3.109)

$$\varsigma'(0,u) - \varsigma(0,u) = \sum_{n=0}^{\infty} \frac{1}{n+1} \sum_{j=0}^{n} \binom{n}{j} (-1)^j (u+j) \log(u+j)$$

$$2\varsigma'(-1,u) - \varsigma(-1,u) = \sum_{n=0}^{\infty} \frac{1}{n+1} \sum_{j=0}^{n} \binom{n}{j} (-1)^j (u+j)^2 \log(u+j)$$

we see that

$$\sum_{k=0}^{\infty} \frac{1}{k+1} \sum_{j=0}^{k} \binom{k}{j} (-1)^j j^2 \log \frac{(u+j)}{(1+j)}$$

$$= 2\varsigma'(-1,u) - \varsigma(-1,u) - 2\varsigma'(-1,1) + \varsigma(-1,1)$$

$$- 2u\varsigma'(0,u) + 2u\varsigma(0,u) + 2\varsigma'(0,1) - 2\varsigma(0,1)$$

$$+ u^2 \psi(u) + \gamma$$

$$= 2\varsigma'(-1,u) + \frac{1}{2}\left(u^2 - u + \frac{1}{6}\right) - 2\varsigma'(-1) - \frac{1}{12}$$

$$- 2u\varsigma'(0,u) + u - 2u^2 + 2\varsigma'(0) + 1$$

$$+ u^2 \psi(u) + \gamma$$

Then using Lerch's identity

$$\varsigma'(0,u) = \log\Gamma(u) - \frac{1}{2}\log(2\pi) \qquad \varsigma'(0,1) = -\frac{1}{2}\log(2\pi)$$

this becomes

$$= 2\varsigma'(-1,u) - 2\varsigma'(-1) - 2u\log\Gamma(u) + (u-1)\log(2\pi) + u^2\psi(u) + \gamma - \frac{1}{2}(u-1)(3u+2)$$

Since from (4.3.84)

$$\sum_{k=0}^{\infty} \frac{1}{k+1} \sum_{j=0}^{k} \binom{k}{j} (-1)^j j^2 \log \frac{(u+j)}{(1+j)} = (u-1)\log(2\pi) - 2\log G(1+u) + u^2\psi(u)$$

$$+\gamma - \frac{1}{2}(u-1)(3u+2)$$

we have therefore shown that

(4.3.126)        $$\log G(u+1) - u\log\Gamma(u) = \varsigma'(-1) - \varsigma'(-1,u)$$

This functional equation was derived by Vardi in 1988 and also by Gosper in 1997 (see [6a]).

Letting $u \to 1+u$ in (4.3.126) we get

$$\log G(1+1+u) - (1+u)\log\Gamma(1+u) = \varsigma'(-1) - \varsigma'(-1,1+u)$$

and we have

$$\log G(1+1+u) = \log G(1+u) + \log\Gamma(1+u)$$

We then see that

(4.3.126i)        $$\log G(1+u) - u\log\Gamma(1+u) = \varsigma'(-1) - \varsigma'(-1,1+u)$$

and upon subtraction we see that

$$\varsigma'(-1,1+u) - \varsigma'(-1,u) = u\log u$$

Differentiating this, by reference to (4.3.126di), we obtain

$$\varsigma'(0,1+u) - \varsigma(0,1+u) - \varsigma'(0,u) + \varsigma(0,u) = 1 + \log u$$

which, using (4.3.110a), concurs with (4.3.121a)

$$\log u = \varsigma'(0,u+1) - \varsigma'(0,u)$$

Letting $u \to u-1$ in (4.3.126) we get with $u-1 > 0$

$$\log G(u) - (u-1)\log\Gamma(u-1) = \varsigma'(-1) - \varsigma'(-1,u-1)$$

and subtraction results in

$$(u-1)\log\frac{\Gamma(u)}{\Gamma(u-1)} = \varsigma'(-1,u-1) - \varsigma'(-1,u)$$

With $u = 2$ we see that $\varsigma'(-1,2) = \varsigma'(-1,1) = \varsigma'(-1)$ and

$$\lim_{u \to 1}(u-1)\log\frac{\Gamma(u)}{\Gamma(u-1)} = \varsigma'(-1,1) - \varsigma'(-1,0)$$

and hence we see that $\varsigma'(-1,1) = \varsigma'(-1,0)$.

It will be noted from (4.3.84) that we need $\log(u+j)$ to exist and hence we require $u+j > 0$; it therefore appeared to me that (4.3.126) might be valid for negative values of $u$ provided $u+j > 0$, and in particular $u+1 > 0$. Hence $u = -1/2$ initially seemed a reasonable candidate.

(4.3.126a) $\qquad \log G\left(\frac{1}{2}\right) + \frac{1}{2}\log\Gamma\left(-\frac{1}{2}\right) = \varsigma'(-1) - \varsigma'\left(-1,-\frac{1}{2}\right)$ (?)

However, we know from [126, p.4] that $\Gamma\left(-\frac{1}{2}\right) = -\frac{1}{2}\sqrt{\pi}$ and hence we end up with the logarithm of a negative number.

As another example we have

(4.3.126b) $\qquad \log G(3/2) - \frac{1}{4}\log\pi = \varsigma'(-1) - \varsigma'(-1,1/2)$

where $G(3/2) = G(1/2)\Gamma(1/2)$ may be evaluated in terms of the Glaisher-Kinkelin constant. Using (4.3.170d)

$$\varsigma'\left(-1,\frac{1}{2}\right) = -\frac{B_2\log 2}{4} - \frac{1}{2}\varsigma'(-1) = -\frac{1}{24}\log 2 - \frac{1}{2}\varsigma'(-1)$$

we easily determine that [126, p.26]

(4.3.126bi) $\qquad \log G(1/2) = -\frac{1}{4}\log\pi + \frac{3}{2}\varsigma'(-1) + \frac{1}{24}\log 2$

We have

$$\log A = \frac{1}{12} - \varsigma'(-1)$$

where $A$ is the Glaisher-Kinkelin constant and from [126, p.25] the numerical value is reported as approximately

$$A = 1.282427130...$$

The identity (4.3.126) also implies that

(4.3.126c) $\qquad \lim_{u \to 1} \varsigma'(-1, u) = \varsigma'(-1, 1) = \varsigma'(-1)$ and $\lim_{u \to 0} \varsigma'(-1, u) = \varsigma'(-1, 0) = \varsigma'(-1)$

This is considered further in (4.3.132).

In this regard, it may be noted Apostol [14b] has shown that for $-1 < s < 0$

(4.3.126d) $\qquad\qquad \lim_{u \to 0+} \varsigma(s, u) = \varsigma(s)$

Differentiating (4.3.126) we obtain

$$\frac{G'(u+1)}{G(u+1)} - u\psi(u) - \log \Gamma(u) = -\frac{d}{du} \varsigma'(-1, u)$$

We have

$$\frac{\partial}{\partial u} \frac{\partial}{\partial s} \varsigma(s, u) = \frac{\partial}{\partial s} \frac{\partial}{\partial u} \varsigma(s, u) = -\frac{\partial}{\partial s} [s\varsigma(s+1, u)]$$

$$= -\varsigma(s+1, u) - s\frac{\partial}{\partial s} \varsigma(s+1, u)$$

and therefore we get

(4.3.126di) $\qquad -\frac{d}{du} \varsigma'(-1, u) = -\varsigma'(0, u) + \varsigma(0, u)$

Using Lerch's identity, we then see that

$$\frac{G'(u+1)}{G(u+1)} - u\psi(u) - \log \Gamma(u) = -\log \Gamma(u) + \frac{1}{2}\log(2\pi) + \varsigma(0, u)$$

and integrating this we obtain

(4.3.126e) $\qquad \log G(x+1) - \int_0^x u\psi(u) = \frac{1}{2}x\log(2\pi) + \frac{1}{2}x(1-x)$

which is equivalent to Alexeiewsky's theorem (4.3.85)

$\qquad\qquad\qquad\qquad\qquad\qquad\qquad\qquad\qquad\qquad$ □

We now revisit (4.3.74)

$$\psi(u) = \sum_{n=0}^{\infty} \frac{1}{n+1} \sum_{k=0}^{n} \binom{n}{k} (-1)^j \log(u+k)$$

where multiplication by $u$ and integration gives us

$$\int_0^t u\psi(u)du = \sum_{n=0}^{\infty} \frac{1}{n+1} \sum_{k=0}^{n} \binom{n}{k} (-1)^j \int_0^t u \log(u+k)du$$

With integration by parts we have

$$\int_0^t u\psi(u)\,du = u\log\Gamma(u)\Big|_0^t - \int_0^t \log\Gamma(u)\,du$$

By L'Hôpital's rule we see that

$$\lim_{u \to 0}\big[u\log\Gamma(u)\big] = \lim_{u \to 0}\frac{\log\Gamma(u)}{1/u} = -\lim_{u \to 0}\frac{\psi(u)}{1/u^2}$$

From (E.14) we have

$$\psi(u) = -\gamma - \frac{1}{u} + \sum_{k=1}^{\infty}\left(\frac{1}{k} - \frac{1}{u+k}\right)$$

and hence we see that $\lim_{u \to 0}\big[u^2\psi(u)\big] = 0$ and this therefore implies that $\lim_{u \to 0}\big[u\log\Gamma(u)\big] = 0$. This may also be seen more directly by considering the form

$$u\log\Gamma(u+1) = u\log u + u\log\Gamma(u)$$

Accordingly we have

$$\int_0^t u\psi(u)\,du = t\log\Gamma(t) - \int_0^t \log\Gamma(u)\,du$$

An easy integration by parts gives us

$$\int_0^t u\,\log(u+k)du = \frac{1}{4}\Big[(2k-u) - 2(k^2-u^2)\log(u+k)\Big]\Big|_0^t$$

and a little algebra gives us

$$= \frac{1}{4}\Big[ (2k-u) - 2[(k+u)^2 - 2u(k+u)]\log(u+k) \Big]\Big|_0^t$$

Hence we have from (4.3.108) and (4.3.109)

$$\int_0^t u\psi(u)du = F(t) - F(0)$$

where

$$F(t) = \frac{1}{4}\Big( -1 - t - 2[2\varsigma'(-1,t) - \varsigma(-1,t)] + 4t[\varsigma'(0,t) - \varsigma(0,t)] \Big)$$

$$F(0) = \frac{1}{4}\Big( -1 - 2[2\varsigma'(-1,0) - \varsigma(-1,0)] \Big)$$

Using (4.3.126c) we get

$$F(0) = \frac{1}{4}\Big( -1 - 2[2\varsigma'(-1) - \varsigma(-1,0)] \Big)$$

$$F(t) - F(0) = \frac{1}{4}\Big( -t + 4[\varsigma'(-1) - \varsigma'(-1,t)] + 2[\varsigma(-1,t) - \varsigma(-1,0)] + 4t[\varsigma'(0,t) - \varsigma(0,t)] \Big)$$

and hence

$$\int_0^t u\psi(u)du = \frac{1}{4}\Big( -t + 4[\varsigma'(-1) - \varsigma'(-1,t)] - \varsigma(-1,t)] + 4t[\varsigma'(0,t) - \varsigma(0,t)] \Big)$$

Then using (4.3.110a), (4.3.116) and (4.3.126) we get another derivation of Alexeiewsky's theorem (4.3.85).

$$\int_0^t u\psi(u)du = \log G(u+1) - \frac{1}{2}\log(2\pi) - \frac{1}{2}u(u-1)$$

□

We now refer back to (4.3.75) and (4.3.84) as shown below

$$\log \Gamma(u) = \sum_{n=0}^{\infty} \frac{1}{n+1}\sum_{k=0}^{n}\binom{n}{k}(-1)^k k \log \frac{(u+k)}{(1+k)} + u\psi(u) + \gamma + 1 - u$$

$$2\log G(1+u) = (u-1)\log(2\pi) - \sum_{k=0}^{\infty}\frac{1}{k+1}\sum_{j=0}^{k}\binom{k}{j}(-1)^j j^2 \log\frac{(u+j)}{(1+j)}$$

$$+ u^2\psi(u) + \gamma - \frac{1}{2}(u-1)(3u+2)$$

Letting $u = 0$ in (4.3.84) we deduce that

$$\sum_{n=0}^{\infty}\frac{1}{n+1}\sum_{k=0}^{n}\binom{n}{k}(-1)^k k^2 \log\left[\frac{k}{1+k}\right] = \gamma + 1 - \log(2\pi)$$

which we also saw in (4.3.77a).

$\square$

A slightly different proof of (4.3.126) is shown below. Hasse's formula gives us

$$\varsigma'(s) = -\frac{1}{s-1}\sum_{n=0}^{\infty}\frac{1}{n+1}\sum_{k=0}^{n}\binom{n}{k}\frac{(-1)^k \log(k+1)}{(k+1)^{s-1}} - \frac{1}{(s-1)^2}\sum_{n=0}^{\infty}\frac{1}{n+1}\sum_{k=0}^{n}\binom{n}{k}\frac{(-1)^k}{(k+1)^{s-1}}$$

$$\varsigma'(-1) = \frac{1}{2}\sum_{n=0}^{\infty}\frac{1}{n+1}\sum_{k=0}^{n}\binom{n}{k}(-1)^k(k+1)^2 \log(k+1) - \frac{1}{4}\sum_{n=0}^{\infty}\frac{1}{n+1}\sum_{k=0}^{n}\binom{n}{k}(-1)^k(k+1)^2$$

We will see in (F.21) of Volume VI that (see also Guillera and Sondow's paper [75aa])

$$\sum_{n=0}^{\infty}\frac{1}{n+1}\sum_{k=0}^{n}\binom{n}{k}(-1)^k(k+1)^2 = B_2 = \frac{1}{6}$$

which can also be seen from

$$\sum_{n=0}^{\infty}\frac{1}{n+1}\sum_{k=0}^{n}\binom{n}{k}(-1)^k(k+1)^2 = -2\varsigma(-1,1)$$

We have

$$\sum_{n=0}^{\infty}\frac{1}{n+1}\sum_{k=0}^{n}\binom{n}{k}(-1)^k(k+1)^2 \log(k+1) =$$

$$\sum_{n=0}^{\infty}\frac{1}{n+1}\sum_{k=0}^{n}\binom{n}{k}(-1)^k k^2 \log(k+1) + 2\sum_{n=0}^{\infty}\frac{1}{n+1}\sum_{k=0}^{n}\binom{n}{k}(-1)^k k \log(k+1)$$

$$+\sum_{n=0}^{\infty}\frac{1}{n+1}\sum_{k=0}^{n}\binom{n}{k}(-1)^k\log(k+1)$$

Using (4.3.115a) this then gives us

$$2\zeta'(-1)=\sum_{n=0}^{\infty}\frac{1}{n+1}\sum_{k=0}^{n}\binom{n}{k}(-1)^k k^2\log(k+1)-\frac{1}{12}+2\left[\gamma+\frac{1}{2}-\frac{1}{2}\log(2\pi)\right]-\gamma$$

and therefore we get

$$\sum_{n=0}^{\infty}\frac{1}{n+1}\sum_{k=0}^{n}\binom{n}{k}(-1)^k k^2\log(k+1)=2\zeta'(-1)-\frac{11}{12}-\gamma+\log(2\pi)$$

Hence using (4.3.84) we obtain

$$2\log G(1+u)=(u-1)\log(2\pi)-\sum_{n=0}^{\infty}\frac{1}{n+1}\sum_{k=0}^{n}\binom{n}{k}(-1)^k k^2\log(u+k)+u^2\psi(u)$$

$$-\frac{1}{2}(u-1)(3u+2)+2\zeta'(-1)-\frac{11}{12}+\log(2\pi)$$

$$=-\sum_{n=0}^{\infty}\frac{1}{n+1}\sum_{k=0}^{n}\binom{n}{k}(-1)^k(u+k)^2\log(u+k)+2u\sum_{n=0}^{\infty}\frac{1}{n+1}\sum_{k=0}^{n}\binom{n}{k}(-1)^k k\log(u+k)$$

$$+2u^2\psi(u)+(u-1)\log(2\pi)-\frac{1}{2}(u-1)(3u+2)+2\zeta'(-1)-\frac{11}{12}+\log(2\pi)$$

where we have used (4.3.74) for $u^2\psi(u)$. We now employ (4.3.119) for

$$2u\sum_{n=0}^{\infty}\frac{1}{n+1}\sum_{k=0}^{n}\binom{n}{k}(-1)^k k\log(u+k)=2u\log\Gamma(u)-2u^2\psi(u)-u+2u^2-u\log(2\pi)$$

and obtain

$$2\log G(1+u)=$$

$$-\sum_{n=0}^{\infty}\frac{1}{n+1}\sum_{k=0}^{n}\binom{n}{k}(-1)^k(u+k)^2\log(u+k)+2u\log\Gamma(u)-2u^2\psi(u)-u+2u^2-u\log(2\pi)$$

$$+2u^2\psi(u)+(u-1)\log(2\pi)-\frac{1}{2}(u-1)(3u+2)+2\zeta'(-1)-\frac{11}{12}+\log(2\pi)$$

$$= -\sum_{n=0}^{\infty} \frac{1}{n+1} \sum_{k=0}^{n} \binom{n}{k} (-1)^k (u+k)^2 \log(u+k) + 2u \log \Gamma(u) - u + 2u^2$$

$$-\frac{1}{2}(u-1)(3u+2) + 2\varsigma'(-1) - \frac{11}{12}$$

Finally, from (4.3.109) we have

$$\sum_{n=0}^{\infty} \frac{1}{n+1} \sum_{k=0}^{n} \binom{n}{k} (-1)^k (u+k)^2 \log(u+k) = 2\varsigma'(-1,u) - \varsigma(-1,u)$$

$$= 2\varsigma'(-1,u) + \frac{1}{2}u^2 - \frac{1}{2}u + \frac{1}{12}$$

and we therefore obtain (4.3.126)

$$\log G(u+1) - u \log \Gamma(u) = \varsigma'(-1) - \varsigma'(-1,u)$$

We recall (4.3.109)

$$\varsigma'(-1,u) = \frac{1}{2}\varsigma(-1,u) + \frac{1}{2} \sum_{n=0}^{\infty} \frac{1}{n+1} \sum_{k=0}^{n} \binom{n}{k} (-1)^k (u+k)^2 \log(k+u)$$

and differentiation results in

$$\frac{d}{du}\varsigma'(-1,u) =$$

$$\frac{1}{4} - \frac{1}{2}u + \frac{1}{2} \sum_{n=0}^{\infty} \frac{1}{n+1} \sum_{k=0}^{n} \binom{n}{k} (-1)^k (u+k) + \sum_{n=0}^{\infty} \frac{1}{n+1} \sum_{k=0}^{n} \binom{n}{k} (-1)^k (u+k) \log(k+u)$$

$$= \frac{1}{4} - \frac{1}{2}u + \frac{1}{2} \sum_{n=0}^{\infty} \frac{1}{n+1} \sum_{k=0}^{n} \binom{n}{k} (-1)^k (u+k) + \sum_{n=0}^{\infty} \frac{1}{n+1} \sum_{k=0}^{n} \binom{n}{k} (-1)^k (u+k) \log(k+u)$$

Using (4.3.119) we have

$$\sum_{n=0}^{\infty} \frac{1}{n+1} \sum_{k=0}^{n} \binom{n}{k} (-1)^k (u+k) \log(u+k) = \log \Gamma(u) - \frac{1}{2} + u - \frac{1}{2}\log(2\pi)$$

and we obtain

$$\frac{d}{du}\varsigma'(-1,u) = \log\Gamma(u) - \frac{1}{2} + u - \frac{1}{2}\log(2\pi)$$

From (4.3.126di) we have

$$\frac{d}{du}\varsigma'(-1,u) = \varsigma'(0,u) - \varsigma(0,u)$$

and hence we obtain Lerch's identity (4.3.116)

$$\varsigma'(0,u) = \log\Gamma(u) - \frac{1}{2}\log(2\pi)$$

A further differentiation results in

$$\frac{d^2}{du^2}\varsigma'(-1,u) =$$

$$-\frac{1}{2} + \frac{1}{2}\sum_{n=0}^{\infty}\frac{1}{n+1}\sum_{k=0}^{n}\binom{n}{k}(-1)^k + \sum_{n=0}^{\infty}\frac{1}{n+1}\sum_{k=0}^{n}\binom{n}{k}(-1)^k + \sum_{n=0}^{\infty}\frac{1}{n+1}\sum_{k=0}^{n}\binom{n}{k}(-1)^k\log(k+u)$$

$$= 1 + \sum_{n=0}^{\infty}\frac{1}{n+1}\sum_{k=0}^{n}\binom{n}{k}(-1)^k\log(k+u)$$

and hence we get

(4.3.126f) $\qquad \dfrac{d^2}{du^2}\varsigma'(-1,u) = 1 + \psi(u)$

Differentiating (4.3.84) we get

$$2\frac{d}{du}\log G(1+u) =$$

$$-\sum_{k=0}^{\infty}\frac{1}{k+1}\sum_{j=0}^{k}\binom{k}{j}(-1)^j\frac{j^2}{u+j} + \frac{d}{du}\left[(u-1)\log(2\pi) + u^2\psi(u) + \gamma - \frac{1}{2}(u-1)(3u+2)\right]$$

Since $\dfrac{j^2}{u+j} = j - u + \dfrac{u^2}{u+j}$ we have

$$\sum_{k=0}^{\infty} \frac{1}{k+1} \sum_{j=0}^{k} \binom{k}{j} (-1)^j \frac{j^2}{u+j} = \sum_{k=0}^{\infty} \frac{1}{k+1} \sum_{j=0}^{k} \binom{k}{j} (-1)^j j - u \sum_{k=0}^{\infty} \frac{1}{k+1} \sum_{j=0}^{k} \binom{k}{j} (-1)^j$$

$$+ u^2 \sum_{k=0}^{\infty} \frac{1}{k+1} \sum_{j=0}^{k} \binom{k}{j} (-1)^j \frac{1}{u+j}$$

$$= -\frac{1}{2} - u + u^2 \psi'(u)$$

We then obtain

$$2 \frac{d}{du} \log G(1+u) =$$

$$\frac{1}{2} + u - u^2 \psi'(u) + \frac{d}{du} \left[ (u-1) \log(2\pi) + u^2 \psi(u) + \gamma - \frac{1}{2}(u-1)(3u+2) \right]$$

and upon integrating this simpler expression we obtain

$$2 \log G(1+u) =$$

$$-\frac{1}{2} u - \frac{1}{2} u^2 + \int_0^u x^2 \psi'(x) dx + (u-1) \log(2\pi) + u^2 \psi(u) + \gamma - \frac{1}{2}(u-1)(3u+2) + [\log(2\pi) - 1 - \gamma]$$

$$= -u - 2u^2 + \int_0^u x^2 \psi'(x) dx + u \log(2\pi) + u^2 \psi(u)$$

By integration by parts we have

$$\int_0^u x^2 \psi'(x) dx = u^2 \psi(u) - 2 \int_0^u x \psi(x) dx$$

$$= u^2 \psi(u) - 2u \log \Gamma(u) + 2 \int_0^u \log \Gamma(x) dx$$

and we have now obtained a further derivation of Alexeiewsky's theorem.

We recall from (4.3.109) that

$$(4.3.127) \qquad \varsigma'(-1,u) = \frac{1}{2}\varsigma(-1,u) + \frac{1}{2}\sum_{n=0}^{\infty}\frac{1}{n+1}\sum_{k=0}^{n}\binom{n}{k}(-1)^k(u+k)^2\log(u+k)$$

and, combining this with (4.3.126), we get

(4.3.127a)

$$\log G(u+1) - u\log\Gamma(u) = \varsigma'(-1) - \frac{1}{2}\varsigma(-1,u) - \frac{1}{2}\sum_{n=0}^{\infty}\frac{1}{n+1}\sum_{k=0}^{n}\binom{n}{k}(-1)^k(u+k)^2\log(u+k)$$

We have

$$\varsigma(-1,u) = -\frac{1}{2}B_2(u) = -\frac{1}{2}u^2 + \frac{1}{2}u - \frac{1}{12}$$

and therefore

$$\varsigma(-1,1) = -\frac{1}{12}$$

We accordingly see from (3.11b) in Volume I that

$$\varsigma(-1,0) = -\frac{1}{12} = \varsigma(-1)$$

With $u = 1$ in (4.3.127) we get

$$\varsigma'(-1,1) = \frac{1}{2}\varsigma(-1,1) + \frac{1}{2}\sum_{n=0}^{\infty}\frac{1}{n+1}\sum_{k=0}^{n}\binom{n}{k}(-1)^k(1+k)^2\log(1+k)$$

$$= -\frac{1}{24} + \frac{1}{2}\sum_{n=0}^{\infty}\frac{1}{n+1}\sum_{k=0}^{n}\binom{n}{k}(-1)^k(1+k)^2\log(1+k)$$

From (4.3.126) we note that

$$\log G(2) - \log\Gamma(1) = \varsigma'(-1) - \varsigma'(-1,1)$$

and hence, since $G(2) = 1$, we have as expected

$$\varsigma'(-1,1) = \varsigma'(-1)$$

Accordingly we see that

(4.3.128) $$\varsigma'(-1) = -\frac{1}{24} + \frac{1}{2}\sum_{n=0}^{\infty}\frac{1}{n+1}\sum_{k=0}^{n}\binom{n}{k}(-1)^k(1+k)^2\log(1+k)$$

Guillera and Sondow [75aa] have recently derived an equivalent expression

(4.3.128a) $$\log A = \frac{1}{12} - \varsigma'(-1) = \frac{1}{8} - \frac{1}{2}\sum_{n=0}^{\infty}\frac{1}{n+1}\sum_{k=0}^{n}\binom{n}{k}(-1)^k(1+k)^2\log(1+k)$$

where $A$ is the Glaisher-Kinkelin constant.

With $u = 2$ in (4.3.126) we see that since $G(3) = 1$

$$\log G(3) - 2\log\Gamma(2) = \varsigma'(-1) - \varsigma'(-1,2) = 0$$

and hence we obtain

(4.3.128b) $$\varsigma'(-1) = \varsigma'(-1,2)$$

Hence we get from (4.3.127)

(4.3.128c) $$\sum_{n=0}^{\infty}\frac{1}{n+1}\sum_{k=0}^{n}\binom{n}{k}(-1)^k(2+k)^2\log(2+k) = 2\varsigma'(-1,2) - \varsigma(-1,2)$$

$$= 2\varsigma'(-1) - \varsigma(-1,2) = 2\varsigma'(-1) - \frac{13}{12}$$

$\square$

Differentiating (4.3.127a) gives us

$$\frac{G'(u+1)}{G(u+1)} - u\psi(u) - \log\Gamma(u) = \frac{1}{2}u - \frac{1}{4} - \sum_{n=0}^{\infty}\frac{1}{n+1}\sum_{k=0}^{n}\binom{n}{k}(-1)^k(u+k)\log(u+k)$$

$$-\frac{1}{2}\sum_{n=0}^{\infty}\frac{1}{n+1}\sum_{k=0}^{n}\binom{n}{k}(-1)^k(u+k)$$

and reference to (4.3.108) gives us

$$= \frac{1}{2}u - \frac{1}{4} - [\varsigma'(0,u) - \varsigma(0,u)] - \frac{1}{2}u + \frac{1}{4}$$

$$= -[\varsigma'(0,u) - \varsigma(0,u)]$$

Then using Lerch's identity (4.3.116) we have

$$= -\log \Gamma(u) - \frac{1}{2}\log(2\pi) + \frac{1}{2} - u$$

Hence we see that

$$\frac{G'(u+1)}{G(u+1)} - u\psi(u) = -\frac{1}{2}\log(2\pi) + \frac{1}{2} - u$$

and integrating this gives us yet another derivation of Alexeiewsky's theorem (4.3.85).

ASYMPTOTIC FORMULA FOR $\log G(u+1)$

As mentioned in (4.3.127a) we have

$$\log G(u+1) - u\log \Gamma(u) = \varsigma'(-1) - \frac{1}{2}\varsigma(-1,u) - \frac{1}{2}\sum_{n=0}^{\infty}\frac{1}{n+1}\sum_{k=0}^{n}\binom{n}{k}(-1)^k(u+k)^2\log(u+k)$$

We have

$$\sum_{k=0}^{n}\binom{n}{k}(-1)^k(u+k)^2\log(u+k) = \log\left[\prod_{k=0}^{n}(u+j)^{(-1)^k\binom{n}{k}(u+k)^2}\right]$$

and as $u \to \infty$ we see that

$$\sum_{k=0}^{n}\binom{n}{k}(-1)^k(u+k)^2\log(u+k) \approx \log\left[u^{\sum_{i=0}^{n}\binom{n}{k}(-1)^i(u+k)^2}\right]$$

Using (3.11bi) in Volume I we have

$$\sum_{k=0}^{n}\binom{n}{k}(-1)^k(1+k)^2 = \delta_{n,0} - 3n\delta_{n,1} + n(n-1)\delta_{2,1}$$

and therefore we see that

$$\sum_{k=0}^{n}\binom{n}{k}(-1)^k k^2 = \delta_{n,0} - 3n\delta_{n,1} + n(n-1)\delta_{2,1} - 2\sum_{k=0}^{n}\binom{n}{k}(-1)^k k - \sum_{k=0}^{n}\binom{n}{k}(-1)^k$$

$$= -n\delta_{n,1} + n(n-1)\delta_{2,1}$$

$$\sum_{k=0}^{n}\binom{n}{k}(-1)^k(u+k)^2 = u^2\sum_{k=0}^{n}\binom{n}{k}(-1)^k + 2u\sum_{k=0}^{n}\binom{n}{k}(-1)^k k + \sum_{k=0}^{n}\binom{n}{k}(-1)^k k^2$$

$$= u^2\delta_{n,0} - 2un\delta_{n,1} - n\delta_{n,1} + n(n-1)\delta_{2,1}$$

$$= u^2\delta_{n,0} - (2un+n)\delta_{n,1} + n(n-1)\delta_{2,1}$$

and therefore we have as $u \to \infty$

$$\sum_{n=0}^{\infty}\frac{1}{n+1}\sum_{k=0}^{n}\binom{n}{k}(-1)^k(u+k)\log(u+k) \to \sum_{n=0}^{\infty}\frac{1}{n+1}\Big[u^2\delta_{n,0} - (2un+n)\delta_{n,1} + n(n-1)\delta_{2,1}\Big]\log u$$

$$= \Big[u^2 - u + \frac{1}{6}\Big]\log u$$

Therefore we obtain

$$\log G(u+1) - u\log\Gamma(u) \approx \varsigma'(-1) - \frac{1}{2}\varsigma(-1,u) - \frac{1}{2}\Big[u^2 - u + \frac{1}{6}\Big]\log u$$

and hence employing Stirling's approximation for $\log\Gamma(u)$ we get

$$\log G(u+1) \approx \varsigma'(-1) + \frac{1}{4}\Big[u^2 - u + \frac{1}{6}\Big] - \frac{1}{2}\Big[u^2 - u + \frac{1}{6}\Big]\log u + u\Big[u - \frac{1}{2}\Big]\log u - u^2 + \frac{1}{2}u\log(2\pi)$$

or

(4.3.128d)    $$\log G(u+1) \approx \frac{1}{2}u\log(2\pi) + \varsigma'(-1) - \frac{3}{4}u^2 + \frac{1}{2}\Big[u^2 - \frac{1}{6}\Big]\log u - \frac{1}{4}u + \frac{1}{24}$$

It should however be noted that this is similar to, but does not quite agree with, the result reported by Choi and Srivastava [45ac].

(4.3.128e)    $$\log G(u+1) \approx \frac{1}{2}u\log(2\pi) + \varsigma'(-1) - \frac{3}{4}u^2 + \frac{1}{2}\Big[u^2 - \frac{1}{6}\Big]\log u$$

and hence the analysis here needs to be refined (and/or corrected). Note that Adamchik's result [6a] is also slightly different.

This time round we employ the method shown in (4.3.119iv); we have for $k/u < 1$

$$\log(u+k) - \log u = \log\left(1 + \frac{k}{u}\right) = -\sum_{j=1}^{\infty} (-1)^j \frac{k^j}{ju^j}$$

and we then obtain (it should be noted that there are major issues of convergence associated with the following presentation; rigour may be provided through an appropriate limiting procedure)

$$\sum_{n=0}^{\infty} \frac{1}{n+1} \sum_{k=0}^{n} \binom{n}{k} (-1)^k (u+k)^2 \log(u+k) = \sum_{n=0}^{\infty} \frac{1}{n+1} \sum_{k=0}^{n} \binom{n}{k} (-1)^k (u+k)^2 \left[\log u - \sum_{j=1}^{\infty} (-1)^j \frac{k^j}{ju^j}\right]$$

We have from (A.23)

$$\sum_{n=0}^{\infty} \frac{1}{n+1} \sum_{k=0}^{n} \binom{n}{k} (-1)^k (u+k)^2 \left[\log u - \sum_{j=1}^{\infty} (-1)^j \frac{k^j}{ju^j}\right] =$$

$$B_2(u) \log u - \sum_{n=0}^{\infty} \frac{1}{n+1} \sum_{k=0}^{n} \binom{n}{k} (-1)^k (u+k)^2 \sum_{j=1}^{\infty} (-1)^j \frac{k^j}{ju^j}$$

Rearrangement of the series results in

$$\sum_{n=0}^{\infty} \frac{1}{n+1} \sum_{k=0}^{n} \binom{n}{k} (-1)^k (u+k)^2 \sum_{j=1}^{\infty} (-1)^j \frac{k^j}{ju^j} = \sum_{j=1}^{\infty} (-1)^j \frac{1}{ju^j} \sum_{n=0}^{\infty} \frac{1}{n+1} \sum_{k=0}^{n} \binom{n}{k} (-1)^k (u+k)^2 k^j$$

and we see that

$$\sum_{j=1}^{\infty} (-1)^j \frac{1}{ju^j} \sum_{n=0}^{\infty} \frac{1}{n+1} \sum_{k=0}^{n} \binom{n}{k} (-1)^k (u+k)^2 k^j = u^2 \sum_{j=1}^{\infty} (-1)^j \frac{1}{ju^j} \sum_{n=0}^{\infty} \frac{1}{n+1} \sum_{k=0}^{n} \binom{n}{k} (-1)^k k^j$$

$$+ 2u \sum_{j=1}^{\infty} (-1)^j \frac{1}{ju^j} \sum_{n=0}^{\infty} \frac{1}{n+1} \sum_{k=0}^{n} \binom{n}{k} (-1)^k k^{j+1} + \sum_{j=1}^{\infty} (-1)^j \frac{1}{ju^j} \sum_{n=0}^{\infty} \frac{1}{n+1} \sum_{k=0}^{n} \binom{n}{k} (-1)^k k^{j+2}$$

Another reference to (A.23aa) then results in

$$\sum_{j=1}^{\infty} (-1)^j \frac{1}{ju^j} \sum_{n=0}^{\infty} \frac{1}{n+1} \sum_{k=0}^{n} \binom{n}{k} (-1)^k (u+k)^2 k^j = u^2 \sum_{j=1}^{\infty} \frac{(-1)^j B_j}{ju^j} + 2u \sum_{j=1}^{\infty} \frac{(-1)^j B_{j+1}}{ju^j} + \sum_{j=1}^{\infty} \frac{(-1)^j B_{j+2}}{ju^j}$$

We saw previously that

$$\sum_{j=1}^{\infty} \frac{(-1)^j B_j}{ju^j} = \sum_{m=1}^{\infty} \frac{B_{2m}}{2mu^{2m}} + \frac{1}{2u} = \sum_{m=2}^{\infty} \frac{B_{2m}}{2mu^{2m}} + \frac{1}{2u} + \frac{1}{12u^2}$$

and

$$\sum_{j=1}^{\infty} \frac{(-1)^j B_{j+1}}{ju^j} = -\sum_{m=1}^{\infty} \frac{B_{2m}}{(2m-1)u^{2m-1}} = -\frac{1}{6u} - \sum_{m=2}^{\infty} \frac{B_{2m}}{(2m-1)u^{2m-1}}$$

We also see that

$$\sum_{j=1}^{\infty} \frac{(-1)^j B_{j+2}}{ju^j} = \sum_{m=1}^{\infty} \frac{B_{2m+2}}{2mu^{2m}} - \sum_{m=0}^{\infty} \frac{B_{2m+3}}{(2m+1)u^{2m+1}} = \sum_{m=1}^{\infty} \frac{B_{2m+2}}{2mu^{2m}}$$

$$= \sum_{m=2}^{\infty} \frac{B_{2m}}{(2m-2)u^{2m-2}}$$

We then get

$$u^2 \sum_{j=1}^{\infty} \frac{(-1)^j B_j}{ju^j} + 2u \sum_{j=1}^{\infty} \frac{(-1)^j B_{j+1}}{ju^j} + \sum_{j=1}^{\infty} \frac{(-1)^j B_{j+2}}{ju^j}$$

$$= \sum_{m=2}^{\infty} \frac{B_{2m}}{2mu^{2m-2}} + \frac{1}{2}u + \frac{1}{12} - \frac{1}{3} - 2\sum_{m=2}^{\infty} \frac{B_{2m}}{(2m-1)u^{2m-2}} + \sum_{m=2}^{\infty} \frac{B_{2m}}{(2m-2)u^{2m-2}}$$

Then using (4.3.127a)

$$\log G(u+1) - u\log \Gamma(u) = \varsigma'(-1) - \frac{1}{2}\varsigma(-1,u) - \frac{1}{2}\sum_{n=0}^{\infty} \frac{1}{n+1} \sum_{k=0}^{n} \binom{n}{k}(-1)^k (u+k)^2 \log(u+k)$$

we obtain using (4.3.119v)

$$\log G(u+1) = \frac{1}{2}u\log(2\pi) + u\left[u - \frac{1}{2}\right]\log u - u^2 + \sum_{m=1}^{\infty} \frac{B_{2m}}{2m(2m-1)u^{2m-2}}$$

$$+ \varsigma'(-1) - \frac{1}{2}\varsigma(-1,u) - \frac{1}{2}B_2(u)\log u$$

$$+ \frac{1}{2}\sum_{m=2}^{\infty} \frac{B_{2m}}{2mu^{2m-2}} + \frac{1}{4}u - \frac{1}{8} - \sum_{m=2}^{\infty} \frac{B_{2m}}{(2m-1)u^{2m-2}} + \frac{1}{2}\sum_{m=2}^{\infty} \frac{B_{2m}}{(2m-2)u^{2m-2}}$$

which may be written as

$$\log G(u+1) = \frac{1}{2}u\log(2\pi) + \left[u^2 - \frac{1}{2}u\right]\log u - \left[\frac{1}{2}u^2 - \frac{1}{2}u + \frac{1}{12}\right]\log u - u^2 + \frac{1}{4}u$$

$$-\frac{1}{8} + \varsigma'(-1) + \frac{1}{4}\left[u^2 - u + \frac{1}{6}\right]$$

$$-\sum_{m=2}^{\infty}\frac{B_{2m}}{2mu^{2m-2}} + 2\sum_{m=2}^{\infty}\frac{B_{2m}}{(2m-1)u^{2m-2}} - \sum_{m=1}^{\infty}\frac{B_{2m}}{(2m-2)u^{2m-2}} - \sum_{m=2}^{\infty}\frac{B_{2m}}{2m(2m-1)u^{2m-2}}$$

and simplified to

$$\log G(u+1) = \frac{1}{2}u\log(2\pi) + \frac{1}{2}u^2\log u - \frac{1}{12}\log u - \frac{3}{4}u^2 + \frac{1}{4} + \varsigma'(-1)$$

$$-\sum_{m=2}^{\infty}\frac{B_{2m}}{2mu^{2m-2}} + 2\sum_{m=2}^{\infty}\frac{B_{2m}}{(2m-1)u^{2m-2}} - \sum_{m=1}^{\infty}\frac{B_{2m}}{(2m-2)u^{2m-2}} - \sum_{m=2}^{\infty}\frac{B_{2m}}{2m(2m-1)u^{2m-2}}$$

We see that

$$-\sum_{m=2}^{\infty}\frac{B_{2m}}{2mu^{2m-2}} + 2\sum_{m=2}^{\infty}\frac{B_{2m}}{(2m-1)u^{2m-2}} - \sum_{m=1}^{\infty}\frac{B_{2m}}{(2m-2)u^{2m-2}} - \sum_{m=2}^{\infty}\frac{B_{2m}}{2m(2m-1)u^{2m-2}}$$

$$= \sum_{m=2}^{\infty}\frac{B_{2m}}{u^{2m-2}}\left[-\frac{1}{2m} + \frac{2}{2m-1} - \frac{1}{2m-2} - \frac{1}{2m(2m-1)}\right]$$

$$= -\sum_{m=2}^{\infty}\frac{B_{2m}}{(2m-1)(2m-2)u^{2m-2}}$$

Hence, except for the additional factor of $-1/12$, this agrees with the asymptotic formula discovered by Adamchik [6a] using Hermite's integral for the Hurwitz zeta function.

$$\log G(u+1) = \frac{1}{2}u^2\left(\log u - \frac{3}{2}\right) - \frac{1}{12}\log u - u\varsigma'(0) + \varsigma'(-1) - \sum_{k=1}^{n}\frac{B_{2k+2}}{4k(k+1)u^{2k}} + O\left(\frac{1}{u^{2n+2}}\right)$$

Corrections welcome!

We know from (4.3.119) that

$$\log \Gamma(u) = \sum_{n=0}^{\infty} \frac{1}{n+1} \sum_{k=0}^{n} \binom{n}{k} (-1)^k (u+k) \log(u+k) + \frac{1}{2} - u + \frac{1}{2} \log(2\pi)$$

and hence we have

$$\int_{1}^{x} \log \Gamma(u) du =$$

$$\sum_{n=0}^{\infty} \frac{1}{n+1} \sum_{k=0}^{n} \binom{n}{k} (-1)^k \left[ \int_{1}^{x} (u+k) \log(u+k) du + \frac{1}{2}(x-1) - \frac{1}{2}(x^2-1) + \frac{1}{2}(x-1) \log(2\pi) \right]$$

The following integral is elementary

$$\int_{1}^{x} (u+k) \log(u+k) du = \frac{1}{2}(x+k)^2 \log(x+k) - \frac{1}{2}(1+k)^2 \log(1+k) - \frac{1}{4} x(x+2k) + \frac{1}{4}(1+2k)$$

and, referring to (4.3.109), we see that

$$\int_{1}^{x} \log \Gamma(u) du = \varsigma'(-1,x) - \frac{1}{2}\varsigma(-1,x) - \varsigma'(-1,1) + \frac{1}{2}\varsigma(-1,1) + R(x)$$

where

$$R(x) = \sum_{n=0}^{\infty} \frac{1}{n+1} \sum_{k=0}^{n} \binom{n}{k} (-1)^k \left[ \frac{1}{2}(x-1) - \frac{1}{2}(x^2-1) + \frac{1}{2}(x-1) \log(2\pi) - \frac{1}{4}x(x+2k) + \frac{1}{4}(1+2k) \right]$$

$$= \frac{1}{2}(x-1) - \frac{1}{2}(x^2-1) + \frac{1}{2}(x-1) \log(2\pi) - \frac{1}{4}x(x+2k) + \frac{1}{4}(1+2k)$$

With a little more algebra we find the well-known formula

$$(4.3.129) \qquad \int_{1}^{x} \log \Gamma(u) du = \frac{1}{2}x(1-x) + \frac{1}{2}(x-1)\log(2\pi) + \varsigma'(-1,x) - \varsigma'(-1)$$

and we also have

(4.3.129a) $\qquad \displaystyle\int_0^x \log \Gamma(u)\,du = \frac{1}{2}x(1-x) + \frac{1}{2}x\log(2\pi) + \varsigma'(-1,x) - \varsigma'(-1)$

The above integral was found by Gosper in 1997 (see Adamchik's paper [4]). An entirely different proof, using sine and cosine integrals, appears in (6.117g) in Volume V. The connection between Alexeiewsky's theorem (4.3.85), (4.3.126) and (4.3.129) is obvious. We immediately see that

$$\int_0^1 \log \Gamma(u)\,du = \frac{1}{2}\log(2\pi)$$

and letting $x = 0$ in (4.3.129a) we have

(4.3.129b) $\qquad \varsigma'(-1,0) = \varsigma'(-1)$ .

The integral (4.3.129a) may also be indirectly verified by differentiation using

$$\frac{d}{dx}\varsigma'(-1,x) = \varsigma'(0,x) - \varsigma(0,x)$$

and then follow up with Lerch's identity.

From (4.3.129) we have

$$\int_1^{n+1} \log \Gamma(u)\,du = -\frac{1}{2}n(n+1) + \frac{1}{2}n\log(2\pi) + \varsigma'(-1,n+1) - \varsigma'(-1)$$

and since

$$\int_u^{n+u} \log \Gamma(u)\,du = \sum_{k=0}^{n-1}(u+k)\log(u+k) - nu - \frac{1}{2}n(n-1) + \frac{1}{2}n\log(2\pi)$$

we see that

$$\int_1^{n+1} \log \Gamma(u)\,du = \sum_{k=0}^{n-1}(1+k)\log(1+k) - n - \frac{1}{2}n(n+1) + \frac{1}{2}n\log(2\pi)$$

and hence we obtain the known result

(4.3.129c) $\qquad \displaystyle\sum_{k=1}^n k\log k = \varsigma'(-1,n+1) - \varsigma'(-1)$

When $n = 1$ we see that

(4.3.129d)          $\varsigma'(-1, 2) = \varsigma'(-1)$.

This is only to be expected because $\varsigma(s, 2) = \varsigma(s) - 1$.

Integrating (4.3.87f) we see that

$$\int_0^x \log \Gamma(u) \, du = -x \log x + x - \frac{1}{2} \gamma x^2 + \sum_{k=2}^{\infty} (-1)^k \frac{\varsigma(k)}{k(k+1)} x^{k+1}$$

and equating this with (4.3.129a) gives us

(4.3.129e)

$$\frac{1}{2} x(1-x) + \frac{1}{2} x \log(2\pi) + \varsigma'(-1, x) - \varsigma'(-1) = -x \log x + x - \frac{1}{2} \gamma x^2 + \sum_{k=2}^{\infty} (-1)^k \frac{\varsigma(k)}{k(k+1)} x^{k+1}$$

With $x = 1$ we have [126, p.222]

$$\frac{1}{2} \log(2\pi) = 1 - \frac{1}{2} \gamma + \sum_{k=2}^{\infty} (-1)^k \frac{\varsigma(k)}{k(k+1)}$$

Differentiating (4.3.129e) results in

$$\frac{1}{2} - x + \frac{1}{2} \log(2\pi) + \varsigma'(0, x) - \varsigma(0, x) = -\log x - \gamma x + \sum_{k=2}^{\infty} (-1)^k \frac{\varsigma(k)}{k} x^k$$

and, using (4.3.87f), we end up with Lerch's identity.

From (4.3.87e) and (4.3.87f) we see that

$$\log G(1+x) - x \log \Gamma(x) = x \log x + \frac{1}{2} x \log(2\pi) + \frac{1}{2} \gamma x^2 - \frac{1}{2} x(x+1) - \sum_{k=2}^{\infty} (-1)^k \frac{\varsigma(k)}{k(k+1)} x^{k+1}$$

and looking ahead to (4.3.150a)

$$\int_0^x \log G(u+1) \, du - \int_0^x u \log \Gamma(u) \, du = x \varsigma'(-1) + \frac{1}{12} B_3(x) - \frac{1}{2} \varsigma'(-2, x) + \frac{1}{2} \varsigma'(-2)$$

we see that

$$x\varsigma'(-1) + \frac{1}{12}B_3(x) - \frac{1}{2}\varsigma'(-2,x) + \frac{1}{2}\varsigma'(-2) =$$

$$\frac{1}{2}x^2\log x - \frac{1}{4}x^2 + \frac{1}{4}x^2\log(2\pi) + \frac{1}{6}\gamma x^3 - \frac{1}{6}x^3 - \frac{1}{4}x^2 - \sum_{k=2}^{\infty}(-1)^k\frac{\varsigma(k)}{k(k+1)(k+2)}x^{k+2}$$

## THE VANISHING INTEGRAL

From (4.3.107a) we see that

$$(s-1)\varsigma'(s,u) + \varsigma(s,u) = -\sum_{m=0}^{\infty}\frac{1}{m+1}\sum_{k=0}^{m}\binom{m}{k}(-1)^k\frac{\log(u+k)}{(u+k)^{s-1}}$$

and hence with $s = 1 - n$ we get $\forall n \in \mathbf{Z}\setminus\{0\}$

$$n\varsigma'(1-n,u) - \varsigma(1-n,u) = \sum_{m=0}^{\infty}\frac{1}{m+1}\sum_{k=0}^{m}\binom{m}{k}(-1)^k(u+k)^n\log(u+k)$$

We now wish to integrate the above equation.

In (4.3.110) we mentioned that

$$\varsigma(1-n,u) = -\frac{B_n(u)}{n} \quad, n \geq 1$$

and from Appendix A of Volume VI we see that

$$\int_1^x B_n(u)du = \frac{B_{n+1}(x) - B_{n+1}(1)}{n+1} = \frac{B_{n+1}(x) - B_{n+1}}{n+1}$$

since (A.4) shows that $B_{n+1}(1) = B_{n+1}(0) = B_{n+1}$ .

Integration by parts gives us

$$\int(u+k)^n\log(u+k)du = \frac{(u+k)^{n+1}}{n+1}\log(u+k) - \frac{1}{n+1}\int(u+k)^ndu$$

$$= \frac{(u+k)^{n+1}}{n+1}\log(u+k) - \frac{(u+k)^{n+1}}{(n+1)^2}$$

and hence we have

$$\int_1^x (u+k)^n \log(u+k)\,du = \frac{(x+k)^{n+1}}{n+1}\log(x+k) - \frac{(x+k)^{n+1}}{(n+1)^2} - \frac{(1+k)^{n+1}}{n+1}\log(1+k) + \frac{(1+k)^{n+1}}{(n+1)^2}$$

Therefore we obtain

$$\sum_{m=0}^{\infty}\frac{1}{m+1}\sum_{k=0}^{m}\binom{m}{k}(-1)^k\int_1^x(u+k)^n\log(u+k)\,du = \frac{1}{n+1}\sum_{m=0}^{\infty}\frac{1}{m+1}\sum_{k=0}^{m}\binom{m}{k}(-1)^k(x+k)^{n+1}\log(x+k)$$

$$-\frac{1}{(n+1)^2}\sum_{m=0}^{\infty}\frac{1}{m+1}\sum_{k=0}^{m}\binom{m}{k}(-1)^k(x+k)^{n+1} - \frac{1}{n+1}\sum_{m=0}^{\infty}\frac{1}{m+1}\sum_{k=0}^{m}\binom{m}{k}(-1)^k(1+k)^{n+1}\log(1+k)$$

$$+\frac{1}{(n+1)^2}\sum_{m=0}^{\infty}\frac{1}{m+1}\sum_{k=0}^{m}\binom{m}{k}(-1)^k(1+k)^{n+1}$$

$$= -\frac{1}{n+1}[(-n-1)\varsigma'(-n,x)+\varsigma(-n,x)] + \frac{1}{n+1}[(-n-1)\varsigma'(-n,1)+\varsigma(-n,1)]$$

$$-\frac{1}{(n+1)^2}(-n-1)\varsigma(-n,x) + \frac{1}{(n+1)^2}(-n-1)\varsigma(-n,1)$$

$$= \varsigma'(-n,x) - \frac{1}{n+1}\varsigma(-n,x) - \varsigma'(-n,1) + \frac{1}{n+1}\varsigma(-n,1)$$

$$+\frac{1}{n+1}\varsigma(-n,x) - \frac{1}{n+1}\varsigma(-n,1)$$

$$= \varsigma'(-n,u) - \varsigma'(-n)$$

Therefore we have for $n \geq 1$

$$(4.3.130) \qquad n\int_1^x \varsigma'(1-n,u)\,du = \frac{B_{n+1} - B_{n+1}(x)}{n(n+1)} + \varsigma'(-n,x) - \varsigma'(-n)$$

A similar result was obtained by Adamchik [4] in 1998 by reference to polygamma functions of negative order. Adamchik reported the following integral over the range $[0,x]$

$$(4.3.131) \qquad n\int_0^x \varsigma'(1-n,u)\,du = \frac{B_{n+1} - B_{n+1}(x)}{n(n+1)} + \varsigma'(-n,x) - \varsigma'(-n)$$

(an alternative and more direct proof is given in (4.3.144)).

Letting $x \to 0$ in (4.3.131) then implies that

(4.3.132) $$\lim_{x \to 0}[\varsigma'(-n,x) - \varsigma'(-n)] = 0$$

and we recall that in (4.3.126b) we previously showed that

$$\lim_{x \to 0}\varsigma'(-1,x) = \varsigma'(-1)$$

By letting $s \to -p$ in (4.3.107a) we obtain

(4.3.132a) $$-(p+1)\varsigma'(-p,x) + \varsigma(-p,x) = -\sum_{n=0}^{\infty}\frac{1}{n+1}\sum_{k=0}^{n}\binom{n}{k}(-1)^k(x+k)^{p+1}\log(x+k)$$

and letting $x \to 0$ we get

(4.3.132b) $$-(p+1)\varsigma'(-p,0) + \varsigma(-p,0) = -\sum_{n=0}^{\infty}\frac{1}{n+1}\sum_{k=0}^{n}\binom{n}{k}(-1)^k k^{p+1}\log k$$

Using (4.3.110) this becomes

$$-(p+1)\varsigma'(-p,0) - \frac{B_{p+1}}{p+1} = -\sum_{n=0}^{\infty}\frac{1}{n+1}\sum_{k=0}^{n}\binom{n}{k}(-1)^k k^{p+1}\log k$$

We now use the Hasse formula (3.12) for the Riemann zeta function

$$\varsigma(s) = \frac{1}{s-1}\sum_{n=0}^{\infty}\frac{1}{n+1}\sum_{k=0}^{n}\binom{n}{k}\frac{(-1)^k}{(k+1)^{s-1}}$$

and, in a similar manner, using (3.11c) we ascertain that

(4.3.132c) $$-(p+1)\varsigma'(-p) - \frac{B_{p+1}}{p+1} = -\sum_{n=0}^{\infty}\frac{1}{n+1}\sum_{k=0}^{n}\binom{n}{k}(-1)^k(1+k)^{p+1}\log(1+k)$$

Equation (4.3.132c) may also be obtained directly from (4.3.132a) by letting $x = 1$. This coincides with (4.3.127) when $p = 1$. We have

$$2\varsigma'(-1) + \frac{B_2}{2} = \sum_{n=0}^{\infty}\frac{1}{n+1}\sum_{k=0}^{n}\binom{n}{k}(-1)^k(1+k)^2\log(1+k)$$

$$= \sum_{n=0}^{\infty} \frac{1}{n+1} \sum_{k=0}^{n} \binom{n}{k} (-1)^k (1+2k+k^2) \log(1+k)$$

$$= -\gamma + 2\left[\gamma + \frac{1}{2} - \frac{1}{2}\log(2\pi)\right] + \sum_{n=0}^{\infty} \frac{1}{n+1} \sum_{k=0}^{n} \binom{n}{k} (-1)^k k^2 \log(1+k)$$

and therefore we get (4.3.128)

(4.3.132d) $$\sum_{n=0}^{\infty} \frac{1}{n+1} \sum_{k=0}^{n} \binom{n}{k} (-1)^k k^2 \log(1+k) = 2\varsigma'(-1) + \log(2\pi) - \gamma - \frac{11}{12}$$

We now refer to (4.3.84)

$$2\log G(1+u) = (u-1)\log(2\pi) - \sum_{n=0}^{\infty} \frac{1}{n+1} \sum_{k=0}^{n} \binom{n}{k} (-1)^k k^2 \log \frac{(u+k)}{(1+k)}$$

$$+ u^2 \psi(u) + \gamma - \frac{1}{2}(u-1)(3u+2)$$

and then note that we may write this as

$$2\log G(1+u) = u\log(2\pi) - \sum_{n=0}^{\infty} \frac{1}{n+1} \sum_{k=0}^{n} \binom{n}{k} (-1)^k k^2 \log(u+k)$$

$$+ 2\varsigma'(-1) - \frac{11}{12} + u^2 \psi(u) - \frac{1}{2}(u-1)(3u+2)$$

Alternatively, we have

$$2\log G(1+u) = u\log(2\pi) - \sum_{n=0}^{\infty} \frac{1}{n+1} \sum_{k=0}^{n} \binom{n}{k} (-1)^k (u^2+k^2) \log(u+k)$$

$$+ 2\varsigma'(-1) - \frac{11}{12} + 2u^2 \psi(u) - \frac{1}{2}(u-1)(3u+2)$$

$$= u \log(2\pi) - \sum_{n=0}^{\infty} \frac{1}{n+1} \sum_{k=0}^{n} \binom{n}{k} (-1)^k (u+k)^2 \log(u+k) + 2\sum_{n=0}^{\infty} \frac{1}{n+1} \sum_{k=0}^{n} \binom{n}{k} (-1)^k k \log(u+k)$$

$$+ 2\varsigma'(-1) - \frac{11}{12} + 2u^2\psi(u) - \frac{1}{2}(u-1)(3u+2)$$

Using (4.3.119) we then have

$$= u \log(2\pi) - \sum_{n=0}^{\infty} \frac{1}{n+1} \sum_{k=0}^{n} \binom{n}{k} (-1)^k (u+k)^2 \log(u+k) + 2u \log\Gamma(u) + 2u^2 - u - u\log(2\pi)$$

$$+ 2\varsigma'(-1) - \frac{11}{12} + 2u^2\psi(u) - 2u^2\psi(u) - \frac{1}{2}(u-1)(3u+2)$$

Therefore we obtain

(4.3.132e)
$$2\log G(1+u) = -\sum_{n=0}^{\infty} \frac{1}{n+1} \sum_{k=0}^{n} \binom{n}{k} (-1)^k (u+k)^2 \log(u+k) + 2u \log\Gamma(u) + 2u^2 - u$$

$$+ 2\varsigma'(-1) - \frac{11}{12} - \frac{1}{2}(u-1)(3u+2)$$

Differentiating (4.3.132e) results in

$$2\frac{G'(1+u)}{G(1+u)} = -\sum_{n=0}^{\infty} \frac{1}{n+1} \sum_{k=0}^{n} \binom{n}{k} (-1)^k (u+k) - 2\sum_{n=0}^{\infty} \frac{1}{n+1} \sum_{k=0}^{n} \binom{n}{k} (-1)^k (u+k)\log(u+k)$$

$$+ 2u\psi(u) + 2\log\Gamma(u) + 4u - 1 - 3u + \frac{1}{2}$$

Using (4.3.119) this becomes

$$2\frac{G'(1+u)}{G(1+u)} = -\sum_{n=0}^{\infty} \frac{1}{n+1} \sum_{k=0}^{n} \binom{n}{k} (-1)^k (u+k) + \log(2\pi) + \frac{1}{2} - u + 2u\psi(u)$$

$$= -u + \frac{1}{2} + \log(2\pi) + \frac{1}{2} - u + 2u\psi(u)$$

and integrating this results in

$$\log G(1+x) = \frac{1}{2}x\log(2\pi) + \frac{1}{2}x - \frac{1}{2}x^2 + \int_0^x u\psi(u)\,du$$

Now making reference to (4.3.109) we see that

$$2\log G(1+u) = -2\varsigma'(-1,u) + \varsigma(-1,u) + 2u\log\Gamma(u) + 2u^2 - u + 2\varsigma'(-1) - \frac{11}{12} - \frac{1}{2}(u-1)(3u+2)$$

and this easily simplifies to Gosper's formula

$$\log G(1+u) - u\log\Gamma(u) = \varsigma'(-1) - \varsigma'(-1,u)$$

Therefore, assuming that (4.3.132) $\lim_{x\to 0}[\varsigma'(-p,x) - \varsigma'(-p)] = 0$ is correct, we have the non-trivial result

$$(4.3.133) \qquad \sum_{n=0}^{\infty}\frac{1}{n+1}\sum_{k=0}^{n}\binom{n}{k}(-1)^k k^{p+1}\log k = \sum_{n=0}^{\infty}\frac{1}{n+1}\sum_{k=0}^{n}\binom{n}{k}(-1)^k(1+k)^{p+1}\log(1+k)$$

It is indeed correct for $p=1$ and we therefore have

$$(4.3.133a) \qquad \sum_{n=0}^{\infty}\frac{1}{n+1}\sum_{k=0}^{n}\binom{n}{k}(-1)^k k^2\log k = \sum_{n=0}^{\infty}\frac{1}{n+1}\sum_{k=0}^{n}\binom{n}{k}(-1)^k(1+k)^2\log(1+k)$$

This then gives us

$$\sum_{n=0}^{\infty}\frac{1}{n+1}\sum_{k=0}^{n}\binom{n}{k}(-1)^k k^2\log\frac{k}{1+k} = \sum_{n=0}^{\infty}\frac{1}{n+1}\sum_{k=0}^{n}\binom{n}{k}(-1)^k(1+2k)\log(1+k)$$

$$= -\gamma + 2\left[\frac{1}{2} + \gamma - \frac{1}{2}\log(2\pi)\right]$$

and this concurs with (4.3.77a).

Letting $u=0$ in (4.3.132e) gives us

$$\sum_{n=0}^{\infty}\frac{1}{n+1}\sum_{k=0}^{n}\binom{n}{k}(-1)^k k^2\log k = 2\varsigma'(-1) + \frac{1}{12}$$

in agreement with (4.3.133a) in conjunction with (4.3.109). See also (4.3.137).

By letting $p = 2m$ in (4.3.132c) we obtain

$$(2m+1)\varsigma'(-2m) = \sum_{n=0}^{\infty} \frac{1}{n+1} \sum_{k=0}^{n} \binom{n}{k} (-1)^k (1+k)^{2m+1} \log(1+k)$$

Therefore, using (F.8a)

$$\varsigma'(-2m) = (-1)^m \frac{(2m)!}{2(2\pi)^{2m}} \varsigma(2m+1)$$

we obtain

$$(4.3.133b) \quad \varsigma(2m+1) = (-1)^m \frac{2(2\pi)^{2m}}{(2m+1)!} \sum_{n=0}^{\infty} \frac{1}{n+1} \sum_{k=0}^{n} \binom{n}{k} (-1)^k (k+1)^{2m+1} \log(k+1)$$

and in particular

$$(4.3.133c) \quad \varsigma(3) = -\frac{4}{3}\pi^2 \sum_{n=0}^{\infty} \frac{1}{n+1} \sum_{k=0}^{n} \binom{n}{k} (-1)^k (1+k)^3 \log(1+k)$$

From this we may derive an expression for $\sum_{n=0}^{\infty} \frac{1}{n+1} \sum_{k=0}^{n} \binom{n}{k} (-1)^k k^3 \log(1+k)$ in terms of $\varsigma(3)$ etc.

A different proof of (4.3.133b) is contained in Volume II(b) where we made reference to the Stieltjes constants.

We see from (4.3.110) that

$$(4.3.134) \quad \int_0^1 \varsigma(-n,u)\,du = -\int_0^1 \frac{B_{n+1}(u)}{n+1}\,du = 0 \quad \text{for } n \in \mathbf{N}_o$$

In fact, as noted by Espinosa and Moll [59] and Coffey [45c] for $\mathrm{Re}(s) < 1$ we have

$$(4.3.134a) \quad \int_0^1 \varsigma(s,u)\,du = 0$$

A proof of this appears in Broughan's paper "Vanishing of the integral of the Hurwitz zeta function" [36a]. Hence, using the Hasse representation for $\varsigma(s,u)$, we have for $\mathrm{Re}(s) < 1$

(4.3.135) $$\sum_{n=0}^{\infty}\frac{1}{n+1}\sum_{k=0}^{n}\binom{n}{k}(-1)^k(1+k)^{2-s}=\sum_{n=0}^{\infty}\frac{1}{n+1}\sum_{k=0}^{n}\binom{n}{k}(-1)^k k^{2-s}$$

In fact, we note that (4.3.135) is simply an expression of the identity $\varsigma(s-1,1)=\varsigma(s-1)$ and is therefore prima facie valid for all $s$ except $s=2$. Inspection shows that both of the summations diverge when $s=2$.

Equation (4.3.134a) may also be derived by integrating Hurwitz's formula (4.4.228i) for the Fourier expansion of $\varsigma(p,t)$

$$\varsigma(p,t)=2\Gamma(1-p)\left[\sin\left(\frac{\pi p}{2}\right)\sum_{n=1}^{\infty}\frac{\cos 2n\pi t}{(2\pi n)^{1-p}}+\cos\left(\frac{\pi p}{2}\right)\sum_{n=1}^{\infty}\frac{\sin 2n\pi t}{(2\pi n)^{1-p}}\right]$$

We also have

(4.3.135a) $$\int_{1}^{2}\varsigma(s,u)du=\frac{1}{s-1}$$

and thus we get

(4.3.135b) $$\sum_{n=0}^{\infty}\frac{1}{n+1}\sum_{k=0}^{n}\binom{n}{k}(-1)^k(1+k)^{2-s}=\sum_{n=0}^{\infty}\frac{1}{n+1}\sum_{k=0}^{n}\binom{n}{k}(-1)^k(2+k)^{2-s}$$

Differentiating (4.3.135) we obtain

(4.3.136) $$\sum_{n=0}^{\infty}\frac{1}{n+1}\sum_{k=0}^{n}\binom{n}{k}(-1)^k(1+k)^{2-s}\log(1+k)=\sum_{n=0}^{\infty}\frac{1}{n+1}\sum_{k=0}^{n}\binom{n}{k}(-1)^k k^{2-s}\log k$$

and, letting $s=1-p$ where $p\neq-1$, we may write this as

(4.3.136a) $$\sum_{n=0}^{\infty}\frac{1}{n+1}\sum_{k=0}^{n}\binom{n}{k}(-1)^k(1+k)^{p+1}\log(1+k)=\sum_{n=0}^{\infty}\frac{1}{n+1}\sum_{k=0}^{n}\binom{n}{k}(-1)^k k^{p+1}\log k$$

which is of course simply (4.3.133). See also (4.3.137a) below. Further differentiations of (4.3.136) with respect to $s$ result in

(4.3.136b) $$\sum_{n=0}^{\infty}\frac{1}{n+1}\sum_{k=0}^{n}\binom{n}{k}(-1)^k(1+k)^{2-s}\log^m(1+k)=\sum_{n=0}^{\infty}\frac{1}{n+1}\sum_{k=0}^{n}\binom{n}{k}(-1)^k k^{2-s}\log^m k$$

(4.3.136c) $$\sum_{n=0}^{\infty}\frac{1}{n+1}\sum_{k=0}^{n}\binom{n}{k}(-1)^k(1+k)^{p+1}\log^m(1+k)=\sum_{n=0}^{\infty}\frac{1}{n+1}\sum_{k=0}^{n}\binom{n}{k}(-1)^k k^{p+1}\log^m k$$

Since we have already shown in (4.3.128) that

$$\varsigma'(-1) = -\frac{1}{24} + \frac{1}{2}\sum_{n=0}^{\infty}\frac{1}{n+1}\sum_{k=0}^{n}\binom{n}{k}(-1)^k(1+k)^2\log(1+k)$$

we deduce

(4.3.137) $\qquad 2\varsigma'(-1) + \frac{1}{12} = \sum_{n=0}^{\infty}\frac{1}{n+1}\sum_{k=0}^{n}\binom{n}{k}(-1)^k k^2 \log k$

We note that (4.3.137) is equivalent to (4.3.109b) as shown below

$$\varsigma'(-1,0) = \frac{1}{2}\varsigma(-1,0) + \frac{1}{2}\sum_{n=0}^{\infty}\frac{1}{n+1}\sum_{k=0}^{n}\binom{n}{k}(-1)^k k^2 \log k$$

With $p = 1$ in (4.3.136c) we recover (4.3.106)

$$\sum_{n=0}^{\infty}\frac{1}{n+1}\sum_{k=0}^{n}\binom{n}{k}(-1)^k k^2 \log k = 1 + \gamma - \log(2\pi) + \sum_{n=0}^{\infty}\frac{1}{n+1}\sum_{k=0}^{n}\binom{n}{k}(-1)^k k^2 \log(1+k)$$

By combining (4.3.136) with (4.3.77a) we obtain another derivation of (4.3.112).

Espinosa and Moll [59] also report that for $s \in \mathbf{R}_0^-$

(4.3.138) $\qquad \int_0^{1/2} \varsigma(s,u)\,du = 2\frac{(2^{s-2}-1)}{1-s}\varsigma(s-1)$

Hence, using the Hasse representation for $\varsigma(s,u)$, we have for $s \in \mathbf{R}_0^-$

$$\sum_{n=0}^{\infty}\frac{1}{n+1}\sum_{k=0}^{n}\binom{n}{k}(-1)^k\left(\frac{1}{2}+k\right)^{2-s} - \sum_{n=0}^{\infty}\frac{1}{n+1}\sum_{k=0}^{n}\binom{n}{k}(-1)^k k^{2-s} = 2(2^{s-2}-1)(s-2)\varsigma(s-1)$$

or equivalently

(4.3.139)

$$2^{s-2}\sum_{n=0}^{\infty}\frac{1}{n+1}\sum_{k=0}^{n}\binom{n}{k}(-1)^k(2k+1)^{2-s} - \sum_{n=0}^{\infty}\frac{1}{n+1}\sum_{k=0}^{n}\binom{n}{k}(-1)^k k^{2-s} = 2(2^{s-2}-1)(s-2)\varsigma(s-1)$$

With $s = 0$ we get using $\varsigma(-1) = -\dfrac{1}{12}$

$$\frac{1}{4}\sum_{n=0}^{\infty}\frac{1}{n+1}\sum_{k=0}^{n}\binom{n}{k}(-1)^k(2k+1)^2 - \sum_{n=0}^{\infty}\frac{1}{n+1}\sum_{k=0}^{n}\binom{n}{k}(-1)^k k^2 = 3\varsigma(-1) = -\frac{1}{4}$$

and this is readily verified by expanding the term $(2k+1)^2$ in the expansion. We also see that with $s = -1$

$$\frac{1}{8}\sum_{n=0}^{\infty}\frac{1}{n+1}\sum_{k=0}^{n}\binom{n}{k}(-1)^k(2k+1)^3 - \sum_{n=0}^{\infty}\frac{1}{n+1}\sum_{k=0}^{n}\binom{n}{k}(-1)^k k^3 = 0$$

and more generally with $s = -2m+1$

$$2^{-2m-1}\sum_{n=0}^{\infty}\frac{1}{n+1}\sum_{k=0}^{n}\binom{n}{k}(-1)^k(2k+1)^{2m+1} - \sum_{n=0}^{\infty}\frac{1}{n+1}\sum_{k=0}^{n}\binom{n}{k}(-1)^k k^{2m+1} = 0$$

Differentiating (4.3.139) results in

$$-2^{s-2}\sum_{n=0}^{\infty}\frac{1}{n+1}\sum_{k=0}^{n}\binom{n}{k}(-1)^k(2k+1)^{2-s}\log(2k+1) + 2^{s-2}\log 2\sum_{n=0}^{\infty}\frac{1}{n+1}\sum_{k=0}^{n}\binom{n}{k}(-1)^k(2k+1)^{2-s}$$

$$+\sum_{n=0}^{\infty}\frac{1}{n+1}\sum_{k=0}^{n}\binom{n}{k}(-1)^k k^{2-s}\log k$$

$$= 2(2^{s-2}-1)(s-2)\varsigma'(s-1) + 2(2^{s-2}-1)\varsigma(s-1) + 2^{s-1}\log 2(s-2)\varsigma(s-1)$$

and $s = 0$ gives us

$$-\frac{1}{4}\sum_{n=0}^{\infty}\frac{1}{n+1}\sum_{k=0}^{n}\binom{n}{k}(-1)^k(2k+1)^2\log(2k+1) + \frac{1}{4}\log 2\sum_{n=0}^{\infty}\frac{1}{n+1}\sum_{k=0}^{n}\binom{n}{k}(-1)^k(2k+1)^2$$

$$+\sum_{n=0}^{\infty}\frac{1}{n+1}\sum_{k=0}^{n}\binom{n}{k}(-1)^k k^2\log k = 3\varsigma'(-1) - \left(\frac{3}{2} + \log 2\right)\varsigma(-1)$$

$$= 3\varsigma'(-1) + \frac{1}{12}\log 2 + \frac{1}{8}$$

Using (4.3.110) this may be written as

$$-\frac{1}{4}\sum_{n=0}^{\infty}\frac{1}{n+1}\sum_{k=0}^{n}\binom{n}{k}(-1)^k(2k+1)^2\log(2k+1)+\sum_{n=0}^{\infty}\frac{1}{n+1}\sum_{k=0}^{n}\binom{n}{k}(-1)^k k^2\log k$$

$$=3\varsigma'(-1)+\frac{1}{6}\log 2+\frac{1}{8}$$

and using (4.3.109) we obtain

$$\varsigma'\left(-1,\frac{1}{2}\right)=\frac{1}{2}\varsigma\left(-1,\frac{1}{2}\right)+\frac{1}{8}\sum_{n=0}^{\infty}\frac{1}{n+1}\sum_{k=0}^{n}\binom{n}{k}(-1)^k(2k+1)^2\log\frac{2k+1}{2}$$

$$=\frac{1}{2}\varsigma\left(-1,\frac{1}{2}\right)+\frac{1}{8}\sum_{n=0}^{\infty}\frac{1}{n+1}\sum_{k=0}^{n}\binom{n}{k}(-1)^k(2k+1)^2\log(2k+1)-\frac{1}{8}\log 2\sum_{n=0}^{\infty}\frac{1}{n+1}\sum_{k=0}^{n}\binom{n}{k}(-1)^k(2k+1)^2$$

$$=\frac{1}{24}+\frac{1}{8}\sum_{n=0}^{\infty}\frac{1}{n+1}\sum_{k=0}^{n}\binom{n}{k}(-1)^k(2k+1)^2\log(2k+1)+\frac{1}{24}\log 2$$

Hence we have

$$\varsigma'\left(-1,\frac{1}{2}\right)-\frac{1}{48}-\frac{1}{24}\log 2=\frac{1}{8}\sum_{n=0}^{\infty}\frac{1}{n+1}\sum_{k=0}^{n}\binom{n}{k}(-1)^k(2k+1)^2\log(2k+1)$$

and thus

$$-2\varsigma'\left(-1,\frac{1}{2}\right)+\frac{1}{24}+\frac{1}{12}\log 2+\sum_{n=0}^{\infty}\frac{1}{n+1}\sum_{k=0}^{n}\binom{n}{k}(-1)^k k^2\log k$$

$$=3\varsigma'(-1)+\frac{1}{6}\log 2+\frac{1}{8}$$

This gives us

$$3\varsigma'(-1)=-2\varsigma'\left(-1,\frac{1}{2}\right)-\frac{1}{24}-\frac{1}{12}\log 2+\sum_{n=0}^{\infty}\frac{1}{n+1}\sum_{k=0}^{n}\binom{n}{k}(-1)^k k^2\log k$$

We have from (4.3.109b) since $\varsigma'(-1)=\varsigma'(-1,0)$

$$2\varsigma'(-1)=\sum_{n=0}^{\infty}\frac{1}{n+1}\sum_{k=0}^{n}\binom{n}{k}(-1)^k k^2\log k-\frac{1}{12}$$

and hence we see that

$$(4.3.140) \qquad \varsigma'\left(-1, \frac{1}{2}\right) = -\frac{1}{24}\log 2 - \frac{1}{2}\varsigma'(-1)$$

this is in accordance with [103] where Adamchik showed that

$$(4.3.141) \qquad \varsigma'\left(-2k+1, \frac{1}{2}\right) = -\frac{B_{2k}\log 2}{4^k\,k} - \frac{(2^{2k-1}-1)}{2^{2k-1}}\varsigma'(-2k+1)$$

and in particular

$$\varsigma'\left(-1, \frac{1}{2}\right) = -\frac{B_2\log 2}{4} - \frac{1}{2}\varsigma'(-1) = -\frac{1}{24}\log 2 - \frac{1}{2}\varsigma'(-1)$$

$\square$

A more direct method of obtaining Adamchik's result (4.3.130) is set out below. We have

$$(4.3.142) \qquad \frac{\partial}{\partial u}\frac{\partial}{\partial s}\varsigma(s,u) = \frac{\partial}{\partial s}\frac{\partial}{\partial u}\varsigma(s,u) = -\frac{\partial}{\partial s}[s\varsigma(s+1,u)]$$

$$= -\varsigma(s+1,u) - s\frac{\partial}{\partial s}\varsigma(s+1,u)$$

and upon integrating we see that

$$-s\int_0^x \varsigma'(s+1,u)\,du = \int_0^x \frac{\partial}{\partial u}\frac{\partial}{\partial s}\varsigma(s,u)\,du + \int_0^x \varsigma(s+1,u)\,du$$

We therefore get

$$(4.3.143) \qquad -s\int_0^x \varsigma'(s+1,u)\,du = \varsigma'(s,x) - \varsigma'(s,0) + \int_0^x \varsigma(s+1,u)\,du$$

and with $s = -n$ we have

$$(4.3.144) \qquad n\int_0^x \varsigma'(1-n,u)\,du = \varsigma'(-n,x) - \varsigma'(-n,0) + \int_0^x \varsigma(1-n,u)\,du$$

Then using $\varsigma(1-n,u) = -\dfrac{B_n(u)}{n}$ for $n \geq 1$ we obtain (4.3.131) again

$$n \int_0^x \varsigma'(1-n,u)\,du = \frac{B_{n+1} - B_{n+1}(x)}{n(n+1)} + \varsigma'(-n,x) - \varsigma'(-n,0)$$

and from this we immediately see that

(4.3.145) $$\int_0^1 \varsigma'(1-n,u)\,du = 0$$

Similarly we see that

$$\frac{\partial}{\partial u}\frac{\partial^2}{\partial s^2}\varsigma(s,u) = -2\varsigma'(s+1,u) - s\varsigma''(s+1,u)$$

and upon integrating this we see that

$$-s\int_0^x \varsigma''(s+1,u)\,du = \int_0^x \frac{\partial}{\partial u}\frac{\partial^2}{\partial s^2}\varsigma(s,u)\,du + 2\int_0^x \varsigma'(s+1,u)\,du$$

and hence

(4.3.146)

$$-s\int_0^x \varsigma''(s+1,u)\,du = \varsigma''(s,x) - \varsigma''(s,0) - \frac{2}{s}[\varsigma'(s,x) - \varsigma'(s,0)] - \frac{2}{s}\int_0^x \varsigma(s+1,u)\,du$$

With $s = -n$ we have

(4.3.147)

$$n\int_0^x \varsigma''(1-n,u)\,du = \varsigma''(1-n,x) - \varsigma''(1-n,0) + \frac{2}{n}[\varsigma'(1-n,x) - \varsigma'(1-n,0)] + \frac{2[B_{n+1} - B_{n+1}(x)]}{n^2(n+1)}$$

We easily see that

(4.3.148) $$\int_0^1 \varsigma''(1-n,u)\,du = 0$$

Equivalent results for $\int_0^x \varsigma^{(p)}(1-n,u)\,du$ may be easily obtained and it is fairly obvious that we have

(4.3.149)
$$\int_0^1 \varsigma^{(p)}(1-n,u)du = 0$$

We also note from Adamchik's integral (4.3.131) that

$$\int_0^1 \varsigma'(1-n,u)\,du = 0$$

and using (4.3.107a) we see that

$$-n\varsigma'(1-n,u) + \varsigma(1-n,u) = -\sum_{n=0}^{\infty} \frac{1}{n+1} \sum_{k=0}^{n} \binom{n}{k}(-1)^k (u+k)^n \log(u+k)$$

Therefore, referring to (4.3.110) we have $\varsigma(1-n,u) = -\dfrac{B_n(u)}{n}$ and hence

$$\frac{1}{n}\int_0^1 B_n(u)du = \sum_{n=0}^{\infty} \frac{1}{n+1} \sum_{k=0}^{n} \binom{n}{k}(-1)^k \int_0^1 (u+k)^n \log(u+k)\,du$$

Therefore we see that

$$\sum_{n=0}^{\infty} \frac{1}{n+1} \sum_{k=0}^{n} \binom{n}{k}(-1)^k \int_0^1 (u+k)^n \log(u+k)\,du = 0$$

and since

$$\int_0^1 (u+k)^n \log(u+k)du = \frac{(1+k)^{n+1}}{n+1}\log(1+k) - \frac{(1+k)^{n+1}}{(n+1)^2} - \frac{k^{n+1}}{n+1}\log k + \frac{k^{n+1}}{(n+1)^2}$$

we may accordingly show that

(4.3.150)    $$\varsigma'(-n,0) = \varsigma'(-n)\,.$$

We recall from (4.3.126) that

$$\log G(u+1) - u\log\Gamma(u) = \varsigma'(-1) - \varsigma'(-1,u)$$

and, upon integration using Adamchik's integral (4.3.131),we obtain

(4.3.150a) $\displaystyle\int_0^x \log G(u+1)\,du - \int_0^x u \log \Gamma(u)\,du = x\varsigma'(-1) + \frac{1}{12}B_3(x) - \frac{1}{2}\varsigma'(-2,x) + \frac{1}{2}\varsigma'(-2)$

since $\displaystyle 2\int_0^x \varsigma'(-1,u)\,du = \frac{B_3 - B_3(x)}{6} + \varsigma'(-2,x) - \varsigma'(-2)$. With $x=1$ we obtain

$$\int_0^1 \log G(u+1)\,du - \int_0^1 u \log \Gamma(u)\,du = \varsigma'(-1)$$

We see from [126, p.217]

(4.3.151) $\displaystyle \int_0^1 \log G(u+1)\,du = \frac{1}{4}\log(2\pi) + \frac{1}{12} - 2\log A$

and hence we obtain

$$\int_0^1 u \log \Gamma(u)\,du = \frac{1}{4}\log(2\pi) + \frac{1}{12} - 2\log A - \varsigma'(-1)$$

From (6.126) in Volume V we see that

$$\int_0^1 u \log \Gamma(1+u)\,du = \frac{1}{4}\log(2\pi) - \frac{1}{4} - \log A$$

and therefore we have

$$\int_0^1 u \log \Gamma(u)\,du = \frac{1}{4}\log(2\pi) - \log A$$

We then obtain the well-known result

(4.3.152) $\displaystyle \log A = \frac{1}{12} - \varsigma'(-1) = -[\varsigma(-1,0) + \varsigma'(-1)]$

Using (4.3.152) we may write

$$\log G(u+1) = u\varsigma'(0,u) + \frac{1}{2}u\log(2\pi) + \varsigma'(-1) - \varsigma'(-1,u)$$

and this may be written as a Lerch-type identity

(4.3.152a)  $\log G(u+1) = u[\varsigma'(0,u) + \log A_0] - [\varsigma(-1,u) + \log A_1] - \varsigma'(-1,u)$

where we have from [45]

$$\log A_n = \frac{B_{n+1}H_n}{n+1} - \varsigma'(-n)$$

which results in

$$\log A_0 = \frac{1}{2}\log(2\pi)$$

and we have used the notation $A_1 = A$.

More generally, using (4.3.87ai) and (4.3.87c), we see that

$$\int_0^x \log G(1+t)dt - \int_0^x t\log\Gamma(1+t)dt =$$

$$\left(\frac{1}{4} - 2\log A\right)x + \frac{1}{4}\log(2\pi)x^2 - \frac{1}{6}x^3 + (x-1)\log G(1+x) - 2\log\Gamma_3(1+x)$$

$$-\frac{1}{2}\left(\frac{1}{4} - 2\log A\right)x - \frac{1}{2}\left(\frac{1}{2}\log(2\pi) - \frac{1}{4}\right)x^2 + \frac{1}{4}x^3 - \frac{1}{2}x^2\log\Gamma(1+x)$$

$$+\frac{1}{2}\log G(1+x) + \log\Gamma_3(1+x)$$

$$=\left(\frac{1}{8} - \log A\right)x + \frac{1}{8}x^2 + \frac{1}{12}x^3 - \frac{1}{2}x^2\log\Gamma(1+x) + \left(x - \frac{1}{2}\right)\log G(1+x) - \log\Gamma_3(1+x)$$

and we therefore get

$$\int_0^x \log G(1+t)dt - \int_0^x t\log\Gamma(t)dt =$$

$$\left(\frac{1}{8} - \log A\right)x - \frac{1}{8}x^2 + \frac{1}{12}x^3 + \frac{1}{2}x^2\log x - \frac{1}{2}x^2\log\Gamma(1+x) + \left(x - \frac{1}{2}\right)\log G(1+x) - \log\Gamma_3(1+x)$$

Hence we have using (4.3.126)

(4.3.152b) $\qquad x\varsigma'(-1) + \dfrac{1}{12}B_3(x) - \dfrac{1}{2}[\varsigma'(-2,x) - \varsigma'(-2)] =$

$$\left(\frac{1}{8} - \log A\right)x - \frac{1}{8}x^2 + \frac{1}{12}x^3 + \frac{1}{2}x^2\log x - \frac{1}{2}x^2\log\Gamma(1+x) + \left(x - \frac{1}{2}\right)\log G(1+x) - \log\Gamma_3(1+x)$$

With $x = 1$ this reduces to the familiar result $\varsigma'(-1) = \dfrac{1}{12} - \log A$ .

This formula would then enable us to determine, for example, $\log\Gamma_3(3/2)$ and, using (4.3.130) we may also evaluate $\displaystyle\int_0^u \log\Gamma_3(1+x)dx$ (see for example [126, p.207]).

With $x = 1/2$ we get (noting that the term involving $\log G\left(\dfrac{3}{2}\right)$ conveniently cancels)

(4.3.152c) $\qquad \dfrac{1}{2}\varsigma'(-1) + \dfrac{1}{12}B_3\left(\dfrac{1}{2}\right) - \dfrac{1}{2}[\varsigma'\left(-2,\dfrac{1}{2}\right) - \varsigma'(-2)] =$

$$\frac{1}{2}\left(\frac{1}{8} - \log A\right) - \frac{1}{48} - \frac{1}{8}\log 2 - \frac{1}{8}\log\Gamma\left(\frac{3}{2}\right) - \log\Gamma_3\left(\frac{3}{2}\right)$$

and skipping ahead to (4.3.168d) we see that $\varsigma'\left(-2,\dfrac{1}{2}\right) = \dfrac{3\varsigma(3)}{16\pi^2}$ and from (F.8b) we have

$\varsigma'(-2) = -\dfrac{\varsigma(3)}{4\pi^2}$ . This then provides us with

(4.3.152d) $\qquad \log\Gamma_3\left(\dfrac{3}{2}\right) = \dfrac{7\varsigma(3)}{32\pi^2} - \dfrac{1}{16}\log\pi$

We have [126, p.39]

(4.3.152e) $\qquad G_n(x+1) = G_n(x)G_{n-1}(x) \qquad \Gamma_n(x) = \left[G_n(x)\right]^{(-1)^{n-1}}$

and it is easily seen that

$$\log G_n(x+1) = \log G_n(x) + \log G_{n-1}(x)$$

$$\log\Gamma_n(x) = (-1)^{n-1}\log G_n(x)$$

and from this we obtain

(4.3.152f) $$\Gamma_n(x+1) = \frac{\Gamma_n(x)}{\Gamma_{n-1}(x)}$$

Therefore we get

$$\Gamma_3(3/2) = \frac{\Gamma_3(1/2)}{\Gamma_2(1/2)}$$

As shown by Choi, Cho and Srivastava [45ac] we have

(4.3.152g) $$\Gamma_2\left(\frac{1}{2}\right) = 2^{-\frac{1}{24}} \pi^{\frac{1}{4}} e^{-\frac{1}{8}} A^{\frac{1}{2}}$$

(4.3.152h) $$\Gamma_3\left(\frac{1}{2}\right) = 2^{-\frac{1}{24}} \pi^{\frac{3}{16}} e^{-\frac{1}{8}} A^{\frac{1}{2}} \exp\frac{7\varsigma(3)}{32\pi^2}$$

and thus we have

$$\Gamma_3(3/2) = \pi^{-\frac{1}{16}} \exp\frac{7\varsigma(3)}{32\pi^2}$$

in accordance with (4.3.152d) above.

Choi, Cho and Srivastava [45ac] show that

(4.3.152i) $$\log\Gamma_3(1+x) = -\frac{1}{24} + \frac{1}{2}\log A + \frac{\varsigma(3)}{8\pi^2} + \left(\frac{1}{12} - \log A\right)x + \frac{1}{2}(x^2 - x)\log\Gamma(1+x)$$

$$+ \left(x - \frac{1}{2}\right)\varsigma'(-1,1+x) + \left(x - \frac{1}{2}\right)\varsigma'(-1,1+x) + \frac{1}{2}\varsigma'(-2,1+x)$$

and using (4.3.126i) it is easily seen that both expressions are equivalent.

Let us now integrate (4.3.84a) over the interval $[1,x]$ to obtain

$$\sum_{k=0}^{\infty}\frac{1}{k+1}\sum_{j=0}^{k}\binom{k}{j}(-1)^{j}j^{2}(x+j)\log(x+j)-\sum_{k=0}^{\infty}\frac{1}{k+1}\sum_{j=0}^{k}\binom{k}{j}(-1)^{j}j^{2}(1+j)\log(1+j)$$

$$-(x-1)\sum_{k=0}^{\infty}\frac{1}{k+1}\sum_{j=0}^{k}\binom{k}{j}(-1)^{j}j^{2}-(x-1)\sum_{k=0}^{\infty}\frac{1}{k+1}\sum_{j=0}^{k}\binom{k}{j}(-1)^{j}j^{2}\log(1+j)=$$

$$\left(\frac{1}{2}x^{2}-x+\frac{1}{2}\right)\log(2\pi)-2\int_{1}^{x}\log G(1+u)\,du+x^{2}\log\Gamma(x)-2\int_{1}^{x}u\log\Gamma(u)\,du$$

$$+\gamma(x-1)-\frac{1}{2}x^{3}+\frac{1}{4}x^{2}+x-\frac{3}{4}$$

Having regard to the first term of the above equation, we have by simple algebra

$$\sum_{k=0}^{\infty}\frac{1}{k+1}\sum_{j=0}^{k}\binom{k}{j}(-1)^{j}j^{2}(x+j)\log(x+j)=\sum_{k=0}^{\infty}\frac{1}{k+1}\sum_{j=0}^{k}\binom{k}{j}(-1)^{j}(x+j)^{3}\log(x+j)$$

$$-2x\sum_{k=0}^{\infty}\frac{1}{k+1}\sum_{j=0}^{k}\binom{k}{j}(-1)^{j}(x+j)^{2}\log(x+j)+x^{2}\sum_{k=0}^{\infty}\frac{1}{k+1}\sum_{j=0}^{k}\binom{k}{j}(-1)^{j}(x+j)\log(x+j)$$

Referring to (4.3.108) and (4.3.109) we see that

$$\varsigma'(0,x)-\varsigma(0,x)=\sum_{k=0}^{\infty}\frac{1}{k+1}\sum_{j=0}^{k}\binom{k}{j}(-1)^{j}(x+j)\log(x+j)$$

$$2\varsigma'(-1,x)-\varsigma(-1,x)=\sum_{k=0}^{\infty}\frac{1}{k+1}\sum_{j=0}^{k}\binom{k}{j}(-1)^{j}(x+j)^{2}\log(x+j)$$

and we also have

$$3\varsigma'(-2,x)-\varsigma(-2,x)=\sum_{k=0}^{\infty}\frac{1}{k+1}\sum_{j=0}^{k}\binom{k}{j}(-1)^{j}(x+j)^{3}\log(x+j)$$

Hence we get

$$\sum_{k=0}^{\infty}\frac{1}{k+1}\sum_{j=0}^{k}\binom{k}{j}(-1)^j j^2(x+j)\log(x+j) = 3\varsigma'(-2,x) - \varsigma(-2,x) - 2x[2\varsigma'(-1,x) - \varsigma(-1,x)]$$

$$+x^2[\varsigma'(0,x) - \varsigma(0,x)]$$

and letting $x = 1$ in the above equation results in

$$\sum_{k=0}^{\infty}\frac{1}{k+1}\sum_{j=0}^{k}\binom{k}{j}(-1)^j j^2(1+j)\log(1+j) = 3\varsigma'(-2) - \varsigma(-2) - 2[2\varsigma'(-1) - \varsigma(-1)]$$

$$+[\varsigma'(0) - \varsigma(0)]$$

$$= -\frac{3\varsigma(3)}{4\pi^2} - 4\varsigma'(-1) - \frac{1}{6} - \frac{1}{2}\log(2\pi) + \frac{1}{2}$$

$$= -\frac{3\varsigma(3)}{4\pi^2} - 4\varsigma'(-1) - \frac{1}{2}\log(2\pi) + \frac{1}{3}$$

We also see that

$$\sum_{k=0}^{\infty}\frac{1}{k+1}\sum_{j=0}^{k}\binom{k}{j}(-1)^j j^2\log(1+j) = \sum_{k=0}^{\infty}\frac{1}{k+1}\sum_{j=0}^{k}\binom{k}{j}(-1)^j(j+1)^2\log(1+j)$$

$$-2\sum_{k=0}^{\infty}\frac{1}{k+1}\sum_{j=0}^{k}\binom{k}{j}(-1)^j(j+1)\log(1+j)$$

$$+\sum_{k=0}^{\infty}\frac{1}{k+1}\sum_{j=0}^{k}\binom{k}{j}(-1)^j\log(1+j)$$

$$= 2\varsigma'(-1) - \varsigma(-1) - 2[\varsigma'(0) - \varsigma(0)] - \gamma$$

Hence we get (as previously seen in (4.3.132d))

$$\sum_{k=0}^{\infty}\frac{1}{k+1}\sum_{j=0}^{k}\binom{k}{j}(-1)^j j^2\log(1+j) = 2\varsigma'(-1) + \log(2\pi) - \frac{11}{12} - \gamma$$

From (F.20) in Volume VI and (4.3.110) we see that

$$\sum_{k=0}^{\infty}\frac{1}{k+1}\sum_{j=0}^{k}\binom{k}{j}(-1)^j j^m = B_m(0) = B_m$$

and hence

$$\sum_{k=0}^{\infty}\frac{1}{k+1}\sum_{j=0}^{k}\binom{k}{j}(-1)^j j^2 = B_2 = \frac{1}{6}$$

Accordingly we obtain

$$3\varsigma'(-2,x) - \varsigma(-2,x) - 2x[2\varsigma'(-1,x) - \varsigma(-1,x)] + x^2[\varsigma'(0,x) - \varsigma(0,x)]$$

$$+\frac{3\varsigma(3)}{4\pi^2} + 4\varsigma'(-1) + \frac{1}{2}\log(2\pi) - \frac{1}{3} - \frac{1}{6}(x-1) - (x-1)\left[2\varsigma'(-1) + \log(2\pi) - \frac{11}{12} - \gamma\right]$$

$$=\left(\frac{1}{2}x^2 - x + \frac{1}{2}\right)\log(2\pi) - 2\int_1^x \log G(1+u)\,du + x^2 \log\Gamma(x) - 2\int_1^x u\log\Gamma(u)\,du$$

$$+\gamma(x-1) - \frac{1}{2}x^3 + \frac{1}{4}x^2 + x - \frac{3}{4}$$

From (4.3.87ai) we see that

$$2\int_0^x u\log\Gamma(1+u)\,du =$$

$$\left(\frac{1}{4} - 2\log A\right)x + \left(\frac{1}{2}\log(2\pi) - \frac{1}{4}\right)x^2 - \frac{1}{2}x^3 + x^2\log\Gamma(1+x) - \log G(1+x) - 2\log\Gamma_3(1+x)$$

and therefore we have

$$2\int_0^1 u\log\Gamma(1+u)\,du = \frac{1}{2}\log(2\pi) - 2\log A - \frac{1}{2}$$

Since $2\int_1^x u\log\Gamma(1+u)\,du = 2\int_0^x u\log\Gamma(1+u)\,du - 2\int_0^1 u\log\Gamma(1+u)\,du$ we have

$$2\int_1^x u\log\Gamma(1+u)\,du =$$

$$\left(\frac{1}{4} - 2\log A\right)x + \left(\frac{1}{2}\log(2\pi) - \frac{1}{4}\right)x^2 - \frac{1}{2}x^3 + x^2\log\Gamma(1+x) - \log G(1+x) - 2\log\Gamma_3(1+x)$$

$$-\frac{1}{2}\log(2\pi) + 2\log A + \frac{1}{2}$$

It is seen that

$$2\int_1^x u\log\Gamma(u)\,du = 2\int_1^x u\log\Gamma(u+1)\,du - 2\int_1^x u\log u\,du$$

$$= 2\int_1^x u\log\Gamma(u+1)\,du - x^2\log x + \frac{1}{2}x^2$$

$$= \left(\frac{1}{4} - 2\log A\right)x + \left(\frac{1}{2}\log(2\pi) - \frac{1}{4}\right)x^2 - \frac{1}{2}x^3 + x^2\log\Gamma(1+x) - \log G(1+x) - 2\log\Gamma_3(1+x)$$

$$-\frac{1}{2}\log(2\pi) + 2\log A + \frac{1}{2} - x^2\log x - \frac{1}{2}x^2$$

We then obtain

$$3\varsigma'(-2,x) - \varsigma(-2,x) - 2x[2\varsigma'(-1,x) - \varsigma(-1,x)] + x^2[\varsigma'(0,x) - \varsigma(0,x)]$$

$$+\frac{3\varsigma(3)}{4\pi^2} + 4\varsigma'(-1) + \frac{1}{2}\log(2\pi) - \frac{1}{3} - \frac{1}{6}(x-1) - (x-1)\left[2\varsigma'(-1) + \log(2\pi) - \frac{11}{12} - \gamma\right]$$

$$= \left(\frac{1}{2}x^2 - x + \frac{1}{2}\right)\log(2\pi) - 2\int_1^x \log G(1+u)\,du + x^2\log\Gamma(x)$$

$$+ \gamma(x-1) - \frac{1}{2}x^3 + x^2 + x - \frac{3}{2}$$

$$-\left(\frac{1}{4} - 2\log A\right)x - \left(\frac{1}{2}\log(2\pi) - \frac{1}{4}\right)x^2 + \frac{1}{2}x^3 - x^2\log\Gamma(1+x) + \log G(1+x) + 2\log\Gamma_3(1+x)$$

$$+\frac{1}{2}\log(2\pi) - 2\log A - \frac{1}{2} + x^2\log x - \frac{1}{2}x^2$$

This then gives us

$$2\int_1^x \log G(1+u)\,du =$$

$$-3\varsigma'(-2,x)+\varsigma(-2,x)+2x[2\varsigma'(-1,x)-\varsigma(-1,x)]-x^2[\varsigma'(0,x)-\varsigma(0,x)]$$

$$-\frac{3\varsigma(3)}{4\pi^2}-4\varsigma'(-1)-\frac{1}{2}\log(2\pi)+\frac{1}{3}+\frac{1}{6}(x-1)+(x-1)\left[2\varsigma'(-1)+\log(2\pi)-\frac{11}{12}-\gamma\right]$$

$$+\left(\frac{1}{2}x^2-x+\frac{1}{2}\right)\log(2\pi)+x^2\log\Gamma(x)+\gamma(x-1)-\frac{1}{2}x^3+x^2+x-\frac{3}{2}$$

$$-\left(\frac{1}{4}-2\log A\right)x-\left(\frac{1}{2}\log(2\pi)-\frac{1}{4}\right)x^2+\frac{1}{2}x^3-x^2\log\Gamma(1+x)+\log G(1+x)+2\log\Gamma_3(1+x)$$

$$+\frac{1}{2}\log(2\pi)-2\log A-\frac{1}{2}+x^2\log x-\frac{1}{2}x^2$$

Using (4.3.152i)

$$\log\Gamma_3(1+x)=-\frac{1}{24}+\frac{1}{2}\log A+\frac{\varsigma(3)}{8\pi^2}+\left(\frac{1}{12}-\log A\right)x+\frac{1}{2}(x^2-x)\log\Gamma(1+x)$$

$$+\left(x-\frac{1}{2}\right)\varsigma'(-1,1+x)+\left(x-\frac{1}{2}\right)\varsigma'(-1,1+x)+\frac{1}{2}\varsigma'(-2,1+x)$$

and (4.3.180a), we may show that this is equivalent to the result obtained by Choi, Cho and Srivastava [45ac] in 2004

$$-\int_0^x\log G(1+u)\,du=\int_0^x\log\Gamma_2(1+u)\,du=\frac{\varsigma(3)}{4\pi^2}-\left[\frac{1}{12}+\varsigma'(-1)\right]x-\frac{1}{4}x^2\log(2\pi)$$

$$+\frac{1}{6}x^3-x\varsigma'(-1,1+x)+\varsigma'(-2,1+x)$$

$$\text{EVALUATION OF }\int_0^x\pi u^n\cot\pi u\,du$$

In this section we shall show how the above integrals are closely related to the derivatives of the Hurwitz zeta function.

Using Euler's reflection formula for the gamma function for $0<u<1$

$$\Gamma(u)\Gamma(1-u)=\frac{\pi}{\sin\pi u}$$

we have

$$\log \Gamma(u) + \log \Gamma(1-u) = \log \pi - \log \sin \pi u$$

and differentiation results in the well-known formula

$$\psi(1-u) - \psi(u) = \pi \cot \pi u$$

Using the logarithmic expansion for the digamma function (4.3.74), we have upon multiplication by $u$

(4.3.154) $$u \sum_{k=0}^{\infty} \frac{1}{k+1} \sum_{j=0}^{k} \binom{k}{j} (-1)^j \log \frac{j+1-u}{j+u} = \pi u \cot \pi u$$

and integration then gives us

(4.3.155) $$\sum_{k=0}^{\infty} \frac{1}{k+1} \sum_{j=0}^{k} \binom{k}{j} (-1)^j \int_0^x u \log \frac{j+1-u}{j+u} \, du = \int_0^x \pi u \cot \pi u \, du$$

Integration by parts gives us

$$2\int_0^x u \log \frac{j+1-u}{j+u} \, du = \left[ u^2 - (j+1)^2 \right] \log(j+1-u) - (x^2 - j^2) \log(u+j) - (2j+1)u \Big|_0^x$$

and we have

$$2\int_0^x \pi u \cot \pi u \, du = F_1(x) - F_1(0)$$

where

$$F_1(x) =$$

$$\sum_{k=0}^{\infty} \frac{1}{k+1} \sum_{j=0}^{k} \binom{k}{j} (-1)^j \left( \left[ x^2 - (j+1)^2 \right] \log(j+1-x) - (x^2 - j^2) \log(x+j) - (2j+1)x \right)$$

Since $j^p = [(j+b) - b]^p$ we have

$$j^2 = (j+b)^2 - 2b(j+b) + b^2$$

$$(j+1)^2 = (j+1-b)^2 + 2b(j+1-b) + b^2$$

$F_1(x)$ may be written as

$$-\sum_{k=0}^{\infty}\frac{1}{k+1}\sum_{j=0}^{k}\binom{k}{j}(-1)^j\left\{(j+1-x)^2+2x(j+1-x)\right\}\log(j+1-x)$$

$$+\sum_{k=0}^{\infty}\frac{1}{k+1}\sum_{j=0}^{k}\binom{k}{j}(-1)^j\left\{(j+x)^2-2x(j+x)\right\}\log(j+x)$$

$$-\sum_{k=0}^{\infty}\frac{1}{k+1}\sum_{j=0}^{k}\binom{k}{j}(-1)^j(2j+1)x$$

Referring to (4.3.108) and (4.3.109) we see that

$$\varsigma'(0,x)-\varsigma(0,x)=\sum_{n=0}^{\infty}\frac{1}{n+1}\sum_{j=0}^{n}\binom{n}{j}(-1)^j(j+x)\log(j+x)$$

$$2\varsigma'(-1,x)-\varsigma(-1,x)=\sum_{n=0}^{\infty}\frac{1}{n+1}\sum_{j=0}^{n}\binom{n}{j}(-1)^j(j+x)^2\log(j+x)$$

We also have equivalent results when $x\to 1-x$. We also see from (4.3.76e) that

$$\sum_{k=0}^{\infty}\frac{1}{k+1}\sum_{j=0}^{k}\binom{k}{j}(-1)^j(2j+1)=0$$

Hence we obtain

$$F_1(x)=-2\varsigma'(-1,1-x)+\varsigma(-1,1-x)-2x\varsigma'(0,1-x)+2x\varsigma(0,1-x)$$

$$+2\varsigma'(-1,x)-\varsigma(-1,x)-2x\varsigma'(0,x)+2x\varsigma(0,x)$$

and with $x=0$ we have

$$F_1(0)=-2\varsigma'(-1,1)+\varsigma(-1,1)+2\varsigma'(-1,0)-\varsigma(-1,0)=0$$

where we have used (i) (4.3.112) $\varsigma'(-n,0)=\varsigma'(-n)$, (ii) Lerch's identity (4.3.116) to see that $\lim_{x\to 0}x\varsigma'(0,x)=0$ and (iii) (4.3.110) to see that $\varsigma(-1,1)=\varsigma(-1,0)=-\frac{1}{12}$.

We have therefore shown that

$(4.3.156)$ $\quad I_1 = \int\limits_0^x \pi u \cot \pi u \, du = \varsigma'(-1, x) - \varsigma'(-1, 1-x) - x[\varsigma'(0, x) + \varsigma'(0, 1-x)]$

$$+ x[\varsigma(0, x) + \varsigma(0, 1-x)] - \frac{1}{2}[\varsigma(-1, x) - \varsigma(-1, 1-x)]$$

We note from Lerch's identity $(4.3.116)$ that for $0 < x < 1$

$(4.3.156a)$ $\quad \varsigma'(0, x) + \varsigma'(0, 1-x) = \log[\Gamma(x)\Gamma(1-x)] - \log(2\pi) = -\log(2\sin \pi x)$

We see from $(4.3.110)$ that

$$\varsigma(-m, u) = -\frac{B_{m+1}(u)}{m+1} \quad \text{for } m \in \mathbf{N}_o$$

and from $(A.14)$ we have

$$B_{m+1}(1-x) = (-1)^{m+1} B_{m+1}(x)$$

and we therefore see that

$(4.3.157)$ $\quad \varsigma(-m, x) + (-1)^m \varsigma(-m, 1-x) = 0$

Specifically we have

$(4.3.157a)$ $\qquad \varsigma(0, x) + \varsigma(0, 1-x) = 0 \text{ and } \varsigma(-1, x) - \varsigma(-1, 1-x) = 0$

and

$(4.3.157b)$ $\qquad \varsigma\left(-2m, \frac{1}{2}\right) = 0$

We therefore obtain for $0 \le x < 1$

$(4.3.158)$ $\qquad I_1 = \int\limits_0^x \pi u \cot \pi u \, du = \varsigma'(-1, x) - \varsigma'(-1, 1-x) + x\log(2\sin \pi x)$

This integral is most easily verified by differentiating both sides where we get for derivative of the right-hand side

$$I_1'(x) = \frac{d}{dx}\varsigma'(-1, x) - \frac{d}{dx}\varsigma'(-1, 1-x) + \log(2\sin \pi x) + \pi x \cot \pi x$$

We have already seen from (4.3.142) that

$$\frac{d}{dx}\varsigma'(-1,x) = \varsigma'(0,x) - \varsigma(0,x) \text{ and } \frac{d}{dx}\varsigma'(-1,1-x) = -\varsigma'(0,1-x) + \varsigma(0,1-x)$$

and hence

$$I_1'(x) = \varsigma'(0,x) - \varsigma(0,x) + \varsigma'(0,1-x) - \varsigma(0,1-x) + \log(2\sin\pi x) + \pi x\cot\pi x$$

Then, using (4.3.156a) and (4.3.157a) we see that

$$I_1'(x) = \pi x\cot\pi x$$

and this is obviously the derivative of the left-hand side of (4.3.158). Then, using an "easy" integral (say, $x = 1/2$ or $x = 0$, we may complete the verification.

We now refer to equation (6.5a) in Volume V: this enables us to write for $0 \le x < 1$

$$\int_0^x \pi u\cot\pi u\,du = 2\pi\sum_{n=1}^{\infty}\int_0^x u\sin 2\pi u\,du$$

$$= 2\pi\sum_{n=1}^{\infty}\left[\frac{\sin 2\pi nx}{4\pi^2 n^2} - \frac{x\cos 2\pi nx}{2\pi n}\right]$$

Therefore we obtain

$$\varsigma'(-1,x) - \varsigma'(-1,1-x) + x\log(2\sin\pi x) = 2\pi\sum_{n=1}^{\infty}\left[\frac{\sin 2\pi nx}{4\pi^2 n^2} - \frac{x\cos 2\pi nx}{2\pi n}\right]$$

and, using (7.8) we see that

$$\varsigma'(-1,x) - \varsigma'(-1,1-x) = \frac{1}{2\pi}\sum_{n=1}^{\infty}\frac{\sin 2\pi nx}{n^2}$$

Referring to (6.48) in Volume V we see that

$$\pi u\cot\pi u = -2\sum_{n=0}^{\infty}\varsigma(2n)u^{2n} \qquad , (|u| < 1)$$

(where $\varsigma(0) = -1/2$) and substituting this in (4.3.158) we obtain

(4.3.158i)
$$-2\sum_{n=0}^{\infty}\frac{\varsigma(2n)}{2n+1}x^{2n+1} = \varsigma'(-1,x) - \varsigma'(-1,1-x) + x\log(2\sin\pi x)$$

and comparing this with (6.55) in Volume V

(4.3.158ii)
$$2\sum_{n=0}^{\infty}\frac{\varsigma(2n)}{2n+1}x^{2n+1} = -\log(2\pi)x + \log\frac{G(1+x)}{G(1-x)} \quad , \; (|x|<1)$$

we see that (as is also shown in (4.3.162) below)

$$\varsigma'(-1,x) - \varsigma'(-1,1-x) + x\log(2\sin\pi x) = x\log(2\pi) - \log\frac{G(1+x)}{G(1-x)}$$

which is equivalent to Kinkelin's formula (4.3.87).

Integration of (4.3.158) by parts gives us

$$\int_0^x \pi u \cot \pi u \, du = x\log\sin\pi x - \int_0^x \log\sin\pi u \, du$$

and we then obtain

(4.3.158a)
$$\int_0^x \log[2\sin\pi u]du = -[\varsigma'(-1,x) - \varsigma'(-1,1-x)]$$

and this may be written as

(4.3.158b)
$$\int_0^{\pi x} \log[2\sin t]dt = -\pi[\varsigma'(-1,x) - \varsigma'(-1,1-x)]$$

or as

(4.3.158c)
$$\int_0^{\frac{p}{q}\pi} \log(2\sin t)dt = -\pi\left[\varsigma'\left(-1,\frac{p}{q}\right) - \varsigma'\left(-1,1-\frac{p}{q}\right)\right]$$

With $x = 1/2$ in (4.3.158) we immediately obtain

$$\int_0^{\frac{1}{2}} \pi u \cot \pi u \, du = \frac{1}{2}\log 2$$

We see that

$$\int_0^{1/2} \pi u \cot \pi u \, du = \frac{1}{\pi} \int_0^{\pi/2} x \cot x \, dx$$

and integration by parts results in

$$= -\frac{1}{\pi} \int_0^{\pi/2} \log \sin x \, dx = \frac{1}{2} \log 2$$

where we have used Euler's integral (3.2) from Volume I. We see that this is in agreement with (4.3.156).

A more interesting result is obtained by letting $x = 1/4$ in (4.3.156) when we get

(4.3.159) $$\int_0^{1/4} \pi u \cot \pi u \, du = \varsigma'\left(-1, \frac{1}{4}\right) - \varsigma'\left(-1, \frac{3}{4}\right) + \frac{1}{4} \log 2$$

We also have

$$\int_0^{1/4} \pi u \cot \pi u \, du = -\frac{1}{\pi} \int_0^{\pi/4} \log \sin x \, dx$$

and from Bradley's catalogue [33] of formulae for Catalan's constant $G$ we note that

(4.3.160) $$\int_0^{\pi/4} \log \sin x \, dx = -\frac{\pi}{4} \log 2 - \frac{1}{2} G$$

where $G = \sum_{n=0}^{\infty} \frac{(-1)^n}{(2n+1)^2} = 0.915965....$ It may be noted that (4.3.160) is derived in an alternative manner in (6.107d) in Volume V. Therefore we get

(4.3.160a) $$\int_0^{1/4} \pi u \cot \pi u \, du = \frac{1}{4} \log 2 + \frac{G}{2\pi}$$

Comparing the two results we see that

(4.3.160b)
$$\varsigma'\left(-1,\frac{1}{4}\right) - \varsigma'\left(-1,\frac{3}{4}\right) = \frac{G}{2\pi}$$

The above result was derived by Adamchik [2a] in 1997 using Lerch's transformation formula. Adamchik [2a] also showed that

(4.3.161a)
$$\varsigma'\left(-1,\frac{1}{4}\right) = \frac{G}{4\pi} - \frac{1}{8}\varsigma'(-1)$$

(4.3.161b)
$$\varsigma'\left(-1,\frac{3}{4}\right) = -\frac{G}{4\pi} - \frac{1}{8}\varsigma'(-1)$$

In the same paper, using Rademacher's formula [14, p.261] for the Hurwitz zeta function, Adamchik [2a] also showed that

(4.3.161c)
$$\varsigma'\left(1,\frac{1}{4}\right) - \varsigma'\left(1,\frac{3}{4}\right) = \pi\left[\gamma + 4\log 2 + 3\log\pi - 4\log\Gamma\left(\frac{1}{4}\right)\right]$$

As seen above, we have

$$I_1 = \int_0^x \pi u \cot \pi u \, du = \varsigma'(-1,x) - \varsigma'(-1,1-x) + x\log(2\sin\pi x)$$

and we recall Kinkelin's formula (4.3.87) of 1860

$$\int_0^x \pi u \cot \pi u \, du = x\log(2\pi) - \log\frac{G(1+x)}{G(1-x)}$$

which gives us

(4.3.162)
$$\log\frac{G(1-x)}{G(1+x)} = \varsigma'(-1,x) - \varsigma'(-1,1-x) + x\log\frac{\sin\pi x}{\pi}$$

We also recall (4.3.131)

(4.3.163)
$$n\int_0^u \varsigma'(1-n,x)\,dx = \frac{B_{n+1} - B_{n+1}(u)}{n(n+1)} + \varsigma'(-n,u) - \varsigma'(-n)$$

which gives us for $n = 2$

$$2\int_0^u \varsigma'(-1,x)\,dx = -\frac{1}{6}B_3(u) + \varsigma'(-2,u) - \varsigma'(-2)$$

Letting $x \to 1-x$ we get

$$-2\int_1^{1-u} \varsigma'(-1,1-x)\,dx = -\frac{1}{6}B_3(u) + \varsigma'(-2,u) - \varsigma'(-2)$$

and using (4.3.131) we see that

$$\int_0^1 \varsigma'(-1,x)\,dx = \int_0^1 \varsigma'(-1,1-x)\,dx = 0$$

Hence we obtain

$$-2\int_0^{1-u} \varsigma'(-1,1-x)\,dx = -\frac{1}{6}B_3(u) + \varsigma'(-2,u) - \varsigma'(-2)$$

which is equivalent to

$$-2\int_0^u \varsigma'(-1,1-x)\,dx = -\frac{1}{6}B_3(1-u) + \varsigma'(-2,1-u) - \varsigma'(-2)$$

and using (A.14) from Volume VI this becomes

(4.3.163a) $\quad -2\int_0^u \varsigma'(-1,1-x)\,dx = \frac{1}{6}B_3(u) + \varsigma'(-2,1-u) - \varsigma'(-2)$

Hence we see that

(4.3.163b) $\quad \int_0^u [\varsigma'(-1,x) - \varsigma'(-1,1-x)]dx = \frac{1}{2}[\varsigma'(-2,u) + \varsigma'(-2,1-u)] - \varsigma'(-2)$

We therefore obtain

(4.3.164) $\quad \int_0^y dx \int_0^x \pi u \cot \pi u\,du = \frac{1}{2}[\varsigma'(-2,y) + \varsigma'(-2,1-y)] - \varsigma'(-2) + \int_0^y x\log(2\sin \pi x)dx$

Integration by parts gives us

$$\int\limits_0^x \pi u \cot \pi u \, du = x \log \sin \pi x - \int\limits_0^x \log \sin \pi u \, du$$

We note from (7.8) in Volume V that for $u \in (0,1)$

$$\log \sin \pi u = -\log 2 - \sum_{n=1}^{\infty} \frac{\cos 2n\pi u}{n}$$

and hence, as in Tolstov's book [130, p.148], we have for $x \in [0,1]$

$$\int\limits_0^x \log \sin \pi u \, du = -x \log 2 - \frac{1}{2\pi} \sum_{n=1}^{\infty} \frac{\sin 2n\pi x}{n^2}$$

Therefore we have

$$\int\limits_0^y dx \int\limits_0^x \pi u \cot \pi u \, du = \int\limits_0^y x \log \sin \pi x \, dx + \frac{1}{2} y^2 \log 2 - \frac{1}{(2\pi)^2} \sum_{n=1}^{\infty} \frac{\cos 2n\pi y}{n^3} + \frac{\varsigma(3)}{(2\pi)^2}$$

and comparing this with (4.3.164) gives us Adamchik's result [2a] (see (4.3.167) below)

(4.3.165)    $$\frac{1}{2}\left[\varsigma'(-2,y) + \varsigma'(-2,1-y)\right] = -\frac{1}{(2\pi)^2} \sum_{n=1}^{\infty} \frac{\cos 2n\pi y}{n^3} + \frac{\varsigma(3)}{(2\pi)^2} + \varsigma'(-2)$$

Using the functional equation for the Riemann zeta function we showed in (F.8b) in Volume VI that

$$\varsigma'(-2) = -\frac{\varsigma(3)}{4\pi^2}$$

and hence we get

(4.3.166)    $$\varsigma'(-2,y) + \varsigma'(-2,1-y) = -\frac{2}{(2\pi)^2} \sum_{n=1}^{\infty} \frac{\cos 2n\pi y}{n^3}$$

Note that this identity is also valid at the end points $y = 0,1$. A different proof is obtained in (4.4.229l) in Volume IV. Letting $y = 1/2$ gives us

$$\varsigma'\left(-2, \frac{1}{2}\right) = \frac{3\varsigma(3)}{16\pi^2}$$

We now integrate (4.3.166) using (4.3.163), and noting that

(4.3.166a) $\quad -n\int\limits_0^u \varsigma'(1-n,1-y)\,dy = \dfrac{B_{n+1}-B_{n+1}(1-u)}{n(n+1)} + \varsigma'(-n,1-u) - \varsigma'(-n)$

we easily obtain

$$\int\limits_0^x [\varsigma'(-2,y)+\varsigma'(-2,1-y)]\,dy = \frac{1}{3}[\varsigma'(-3,x)-\varsigma'(-3,1-x)]$$

Hence we see that

$$\varsigma'(-3,x)-\varsigma'(-3,1-x) = \frac{3!}{(2\pi)^3}\sum_{n=1}^\infty \frac{\sin 2n\pi y}{n^4}$$

Continuing in this vein, we thereby derive Adamchik's result [2a]

(4.3.167) $\quad \varsigma'(-2n,x)+\varsigma'(-2n,1-x) = (-1)^n \dfrac{(2n)!}{(2\pi)^{2n}}\mathrm{Cl}_{2n+1}(2\pi x)$

$$\varsigma'(-2n-1,x)-\varsigma'(-2n-1,1-x) = \frac{(2n+1)!}{(2\pi)^{2n+1}}\mathrm{Cl}_{2n+2}(2\pi x)$$

where $\mathrm{Cl}_n(x)$ is the Clausen function defined by

$$\mathrm{Cl}_n(x) = \begin{cases} \mathrm{Re}\,Li_n(e^{-ix}), & n \text{ is odd} \\ -\mathrm{Im}\,Li_n(e^{-ix}), & n \text{ is even} \end{cases}$$

It may also be noted that Ramanujan [21, Part I, p.261] showed that

$$\frac{1}{2}\int\limits_0^x u^n \cot\left(\frac{u}{2}\right)du = \cos\left(\frac{n\pi}{2}\right)n!\,\zeta(n+1) - \sum_{j=0}^n (-1)^{j(j+1)/2}\frac{\Gamma(n+1)}{\Gamma(n+1-j)}x^{n-j}\mathrm{Cl}_{j+1}(x)$$

and also see the paper by Srivastava et al [125aa, Eq (5.12)]. Ramanujan's formula was recently employed by Muzaffar [103ac]. With $x = \pi$ we have

$$2^n \int\limits_0^{\pi/2} x^n \cot x\,dx = \cos\left(\frac{n\pi}{2}\right)n!\,\varsigma(n+1) - \sum_{j=0}^n (-1)^{j(j+1)/2}\frac{\Gamma(n+1)}{\Gamma(n+1-j)}\pi^{n-j}\mathrm{Cl}_{j+1}(\pi)$$

In passing, we mention that G&R [74, p.428] contains the integral

$$\int_0^{\pi/2} x^n \cot x \, dx = \left(\frac{\pi}{2}\right)^n \left[\frac{1}{n} - 2\sum_{k=1}^{\infty} \frac{\varsigma(2k)}{4^k(n+2k)}\right]$$

and we therefore have

$$\cos\left(\frac{n\pi}{2}\right) n! \varsigma(n+1) - \sum_{j=0}^{n} (-1)^{j(j+1)/2} \frac{\Gamma(n+1)}{\Gamma(n+1-j)} \pi^{n-j} \mathrm{Cl}_{j+1}(\pi) = \pi^n \left[\frac{1}{n} - 2\sum_{k=1}^{\infty} \frac{\varsigma(2k)}{4^k(n+2k)}\right]$$

With $n = 2$ this becomes

$$-2\varsigma(3) - \sum_{j=0}^{2} (-1)^{j(j+1)/2} \frac{\Gamma(n+1)}{\Gamma(n+1-j)} \pi^{n-j} \mathrm{Cl}_{j+1}(\pi) = \pi^2 \left[\frac{1}{2} - \sum_{k=1}^{\infty} \frac{\varsigma(2k)}{4^k(1+k)}\right]$$

$$-2\varsigma(3) - \pi^2 \mathrm{Cl}_1(\pi) - 2\pi \mathrm{Cl}_2(\pi) - 2\mathrm{Cl}_3(\pi) = \pi^2 \left[\frac{1}{2} - \sum_{k=1}^{\infty} \frac{\varsigma(2k)}{4^k(1+k)}\right]$$

$$-2\varsigma(3) + \pi^2 \log 2 - 2\mathrm{Cl}_3(\pi) = \pi^2 \left[\frac{1}{2} - \sum_{k=1}^{\infty} \frac{\varsigma(2k)}{4^k(1+k)}\right]$$

and we end up with

$$\pi^2 \log 2 - \varsigma(3) = \pi^2 \left[\frac{1}{2} - \sum_{k=1}^{\infty} \frac{\varsigma(2k)}{4^k(1+k)}\right]$$

However, this does not agree with (6.79).

The above analysis may be easily extended as follows. Upon multiplication of (4.3.154) by $u$ we see that

$$u^2 \sum_{k=0}^{\infty} \frac{1}{k+1} \sum_{j=0}^{k} \binom{k}{j} (-1)^j \log \frac{j+1-u}{j+u} = \pi u^2 \cot \pi u$$

and integration then gives us

$$\sum_{k=0}^{\infty} \frac{1}{k+1} \sum_{j=0}^{k} \binom{k}{j} (-1)^j \int_0^x u^2 \log \frac{j+1-u}{j+u} \, du = \int_0^x \pi u^2 \cot \pi u \, du$$

Integration by parts results in

$$3\int_0^x u^2 \log\frac{j+1-u}{j+u}\,du = [u^3 - (j+1)^3]\log(j+1-u) - (u^3 + j^3)\log(u+j) - \frac{1}{2}(2j+1)(u^2 + 2u)\Big|_0^x$$

and we have

$$3\int_0^x \pi u^2 \cot \pi u\,du = F_2(x) - F_2(0)$$

where $F_2(x)$ is defined below.

As before, since $j^p = [(j+b) - b]^p$ we have

$$j^3 = (j+b)^3 - 3b(j+b)^2 + 3b^2(j+b) - b^3$$

$$(j+1)^3 = (j+1-b)^3 + 3b(j+1-b)^2 + 3b^2(j+1-b) + b^3$$

Hence we have

$$F_2(x) =$$

$$\sum_{k=0}^{\infty}\frac{1}{k+1}\sum_{j=0}^{k}\binom{k}{j}(-1)^j\left([x^3 - (j+1)^3]\log(j+1-x) - (x^3 + j^3)\log(x+j) - \frac{1}{2}(2j+1)(x^2 + 2x)\right)$$

and this may be written as

$$-\sum_{k=0}^{\infty}\frac{1}{k+1}\sum_{j=0}^{k}\binom{k}{j}(-1)^j\left\{(j+1-x)^3 + 3x(j+1-x)^2 + 3x^2(j+1-x)\right\}\log(j+1-x)$$

$$-\sum_{k=0}^{\infty}\frac{1}{k+1}\sum_{j=0}^{k}\binom{k}{j}(-1)^j\left\{(j+x)^3 - 3x(j+x)^2 + 3x^2(j+x)\right\}\log(j+x)$$

$$-\frac{1}{2}\sum_{k=0}^{\infty}\frac{1}{k+1}\sum_{j=0}^{k}\binom{k}{j}(-1)^j(2j+1)(x^2 + 2x)$$

As with (4.3.108) and (4.3.109), we also have

$$3\varsigma'(-2, x) - \varsigma(-2, x) = \sum_{n=0}^{\infty}\frac{1}{n+1}\sum_{j=0}^{n}\binom{n}{j}(-1)^j(j+x)^3 \log(j+x)$$

We also have an equivalent result when $x \to 1 - x$.

Hence we obtain

$F_2(x) =$

$-3\varsigma'(-2,1-x) + \varsigma(-2,1-x) - 6x\varsigma'(-1,1-x) + 3x\varsigma(-1,1-x) - 3x^2\varsigma'(0,1-x) + 3x^2\varsigma(0,1-x)$

$-3\varsigma'(-2,x) + \varsigma(-2,x) + 6x\varsigma'(-1,-x) - 3x\varsigma(-1,x) - 3x^2\varsigma'(0,x) + 3x^2\varsigma(0,x)$

and with $x = 0$ we have

$F_2(0) = -3\varsigma'(-2,1) + \varsigma(-2,1) - 3\varsigma'(-2,0) + \varsigma(-2,0)$

$= -6\varsigma'(-2,1) = -6\varsigma'(-2)$

where we have again used (i) (4.3.112) $\varsigma'(-n,0) = \varsigma'(-n)$, (ii) Lerch's identity (4.3.116) to see that $\lim_{x \to 0} x^2\varsigma'(0,x) = 0$ and (iii) (4.3.110) to see that $\varsigma(-2,0) = 0$.

Using (F.8b) $\varsigma'(-2) = -\dfrac{\varsigma(3)}{4\pi^2}$ we therefore have $F(0) = \dfrac{3\varsigma(3)}{2\pi^2}$.

Hence we get

(4.3.168) $\displaystyle\int_0^x \pi u^2 \cot \pi u \, du =$

$-[\varsigma'(-2,x) + \varsigma'(-2,1-x)] + 2x[\varsigma'(-1,x) - \varsigma'(-1,1-x)] - x^2[\varsigma'(0,x) + \varsigma'(0,1-x)]$

$+ \dfrac{1}{3}[\varsigma(-2,x) + \varsigma(-2,1-x)] - x[\varsigma(-1,x) - \varsigma(-1,1-x)] + x^2[\varsigma(0,x) - \varsigma(0,1-x)] - \dfrac{\varsigma(3)}{2\pi^2}$

The above integral may also be easily verified by differentiation.

As seen above, this may be simplified to

(4.3.168a) $\displaystyle\int_0^x \pi u^2 \cot \pi u \, du =$

$-[\varsigma'(-2,x) + \varsigma'(-2,1-x)] + 2x[\varsigma'(-1,x) - \varsigma'(-1,1-x)] + x^2 \log(2\sin \pi x) - \dfrac{\varsigma(3)}{2\pi^2}$

This may also be expressed as

(4.3.168b) $\qquad \pi \int\limits_0^x (u^2 - 2xu)\cot \pi u \, du = -[\varsigma'(-2,x) + \varsigma'(-2,1-x)] - x^2 \log(2\sin \pi x) - \dfrac{\varsigma(3)}{2\pi^2}$

With $x = 1/2$ in (4.3.168) we see that the terms involving $\varsigma'\left(-1, \dfrac{1}{2}\right)$ conveniently cancel and we obtain

(4.3.168c) $\qquad \int\limits_0^{1/2} \pi u^2 \cot \pi u \, du = -2\varsigma'\left(-2, \dfrac{1}{2}\right) + \dfrac{1}{4}\log 2 - \dfrac{\varsigma(3)}{2\pi^2}$

We also have for $x = 1/2$

$$\int\limits_0^{1/2} \pi u^2 \cot \pi u \, du = \dfrac{1}{\pi^2} \int\limits_0^{\pi/2} x^2 \cot x \, dx$$

and integration by parts results in

$$= -\dfrac{2}{\pi^2} \int\limits_0^{\pi/2} x \log \sin x \, dx$$

$$= \dfrac{1}{4}\log 2 - \dfrac{7\varsigma(3)}{8\pi^2}$$

where we have used Euler's integral (1.11) and (1.12) from Volume I. We therefore see that

(4.3.168d) $\qquad \varsigma'\left(-2, \dfrac{1}{2}\right) = \dfrac{3\varsigma(3)}{16\pi^2}$

Adamchik [4] also showed in 1998 that for positive integers $n$ and $0 < x < 1$

(4.3.170) $\qquad \varsigma'\left(-n,x\right) + (-1)^n \varsigma'\left(-n,1-x\right) = \pi i \dfrac{B_{n+1}(x)}{n+1} + e^{-\frac{in\pi}{2}} \dfrac{n!}{(2\pi)^n} Li_{n+1}(e^{2\pi i x})$

and for $x = 1/2$ we have

(4.3.170a) $\qquad \varsigma'\left(-n,\dfrac{1}{2}\right) + (-1)^n \varsigma'\left(-n,\dfrac{1}{2}\right) = \pi i \dfrac{B_{n+1}(1/2)}{n+1} + e^{-\frac{in\pi}{2}} \dfrac{n!}{(2\pi)^n} Li_{n+1}(-1)$

and with $n = 2m$ we get

(4.3.170b) $\qquad 2\varsigma'\left(-2m,\dfrac{1}{2}\right) = (-1)^m \dfrac{(2m)!}{(2\pi)^{2m}} Li_{2m+1}(-1) = (-1)^{m+1} \dfrac{(2m)!}{(2\pi)^{2m}}(1-2^{-2m})\varsigma(2m+1)$

We therefore see that $\varsigma'\left(-2,\dfrac{1}{2}\right) = \dfrac{3\varsigma(3)}{16\pi^2}$ in agreement with (4.3.168d).

Earlier in 1998 Miller and Adamchik [103] proved that

(4.3.170c) $\qquad \varsigma'\left(-2k+1,\dfrac{1}{2}\right) = -\dfrac{B_{2k}\log 2}{k2^{2k}} - \dfrac{(2^{2k-1}-1)\varsigma'(-2k+1)}{2^{2k-1}}$

and with $k=1$ we get

(4.3.170d) $\qquad \varsigma'\left(-1,\dfrac{1}{2}\right) = -\dfrac{B_2\log 2}{4} - \dfrac{1}{2}\varsigma'(-1) = -\dfrac{1}{24}\log 2 - \dfrac{1}{2}\varsigma'(-1)$

We note that (4.3.170d) may also be easily found by differentiating the identity (4.3.123)

(4.3.170e) $\qquad \varsigma\left(s,\dfrac{1}{2}\right) = (2^s-1)\varsigma(s)$

With $x=1/4$ in (4.3.168) we have

(4.3.171)
$$\int_0^{1/4} \pi x^2 \cot \pi x\, dx = -\left[\varsigma'\left(-2,\dfrac{1}{4}\right)+\varsigma'\left(-2,\dfrac{3}{4}\right)\right] + \dfrac{1}{2}\left[\varsigma'\left(-1,\dfrac{1}{4}\right)+\varsigma'\left(-1,\dfrac{3}{4}\right)\right] + \dfrac{1}{32}\log 2 - \dfrac{\varsigma(3)}{2\pi^2}$$

$$= -\left[\varsigma'\left(-2,\dfrac{1}{4}\right)+\varsigma'\left(-2,\dfrac{3}{4}\right)\right] + \dfrac{G}{4\pi} + \dfrac{1}{32}\log 2 - \dfrac{\varsigma(3)}{2\pi^2}$$

We note from (6.69o) in Volume V (and see also [45]) that

$$\int_0^{\pi/4} x^2 \cot x\, dx = -\dfrac{71}{128}\varsigma(3) + \dfrac{\pi G}{4} + \dfrac{\pi^2}{32}\log 2$$

and hence by an elementary change of variables we get

$$\int_0^{1/4} \pi x^2 \cot \pi x\, dx = -\dfrac{71}{128\pi^2}\varsigma(3) + \dfrac{G}{4\pi} + \dfrac{1}{32}\log 2$$

We therefore deduce that

$$(4.3.172) \qquad \varsigma'\left(-2,\frac{1}{4}\right) + \varsigma'\left(-2,\frac{3}{4}\right) = \frac{3\varsigma(3)}{64\pi^2}$$

which was also originally derived by Adamchik [2a]. We also note from (4.3.168d) that

$$(4.3.173) \qquad \varsigma'\left(-2,\frac{1}{4}\right) + \varsigma'\left(-2,\frac{3}{4}\right) = \varsigma'\left(-2,\frac{1}{2}\right)$$

and this is a direct consequence of (i) Adamchik's result [2a]

$$(4.3.174) \qquad \varsigma'\left(-2n,x\right) + \varsigma'\left(-2n,1-x\right) = (-1)^n \frac{(2n)!}{(2\pi)^{2n}} \mathrm{Cl}_{2n+1}(2\pi x)$$

$$\varsigma'\left(-2n-1,x\right) - \varsigma'\left(-2n-1,1-x\right) = \frac{(2n+1)!}{(2\pi)^{2n+1}} \mathrm{Cl}_{2n+2}(2\pi x)$$

and (ii) the formula [45, equation (5.11)] (where I have corrected a misprint)

$$(4.3.175) \qquad \mathrm{Cl}_{2n+1}(\pi/2) = \frac{1-2^{2n}}{2^{4n+1}} \varsigma(2n+1)$$

This then gives us for $x = 1/4$

$$(4.3.176)$$

$$\varsigma'\left(-2n,\frac{1}{4}\right) + \varsigma'\left(-2n,\frac{3}{4}\right) = (-1)^n \frac{(2n)!}{(2\pi)^{2n}} \mathrm{Cl}_{2n+1}(\pi/2) = (-1)^n \frac{1-2^{2n}}{2^{4n+1}} \frac{(2n)!}{(2\pi)^{2n}} \varsigma(2n+1)$$

and, having regard to equation (2.13) in [45], we also have in terms of the generalised Glaisher-Kinkelin constants $A_{2n}$

$$(4.3.177) \qquad \varsigma'\left(-2n,\frac{1}{4}\right) + \varsigma'\left(-2n,\frac{3}{4}\right) = \frac{2^{2n}-1}{2^{4n+2}} \log A_{2n}$$

We also have

$$(4.3.178) \qquad \varsigma'\left(-1,\frac{1}{4}\right) - \varsigma'\left(-1,\frac{3}{4}\right) = \frac{1}{2\pi} \mathrm{Cl}_2(\pi/2) = \frac{G}{2\pi}$$

It is clear that the above analysis may be extended to $\int_0^x \pi u^n \cot \pi u \, du$ and the following could also be tackled

$$(4.3.179) \qquad \sum_{k=0}^{\infty} \frac{1}{k+1} \sum_{j=0}^{k} \binom{k}{j} (-1)^j \int_0^x B_n(u) \log \frac{j+1-u}{j+u} \, du = \int_0^x \pi B_n(u) \cot \pi u \, du$$

Substituting (6.48) in (4.3.168a) we get

$$-\sum_{n=0}^{\infty} \frac{\varsigma(2n)}{n+1} x^{2n+2} =$$

$$-[\varsigma'(-2,x) + \varsigma'(-2,1-x)] + 2x[\varsigma'(-1,x) - \varsigma'(-1,1-x)] + x^2 \log(2\sin \pi x) - \frac{\varsigma(3)}{2\pi^2}$$

and with $x = 1/2$ we get

$$\sum_{n=0}^{\infty} \frac{\varsigma(2n)}{(n+1)2^{2n}} = 8\varsigma'(-2,1/2) - \log 2 + \frac{2\varsigma(3)}{\pi^2}$$

Comparing this with (6.79) in Volume V

$$\sum_{n=0}^{\infty} \frac{\varsigma(2n)}{(n+1)2^{2n}} = \frac{7}{2\pi^2} \varsigma(3) - \log 2$$

we see that $\varsigma'(-2,1/2) = \dfrac{3\varsigma(3)}{16\pi^2}$ in accordance with (4.3.168d).

Integrating (4.3.158i) in conjunction with (4.3.168b) results in

$$-\sum_{n=0}^{\infty} \frac{\varsigma(2n)}{(2n+1)(n+1)} t^{2n+2} = \frac{1}{2}[\varsigma'(-2,t) + \varsigma'(-2,1-t)] - \varsigma'(-2) + \int_0^t x \log(2\sin \pi x) dx$$

and with $t = 1/2$ we get

$$-\sum_{n=0}^{\infty} \frac{\varsigma(2n)}{(2n+1)(n+1)2^{2n+2}} = \varsigma'\left(-2,\frac{1}{2}\right) - \varsigma'(-2) + \int_0^{1/2} x \log(2\sin \pi x) dx$$

Using (4.3.168d) $\varsigma'\left(-2,\dfrac{1}{2}\right) = \dfrac{3\varsigma(3)}{16\pi^2}$ and (F.8b) $\varsigma'(-2) = -\dfrac{\varsigma(3)}{4\pi^2}$ this becomes

$$-\sum_{n=0}^{\infty}\frac{\varsigma(2n)}{(2n+1)(n+1)2^{2n+2}}=\frac{7\varsigma(3)}{16\pi^2}+\int_0^{1/2}x\log(2\sin\pi x)dx$$

and from Euler's integral (1.11) we have

$$\int_0^{1/2}x\log(2\sin\pi x)dx=\frac{1}{8}\log 2+\frac{7\varsigma(3)}{16\pi^2}-\frac{1}{8}\log 2$$

We then deduce that

$$-\sum_{n=0}^{\infty}\frac{\varsigma(2n)}{(2n+1)(n+1)2^{2n+2}}=\frac{7\varsigma(3)}{8\pi^2}$$

and this is a well-known formula (which is also derived in (6.71) in Volume V).

$\square$

Adamchik [6c] has shown that for $\mathrm{Re}(x)>0$

(4.3.180) $\qquad \varsigma'(-n,x)-\varsigma'(-n)=(-1)^n\sum_{k=0}^{n}k!\,Q_{k,n}(x)\log\Gamma_{k+1}(x)$

where the polynomials $Q_{k,n}(x)$ are defined by

$$Q_{k,n}(x)=\sum_{j=k}^{n}(1-x)^{n-j}\binom{n}{j}\begin{Bmatrix}j\\k\end{Bmatrix}$$

and $\begin{Bmatrix}j\\k\end{Bmatrix}$ are the Stirling subset numbers defined by

$$\begin{Bmatrix}j\\k\end{Bmatrix}=k\begin{Bmatrix}n-1\\k\end{Bmatrix}+\begin{Bmatrix}n-1\\k-1\end{Bmatrix},\quad \begin{Bmatrix}n\\0\end{Bmatrix}=\begin{cases}1,&n=0\\0,&n\neq 0\end{cases}$$

Particular cases are

(4.3.180a) $\qquad \varsigma'(-1,x)-\varsigma'(-1)=x\log\Gamma(x)-\log G(x+1)$

(4.3.180b) $\qquad \varsigma'(-2,x)-\varsigma'(-2)=2\log\Gamma_3(x)+(3-2x)\log G(x)-(1-x)^2\log\Gamma(x)$

Letting $x\to 1-x$ we see that

$$\varsigma'(-1, 1-x) - \varsigma'(-1) = (1-x)\log\Gamma(1-x) - \log\Gamma(1-x) - \log G(1-x)$$

$$= -x\log\Gamma(1-x) - \log G(1-x)$$

$$\varsigma'(-2, 1-x) - \varsigma'(-2) = 2\log\Gamma_3(1-x) + (1+2x)\log G(1-x) - x^2\log\Gamma(1-x)$$

Hence we have

$$\varsigma'(-1, x) - \varsigma'(-1, 1-x) = x\log[\Gamma(x)\Gamma(1-x)] + \log\frac{G(1-x)}{G(1+x)}$$

which we have already seen above in (4.3.162).

We note from [45ac] that Choi, Cho and Srivastava have shown that

$$\log\Gamma_2(1+x) = -\frac{1}{12} + \log A - x\log\Gamma(1+x) + \varsigma'(-1, 1+x)$$

and, since $\Gamma_2(1+x) = 1/G(1+x)$, we have

$$\log G(1+x) = \frac{1}{12} - \log A + x\log\Gamma(1+x) - \varsigma'(-1, 1+x)$$

$$= \varsigma'(-1) + x\log\Gamma(1+x) - \varsigma'(-1, 1+x)$$

AN OBSERVATION BY GLASSER

Part of the following is based on an observation made by Glasser [70ab] in 1966. Let us consider the integral

$$I = \int_0^1 f(x)\psi(x)dx$$

where $f(x) = -f(1-x)$ and $f(x)$ is selected so that the integral converges. We then have

$$I = \int_0^1 f(1-t)\psi(1-t)dt = -\int_0^1 f(t)\psi(1-t)dt$$

and hence we see that

$$2I = \int_0^1 f(x)[\psi(x) - \psi(1-x)]dx$$

Therefore, since $\psi(x) - \psi(1-x) = -\pi \cot \pi x$, we have

(4.3.181) $$\int_0^1 f(x)\psi(x)dx = -\frac{\pi}{2}\int_0^1 f(x)\cot \pi x\, dx = -\pi \int_0^{1/2} f(x)\cot \pi x\, dx$$

Some functions which satisfy the condition $f(x) = -f(1-x)$ are listed below:

$$B_{2n+1}(x)\,,\; B_{2n}(x)\cos(2m+1)\pi x\,,\; B_{2n+1}(x)\sin 2m\pi x \;\text{ and }\; h\big(x(1-x)\big)\cos(2m+1)\pi x$$

since $B_n(x) = (-1)^n B_n(1-x)$ where $n$ and $m$ are positive integers. Glasser [70ab] used the function $f(x) = x(1-x)\cos \pi x$ to show that

(4.3.182) $$\int_0^1 \psi(x)x(1-x)\cos \pi x\, dx = \frac{1}{\pi^2}\left[2 - \frac{7}{2}\varsigma(3)\right]$$

Similarly, if $g(x) = g(1-x)$, we obtain

$$\int_0^1 g(x)[\psi(x) - \psi(1-x)]dx = \int_0^1 g(x)\cot \pi x\, dx = 0$$

Choi and Srivastava have evaluated a number of integrals of the form $\int_0^z x^n \psi(x+a)dx$ in [45ab] and, in the particular case, $a = z = 1$ we have

(4.3.183) $$\int_0^1 x\,\psi(x+1)dx = 1 - \frac{1}{2}\log(2\pi)$$

$$\int_0^1 x^2\psi(x+1)dx = \frac{1}{2} - \frac{1}{2}\log(2\pi) + 2\log A$$

$$\int_0^1 x^3\psi(x+1)dx = \frac{1}{3} - \frac{7}{12}\log(2\pi) + 3\log A - 2\int_0^1 \log \Gamma_3(x+1)dx$$

Since

$$\int_0^1 x^n \psi(x)\,dx = -\frac{1}{n} + \int_0^1 x^n \psi(x+1)\,dx$$

we immediately see that

(4.3.183a) $$\int_0^1 x\,\psi(x)\,dx = -\frac{1}{2}\log(2\pi)$$

$$\int_0^1 x^2 \psi(x)\,dx = -\frac{1}{2}\log(2\pi) + 2\log A$$

$$\int_0^1 x^3 \psi(x)\,dx = -\frac{7}{12}\log(2\pi) + 3\log A - 2\int_0^1 \log \Gamma_3(x+1)\,dx$$

We have from Appendix A in Volume VI

$$B_3(x) = x^3 - \frac{3}{2}x^2 + \frac{1}{2}x$$

and determine that

(4.3.183b) $$\int_0^1 B_3(x)\psi(x)\,dx = -\frac{1}{12}\log(2\pi) - 2\int_0^1 \log \Gamma_3(x+1)\,dx$$

We will see the following identity for suitably differentiable functions $p(x)$ in (6.5a) in Volume V

$$\frac{1}{2}\int_a^b p(x)\cot(\alpha x/2)\,dx = \sum_{n=0}^{\infty} \int_a^b p(x)\sin \alpha nx\,dx$$

The above equation is valid provided (i) $\sin(\alpha x/2) \neq 0 \ \forall \ x \in [a,b]$ or, alternatively, (ii) if $\sin(\alpha a/2) = 0$ then $p(a) = 0$ also. Therefore, provided $p(0) = p(1) = 0$, we have

$$\frac{1}{2}\int_0^1 p(x)\cot \pi x\,dx = \sum_{n=0}^{\infty} \int_0^1 p(x)\sin 2n\pi x\,dx$$

We now let $p(x) = B_3(x)$ and deduce that

$$\frac{1}{2}\int_0^1 B_3(x)\cot\pi x\,dx = \sum_{n=0}^{\infty}\int_0^1 B_3(x)\sin 2n\pi x\,dx$$

Using integration by parts it is easily found that

$$\int_0^1 B_3(x)\sin 2n\pi x\,dx = -B_3(x)\frac{\cos 2n\pi x}{2n\pi}\Bigg|_0^1 + \frac{3}{2n\pi}\int_0^1 B_2(x)\cos 2n\pi x\,dx$$

$$= \frac{3}{2n\pi}\int_0^1 B_2(x)\cos 2n\pi x\,dx$$

$$\int_0^1 B_2(x)\cos 2n\pi x\,dx = B_2(x)\frac{\sin 2n\pi x}{2n\pi}\Bigg|_0^1 - \frac{2}{2n\pi}\int_0^1 B_1(x)\sin 2n\pi x\,dx$$

$$= -\frac{2}{2n\pi}\int_0^1 B_1(x)\sin 2n\pi x\,dx$$

$$\int_0^1 B_1(x)\sin 2n\pi x\,dx = B_1(x)\frac{\cos 2n\pi x}{2n\pi}\Bigg|_0^1 + \frac{1}{2n\pi}\int_0^1 B_0(x)\cos 2n\pi x\,dx$$

$$= -\frac{1}{2n\pi}$$

and hence we get

$$\int_0^1 B_3(x)\sin 2n\pi x\,dx = \frac{3}{4n^3\pi^3}$$

which then gives us

$$\sum_{n=1}^{\infty}\int_0^1 B_3(x)\sin 2n\pi x\,dx = \frac{3\varsigma(3)}{4\pi^3}$$

Therefore we see that

$$\int_0^1 B_3(x)\cot\pi x\,dx = \frac{3\varsigma(3)}{2\pi^3}$$

Reference to (4.3.181) then results in

$$\int_0^1 B_3(x)\psi(x)dx = -\frac{3\varsigma(3)}{4\pi^2}$$

and comparing this with (4.3.183b) we see that

$$(4.3.184) \qquad \int_0^1 \log\Gamma_3(x+1)dx = -\frac{1}{24}\log(2\pi) + \frac{3\varsigma(3)}{8\pi^2}$$

This integral was also given by Choi and Srivastava [45ab].

Donal F. Connon
Elmhurst
Dundle Road
Matfield
Kent TN12 7HD
dconnon@btopenworld.com